\documentclass[journal]{IEEEtran}
\usepackage[notcite,notref]{}
\usepackage{amsmath,graphicx}
\usepackage{amssymb}
\usepackage{graphicx}
\usepackage{epstopdf}
\usepackage{tabularx}
\usepackage{color}
\usepackage{times}
\usepackage{epsfig}
\usepackage{array}
\usepackage{threeparttable}
\usepackage{algorithm}
\usepackage{color}
\usepackage{lineno}
\usepackage{bm}
\usepackage{subfigure}
\usepackage{cases}
\usepackage{epsfig,amssymb,subfigure,version}
\usepackage{amssymb,version,graphicx,fancybox,mathrsfs}
\usepackage{graphicx,mathrsfs,fancyhdr,amsthm}
\usepackage{subfigure,slashbox,rotating,version,fancybox,pifont,booktabs,epsfig}
\usepackage{multicol,url,multirow}
\usepackage{setspace}
\usepackage[colorlinks,urlcolor=blue]{hyperref}
\usepackage[normalem]{ulem}

\usepackage[framemethod=tikz]{mdframed}
\usepackage{lipsum}

\newtheorem{thm}{\bf Theorem}
\newtheorem{rem}{ Remark}[section]
\graphicspath{{images/}}

\newtheorem{assump}{\bf Assumption}
\newtheorem{lemma}{\bf Lemma}[section]

\title
{
 Denoising Poisson Phaseless  Measurements  via Orthogonal Dictionary Learning
}

\author{Huibin Chang, Stefano Marchesini
\thanks{H. Chang is  with School of Mathematical Sciences, Tianjin Normal University, Tianjin, China,{\tt E-mail: changhuibin@gmail.com}.
 He is currently a visiting scholar of the Center for Applied Mathematics for Energy Research Applications at Lawrence Berkeley National Laboratory.
 }
\thanks{S. Marchesini is with Computational Research Division, Lawrence Berkeley National Laboratory, Berkeley, CA, 94720 USA, {\tt E-mail: smarchesini@lbl.gov.}}
}

\begin{document}
\maketitle
\begin{abstract}
Phaseless diffraction measurements  recorded by CCD detector  are often affected by Poisson noise. In this paper, we propose a dictionary learning model by employing patches based  sparsity to denoise Poisson phaseless measurement. The model  consists of three terms: (i) A representation term by an orthogonal dictionary, (ii) an $L^0$ pseudo norm of coefficient matrix, and (iii) a Kullback-Leibler divergence  to fit phaseless Poisson data. Fast Alternating Minimization Method (AMM) and Proximal Alternating Linearized Minimization (PALM) are adopted  to solve the established model with convergence guarantee, and { especially global convergence for PALM is derived.} The subproblems for  two algorithms have fast solvers, and indeed, the solutions for the sparse coding and dictionary updating both have closed forms  due to the orthogonality of learned dictionaries.
Numerical experiments for phase retrieval using coded diffraction and ptychographic patterns are performed to show the efficiency and robustness of  proposed methods, which,  by preserving texture features, produce visually and quantitatively  improved denoised images  { compared with other phase retrieval algorithms without regularization and local sparsity  promoting algorithms.}



\end{abstract}

\begin{IEEEkeywords}
Phase Retrieval; Poisson Noise; Learned orthogonal dictionary; Alternating Minimization Method (AMM); Proximal Alternating Linearized Minimization (PALM); Diffraction imaging
\end{IEEEkeywords}

\section{Introduction}\label{intro}
\IEEEPARstart{I}{n} lensless imaging techniques,  only the magnitudes of Fourier transformation can be recorded due to physical limitation of detector technologies, and the central computational task, known as  ``Phase Retrieval'',  is to recover  underling images from phaseless measurements. It is a challenging inverse quadratic problem. For noiseless case,   related multivariable quadratic systems should be computed, the  complexity of which is possibly NP-complete \cite{ben2001lectures}.  In practice, the measurements are possibly contaminated by noise, and one shall  seek an approximation solution to the reformulated noncovex optimization problem by  maximum a posterior (MAP) of noise.

Numerically,  researchers have developed various iterative algorithms involving with a nonconvex constraint set of magnitude, such as an alternating projection algorithm \cite{Gerchberg1972,Fienup1982} for the classical phase retrieval problem{\footnote{Classical phase retrieval refers to retrieval phase from Fourier measurements}}.  Several variant heuristic algorithms popular in the optics community including \cite{Bauschke2002,Bauschke2003,Elser2003,Luke2005} were proposed to solve the classical phase retrieval problems, and one can  refer to  \cite{marchesini2007invited,wen2012} and the reference therein. These methods are based on projections onto a nonconvex set, and  it is therefore very difficult to describe and prove the convergence mathematically. Consequently, more recent work has focused on designing fast algorithms with convergence guarantee for such nonconvex minimization problem, \emph{e.g.} \cite{candes2015phase,chen2015solving,netrapalli2015phase,marchesini2015alternating,maglobally,noll2016local,chen2015fourier}, \emph{etc}, with random-masking or oversampling measurements. { In the optics community, in order to increase the imaging resolutions and robustness, a popular approach, known as ptychography \cite{rodenburg2004phase}, uses a sequence of phaseless diffraction measurements recorded while translating the specimen with respect to a constant   illumination mask.}
In order to prevent trapping into a local minimum for solving conventional nonconvex projection or minimizing methods, another scheme is to model phase retrieval by convex methods including PhaseLift \cite{candes2013} and PhaseCut \cite{Waldspurger2012}  based on semi-definite programming(SDP), as well as PhaseMax \cite{goldstein2016phasemax,bahmani2016phase} operating in the original signal space with much less computational cost compared with SDP based convex methods.

To the best of our knowledge, an early work for phase retrieval by employing  sparse \emph{prior} of underlining images is  the \textsc{shrink-wrap} algorithm \cite{marchesini2003x}, which was established by iteratively shrinking the size of support set for unknown objects. It has been successfully applied to  several ground breaking X-ray Free Electron Lasers based  experiments in \cite{chapman2006femtosecond}, which also demonstrates that sparse \emph{prior} should be considered for modeling phase retrieval problem in order to increase the robustness and accuracy. The $L^1$ norm based  phase retrieval method with a compressibility constraint was proposed in \cite{moravec2007compressive,yang2013robust}.  SDP-based convex  methods combined with $L^1$ regularization of the lifting matrix were proposed to solve a sparse phase retrieval problem in \cite{ohlsson2012cprl,li2013sparse}. With \emph{prior} knowledge of exact sparsity level(number of nonzeros), a sparse Fieup method \cite{Fienup1982} was proposed in  \cite{mukherjee2014fienup}.  A \emph{generic} two-step iterative scheme   was proposed in \cite{shechtman2014gespar}, which consisted of updating support set, and using damped Gauss-Newton algorithm to solve a nonconvex optimization problem with updated support set. A  two-stage solution technique was presented in \cite{iwen2015robust} by first using arbitrary phase retrieval technique to recover sampling measurements and then compressive sensing method  to recover the unknown signals.  If the transform domain of the unknown objects are corrupted by additive noise, a probabilistic method  based on the generalized approximate message passing algorithm was given in \cite{schniter2015compressive}. Hard thresholding based truncated gradient iterations in  \cite{wang2016sparse} were employed in order to refine  a sparse orthogonality-promoting initialization. An $L^0$ regularized variational model was established for sparse phase retrieval problem accompanied by an efficient algorithms based on  adaptive-step alternating direction method of multipliers in \cite{duan20160}. Under the assumption of objects possessing a sparse representation in the transform domain, shearlet and total variation based regularization methods were considered in \cite{loock2014phase,chang2016phase}.
Dictionary learning methods were proposed to reconstruct the image in \cite{tillmann2016dolphin,qiu2016undersampled} from Gaussian noised measurements with outliers, where   Gaussian distributions only provide a limited approximation in most real applications.


Roughly speaking, the intensity  is recorded by a detector for a limited time, which is usually affected by Poisson noise. For linear measurements, in order to remove the Poisson noise, one can  adopt a variance stabilizing transformation such as the Anscombe root transformation in \cite{anscombe1948transformation,makitalo2011optimal}, where a pointwise nonlinear transformation was introduced such that one just needs to deal with Gaussian  data. When the peak level decreases, this transformation is less efficient. One should directly tackle it and total variation  based regularization based method was proposed in \cite{le2007variational},  {where Kullback-Leibler (KL) divergence was derived based on the Bayesian framework in the case of Poisson noise}. In a similar manner, for phaseless nonlinear measurements of phase retrieval,   variational methods based on  Tikhonov and total variation regularization were proposed to recover images from  Poisson measurements in  \cite{chang2016Total,thibault2012maximum}.
{The reconstructed images shown in \cite{chang2016Total} exhibit sharp edges and clean background for very  noisy and limited data. However, as such method was built upon  the local-sparsity regularization term, there exist visible staircase artifacts, and the smaller repetitive features can not been  preserved particularly  for measurements generated with deterministic illumination in ptychographic phase retrieval problem.} In order to further increase the accuracy and visual  quality of recovery results for practical applications, one possible direction is to  exploit the redundancy between overlapping image patches to boost the sparsity, such as using dictionary learning based denoising method \cite{elad2006image,ma2013dictionary,giryes2014sparsity}, which pursues a sparse representation by a dictionary learned from image patches.

In this paper, we propose to build up a novel dictionary learning model to reconstruct  high quality images from Poisson phaseless measurements, which  can deal with both real and complex-valued images. It consists of three terms, \emph{i.e.} a representation terms by an orthogonal dictionary, an $L^0$ pseudo norm of coefficient matrix, and a Kullback-Leibler divergence to fit phaseless Poisson data. The pseudo $L^0$ norm is used as a sparse regularization of coefficient matrix with both isotropic and anisotropic forms.  Meanwhile, fast algorithms are designed with global convergence guarantee in order to solve the proposed nonconvex nondifferential optimization model. Since an  orthogonal dictionary is employed, subproblems \emph{w.r.t.} the dictionary and coefficients  both have close-form solutions. Numerous experiments demonstrate  for coded diffraction and ptychographic patterns are performed to show the efficiency of our proposed methods, which demonstrate our proposed methods can recover images with sharp edges, clean background compared with the phase retrieval methods without any regularization. Moreover, they are capable  of producing higher quality results by preserving the texture features compared with total variation based method \cite{chang2016Total}. Additionally, it also shows that the proposed method  helps to further increase the quality for reconstructing complex-valued images by unitizing the introduced anisotropic $L^0$ pseudo norm.

The rest of this paper is organized in the following way. The proposed dictionary learning model is given in section \ref{sec1}. We first propose an alternative minimization method in section \ref{sec2} by alternately minimizing the coupled variables  to compute the proposed model. Then a proximal alternating linearized minimization method is further designed in section \ref{sec3} with theoretical  analysis of global convergence. Numerous experiments are performed in section \ref{sec5} to demonstrate the efficiency of the proposed method. We conclude the paper in Section \ref{secCon}.

\section{Proposed Model}\label{sec1}

We consider phase retrieval for a 2-dimensional (2D) image  in a discrete setting, \emph{i.e.}, an underlying object  $u:\Omega=\{0,1,\cdots,n-1\}\rightarrow  \mathbb{C}^{}$ is of size $n$ with $n=n_1\times n_2$, in which   we represent a 2D object with resolution $n_1\times n_2$ in terms of a vector of size $n$  by a lexicographical order.  For classical phase retrieval problem, measured data are only the  magnitudes of the discrete Fourier transform of $u$, \emph{i.e.},
$|\mathcal{F} u|^2,$ where $|\cdot|^2$ denotes the pointwise  square of the absolute value of a vector, $\mathcal{F}:\mathbb{C}^{n}\rightarrow \mathbb{C}^{n}$ denotes the discrete Fourier transform (DFT).
 In this paper, we generalize it to a general phase retrieval problem \cite{candes2013}, \emph{i.e.}
\begin{equation}\label{eqGPR}
{\text{To find ~}u\in \mathbb C^n,~~} s.t.~~|\mathcal A u|^2= b,
\end{equation}
where $\mathcal{A}:\mathbb C^n\rightarrow \mathbb C^m$ is a linear operator in the complex Euclidean space and $b:\tilde\Omega=\{0,1,\cdots,m-1\}\rightarrow\mathbb R_+$. Assume that in the best scenario, photon counting procedures to collect the measurements $b$ are contaminated  by Poisson noise, \emph{i.e.} for noisy measurements $f\in\mathbb R^m_+$,
 the intensity values located at each position are independent and identically distributed (i.i.d.) random variables, and the probability can be written as
 \[
 \mathrm{Pr}(f|b)=\prod\limits_{j\in \tilde\Omega} \dfrac{\exp{(-b(j))(b(j))^{f(j)}}}{f(j)\!},
 \]
 with the ground truth $b\in\mathbb R_+^m$ as the mean value and variance.

 We will first give a brief review of dictionary learning based denoising, and regularized methods for denosing noisy phaseless measurements contaminated by Poisson noise in the following subsection.
 \subsection{Review of Dictionary Learning and Regularized Poisson Denoising for Phase retrieval}
\subsubsection{Dictionary Learning and Denoising}
Given a matrix $Y=[Y(0),Y(1),\cdots,Y(d-1)]\in \mathbb C^{l,d}$, each of whose columns represents a 1-D signal, the target of dictionary learning (or sparse coding)   \cite{tosic2011dictionary} is to find a ``sparse'' linear presentation for the matrix $Y$ by a learned dictionary $D=[D(0), D(1),\cdots,D(c-1)]\in\mathbb C^{l,c}$, \emph{i.e.}
$
Y\approx D\bm \alpha,
$
such that the coefficient matrix $\bm\alpha$ has more zeros (``sparse'') entries.
Obviously, the dictionary $D$ and coefficient matrix $\bm \alpha\in \mathbb C^{c,d}$ are coupled and both needed to be determined in order to derive a sparse representation.
Roughly speaking, it can be reformulated as the following minimization problems
\[
\begin{split}
\min\limits_{D,\bm\alpha} \frac12\|Y-D\bm\alpha\|^2+\tau \|\bm\alpha\|_0, ~~ s.t.~\|D(j)\|=1~\forall 0\leq j\leq c,
\end{split}
\]
where $\|\cdot\|$ denotes the Frobenius norm and $L^2$ norm for a matrix and vector respectively, $\|\cdot\|_0$ denotes the $L^0$ pseudo norm for counting the number of nonzero elements, $\tau$ is a positive constant to balance the data fitting term (first term) and the regularization term (second term), and the constraints guarantee that the dictionary is normalized.

The well known work named K-SVD denoising \cite{elad2006image} is to first produce redundant  signals  by sliding a small  window over an entire image, and then to solve  a sparse coding problem by K-SVD in \cite{aharon2006k}. The K-SVD denoising model  can be established as follows,
\[
\begin{split}
&\min\limits_{u,D,\bm\alpha} \frac12\| R(u)-D\bm\alpha\|^2+\tau \|\bm\alpha\|_0+\frac{\eta}{2}\|u-f\|^2, ~~\\
&\qquad s.t.~\|D(j)\|=1~\forall 0\leq j\leq c,
\end{split}
\]
where one redundantly  selects image patches by the linear selection operator $R(u)=[R_0 u, R_1 u, \cdots, R_{d-1} u]\in \mathcal \mathbb C^{l,d}$ to build totally $d$ atoms of the underlying image $u$ such that it can be sparsely represented by  $D$ with the coefficient matrix $\bm \alpha$,  $f$ is  the noisy image, and $D\in \mathbb C^{l,c}$ is a normalized  dictionary.
The variational methods combining K-SVD were widely applied to different image processing tasks, especially for Poisson denoising \cite{giryes2014sparsity} and deblurring \cite{ma2013dictionary}. In \cite{ma2013dictionary}, the data fitting term in the aforementioned model is replaced with Kullback-Leibler divergence by  maximum a posterior (MAP) estimate of Poisson noise, and additional total variation regularization term was also employed.
\subsubsection{Regularized Poisson Denoising for Phase Retrieval}
 Thibault and Guizar-Sicairos \cite{thibault2012maximum} first proposed to solving the following minimization problem in coherent diffractive imaging with Tikhonov regularization
\[
\min\limits_{u\in \mathbb C^n} {\lambda} \|\nabla u\|^2+\mathcal B(|\mathcal A u|^2,f),
\]
where the Kullback-Leibler divergence following Bayesian framework by MAP of Poisson noise is denoted as follows,
\begin{equation}\label{data}
\mathcal B(h,f):=\frac12\sum\limits_{j\in \tilde\Omega}  (h(j)-f(j)\log h(j)),
\end{equation}
where
$\nabla$ denotes the gradient operator and $\mathcal A$ represents the masked Fourier transformation generated by illumination masks.
Consequently, significant gains in reconstructed accuracy and sensitivity to noise  were achieved compared with the model without any regularization. In order to further increase the quality (\emph{e.g.} sharpness of edges) of recovery images, Chang \emph{et al.}  \cite{chang2016phase} established total variation regularized model for very general phase retrieval tasks
\[
\min\limits_{u}{\lambda} \mathrm{TV}(u)+\mathcal B(|\mathcal A u|^2,f),
\]
with $\mathrm{TV}$ as the total variation, which was successfully applied to denoise phaseless measurements for phase retrieval with coded diffraction, holographic and ptychographic patterns for both the real and complex-valued images.
\subsection{The proposed dictionary learning phase retrieval (DicPR) for Poisson denoising}
In this paper, we consider the noisy phase retrieval problem in which the measurements
$f=\mathrm{Poisson}(|\mathcal A u|^2) \in \mathcal {\mathbb R}^m$ in \eqref{eqGPR} are corrupted by Poisson noise, and establish a minimization problem driven by the patch based sparsity \cite{elad2006image,ma2013dictionary,giryes2014sparsity} of the underlying images, referred to as ``DicPR'', which reads,
\begin{equation}\label{model}
\begin{split}
\min\limits_{u,D,\bm \alpha} &\frac12\|D{\bm\alpha}-R(u)\|^2+\tau \|\bm \alpha\|_0+\eta\mathcal B(|\mathcal A u|^2,f),\\
s.t.&\qquad D^* D=\mathbf{I},
\end{split}
\end{equation}
where $u$ is an underlying image that we want to reconstruct from magnitude data $f$, $D\in \mathbb C^{l,c}$ is a learned orthogonal
 dictionary as \cite{bao2013fast} which combines the sparse coefficients $\bm \alpha$ by the $L^0$ pseudo norm denoted by $\|\cdot\|_0$ to represent the underling image sparsely, $\tau$ and $\eta$ are two parameters to balance the sparsity and data fitting terms and the identical matrix is denoted as $\mathbf I$, $\mathcal B(\cdot,\cdot)$ is defined in \eqref{data}.
We rewrite the model \eqref{model} as
\begin{equation}\label{modelR}
\begin{split}
\min\limits_{u,D,\bm \alpha} \Upsilon(u,D,\bm \alpha):=\mathcal H(u,D,\bm\alpha)+\mathcal F(u)+\mathcal I_{\mathbf K}(D)+\tau \|\bm \alpha\|_0,
\end{split}
\end{equation}
with  \[\mathcal H(u,D,\bm\alpha):=\frac12\|D{\bm\alpha}-R(u)\|^2,\]
\[\mathcal F(u):=\eta\mathcal B(|\mathcal A u|^2,f), \mathbf K=\{D:~D^* D=\mathbf I\},\]
 and related indicator function is introduced as $\mathcal I_{\mathbf K}$.
We introduce two different definitions of $\|\bm\alpha\|_0$ for complex-valued $\bm\alpha$ with  isotropic and anisotropic $L^0$ pseudo norms as follows:
\begin{equation}
\left\{
\begin{aligned}
&\text{Isotropic form:}  \\
&\|\bm\alpha\|^{iso}_0:=\#\big\{(s,t): |\bm\alpha(s,t)|\not=0~\forall 0\leq s\leq c-1,\\
                       &\qquad\qquad\qquad\qquad ~~0\leq t\leq d-1\big\}   \\
&\text{Anisotropic form:}\\
&\|\bm\alpha\|^{aniso}_0:= \|\Re(\bm\alpha)\|^{iso}_0+\|\Im(\bm\alpha)\|^{iso}_0, \\
\end{aligned}
\right.
\label{L0def}
\end{equation}
where $\Re$ and $\Im$ denote the real and complex parts of some complex-valued matrix.
One readily has $\|\bm\alpha\|_0^{iso}=\|\bm\alpha\|_0^{aniso}$ if $\bm \alpha\in \mathbb R^{c,d}$. For complex-valued matrix,
it is quite different, and one can see the advantage by using the anisotropic form in the numerical experiments of this paper.
 {Hereafter, for simplicity $\|\bm\alpha\|_0$ represents the isotropic version pseudo norm if not specified.}
\begin{rem}
In this paper, we propose to learn an  orthogonal dictionary such that the subproblems for solving the proposed model involving with $D$ and $\bm\alpha$ both have closed forms. In   \cite{bao2013fast} it was demonstrated that the denoising results  with  orthogonal dictionary were comparable with nonorthogonal ones. In our experiments, we set $D$ to a square matrix, \emph{i.e.} $l=c$.
\end{rem}
\begin{rem}
Tillmann, Eldar, and Mairal \cite{tillmann2016dolphin} proposed the following dictionary learning based model to denoising Gaussian noised measurements for real-valued images
\[
\begin{split}
&\min\limits_{u\in \mathbb R^n,D,\bm\alpha} \frac12\| R(u)-D\bm\alpha\|^2+\tau \|\bm\alpha\|_1+ \frac{\eta}{4}\left\||\mathcal A u|^2-f\right\|^2, ~~\\
&\qquad s.t.\quad\|D(j)\|=1~\forall 0\leq j\leq c.
\end{split}
\]

Qiu and Palomar \cite{qiu2016undersampled} generalized it to complex-valued images with different data fitting terms for the amplitude  \emph{i.e.}
$$\frac12\left\||\mathcal A u|-\sqrt{f}\right\|^2.$$
\end{rem}
\vskip.1in
The existence of the minimizer for the proposed model \eqref{model} can be readily obtained, and we just list it without proof.
\begin{thm}
There exists a triple $(u^\star,D^\star,\bm \alpha^\star)$ which minimizes the optimization problem \eqref{model},
\emph{i.e.}
\[
(u^\star,D^\star,\bm \alpha^\star)=\arg\min\limits_{u,D,\bm \alpha} \Upsilon (u,D,\bm \alpha).
\]
\end{thm}

Due the the existence of the data fitting term for phase retrieval and $L_0$ regularization of the coefficient matrix, the proposed model is a nonconvex and discontinuous optimization problem, which is very challengeable for designing an algorithm with theoretical convergence guarantee. In the following sections, we first give a simple Alternating Minimization Method (AMM), and only local convergence can be obtained. Then a global convergent algorithm called Proximal Alternating Linearized Minimization method (PALM) will be provided.

\section{Algorithms I:  Alternating Minimization Method (AMM)}\label{sec2}
Since three variables $u,D,\bm\alpha$ of \eqref{model} are coupled, in a natural way we use AMM to solve the problem. If giving the approximation solution $(u^k,D^k,\bm\alpha^k)$, the overall iterative algorithm consists of three steps  \emph{w.r.t.}  the three coupled variables as
\begin{equation}\label{ADM}
\left\{
\begin{aligned}
&u^{k+1}=\arg\min\limits_{u} \Upsilon(u,D^k,\bm \alpha^k),\\
&D^{k+1}=\arg\min\limits_{D} \Upsilon(u^{k+1},D,\bm \alpha^k),\\
&\bm \alpha^{k+1}=\arg\min\limits_{\bm \alpha} \Upsilon(u^{k+1},D^{k+1},\bm \alpha).
\end{aligned}
\right.
\end{equation}
Readily one can know that there exist at least a solution to each subproblem. In the following subsections, we will present how to solve these subproblems. For simplicity we omit all the superscripts of the notations in \eqref{ADM}.
\subsection{Subproblem  \emph{w.r.t.}   variable $u$}
We  concentrate on the first subproblem  \emph{w.r.t.}  variable $u$ in this subsection. By introducing an auxiliary variable $z$,
we have
\begin{equation}\label{sub1}
\begin{split}
\min\limits_{u} \frac12\|Y-R(u)\|^2+\eta\mathcal B(|z|^2,f), ~s.t. ~z=\mathcal A u,
\end{split}
\end{equation}
with $Y=[Y(0),Y(1),\cdots,Y(d-1)]:=D^k{\bm\alpha^k}$.
An equivalent form of \eqref{sub1} is derived as
\begin{equation}\label{sub11}
\begin{split}
\min\limits_{u} \frac12\sum\limits_{t=0}^{d-1}\|Y(t)-R_tu\|^2+\eta\mathcal B(|z|^2,f),
s.t.  ~z=\mathcal A u.
\end{split}
\end{equation}
An alternating direction of multiplier method (ADMM) is adopted to solve the above optimization problem similar to \cite{chang2016phase}. The corresponding augmented Lagrangian reads
\begin{equation}\label{lagrangian1}
\begin{split}
&\max\limits_{\Lambda}\min\limits_{u,z}  \mathcal L_r(u,z;\Lambda):=\frac12\sum\limits_{t=0}^{d-1}\|Y_t-R_tu\|^2\\
&~~+\eta\mathcal B(|z|^2,f)+\Re(\langle z-\mathcal A u,\Lambda \rangle)+\frac{r}{2}\|z-\mathcal A u\|^2,
\end{split}
\end{equation}
where $\Re(\cdot)$ denotes the real part of a complex-valued number, $r$ is a positive parameter, $\langle\cdot\rangle$ denotes the inner product of two vectors. The ADMM is designed as below,
\begin{equation}\label{ADMMU}
\left\{
\begin{aligned}
&u_{j+1}=\arg\min \mathcal L_r(u,z_j;\Lambda_j),\\
&z_{j+1}=\arg\min \mathcal L_r(u_{j+1},z;\Lambda_j),\\
&\Lambda_{j+1}=\Lambda_j+r(z_{j+1}-\mathcal Au_{j+1}),
\end{aligned}
\right.
\end{equation}
starting from the previous iterative solution $(u_j,z_j,\Lambda_j)$.
Following the framework in \cite{chang2016Total}, one can readily obtain that
the closed form of the first subproblems of \eqref{ADMMU} as
	\begin{equation}\label{subsubU}
\begin{split}
	&\left[
	\begin{matrix}
	r{\Re(\mathcal A^* \mathcal A)}+W& -{\Im(\mathcal A^* \mathcal A)}\\
	&\\
	{\Im(\mathcal A^* \mathcal A)}& r{\Re(\mathcal A^* \mathcal A)}+W
	\end{matrix}
	\right]
	\left[\begin{matrix}
	\Re(u_{j+1})\\
	\\
	\Im(u_{j+1})
\end{matrix}\right]\\
&\quad=
	\left[
	\begin{matrix}
	r\Re(\mathcal A^* v_j)+\sum\limits_t R_t^T \Re(Y(t))\\
\\
	r\Im(\mathcal A^* v_j)+\sum\limits_k R_t^T \Im(Y(t))
	\end{matrix}
	\right],
\end{split}
	\end{equation}
	with $v_j=z_j+\Lambda_j/r,$ and $W=\sum\limits_{t=0}^{d-1} R_t^T R_t$.
We can simplify the solution of above subproblem if the matrix $\mathcal A$ involves
Fourier measurements with masks $\{I_k\}_{k=0}^{K-1}$  as
	\begin{equation}\label{cdp}
	\mathcal Au=
	\left[
	\begin{matrix}
	\mathcal F(I_0\circ u)\\
	\mathcal F(I_1\circ u)\\
	\vdots\\
	\mathcal F(I_{K-1}\circ u)
	\end{matrix}
	\right],
	\end{equation}
where $\circ$ denotes the pointwise multiplication, $I_k$ is a masked matrix indexed by $k$, each of which is represented by a vector in $\mathbb C^{n}$ in a lexicographical order.
	Therefore, we have
\begin{equation}\label{AA}
\mathcal A^* \mathcal A=\sum\limits_j I_j^*\circ I_j=\sum\limits_j |I_j|^2,
\end{equation}
which is a real-valued matrix. Finally we derive the solution
\begin{equation}\label{subsubUSimple}
u_{j+1}=\big(r{\mathcal A^* \mathcal A}+W\big)^{-1}\big(r\mathcal A^*v_j+\sum\limits_t R_t^T Y(t)\big),
\end{equation}
if the diagonal matrix ${\mathcal A^* \mathcal A}+W$ is non-singular.

For the second subproblem of \eqref{ADMMU}, we have
\begin{equation}
\begin{split}
\min\limits_{z}  \eta\mathcal B(|z|^2,f)+\frac{r}{2}\|z-\mathcal A u_{j+1}+\Lambda_j/r\|^2.
\end{split}
\end{equation}
With Poisson noised data $f$,
the close form solution  is readily obtained as
\begin{equation}\label{PoisSolver}
z_{j+1}(t)=\dfrac{r|w(t)|+\sqrt{r^2|w(t)|^2+4\eta(\eta+r)f(t)}}{2(\eta+r)}\mathrm{sign}(w(t)),
\end{equation}
by letting $w_j=\mathcal A u_{j+1}-\Lambda_j/r,$
where
$$
\mathrm{sign}(w(t))=\dfrac{w(t)}{|w(t)|}.
$$
In summary, the overall algorithm of the subproblem  \emph{w.r.t.}  $u$ of \eqref{ADM} is listed as follows:
\vskip .1in

\begin{mdframed}[hidealllines=true,backgroundcolor=gray!25,innerleftmargin=1pt,innerrightmargin=1pt,leftmargin=-5pt,rightmargin=-5pt]
\begin{center}
\begin{minipage}{.95\linewidth}
{\hskip .4in \hrule \vskip .05in \hrule\vskip .1in
\centering Algorithm I-I: ADMM for $u-$subproblem of \eqref{ADM}   \vskip .1in
\hrule\vskip .05in}
\begin{enumerate}
\item[1.] Initialization: Set $u_0$, $z_0=\mathcal A u_0, \Lambda_0=0,$ and $j=0,$ and parameter $r$.
\item[2.] Solve $u_{j+1}$ by \eqref{subsubU}. Furthermore, if $\mathcal A$ is generated by \eqref{cdp}, solve it by \eqref{subsubUSimple}.
\item[3.]  Solve $z_{j+1}$  by \eqref{PoisSolver}.
\item [4.] Update multipliers as
$\Lambda_{j+1}=\Lambda_j+r(z_{j+1}-\mathcal A u_{j+1}).$
\item[5.] If  some stopping condition is satisfied, stop the iterations and output the iterative solution $u_{j+1}$; else set $j=j+1$, and goto Step 2.
\end{enumerate}
\hrule\vskip .05in \hrule
\end{minipage}
\end{center}
\end{mdframed}
\vskip .1in
\begin{rem}
One can also use a gradient descent type algorithm with adaptive steps as Wirtinger flow \cite{candes2015phase}, Newton type algorithm \cite{qian2014efficient} to solve it if the number of measurements are sufficient. It seems that ADMM can work pretty well with very few measurements in \cite{chang2016Total}, and  has low computation cost compared to Newton type algorithm.
\end{rem}

\subsection{Subproblem  \emph{w.r.t.}    variables $D$ and $\bm \alpha$}
For the second subproblem  of \eqref{ADM} \emph{w.r.t.} $D$, we have
\begin{equation}\label{subD}
\hat D:=\arg\min\limits_{D} \frac12\|D{\bm\alpha}-R(u)\|^2,~~s.t.~D^* D=\mathbf I.
\end{equation}
The minimizer $\hat D$ has a close form \cite{bao2013fast} as
\begin{equation}\label{solverD}
\hat D=UV^*,
\end{equation}
where $R(u){\bm\alpha}^*=U\Sigma V^*$, with corresponding  eigenvectors $U,V$ and singular values $\Sigma$  by singular value decomposition (SVD).

For the third subproblem of \eqref{ADM} \emph{w.r.t} $\bm \alpha$, we have
\begin{equation*}
\hat{\bm \alpha}:=\arg\min\limits_{\bm \alpha} \frac12\|{\bm\alpha}-D^*R(u)\|^2+\tau \|\bm \alpha\|_0,
\end{equation*}
since $D^*D=\mathbf I.$
For the isotropic $L^0$ norm, it  has a close form solution  known as ``hard thresholding'' as below,
\begin{equation}\label{hard}
{\hat{\bm \alpha}}=\mathrm{Thresh}_\tau(D^*R(u)),
\end{equation}
where the hard thresholding $\mathrm{Thresh}_\tau({{\bm \alpha}})$ is defined as
\begin{equation}\label{hardDef}
\mathrm{Thresh}_\tau({{\bm \alpha}})(s,t)=\left\{
\begin{aligned}
 & \bm\alpha(s,t), \text{if~}|(D^*R(u))(s,t)|\geq \tau,\\
 & 0, \qquad\qquad\qquad\text{otherwise},
\end{aligned}
\right.
\end{equation}
for $0\leq s\leq c-1,0\leq t\leq d-1$.
For the anisotropic case, similarly, one readily obtains the close form solution by separating the real and imaginary parts
\begin{equation}\label{hardAniso}
{\hat{\bm \alpha}}=\mathrm{Thresh}_\tau(\Re(D^*R(u)))+\mathbf{i} \times\mathrm{Thresh}_\tau(\Im(D^*R(u))).
\end{equation}
Therefore, we can ready to list an overall algorithm for our proposed model \eqref{model} as follows:
\begin{mdframed}[hidealllines=true,backgroundcolor=gray!25,innerleftmargin=1pt,innerrightmargin=1pt,leftmargin=-5pt,rightmargin=-5pt]
\begin{center}
\begin{minipage}{.95\linewidth}
{\hskip .4in \hrule \vskip .05in \hrule\vskip .1in
\centering Algorithm I: AMM for ``DicPR'' \eqref{model}   \vskip .1in
\hrule\vskip .05in}
\begin{enumerate}
\item[1.] Initialization:  Initialize $D^0$ by a discrete cosine dictionary, $\bm\alpha^0=0$ and $k=0,$ and parameter $\tau,\eta$.
\item[2.] Solve $u^{k+1}$ by Algorithm I-I with $Y:=D^k\bm\alpha^k$.
\item[3.]  Solve $D^{k+1}$  by \eqref{solverD} where $u:=u^{k+1}$ and $\bm \alpha:=\bm \alpha^k$.
\item [4.] Solve $\bm \alpha^{k+1}$ by \eqref{hard} or \eqref{hardAniso} where $u:=u^{k+1}$ and $D:=D^{k+1}$.
\item[5.] If  some stopping condition is satisfied, stop the iterations and output the iterative solution; else set $k=k+1$, and goto Step 2.
\end{enumerate}
\hrule\vskip .05in \hrule
\end{minipage}
\end{center}
\end{mdframed}
\vskip .1in

\begin{rem}
In order to accelerate the computation, in the first few iterations (within three iterations in the numerical experiments),  we use a truncated version for \eqref{solverD} where only parts (one half in our numerical examples) of the all columns of $U$ and $V$ are selected corresponding to the larger singular values.
 \end{rem}


One can readily  conclude to the decrease of the functional values of iterative sequences, and we list it as below without proof, which can be directly obtained by \eqref{ADM} following Algorithm I.
\begin{lemma}\label{lemDec}
The objective functional values of the iterative sequences generated by Algorithm I  are non-increasing, i.e.
\[
 \Upsilon (u^{k+1},D^{k+1},\bm \alpha^{k+1})\leq  \Upsilon (u^k,D^k,\bm \alpha^k).
\]
\end{lemma}
We give an assumption of  operator $\mathcal A$ as follows.
\begin{assump}\label{assump1}
The sequence $\{u^k\}_{k=0}^\infty$ is bounded if and only if $\{|\mathcal A u|\}_{k=0}^\infty$ is bounded.
\end{assump}
\begin{rem}
In order to guarantee the uniqueness of the solution for phase retrieval, the oversampling or multiple measurements are required. As a result, it make the above assumption available. For example, for the coded diffraction patter and ptychographic pattern, $\mathcal A^*\mathcal A $ is a diagonal matrix as \eqref{AA}. If we assume that the diagonal elements are all non-zeros, Assumption \ref{assump1}  holds. In practise, the matrix  $\mathcal A^*\mathcal A $ is usually invertible.
\end{rem}

\begin{lemma}
The iterative sequence $Z^k:=(u^k,D^k,\bm \alpha^k)$ is bounded under Assumption \ref{assump1}.
\end{lemma}
\begin{IEEEproof}
By Lemma \ref{lemDec}, the sequence $\Upsilon(Z^k)$ is bounded, \emph{i.e.} there exists a positive constant $C$ independent with $k$, such that
$\mathcal H(u^k,D^k,\bm\alpha^k)+\mathcal F(u^k)
+\mathcal I_{\mathbf K}(D^k)+\tau \|\bm \alpha^k\|_0\leq C.$
First we have $\mathcal F(u^k)\leq C$ and readily one has $\{|\mathcal A u^k|\}$ is bounded. By Assumption \ref{assump1}, the boundedness of $\{u^k\}$ is proved.
Since $(D^k)^*D^k=\mathbf I $, it  implies the boundedness of $\{D^k\}.$
Since $\mathcal H(u^k,D^k,\bm\alpha^k)\leq C,$ we have $\|\bm\alpha^k-(D^k)^*R(u^k)\|\leq C$ by the orthogonality of the dictionary, which leads to the boundedness of $\{\bm\alpha^k\}$.
\end{IEEEproof}

\begin{thm}
There exists an accumulative point $\tilde Z$ of $\{Z^k\}$ such that
\[
\lim\limits_{k\rightarrow \infty} \Upsilon (Z^k)=\Upsilon (\tilde Z).
\]
\end{thm}
\begin{IEEEproof}
Since $\{Z^k\}$ is bounded, there exists a subsequence $\{Z^{k_n}\}^\infty_{n=0}\subset\{Z^k\}$ with the limit point $\tilde Z=(\tilde u,\tilde D,\tilde{\bm \alpha}),$
 s.t.
 $
 \lim\limits_{n\rightarrow \infty} Z^{k_n}=\tilde Z.
 $
 Therefore by Lemma \ref{lemDec}, we conclude to the theorem.
\end{IEEEproof}
Here we only derive the local convergence, and more future work should be done in order to investigate the convergence of global minimizer.

\section{Algorithm II: Proximal alternating linearized minimization (PALM)}\label{sec3}
In order to guarantee the global convergence, we propose a proximal alternating linearized minimization method (PALM) to solve the established dictionary learning model \eqref{model} following the multi-block splitting algorithm proposed in \cite{bolte2014proximal}. First we give the derivative of $\mathcal H(u,D,\bm\alpha)$ as below:
\begin{equation}\label{Lip}
\left\{
\begin{aligned}
&\nabla_u \mathcal H(u,D,\bm\alpha)=\sum\limits_{t=0}^{d-1} R_t^T(R_tu-D\bm\alpha(t)),\\
&\nabla_D \mathcal H(u,D,\bm\alpha)=(D\bm\alpha-R(u))\bm\alpha^*,\\
&\nabla_{\bm\alpha}\mathcal H(u,D,\bm\alpha)=D^*(D\bm\alpha-R(u)),
\end{aligned}
\right.
\end{equation}
with $\bm\alpha=[\bm\alpha(0),\bm\alpha(1),\cdots,\bm\alpha(d-1)]$.
The PALM is given as
\begin{eqnarray}
&u^{k+1}&=\arg\min\limits_u \mathcal F(u)+\frac{c^k}{2}\|u-\hat u^k\|^2~\label{eq1-1}\\
&&\qquad\mbox{with~}\hat u^k=u^k-\frac{\nabla_u \mathcal H(u^k,D^k,\bm \alpha^k)}{c^k},\nonumber\\
&D^{k+1}&=\arg\min\limits_D \mathcal I_{\mathbf K}(D)+\frac{d^k}{2}\|D-\hat D^k\|^2~\label{eq1-2}\\
&&\qquad\mbox{with~}\hat D^k=D^k-\frac{\nabla_D \mathcal H(u^{k+1},D^k,\bm \alpha^k)}{d^k},\nonumber\\
&\bm\alpha^{k+1}&=\arg\min\limits_{\bm\alpha} \tau \|\bm \alpha\|_0+\frac{e^k}{2}\|{\bm \alpha}-\hat{\bm\alpha}^k\|^2~\label{eq1-3}\\
&&\qquad\mbox{with~}\hat{\bm\alpha}^k=\bm \alpha^k-\frac{\nabla_{\bm \alpha} \mathcal H(u^{k+1},D^{k+1},\bm \alpha^k)}{e^k},\nonumber
\end{eqnarray}
with three positive steps $c^k,d^k$ and $e^k$.

For $u-$subproblem in \eqref{eq1-1} for PALM, one can readily give the ADMM following Algorithm I-I. Only the solver for $u$ in Step 2 has slightly difference, and we directly list Algorithm II-I as follows.
\vskip .05in
\begin{mdframed}[hidealllines=true,backgroundcolor=gray!25,innerleftmargin=1pt,innerrightmargin=1pt,leftmargin=-5pt,rightmargin=-5pt]
\begin{center}
 \begin{minipage}{.95\linewidth}
{\hskip .2in \hrule \vskip .05in \hrule\vskip .1in
\centering Algorithm II-I: ADMM for $u-$subproblem of \eqref{eq1-1}   \vskip .1in
\hrule\vskip .05in}
\begin{enumerate}
\item[1.] Initialization: Set $u_0$, $z_0=\mathcal A u_0, \Lambda_0=0,$ and $j=0.$
\item[2.] Solve $u_{j+1}$ by
\begin{equation}\label{II-I-subsubUSimple}
\begin{split}
u_{j+1}=&\big(r{\mathcal A^* \mathcal A}+c^k\mathbf I\big)^{-1}\big(r\mathcal A^*v_j+\sum R_t^T Y(t)\\
&\qquad\qquad\qquad\qquad\quad+ (c^k\mathbf I-W)  u^k\big),
\end{split}
\end{equation}
with $v_j=z_j+\Lambda_j/r,$ and $Y=D^k\bm\alpha^k$. 
\item[3.]  Solve $z_{j+1}$  by \eqref{PoisSolver}.
\item [4.] Update multipliers as $
\Lambda_{j+1}=\Lambda_j+r(z_{j+1}-\mathcal A u_{j+1}).
$
\item[5.] If  some stopping condition is satisfied, stop the iterations and output the iterative solution $u_{j+1}$; else set $j=j+1$, and goto Step 2.
\end{enumerate}
\hrule\vskip .05in \hrule
\end{minipage}
\end{center}
\end{mdframed}
For the $D-$subproblem in \eqref{eq1-2} for PALM, we need to solve the following problem as
\begin{equation}\label{eqD}
D^{k+1}=\arg\min\limits_D \|D-\hat D^k\|^2,~s.t. ~D^*D=\mathbf I.
\end{equation}
One can readily derive the closed form for it as
\begin{equation}\label{solverIID}
D^{k+1}=UV^*,
\end{equation}
with
\[
\begin{split}
\hat D^k&=D^k-\frac{1}{d^k}(D^k\bm\alpha^k-R(u^{k+1}))(\bm\alpha^k)^* \\
&=D^k\big(\mathbf I-\frac{1}{d^k}\bm\alpha^k(\bm\alpha^k)^*\big)+\frac{1}{d^k}R(u^{k+1})(\bm\alpha^k)^*\\
&:=U\Lambda V^*
\end{split}
\] as its SVD.

For the $\bm\alpha-$subproblem in \eqref{eq1-3}, one can directly  get the closed form  represented by the  hard thresholding as
\begin{equation}\label{solverIIalpha}
\bm\alpha^{k+1}=\mathrm{Thresh}_{~\tau/e^k}\big( \hat{\bm\alpha}^k\big),
\end{equation}
and
\begin{equation}\label{solverIIalphaAniso}
\bm\alpha^{k+1}=\mathrm{Thresh}_{~\tau/e^k}( \Re(\hat{\bm\alpha}^k))+\mathbf i\times \mathrm{Thresh}_{~\tau/e^k}( \Im(\hat{\bm\alpha}^k)),
\end{equation}
with $\hat{\bm\alpha}^k=(1-\frac{1}{e^k})\bm\alpha^k+\frac{1}{e^k}(D^{k+1})^*R(u^{k+1})$
for the isotropic and anisotropic $L^0$ pseudo norms respectively.

Therefore we can give an overall PALM for our proposed model \eqref{model} as below.
\vskip .05in
\begin{mdframed}[hidealllines=true,backgroundcolor=gray!25,innerleftmargin=1pt,innerrightmargin=1pt,leftmargin=-5pt,rightmargin=-5pt]
\begin{center}
 \begin{minipage}{.95\linewidth}
{\hskip .4in \hrule \vskip .05in \hrule\vskip .1in
\centering Algorithm II: PALM for ``DicPR'' \eqref{model}   \vskip .1in
\hrule\vskip .05in}
\begin{enumerate}
\item[1.] Initialization: Initialize $D^0$ by a discrete cosine dictionary, $\bm\alpha^0=0$ and $k=0,$  parameter $\tau,\eta$ and $(c^k,d^k,e^k)$.
\item[2.] Solve $u^{k+1}$ by Algorithm II-I with $Y:=D^k\bm\alpha^k$.
\item[3.]  Solve $D^{k+1}$  by \eqref{solverIID}.
\item [4.] Solve $\bm \alpha^{k+1}$ by \eqref{solverIIalpha} or \eqref{solverIIalphaAniso}.
\item[5.] If  some stopping condition is satisfied, stop the iterations and output the iterative solution; else set $k=k+1$, and goto Step 2.
\end{enumerate}
\vskip .05in\hrule\vskip .05in \hrule
\end{minipage}
\end{center}
\end{mdframed}
In the following parts, we will provide convergence analysis for PALM. For simplicity, the analysis is conducted for real-valued images, and it can readily  generalized to the complex-valued images, where one only needs to rebuild the related norms and operators of variables on their corresponding real and imaginary parts respectively as in \cite{chang2016Total}.

The existence of the critical point of proposed model is given below.
\begin{lemma}
The critical point of \eqref{sub11} exists with the finite functional value.
\end{lemma}
\begin{IEEEproof}
Following Lemma 5 in \cite{bolte2014proximal} and the lower semi-continuity of $\mathcal F, \mathcal I_{D}$ and $\|\cdot\|_0$, one can readily prove it, and details are omitted.
\end{IEEEproof}

We need the assumption of the boundedness of $u^k$ as follows.
\begin{assump}[Boundedness]\label{assump2}
The sequences $\{u^k\}$ generated by \eqref{eq1-1} is bounded.
\end{assump}
Based on above assumption, the boundedness of the iterative sequences $(u^k,D^k,\bm\alpha^k)$ can be proved, and see details in Lemma \ref{lem2} in the appendix.

\begin{lemma}[Lipschitz-Gradient]\label{lem1}
For the functional $\mathcal H(u,D,\bm \alpha)$, we have
 \begin{itemize}
\item[1.] $\mathcal H$ is differential, and $\inf F>-\infty$.
\item[2.] $\nabla_u\mathcal H,$$\nabla_D \mathcal H$ and $\nabla_{\bm\alpha}\mathcal H$ are Lipschitz continuous respectively with moduli $L_u(D,\bm\alpha),$ $L_D(u,\bm\alpha),$ $L_\alpha(u,D)$, and the Lipschitz constants are bounded, \emph{i.e.}
 there exist three positive constants $\lambda_u^+,$$\lambda_D^+,$$\lambda_{\bm\alpha}^+$, such that
\[
\begin{split}
&\sup L_u(D^k,{\bm\alpha}^k)\leq \lambda_u^+,\\
&\sup L_D(u^k,\bm\alpha^k)\leq \lambda_D^+,\\
&\sup L_{\bm\alpha}(u^k,D^k)\leq \lambda_{\bm\alpha}^+.
\end{split}
\]

\item[3.] The gradient $\nabla \mathcal H(u,D,\bm\alpha) $ is Lipschitz continuous with Lipschitz constant $M$ on a bounded domain $\{(u,D,\bm\alpha):~\|(u,D,\bm\alpha)\|\leq \text{Const}\}.$
\item[4.] $\nabla_D \mathcal H(u,D,\bm\alpha)$ is Lipschitz continuous with Lipschitz constant $\lambda_D^-$  \emph{w.r.t.}  $\bm \alpha$ for  the bounded  sequences.
\end{itemize}
\end{lemma}
\begin{IEEEproof}
The proof for the first item is trivial. By  \eqref{Lip} and Assumption \ref{assump2}, one can readily finish the proof for the left parts based on the boundedness of iterative sequences proved by Lemma \ref{lem2} in the Appendix.
\end{IEEEproof}
For the PALM, we assume the step sizes should be large enough and  the following assumption is needed.
\begin{assump}[Step Sizes]\label{ass2}
\[
\lambda^+:=\min\{c^k-\lambda_u^+,d^k-\lambda_D^+,e^k-\lambda_{\bm\alpha}^+\}>0.
\]
\end{assump}
Then we can consider the global convergence of PALM.
Readily one knows that $\Upsilon(u,D,\bm\alpha)$ is semi-algebra \cite{bolte2014proximal,bao2016dictionary}  with the data term $\mathcal B(\cdot,\cdot)$ in \eqref{data}, and finally we can get the final convergence theorem as follows.
\begin{thm}
 Let Assumption \ref{assump2}, \ref{ass2} hold and $e^k>1/2$. The sequence $\{u^k,D^k,\bm\alpha^k\}$ generated by  Algorithm II globally converges to a critical point of the proposed model \eqref{modelR}.
\end{thm}
\begin{IEEEproof}
Based on Lemma \ref{energyDecay} and Lemma \ref{boundedGrad} in the Appendix, one can finish the proof following Theorem 1 in \cite{bolte2014proximal}.
\end{IEEEproof}
\begin{rem}
Although in this paper, we focus on the orthogonal dictionary, the proposed PALM can also be applied to  the case with non-orthogonal dictionary, and further speedup of PALM should be investigated as in \cite{bao2016dictionary}. We leave them as  future work.
\end{rem}

At the end of this section, we analyze  the computational complexity of Algorithm I and Algorithm II, which have similar cost. Assume that the operator $\mathcal A$ generates  Fourier masked measurements which means fast Fourier transformation can be adopted,  computation complexity  for Algorithm I-I and Algorithm II-I is of $O(n\log(n)T_{in}+ln )$ after $T_{in}$ inner iterations. The complexity for Step 3 and 4 is $O(l^3+l^2d)$. Hence the complexity of Algorithm I and Algorithm II is $O\big((n\log(n)T_{in}+ln +l^3+l^2d)T_{out}\big)$ after $T_{out}$ outer iterations. In our experiments, we set $d\approx n, l\ll n$ and as a result the complexity is
about $O\big(n(\log(n)T_{in}+l^2)T_{out}\big)$.

\section{Numerical experiments}\label{sec5}

All the tests are performed on a laptop with Intel I7-5600U2.6GHZ, and 16GB RAM, and  the codes are implemented in MATLAB. The dictionary is initialized by  discrete cosine transform.  The model parameters $\eta$ and $\tau$ for Algorithm I and Algorithm II are selected by hand. In the inner iterative algorithm e.g. Algorithm I-I and Algorithm II-I, we set defaulted parameter $r=1\times 10^{-3}$, and the defaulted inner iteration number  to five. For the left parameters  for Algorithm II,  one can use dynamic schemes to update the parameters as \cite{bolte2014proximal,bao2016dictionary}. In our experiments, for simplicity, these parameters are set to  be fixed.  We stop  Algorithm I and Algorithm II after a given maximum outer iteration number $T$ to guarantee the convergence, which will be specified in the following subsections. All the other needed parameters  will also be addressed in the following subsections if we do not give or use the defaulted values.

Set image patch size  to $8\times 8$  empirically, and hence $l=c=64.$ If  the patch sizes are too large, the computational cost increases dramatically based on the given complexity analysis. The number of patches  $d=(\sqrt{n}-7)^2$ for square images with $n$ pixels. 
In this paper we have given a  framework for  phase retrieval with arbitrary linear operator $\mathcal A$. However, it it more practical to consider  Fourier  type transformation, and we will show the performance  on Fourier masked measurements involved with two  types of patterns for linear operators $\mathcal A$: Coded diffraction pattern (CDP) with random masks  and ptychographic phase retrieval with deterministic masks generated by zone plate lens. Especially, ptychographic phase retrieval is  a very promising  technique to generate high resolution images with large field of view compared with the traditional diffraction imaging and meanwhile it requires less temporal and spatial coherence, while the related phase retrieval problem is more challengeable.

\begin{figure}
\vskip -.2in
\begin{center}
\subfigure[]{\includegraphics[width=.06\textwidth]{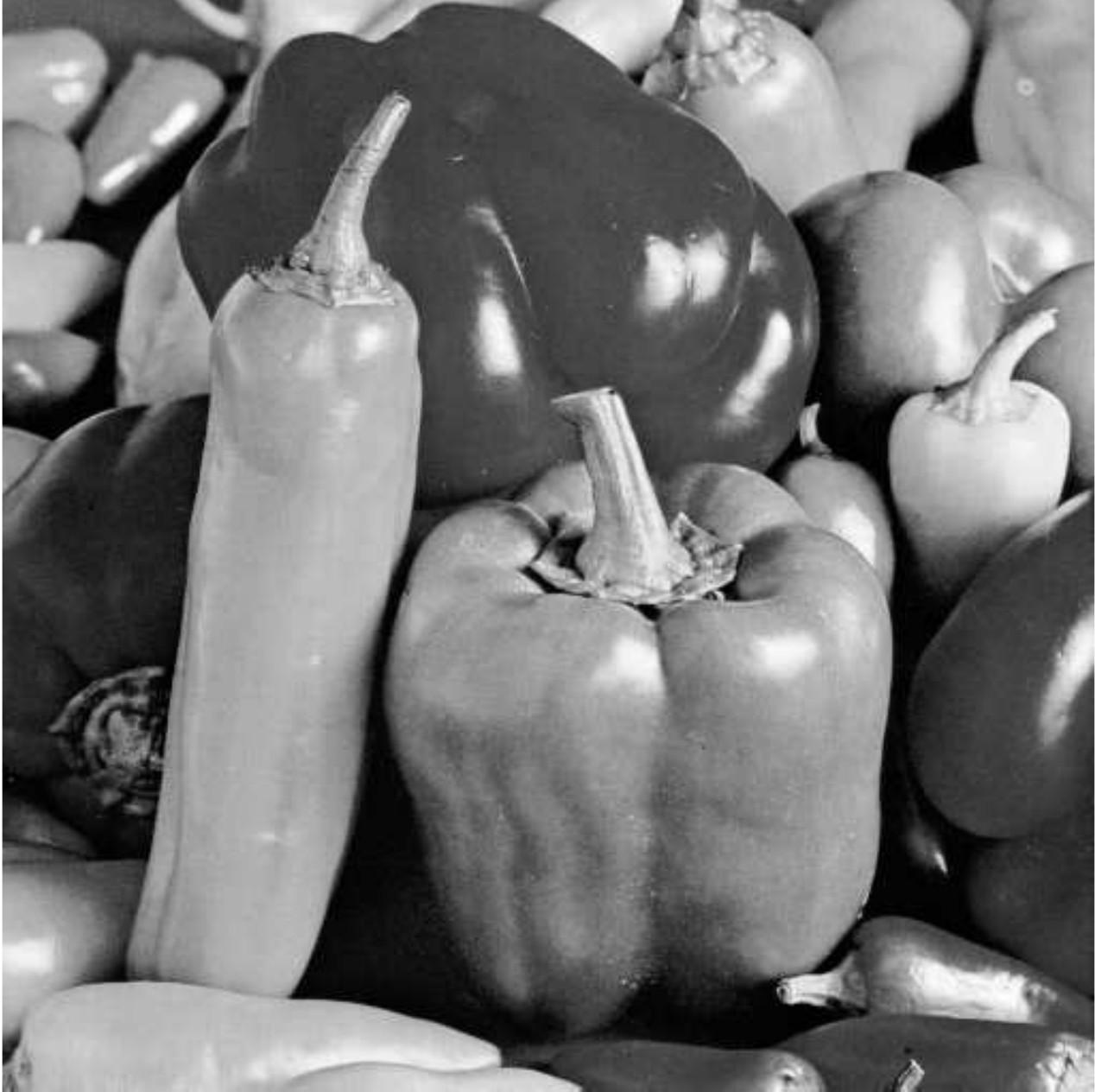}}
\subfigure[]{\includegraphics[width=.06\textwidth]{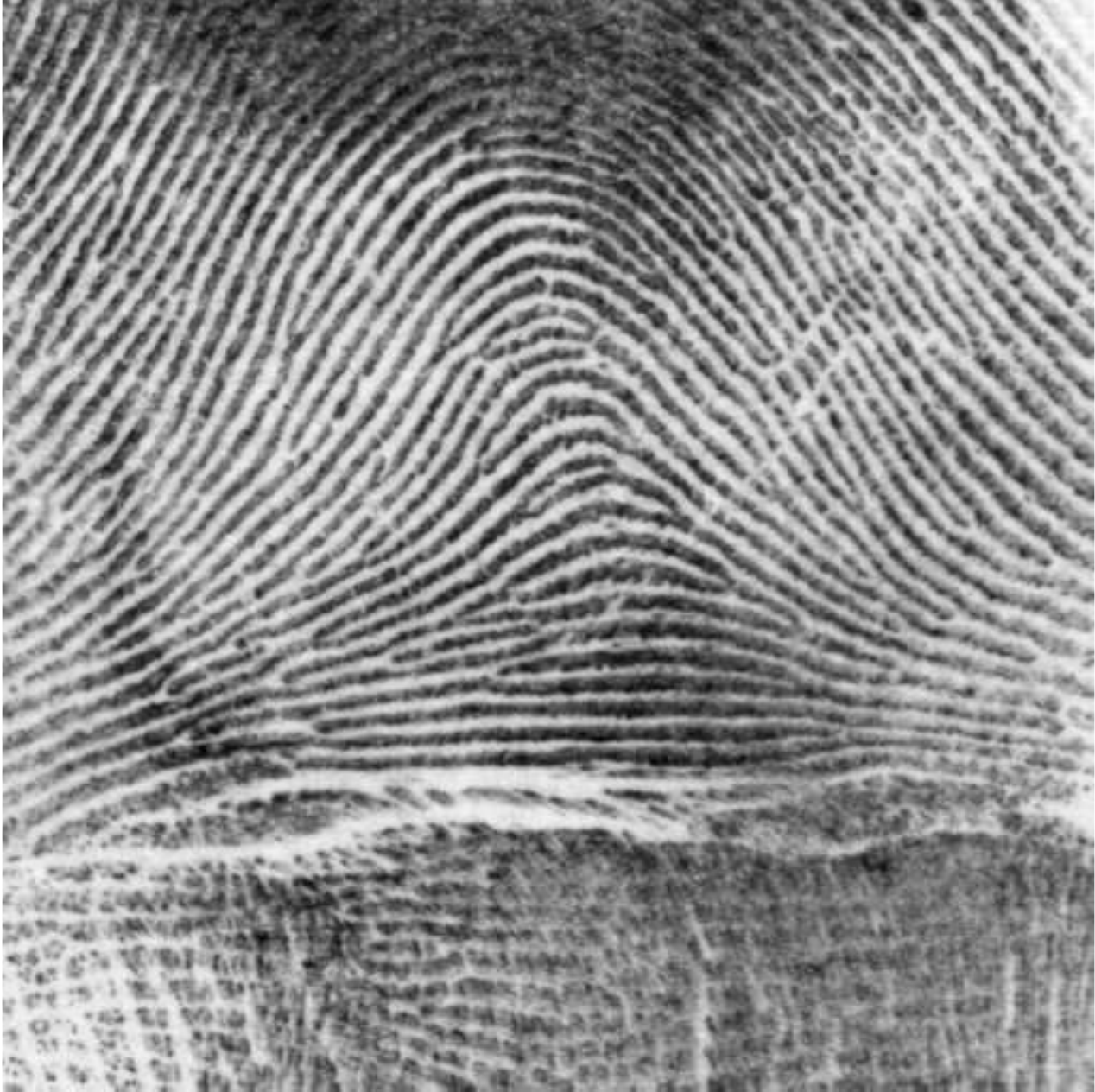}}
\subfigure[]{\includegraphics[width=.06\textwidth]{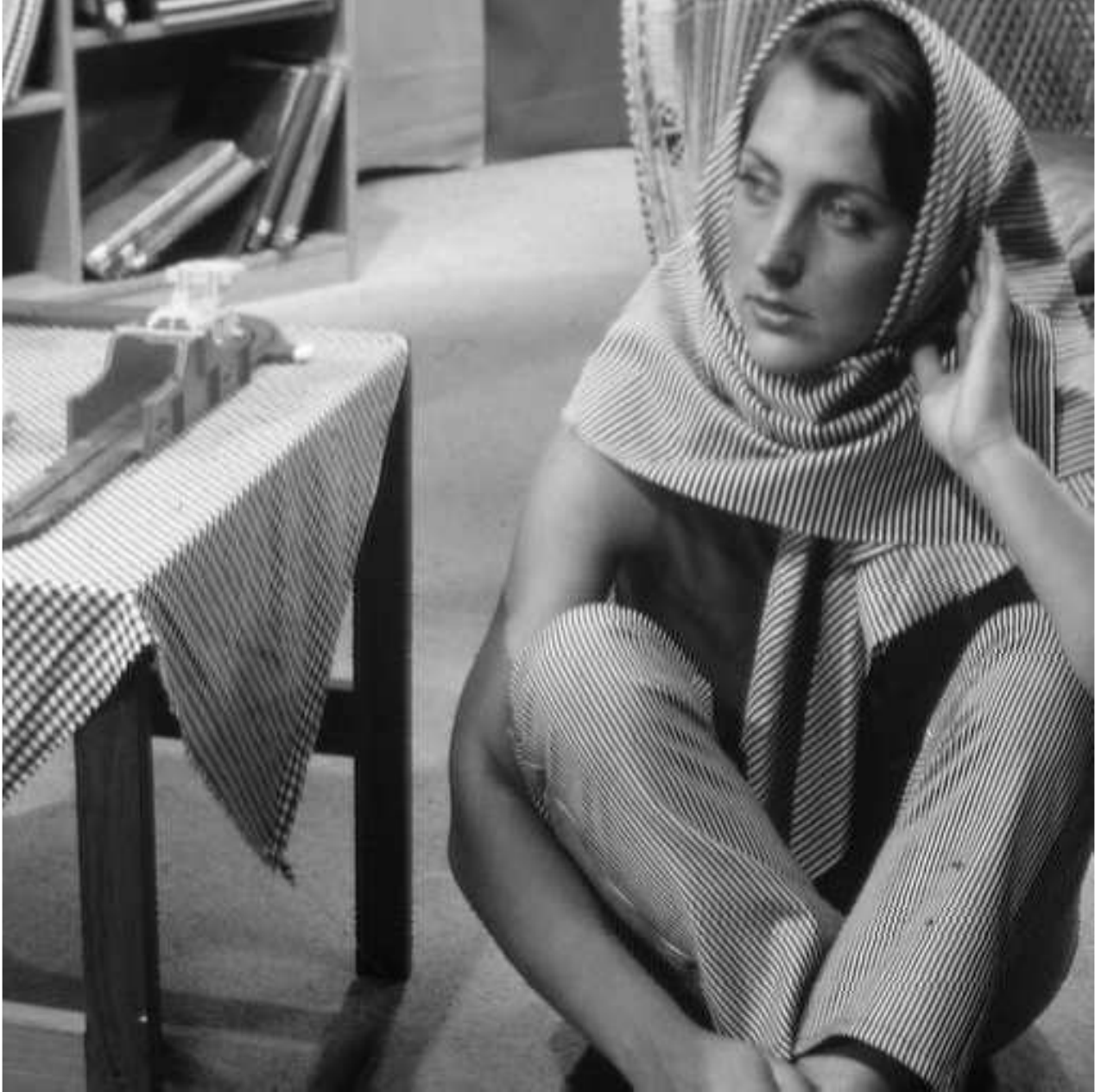}}
\subfigure[]{\includegraphics[width=.06\textwidth]{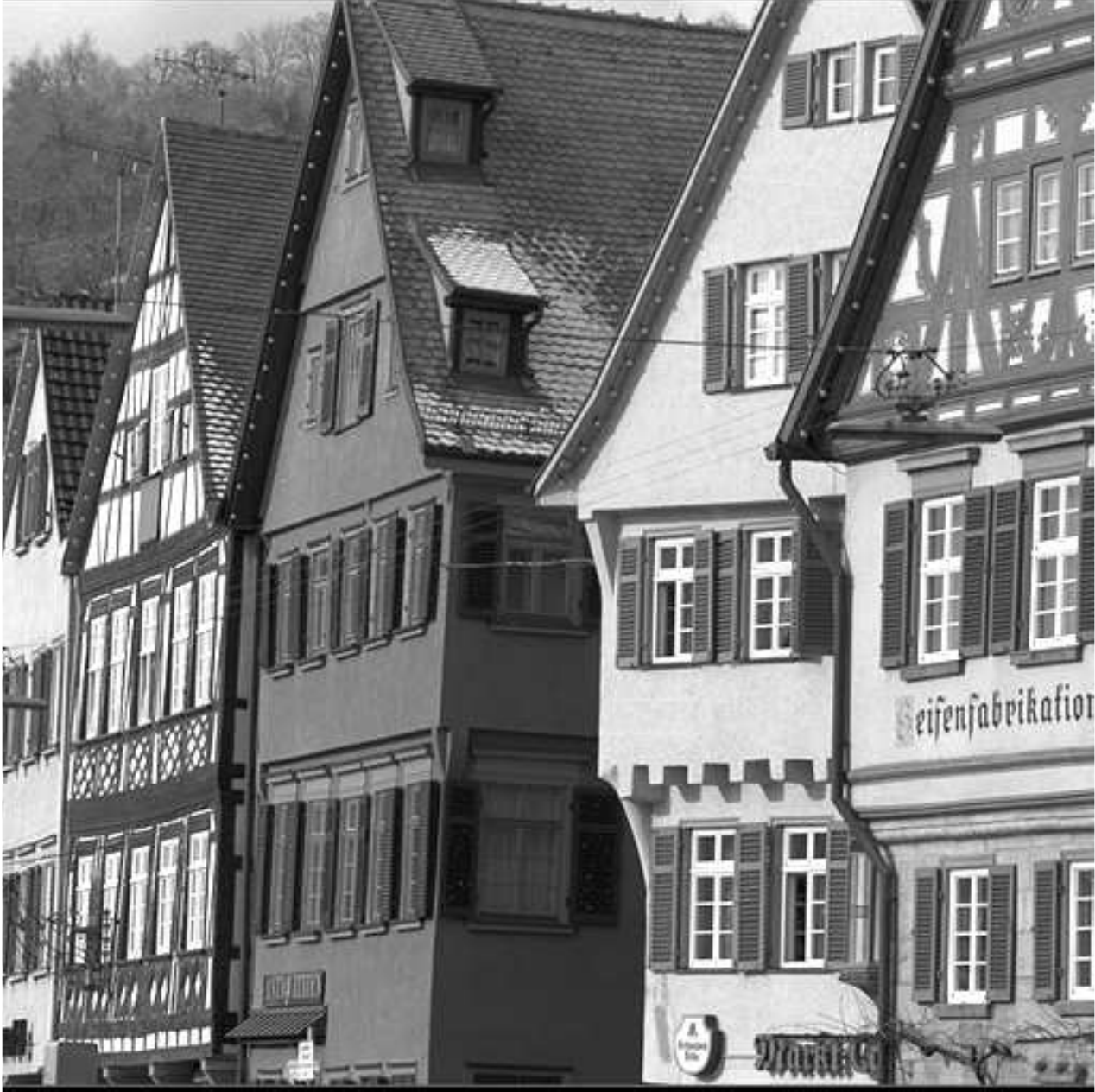}}
\subfigure[]{\includegraphics[width=.06\textwidth]{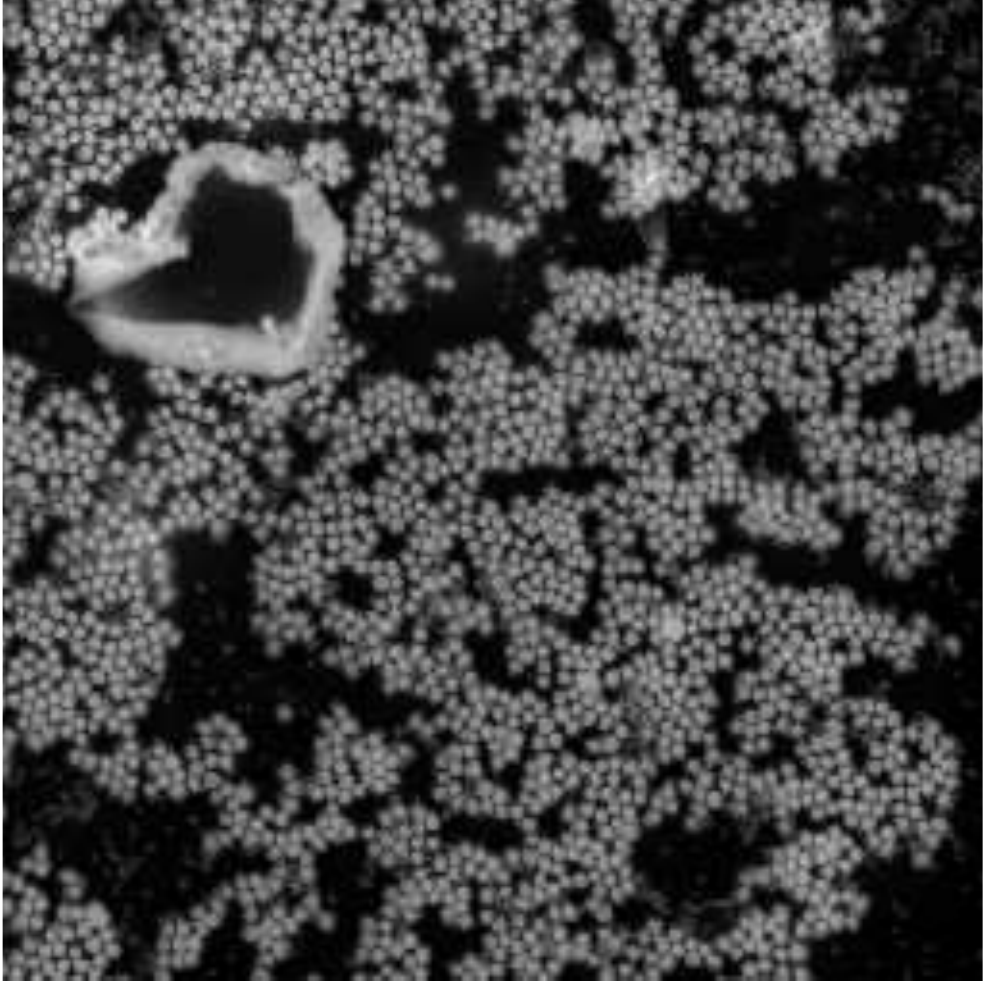}}
\subfigure[]{\includegraphics[width=.06\textwidth]{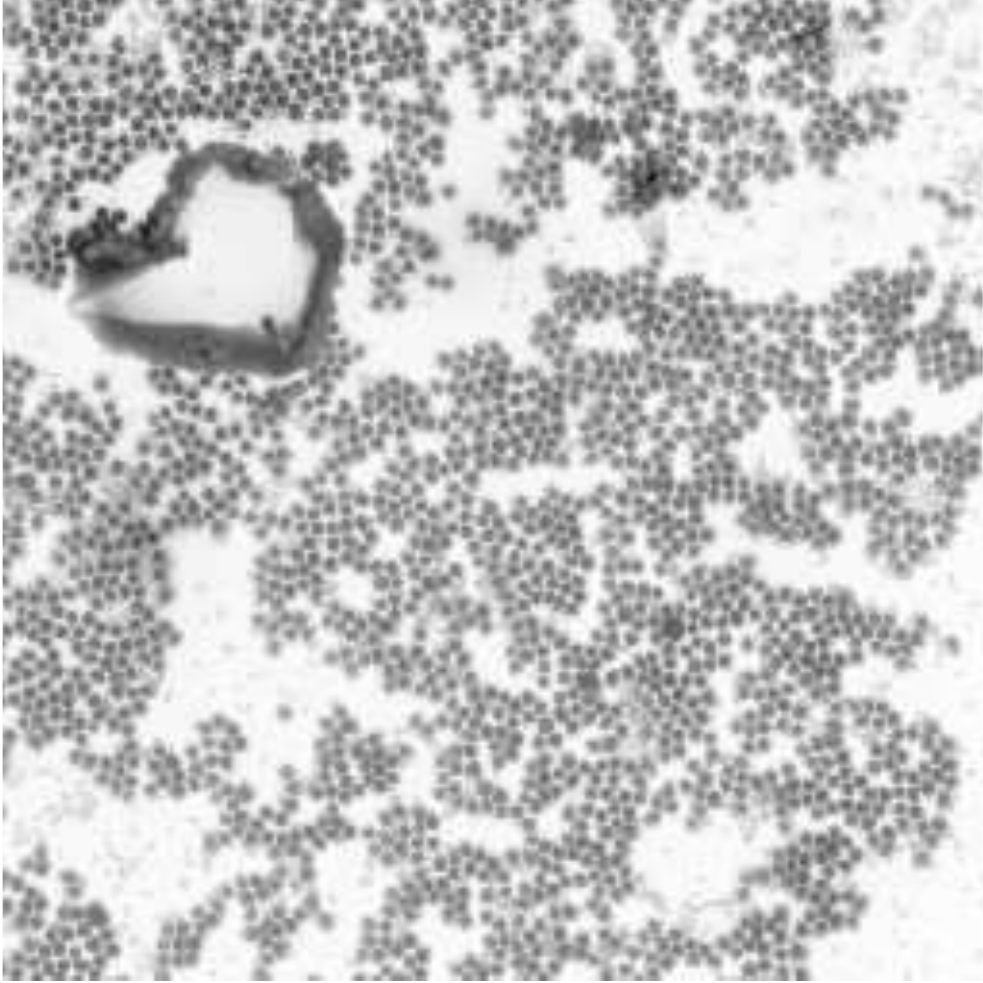}}
\subfigure[]{\includegraphics[width=.06\textwidth]{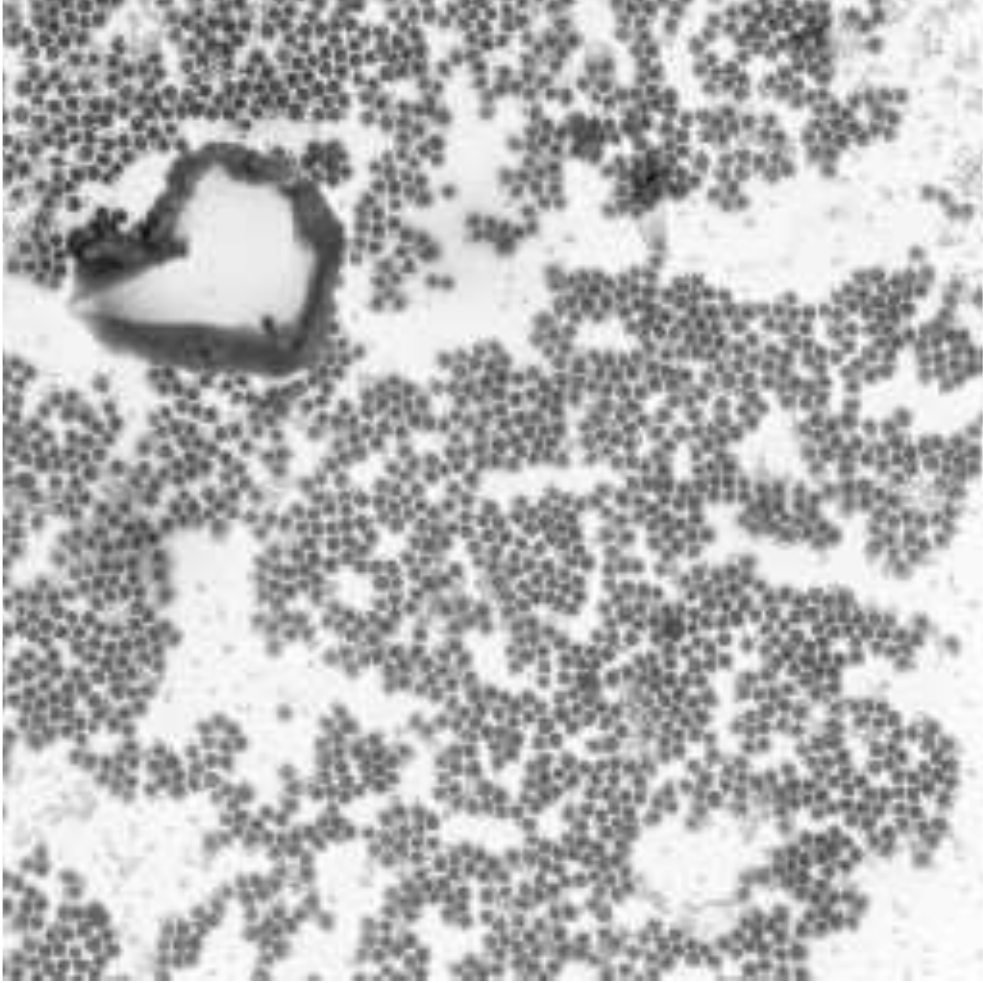}}
\end{center}
\vskip -.1in
\caption{Ground truth images.  First row: Real-valued images of resolution $512\times 512$ in (a): Peppers; (b): Fingerprint; (c): Barbara; (d): House;
Second row: A complex-valued image ``Goldballs'' of resolution $256\times 256$ with magnitude, real and imaginary parts in (e), (f) and (g), respectively.}
\label{ground1}
\vskip -.2in
\end{figure}

 The ground truth images are provided in Fig. \ref{ground1}, where  four real-valued images with resolution $512\times 512$ are put in Fig. \ref{ground1}(a)-(d), and a complex-valued image with resolution $256\times 256$ is put in Fig. \ref{ground1}(e)-(g). Given a ground truth  $u$, the noisy measurement is generated as 
$f(j)=\mathrm{Poisson}(|(\mathcal A u_\delta)(j)|^2)~\forall j\in \tilde\Omega,$
with $u_\delta=\delta u$ at  peak level\footnote{Noise level increases as peak level $\delta$ decreases.} $\delta$.
We measure the quality of the reconstructed image $\tilde u$ by signal-to-noise ratio (SNR)\footnote{The SNR of noise free image is $+\infty$. The reconstructed  image with larger SNR usually means higher visual  quality.}
\[
\mathrm{SNR}(\tilde u,u)=-20\min\limits_{\varsigma\in\{\varsigma\in \mathbb C:~~|\varsigma|=1\}}\log(\|\varsigma\tilde u-u\|/\|\tilde  u\|),
\]
with the ground truth image $u$.
In order to measure the sparsity of coefficient matrix $\bm \alpha\in\mathbb C^{c,d}$ or $\mathbb R^{c,d}$, we introduce the sparsity level \footnote{Smaller sparsity values means  sparser of the data.}
$
S(\bm\alpha)=\dfrac{\|\bm\alpha\|^{iso}_0}{c \times d}\times 100\%.
$

 We compare  our proposed methods with ADMM for phase retrieval algorithm without any regularization \cite{wen2012} and the total variation (TV) based Poisson noise removal method  \cite{chang2016Total}. Denote  ADMM for phase retrieval \cite{wen2012}, TV based Poisson denoising method \cite{chang2016Total},  our proposed  Algorithm I and Algorithm II with isotropic $L^0$ norm by  ``PR'', ``TVPR'', ``ALGI'' and ``ALGII'', respectively. We denote  Algorithm I with anisotropic $L^0$ term for complex-valued images by ``ALGI${}_{aniso}$''. Since our proposed algorithms solve the nonconvex optimization problem,  they are initialized by the results of ``PR'' in order to increase the robustness and convergence speed.

\subsection{CDP for real-valued and complex-valued images}
 For the first pattern,  the octanary CDP is explored, and specifically   each
element of $I_j$ in \eqref{cdp} takes a value randomly among the eight candidates, \emph{i.e.}, $\{\pm \sqrt{2}/2, \pm \sqrt{2}{\mathbf i}/2, \pm \sqrt{3},\pm \sqrt{3}{\mathbf i}\}$. Set  the number of masks $K=2, 4$ as \cite{chang2016Total} for real-valued and complex-valued images, respectively.

\subsubsection{Real-valued Images}
We  will show the performance of our proposed algorithms from noisy measurement on four real-valued images shown in Fig. \ref{ground1}(a)-(d), with  peak level $\delta\in\{5.0\times 10^{-3},1.0\times 10^{-2}\}$. Set $\eta=8.0\times 10^{-6}, \tau=4.5\times 10^{-4}$ in Algorithm I and Algorithm II. Set step sizes $c^k=10, d^k=50,$ and $e^k=1.5, 1.2$ for $\delta=5.0\times 10^{-3},1.0\times 10^{-2}$ respectively in Algorithm II. Both two algorithms stop after $T=100$ iterations. All four images share the same parameters with the others.

\begin{figure}
\vskip -.2in
\begin{center}
\subfigure[]{\includegraphics[width=.11\textwidth]{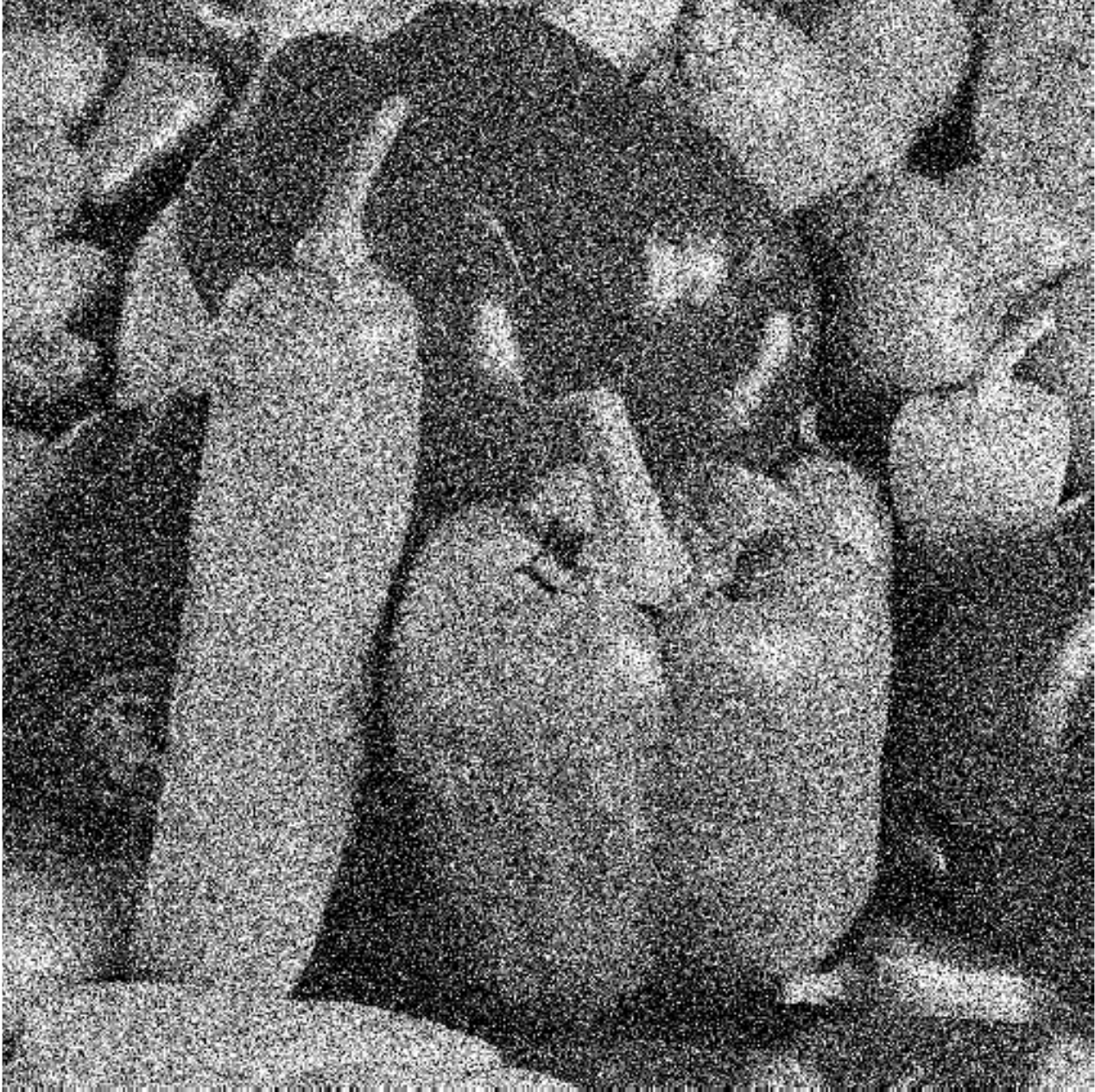}}
\subfigure[]{\includegraphics[width=.11\textwidth]{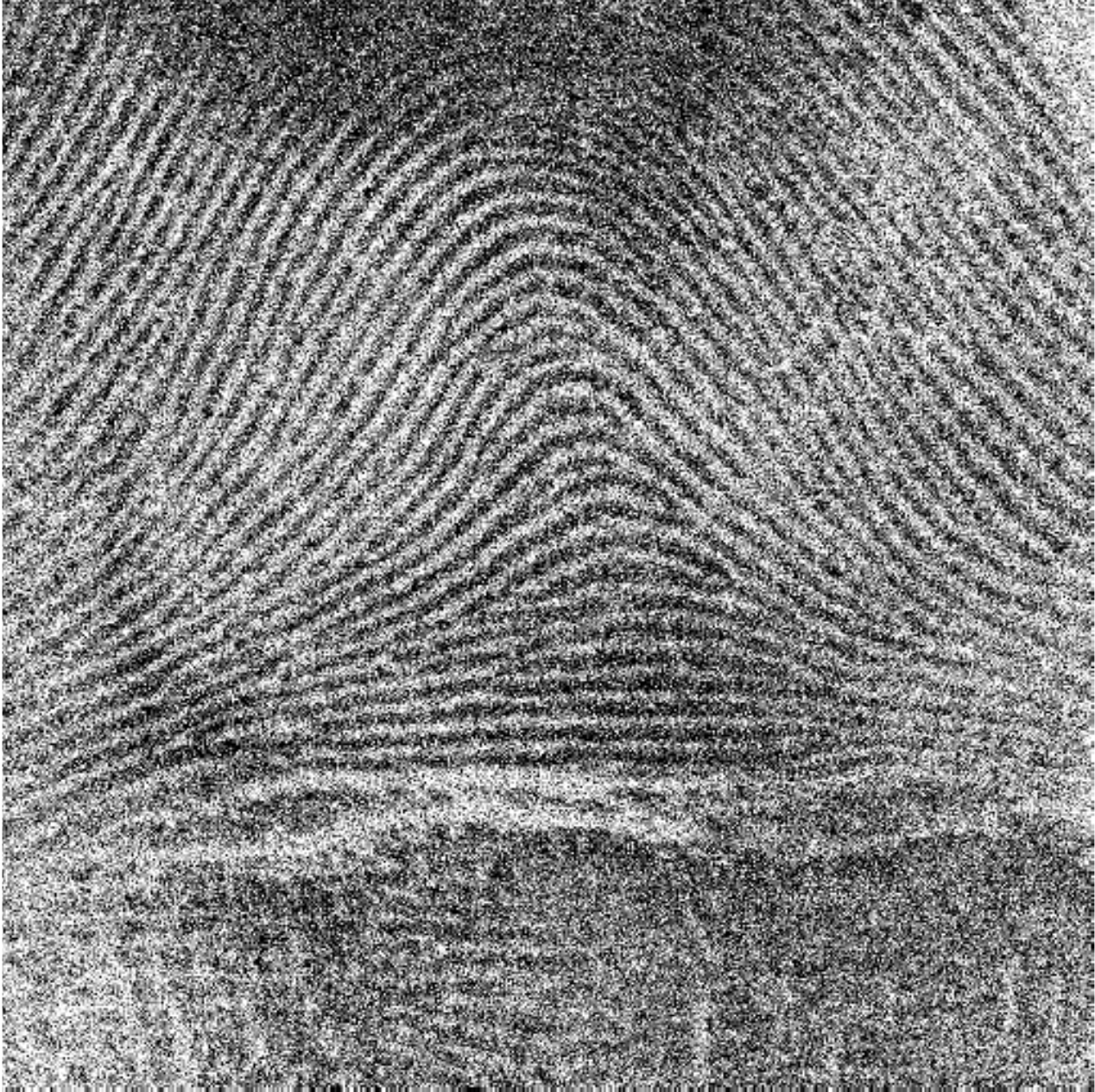}}
\subfigure[]{\includegraphics[width=.11\textwidth]{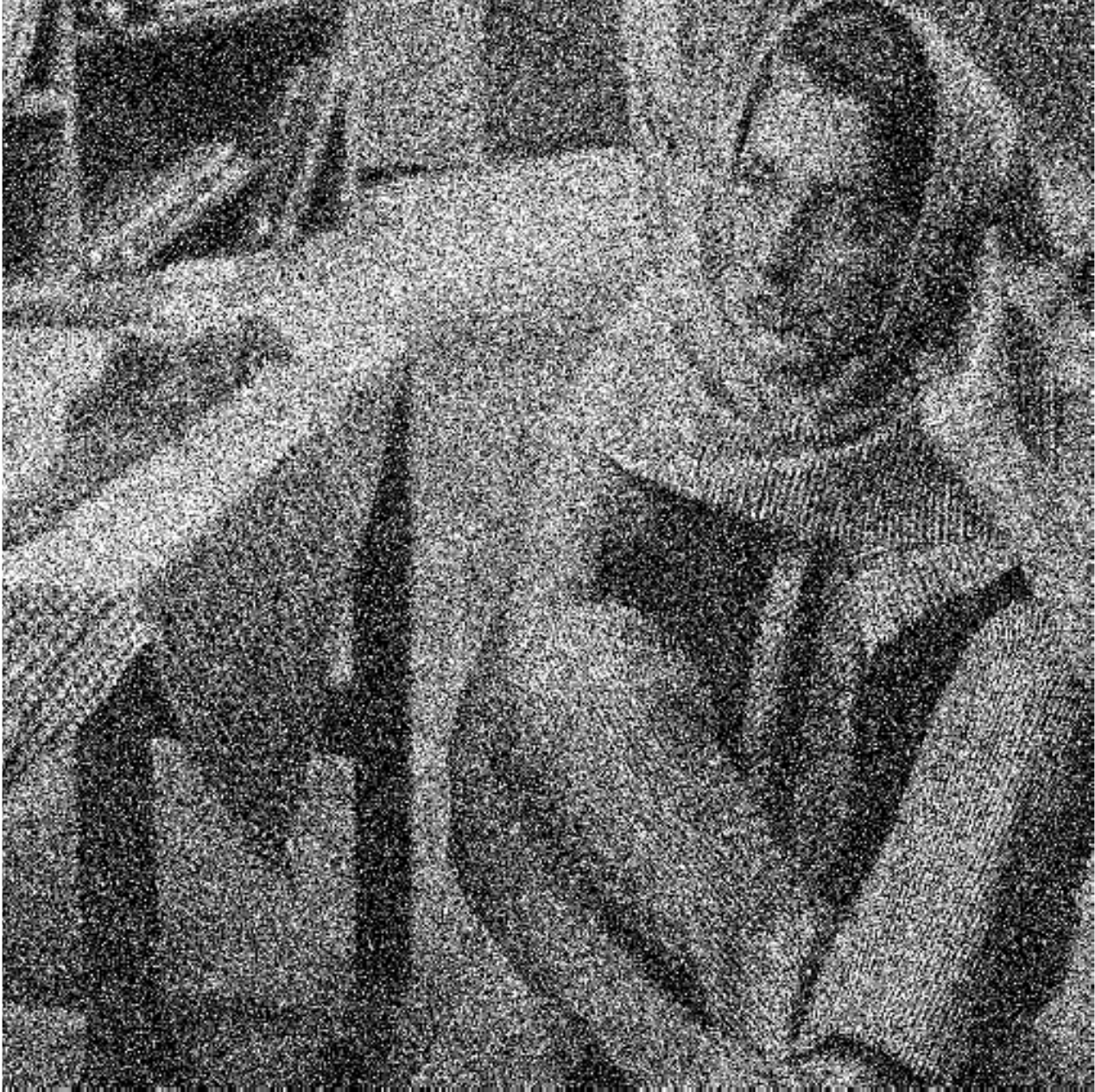}}
\subfigure[]{\includegraphics[width=.11\textwidth]{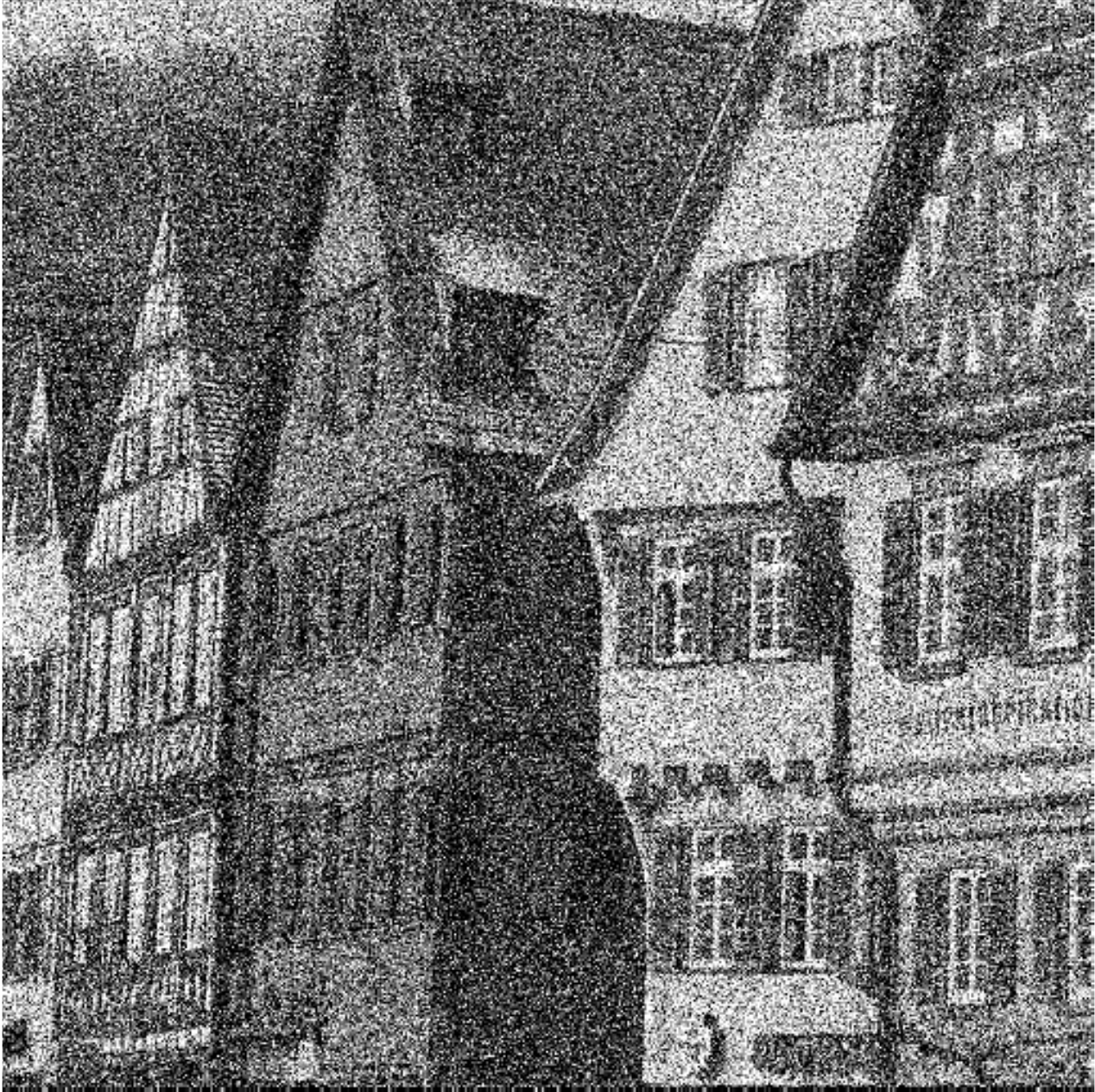}}\\\vskip -.05in
\subfigure[]{\includegraphics[width=.11\textwidth]{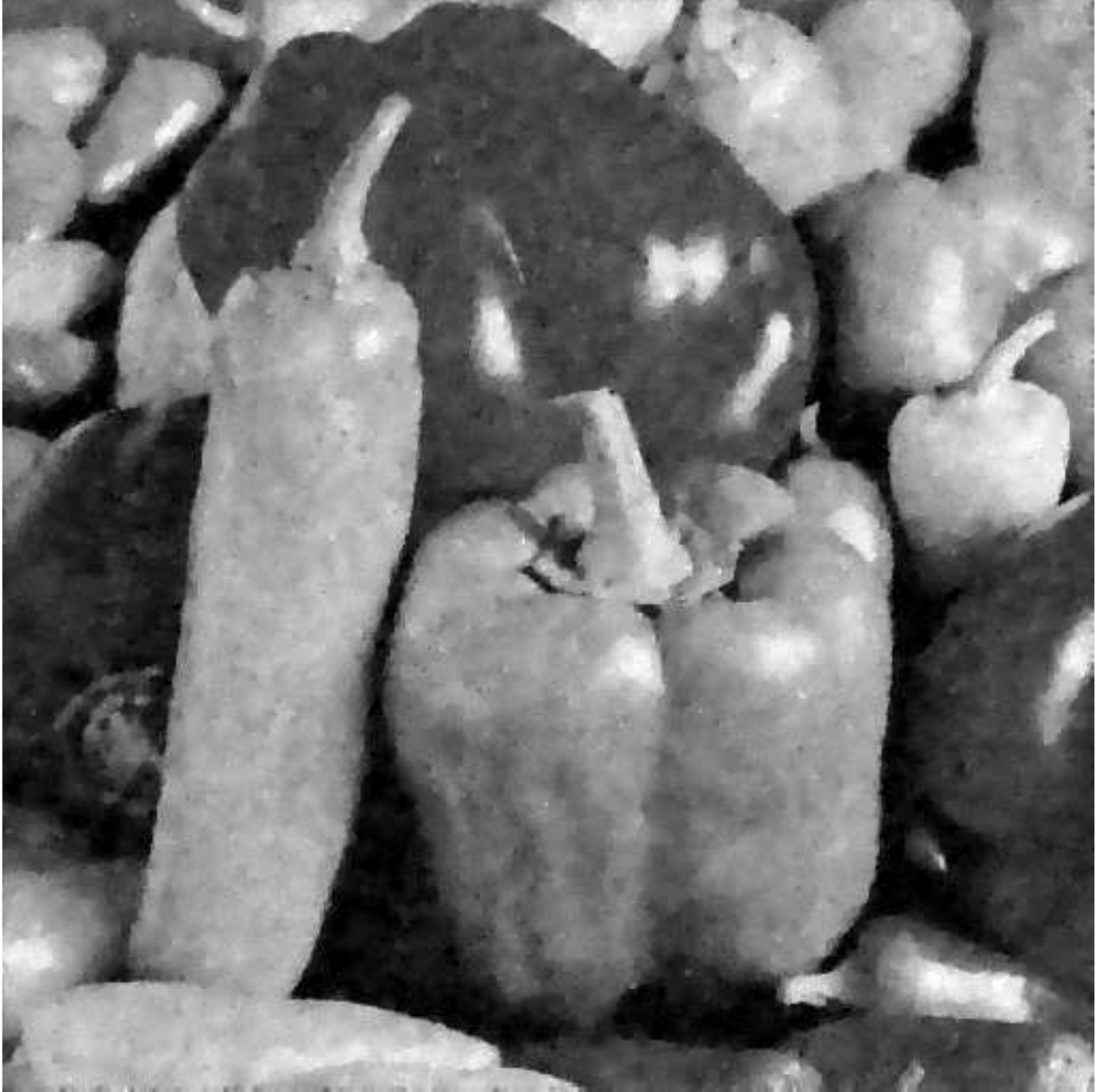}}
\subfigure[]{\includegraphics[width=.11\textwidth]{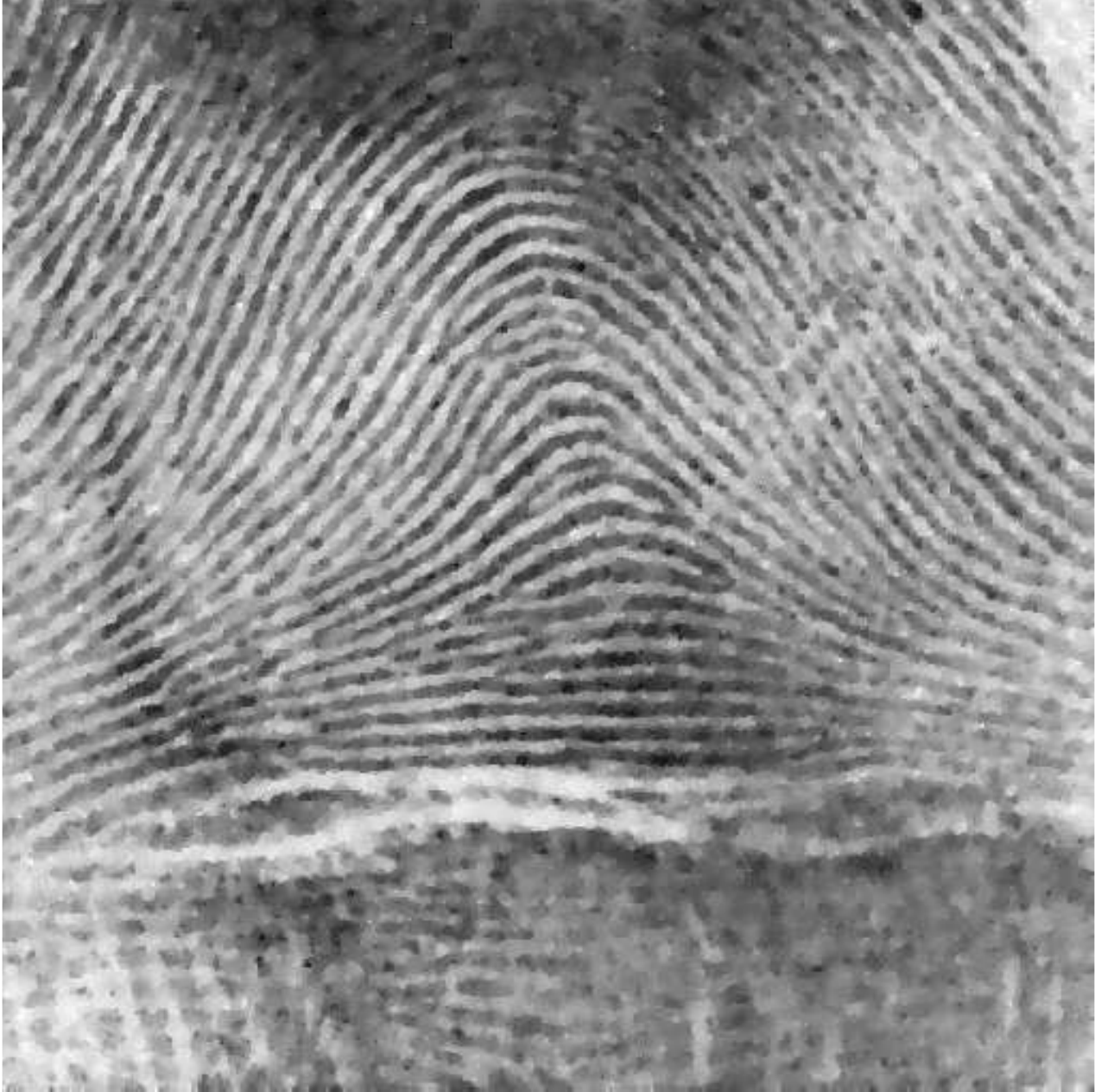}}
\subfigure[]{\includegraphics[width=.11\textwidth]{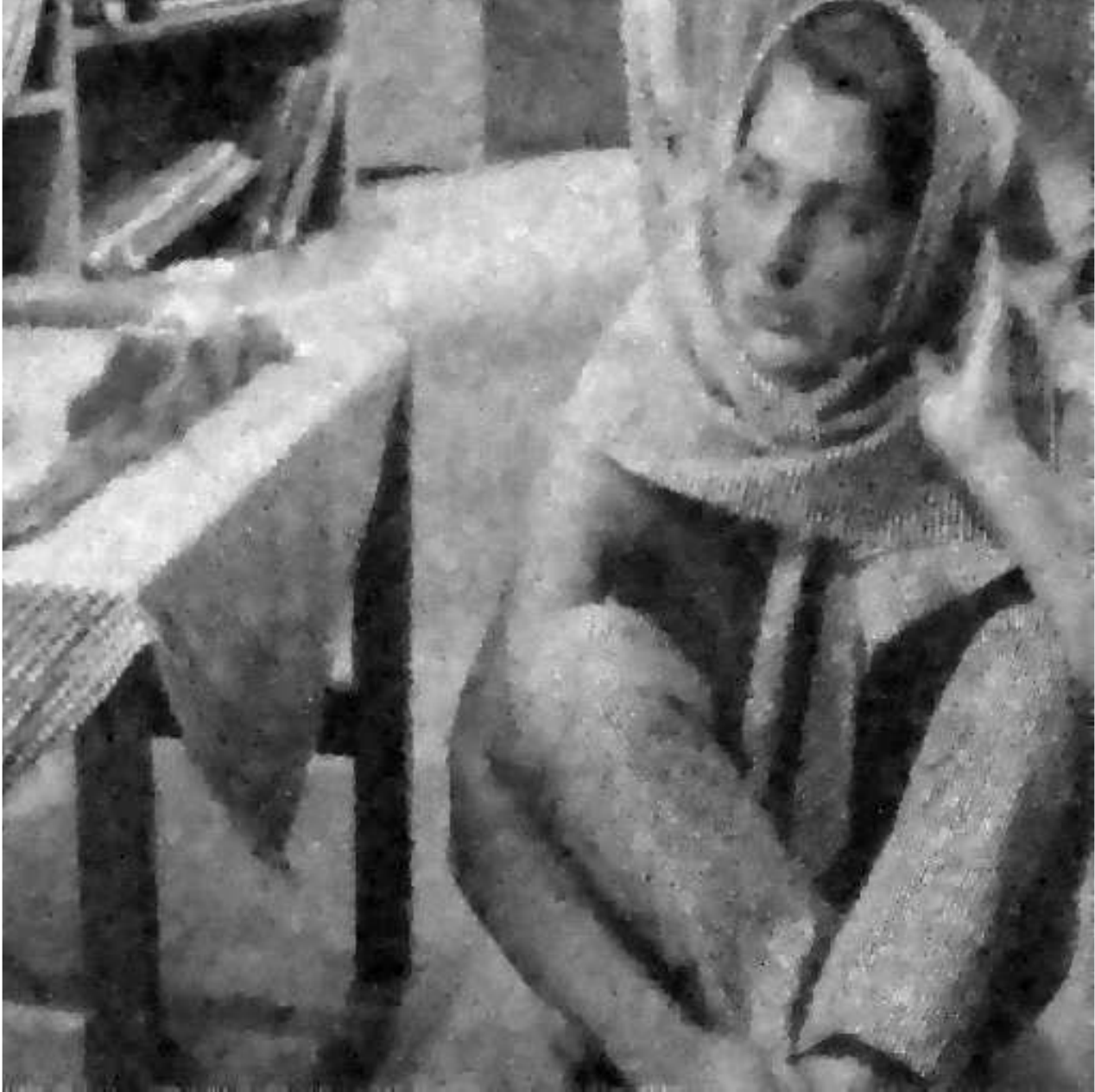}}
\subfigure[]{\includegraphics[width=.11\textwidth]{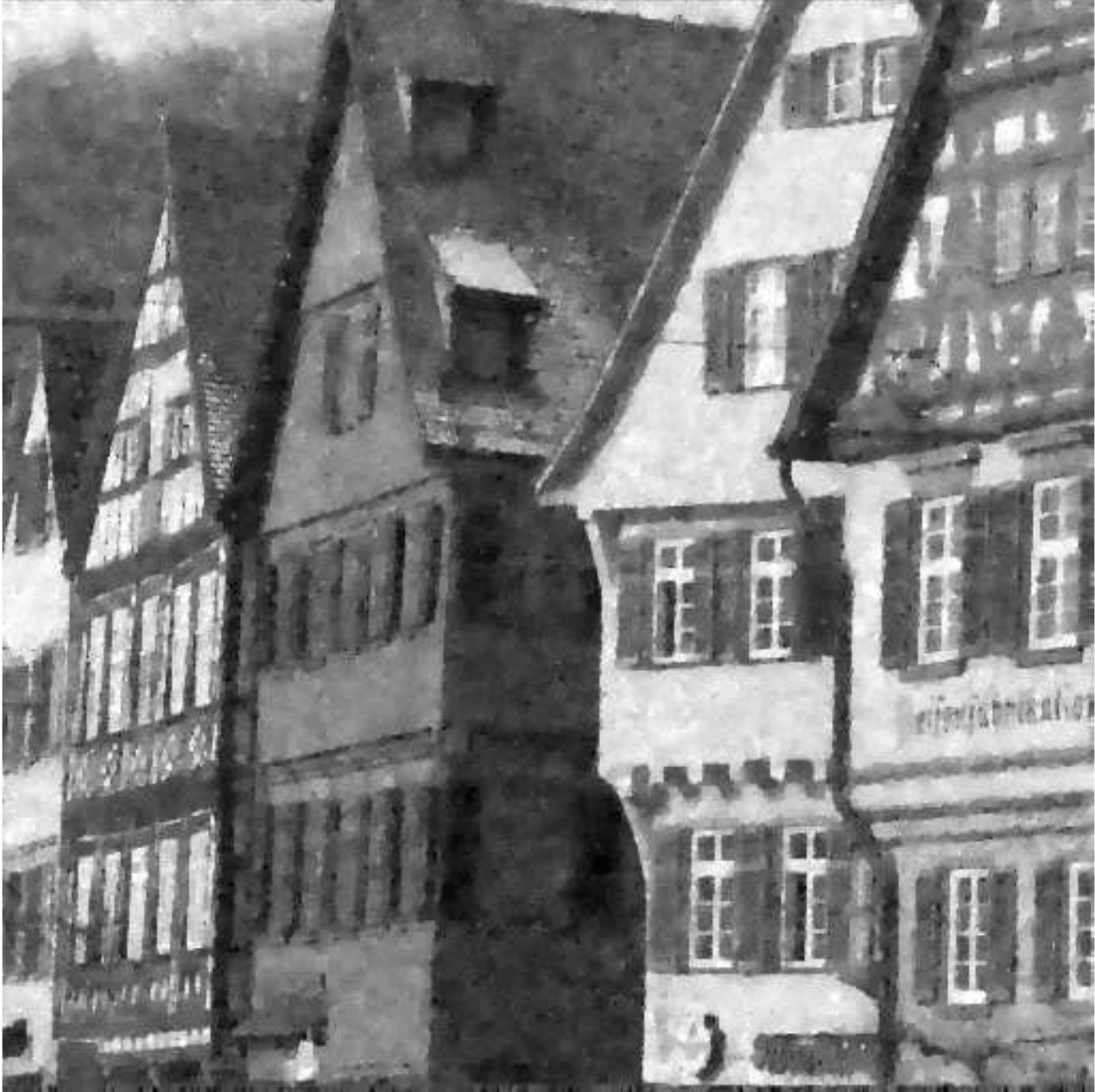}}\\\vskip -.05in
\subfigure[]{\includegraphics[width=.11\textwidth]{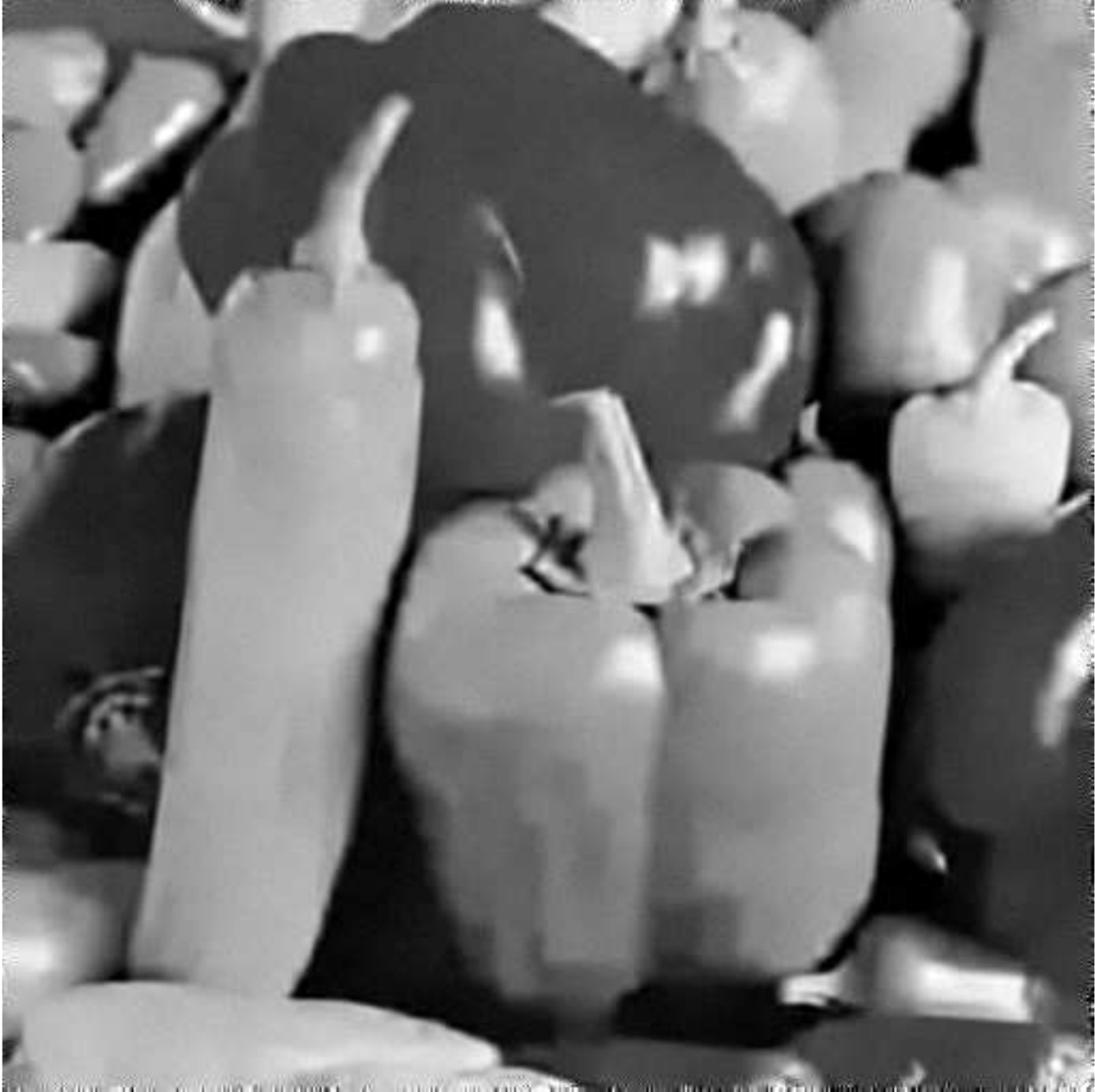}}
\subfigure[]{\includegraphics[width=.11\textwidth]{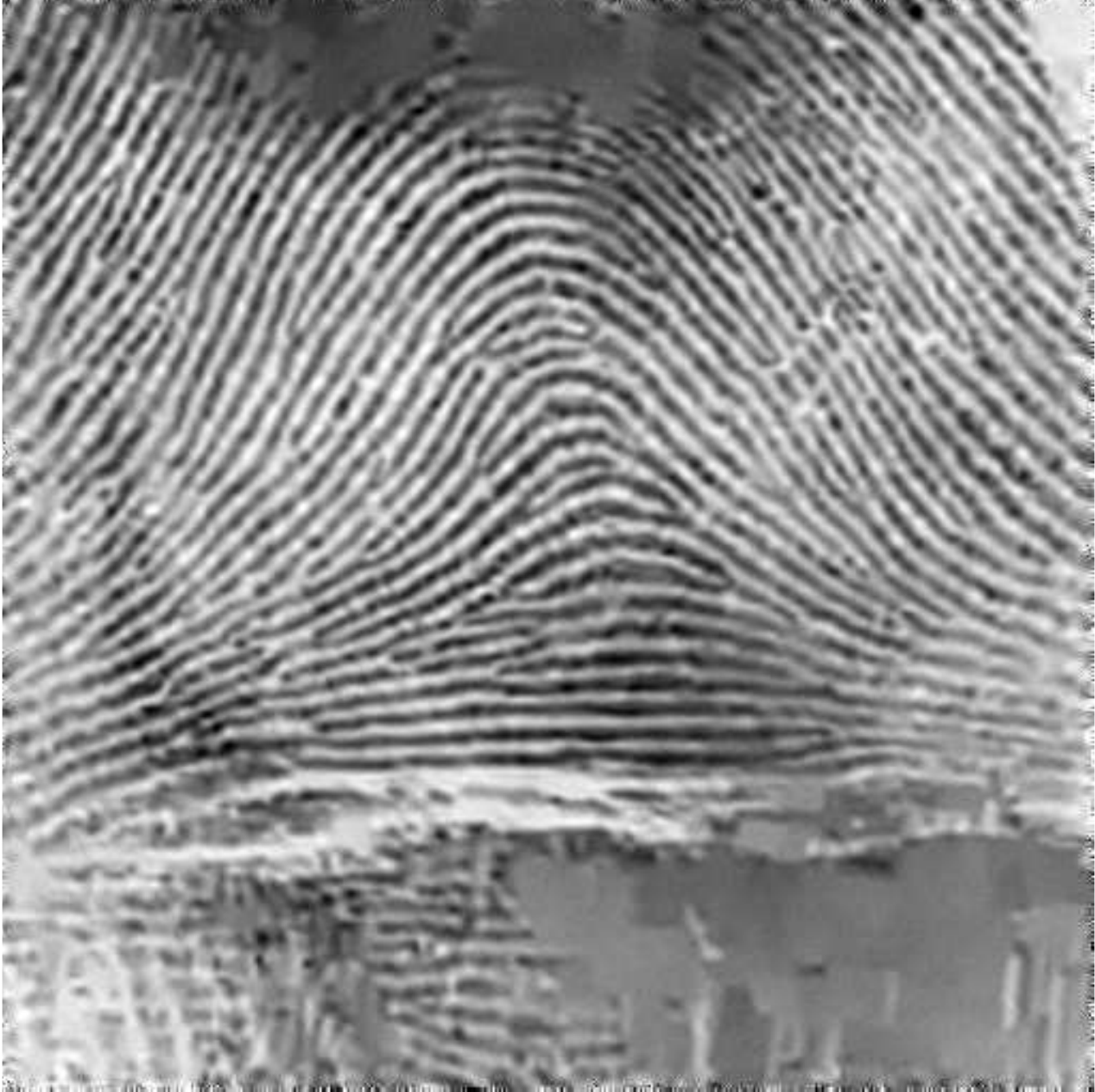}}
\subfigure[]{\includegraphics[width=.11\textwidth]{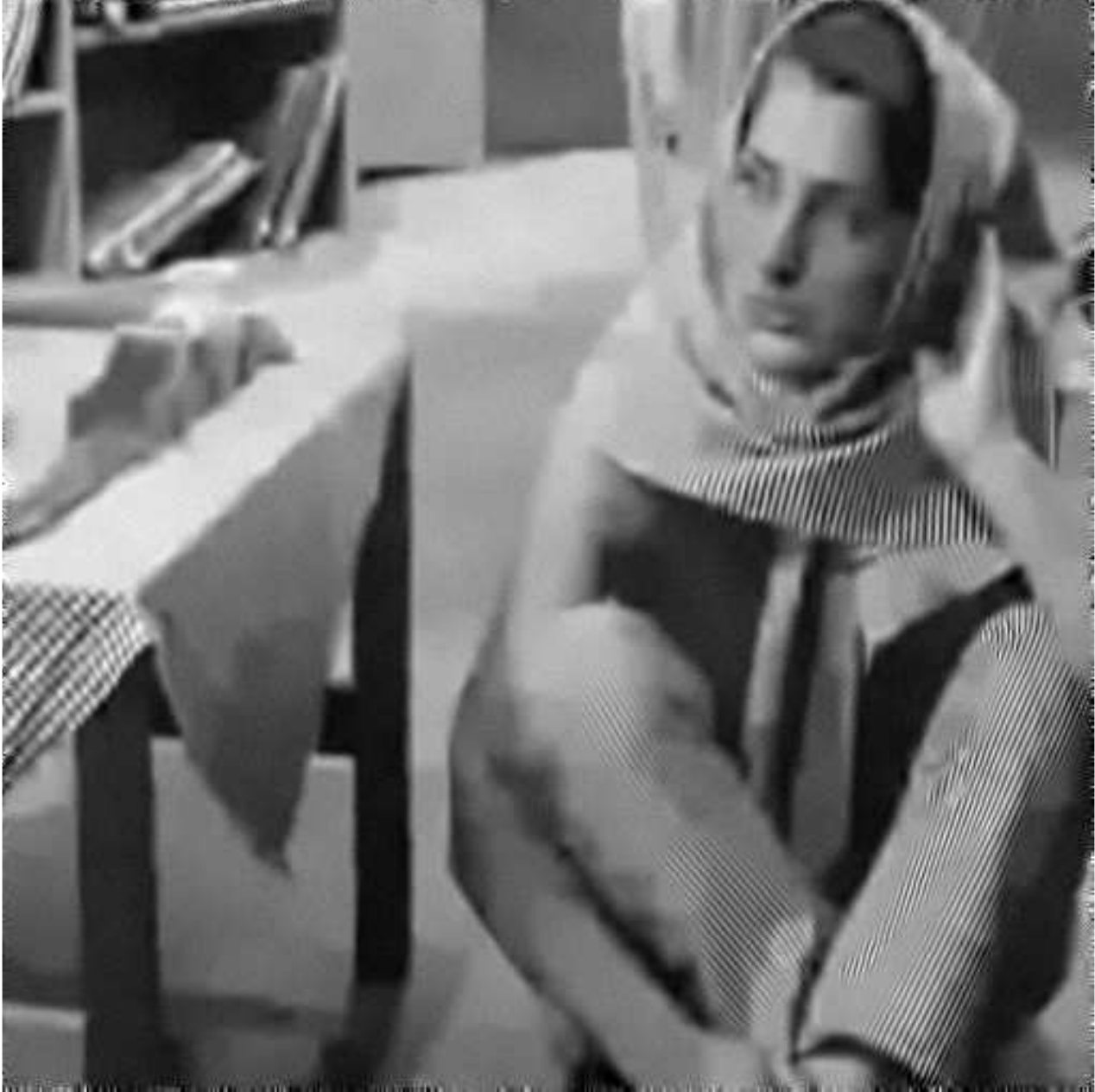}}
\subfigure[]{\includegraphics[width=.11\textwidth]{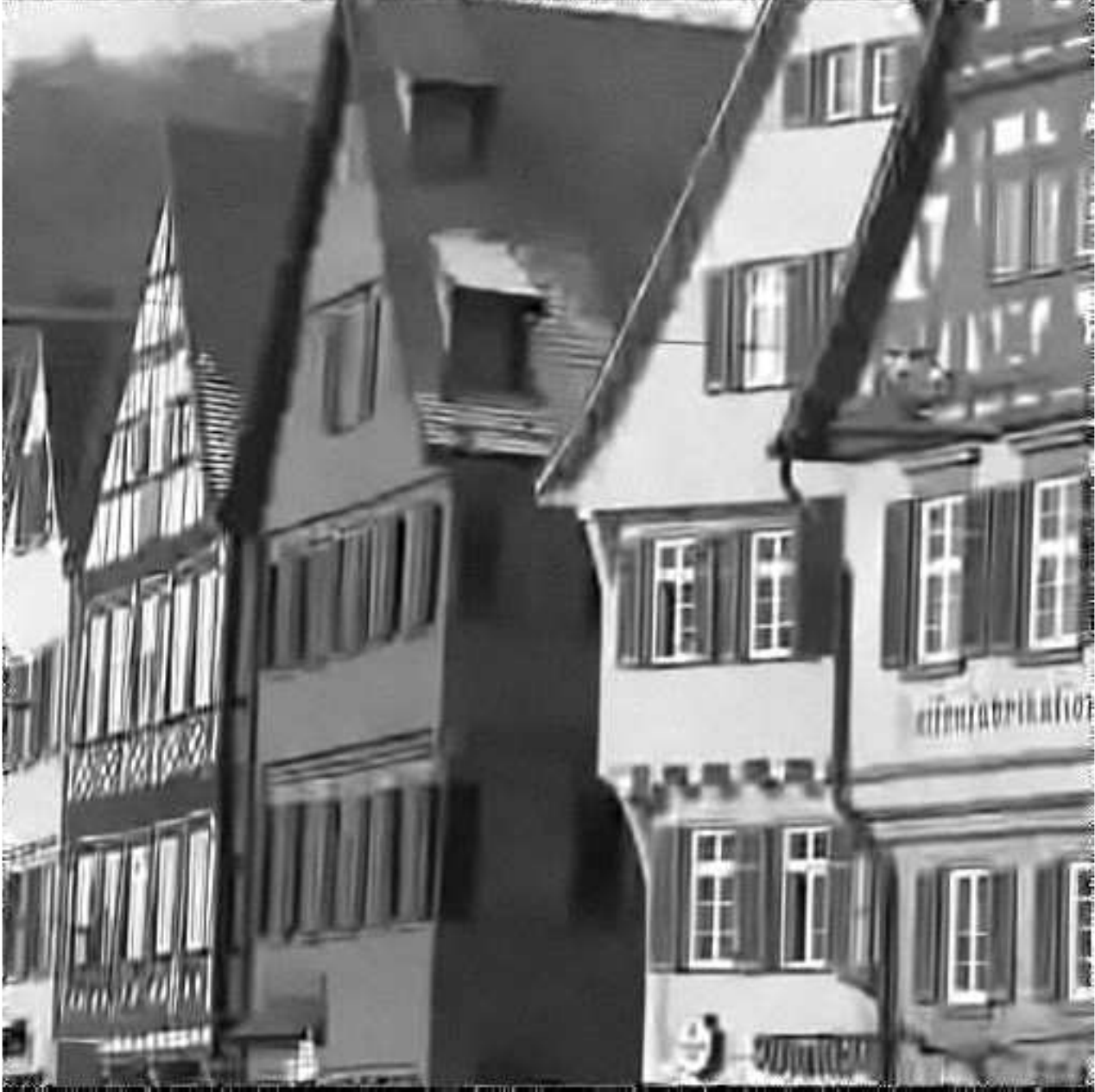}}\\\vskip -.05in
\subfigure[]{\includegraphics[width=.11\textwidth]{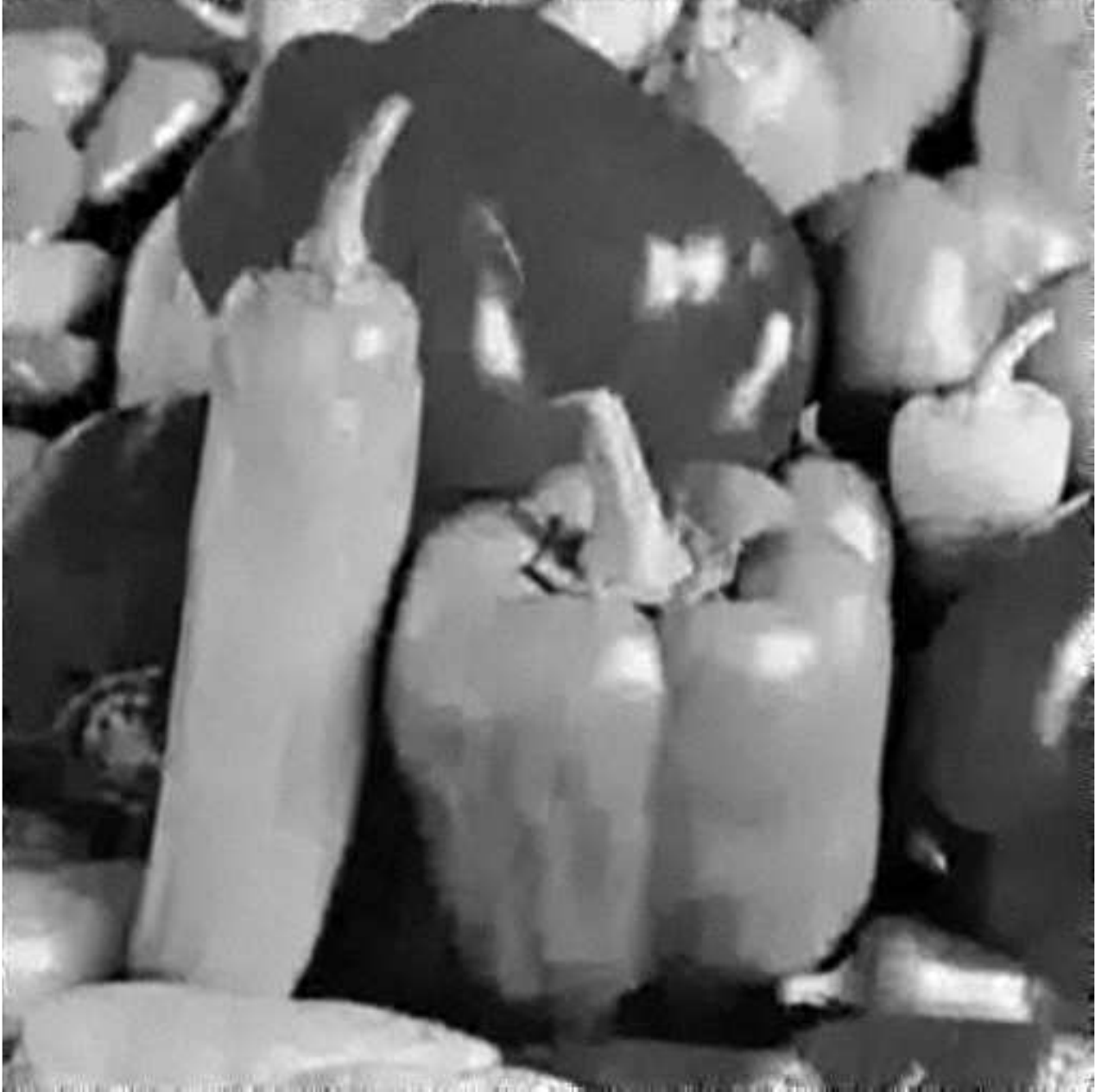}}
\subfigure[]{\includegraphics[width=.11\textwidth]{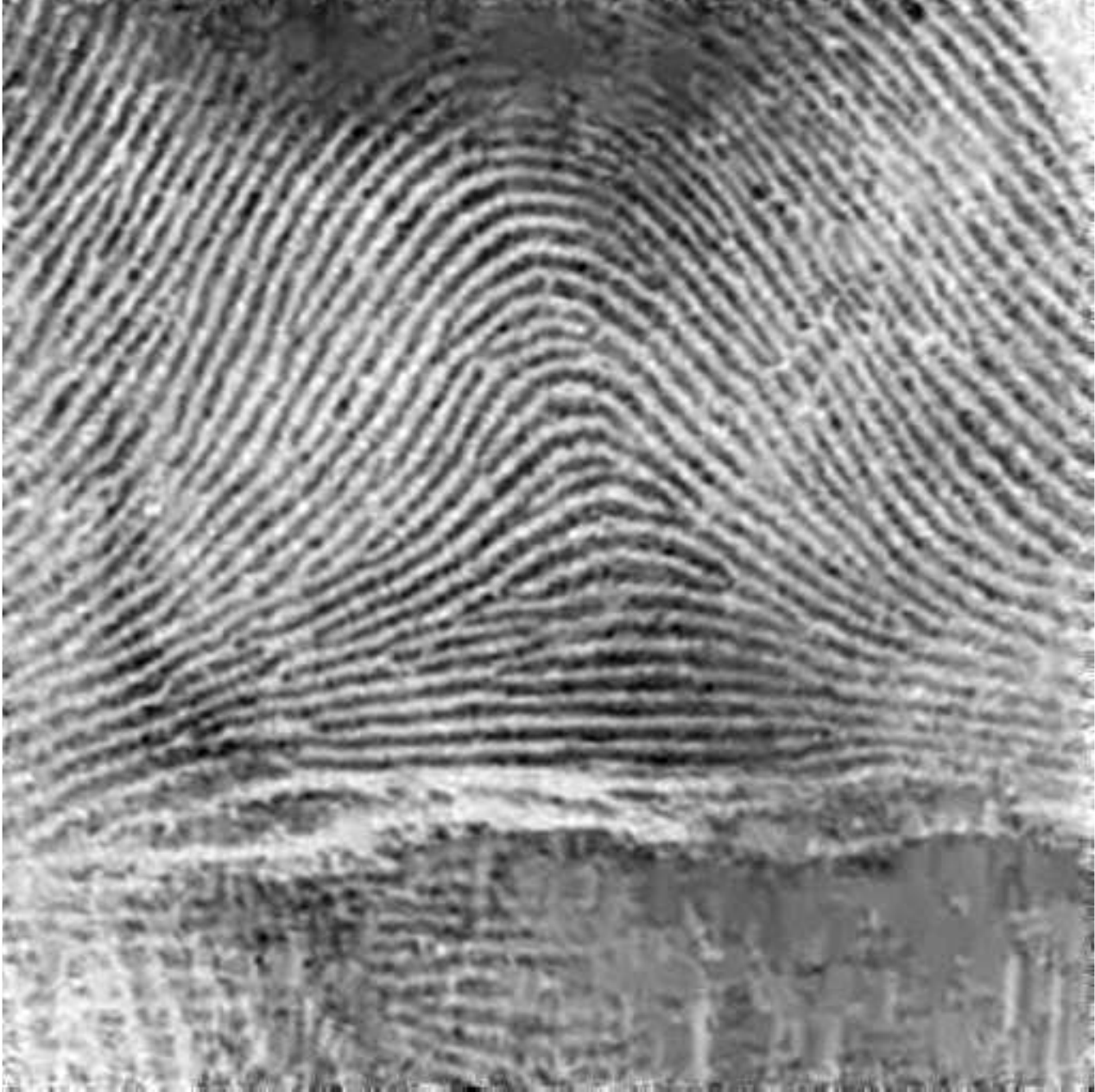}}
\subfigure[]{\includegraphics[width=.11\textwidth]{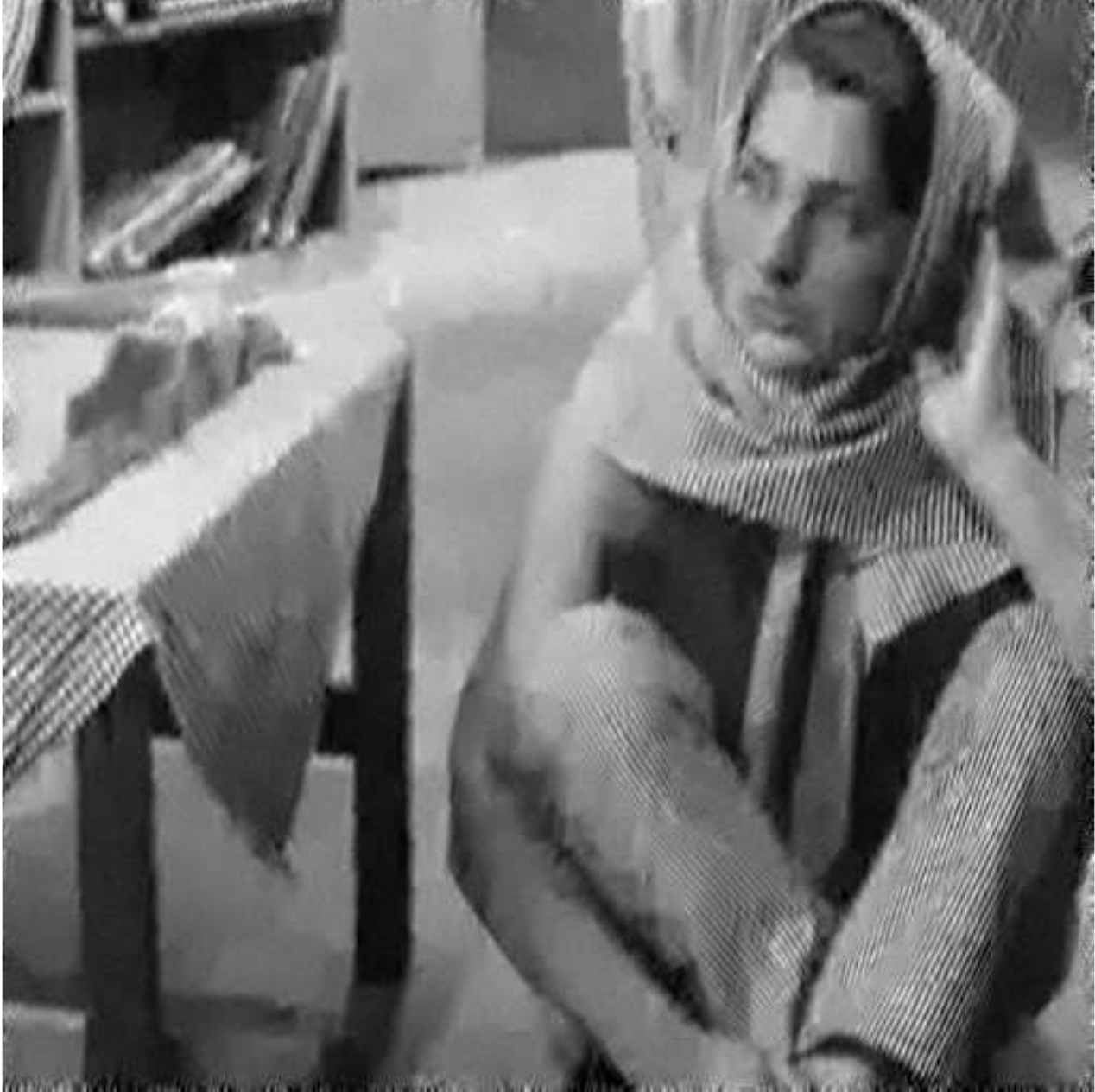}}
\subfigure[]{\includegraphics[width=.11\textwidth]{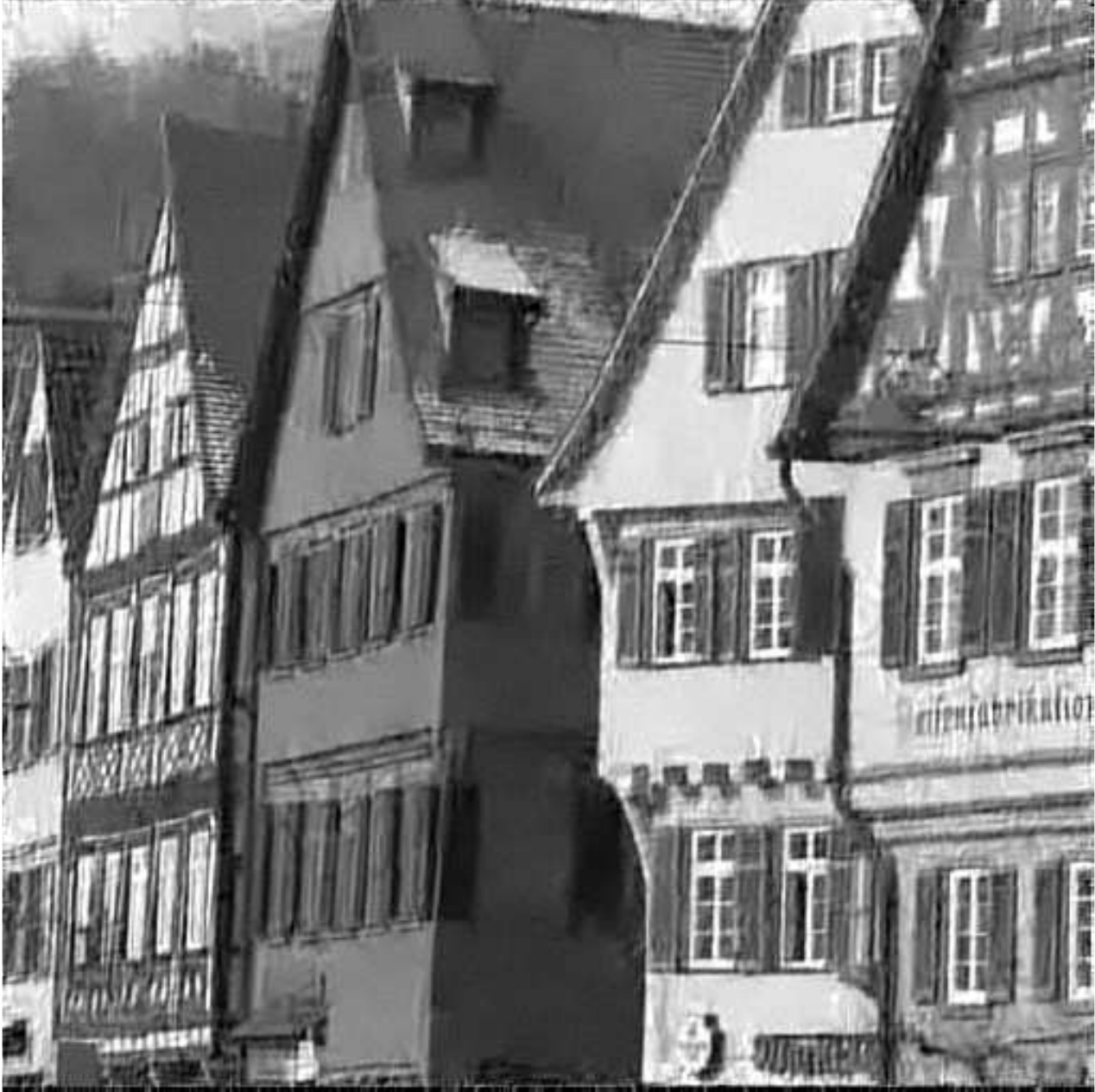}}\vskip -.05in
\end{center}
\vskip -.1in
\caption{CDP with $\delta=5.0\times 10^{-3}$. First row: PR; Second row: TVPR; Third row: ALGI; Fourth row: ALGII.}
\label{cdp1}
\vskip -.1in
\end{figure}

\begin{figure}
\vskip -.3in
\begin{center}
\subfigure[]{\includegraphics[width=.11\textwidth]{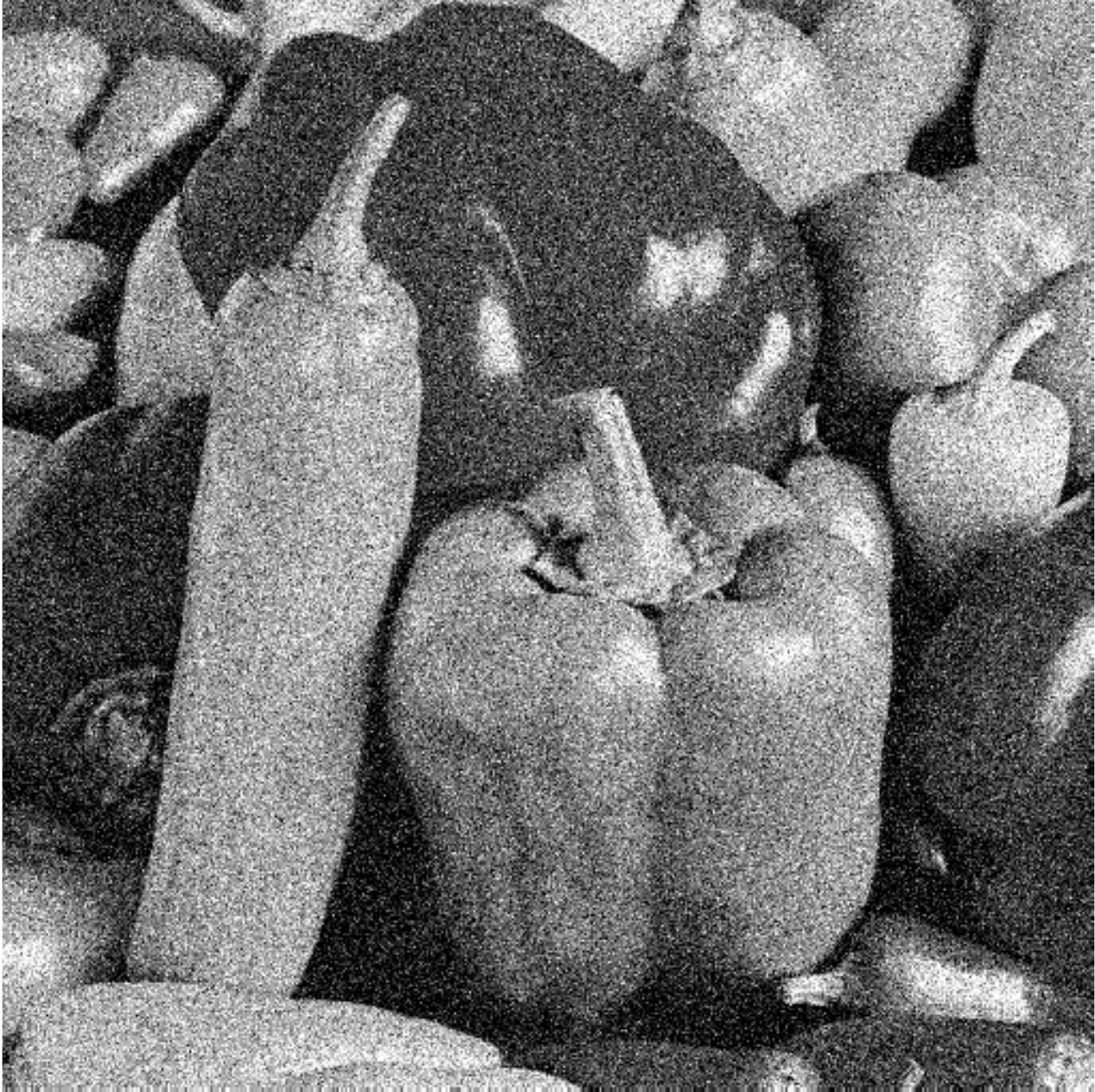}}
\subfigure[]{\includegraphics[width=.11\textwidth]{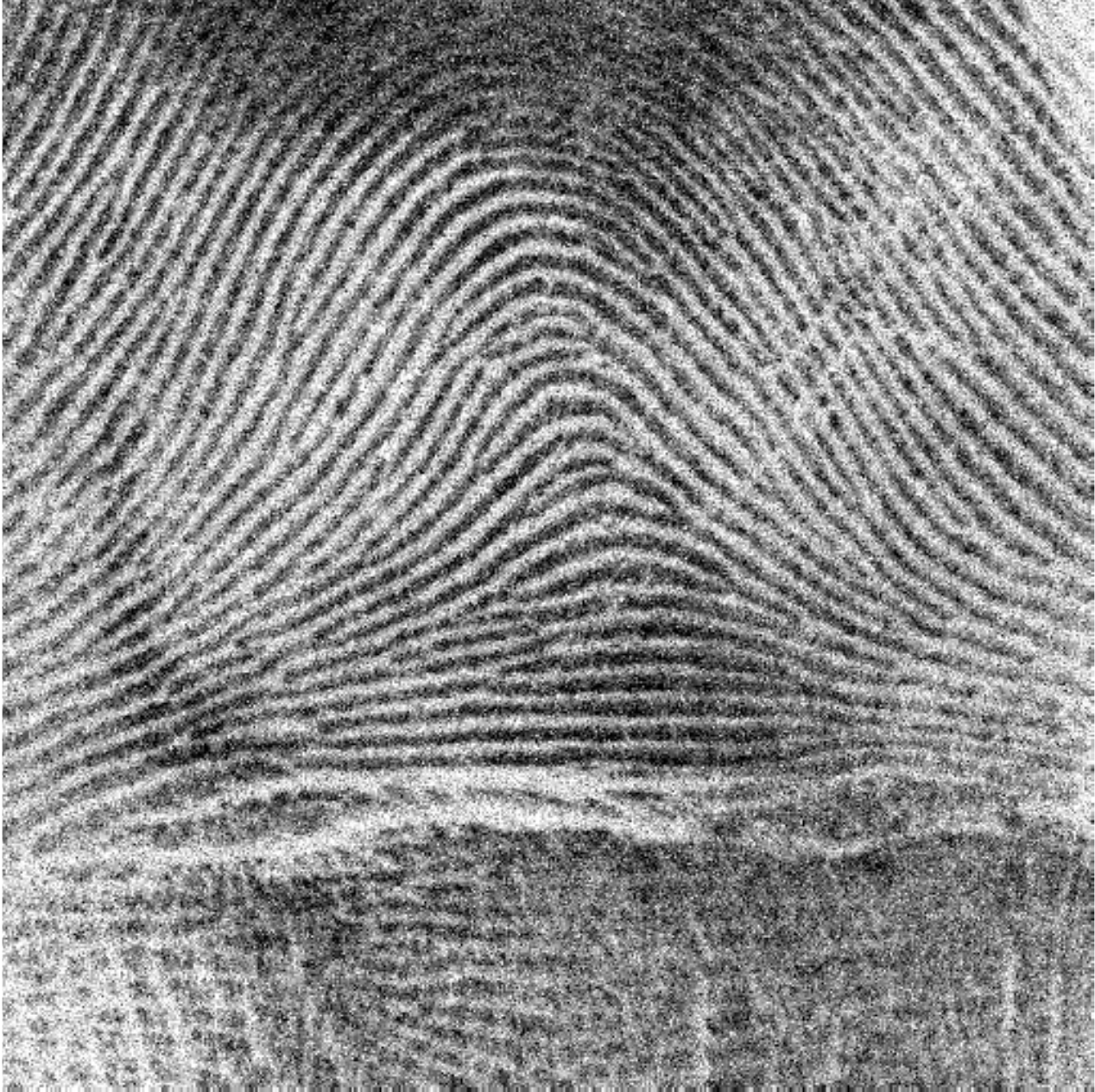}}
\subfigure[]{\includegraphics[width=.11\textwidth]{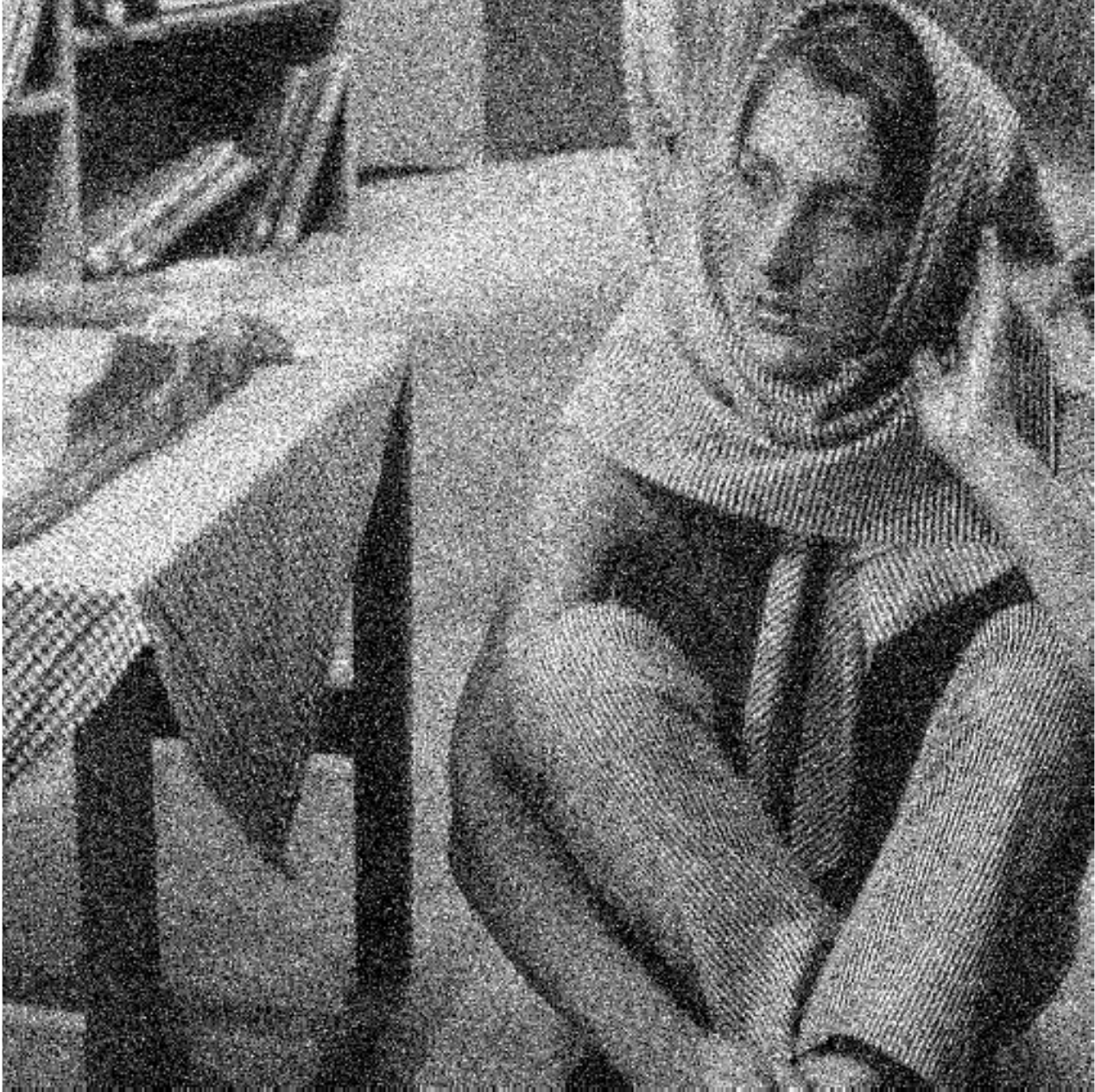}}
\subfigure[]{\includegraphics[width=.11\textwidth]{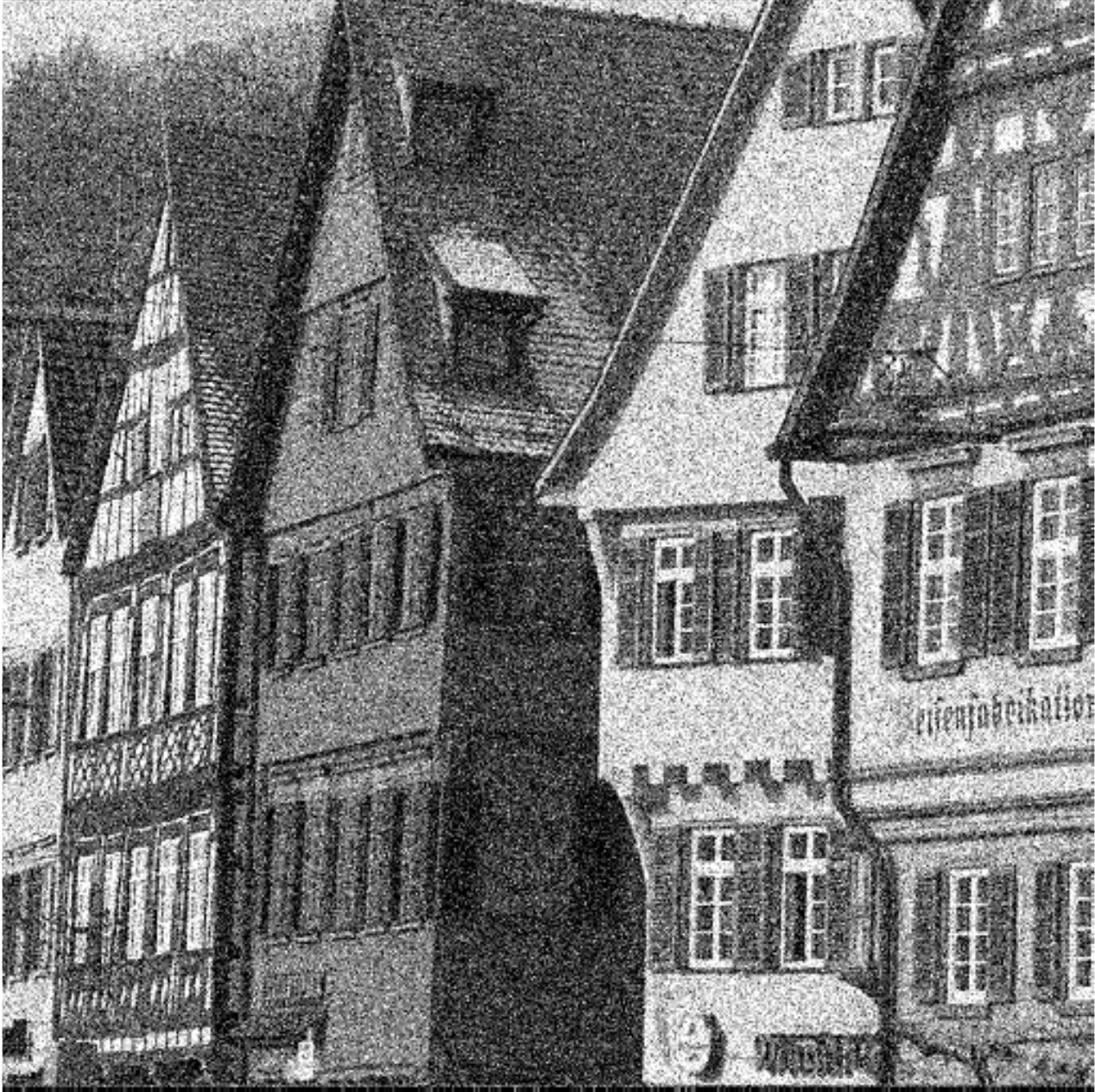}}\\\vskip -.05in
\subfigure[]{\includegraphics[width=.11\textwidth]{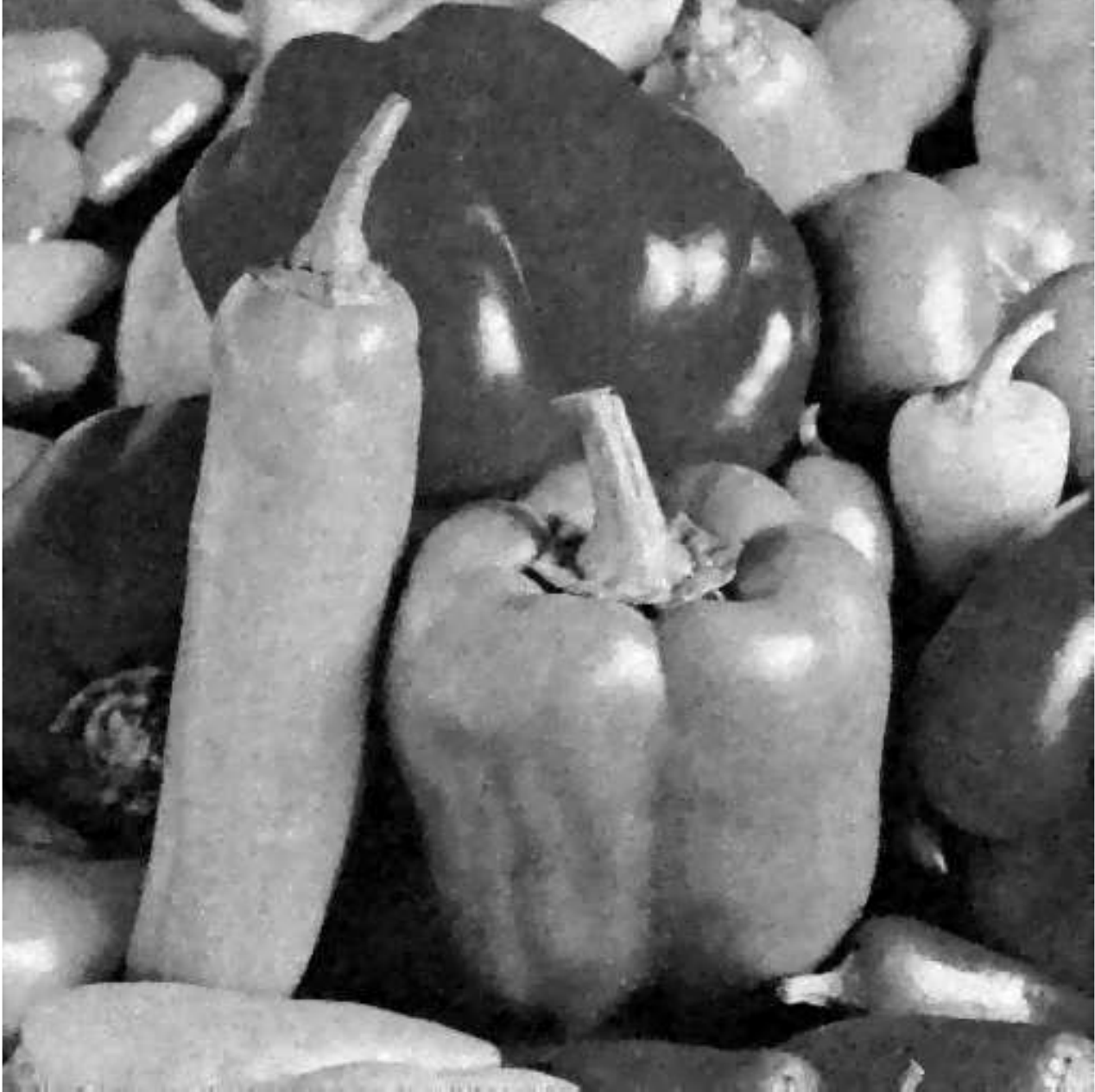}}
\subfigure[]{\includegraphics[width=.11\textwidth]{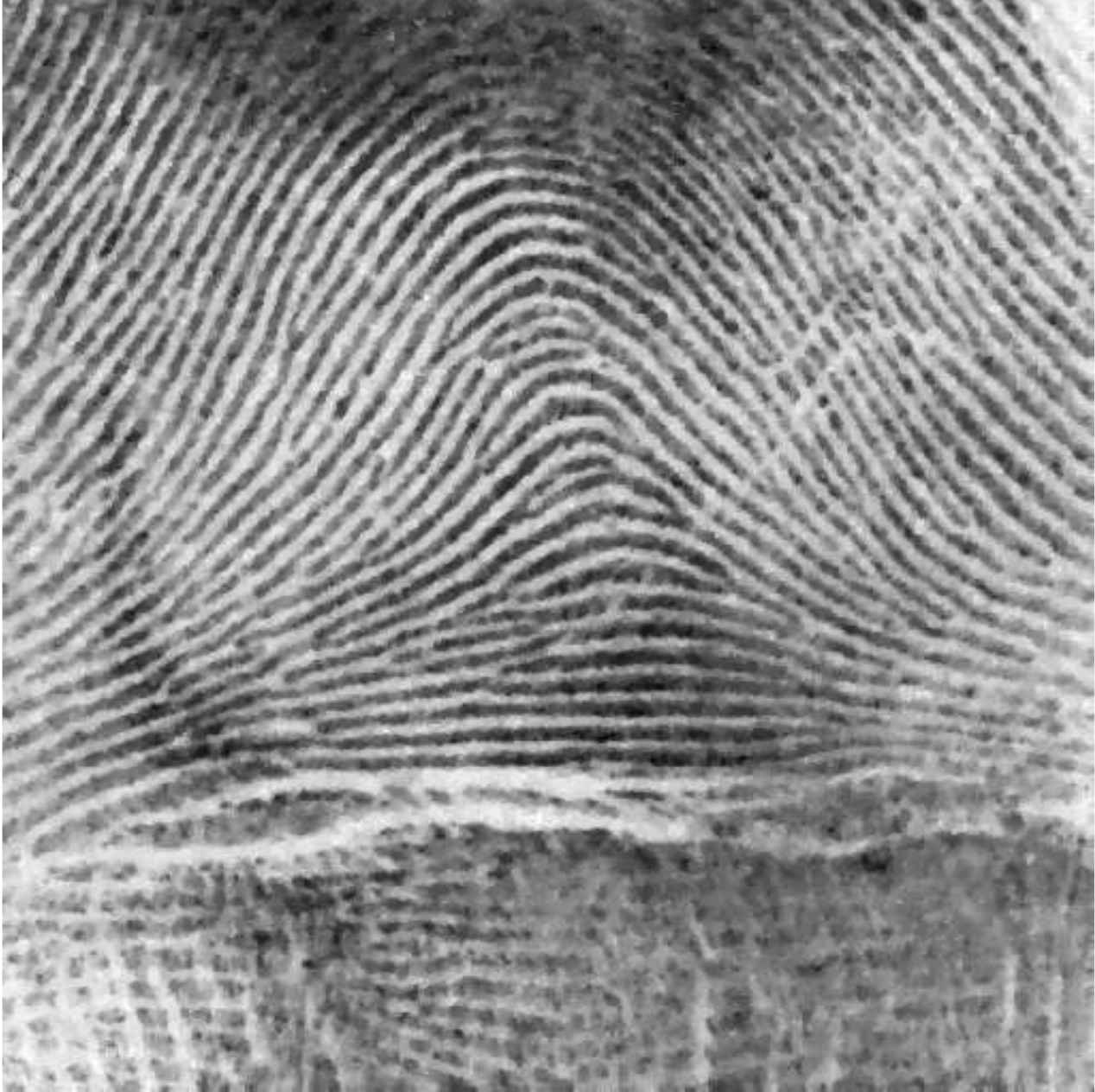}}
\subfigure[]{\includegraphics[width=.11\textwidth]{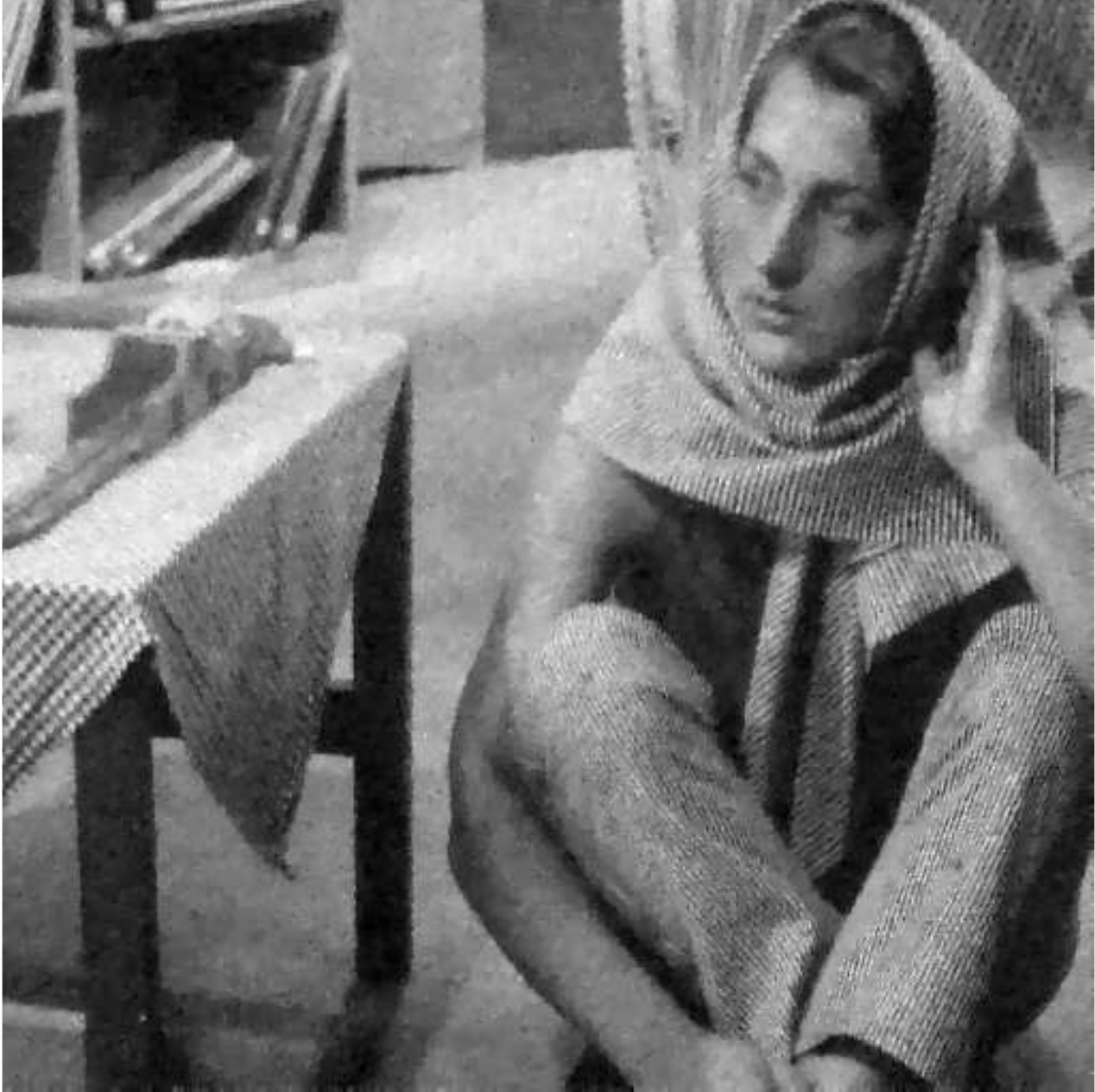}}
\subfigure[]{\includegraphics[width=.11\textwidth]{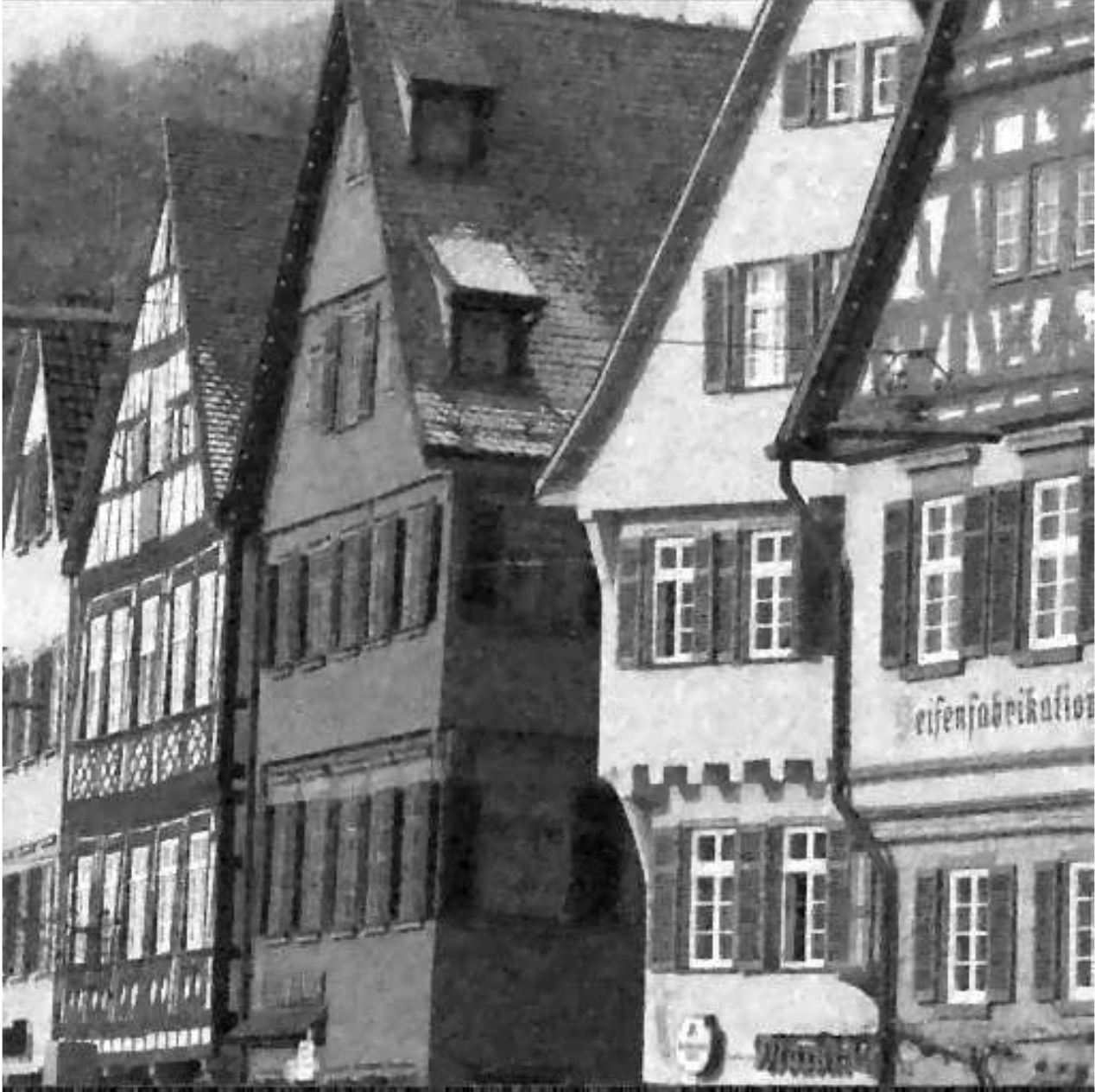}}\\\vskip -.05in
\subfigure[]{\includegraphics[width=.11\textwidth]{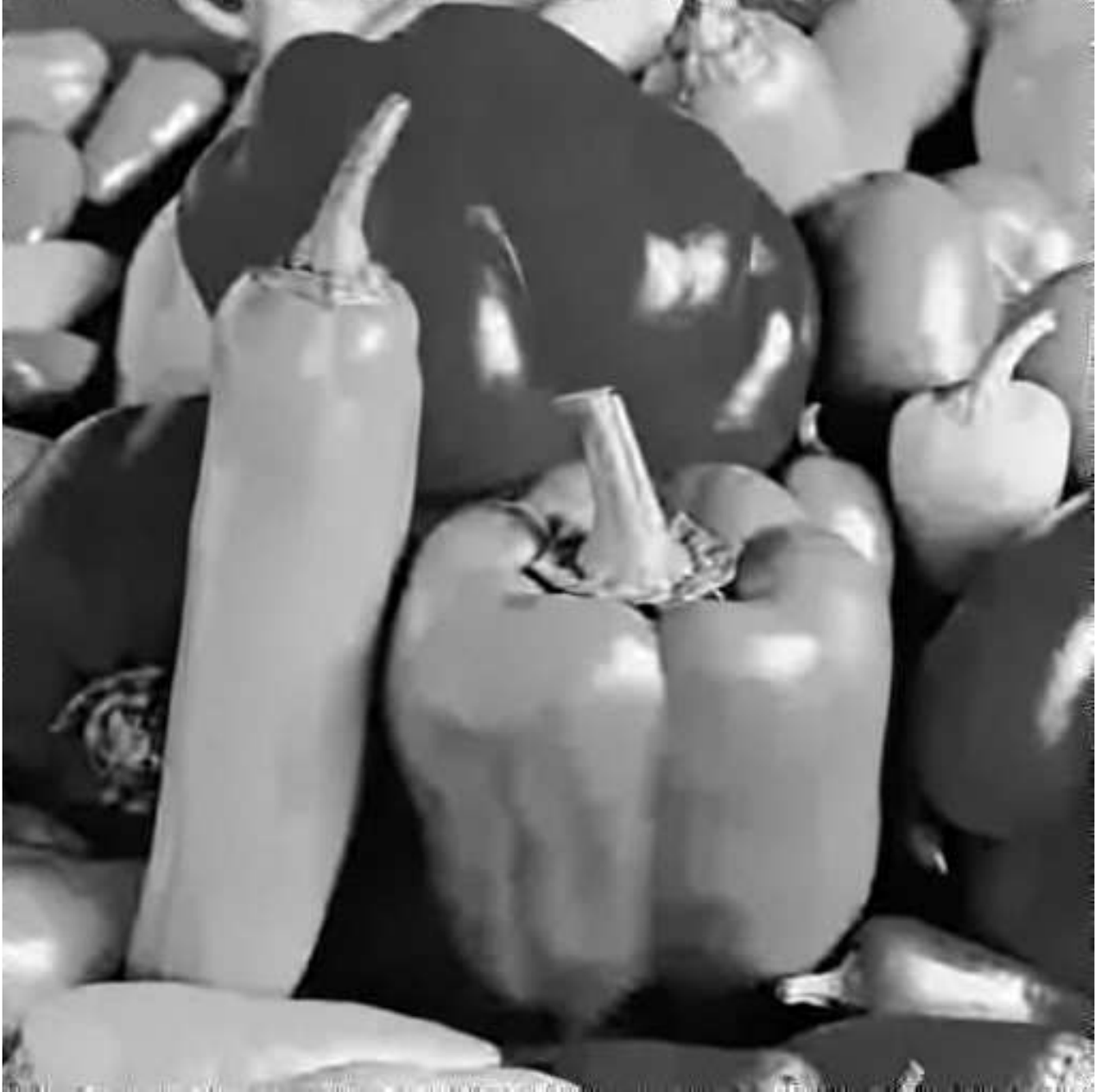}}
\subfigure[]{\includegraphics[width=.11\textwidth]{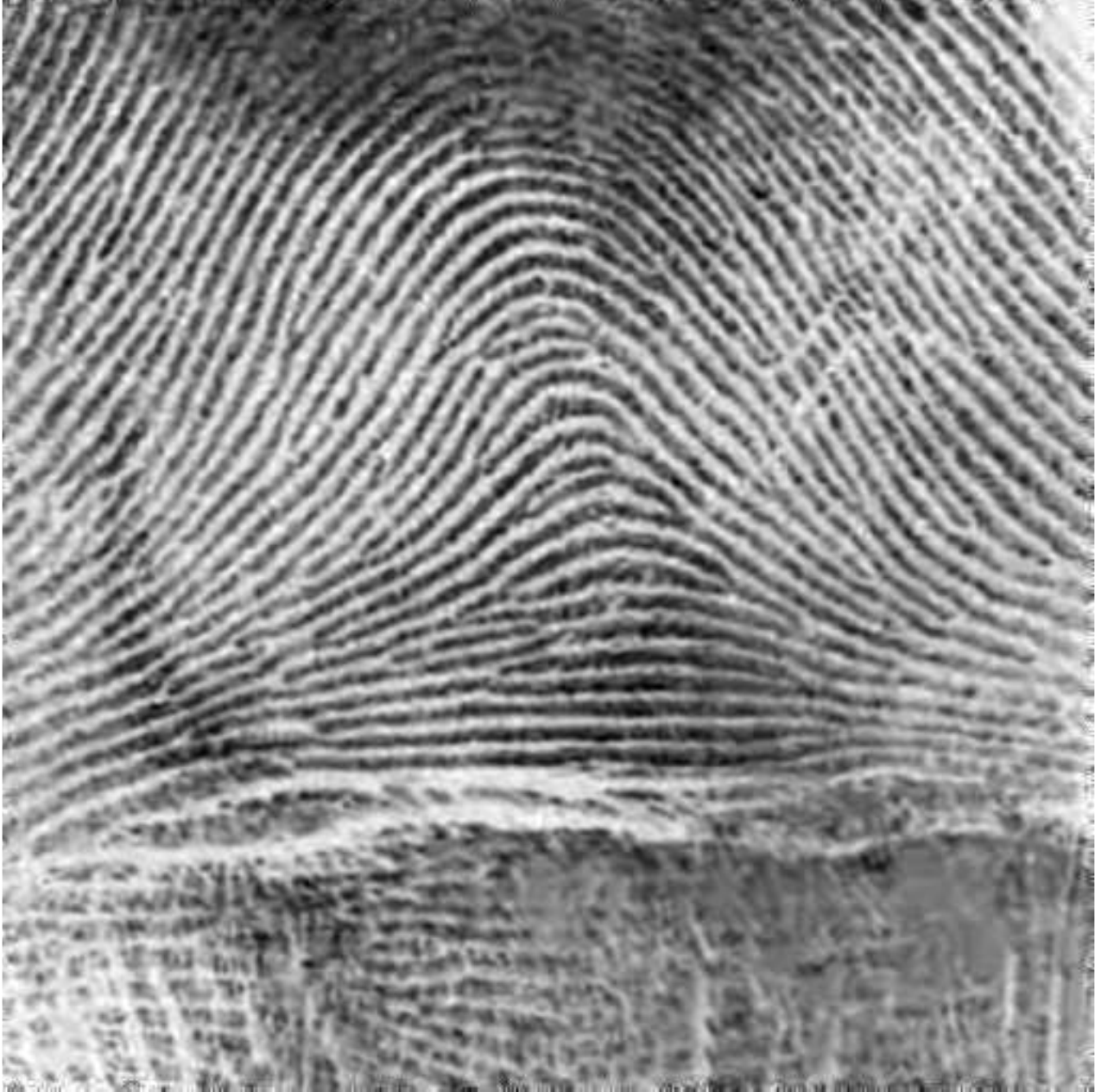}}
\subfigure[]{\includegraphics[width=.11\textwidth]{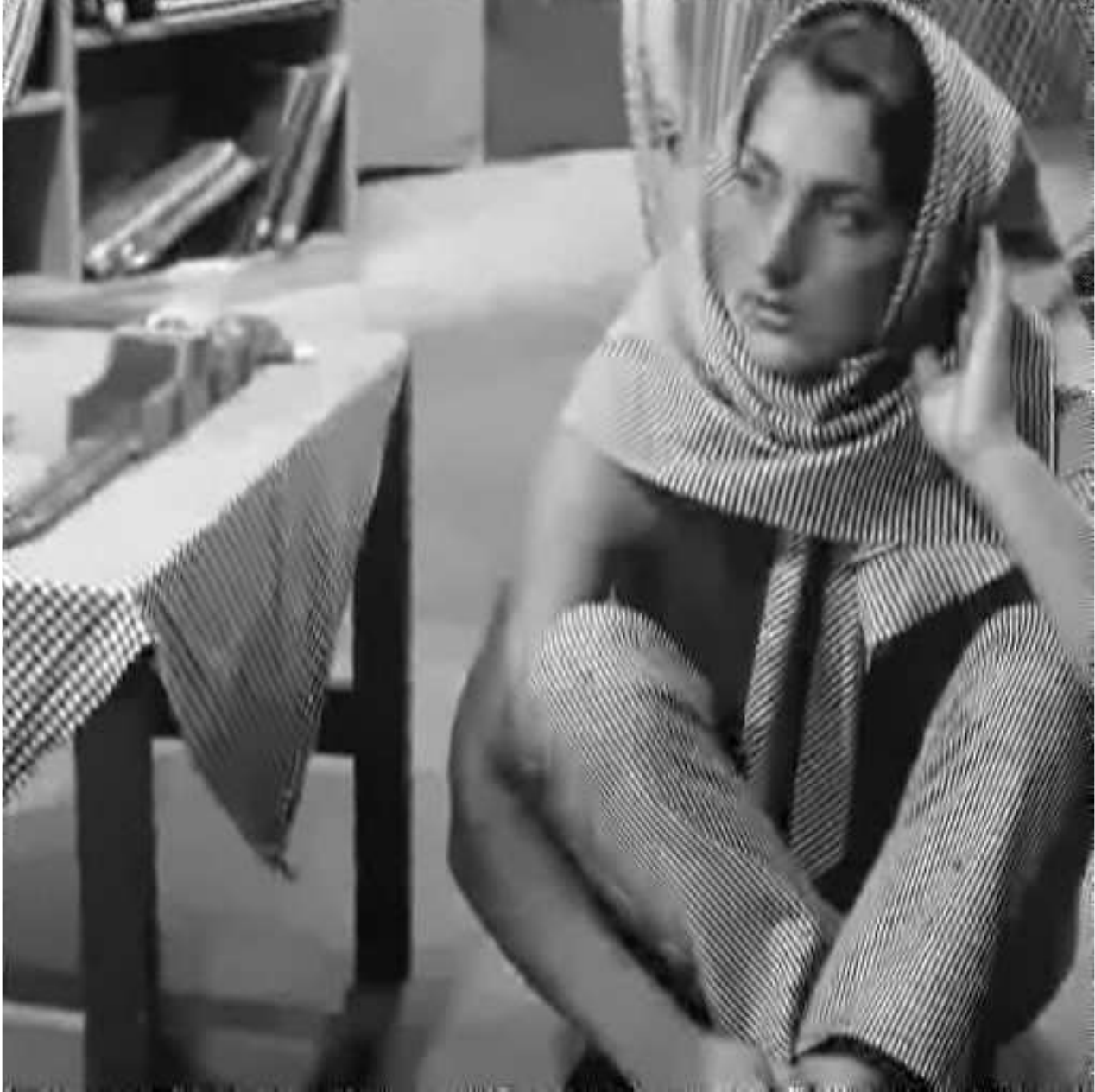}}
\subfigure[]{\includegraphics[width=.11\textwidth]{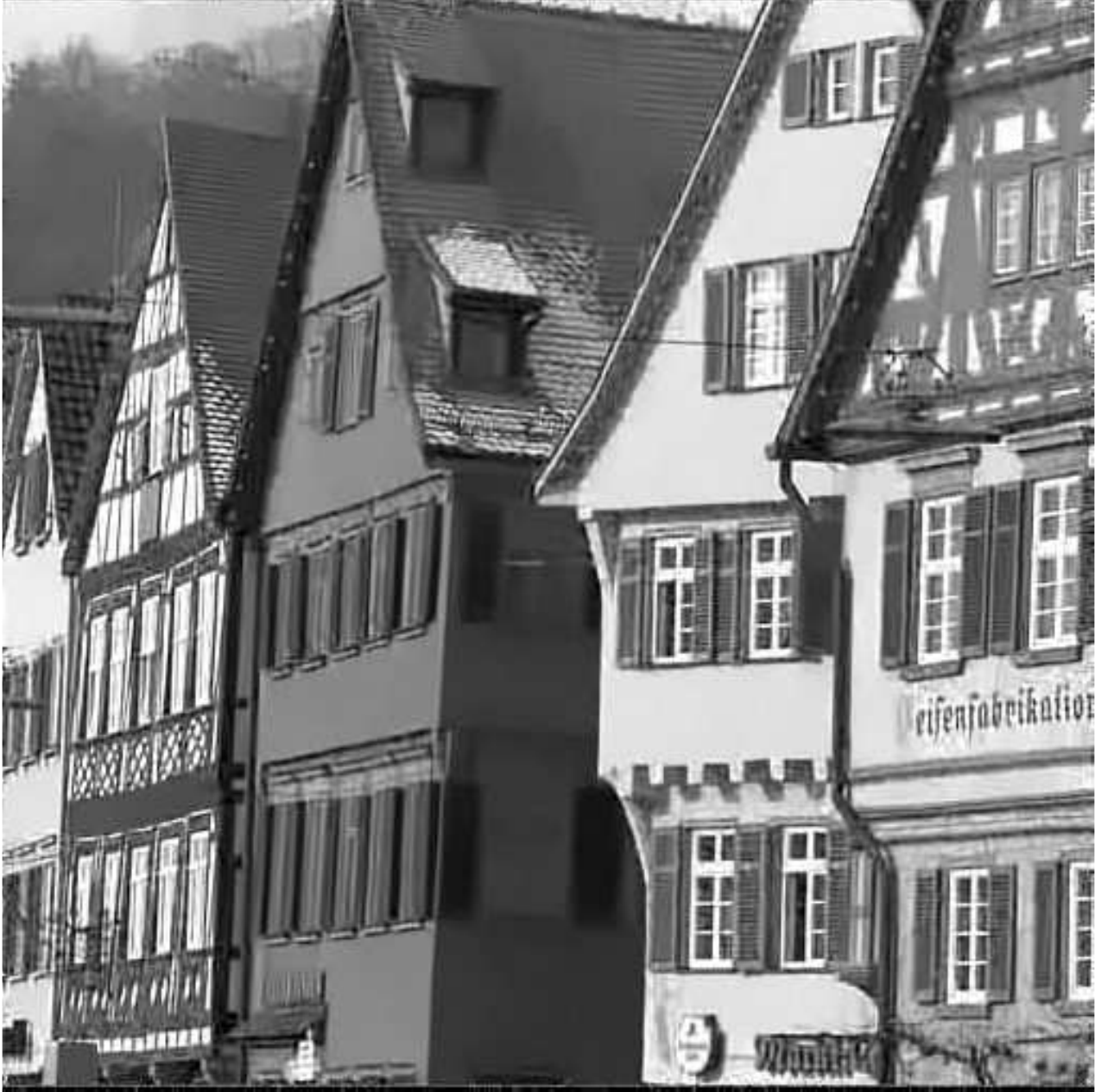}}\\\vskip -.05in
\subfigure[]{\includegraphics[width=.11\textwidth]{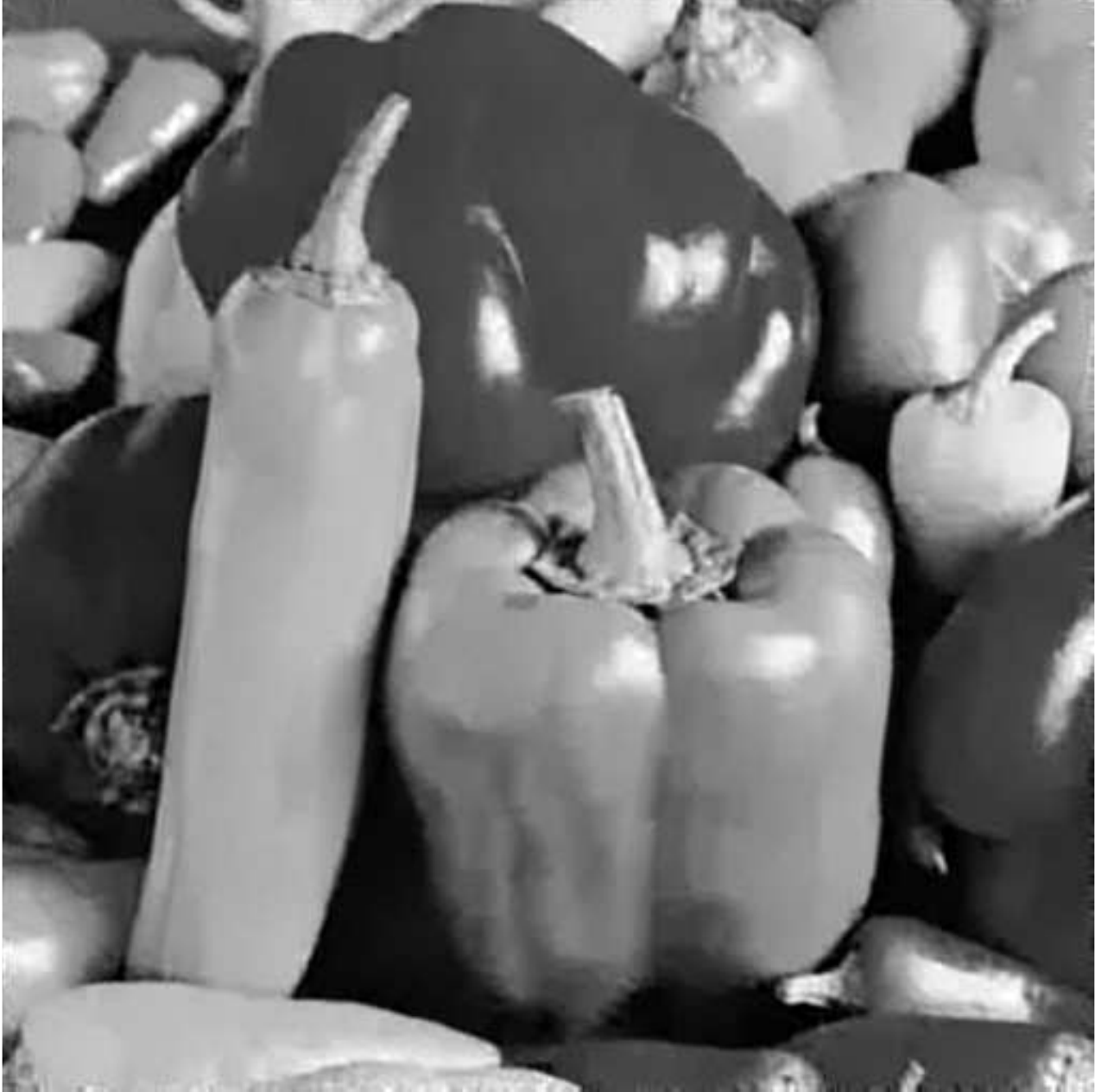}}
\subfigure[]{\includegraphics[width=.11\textwidth]{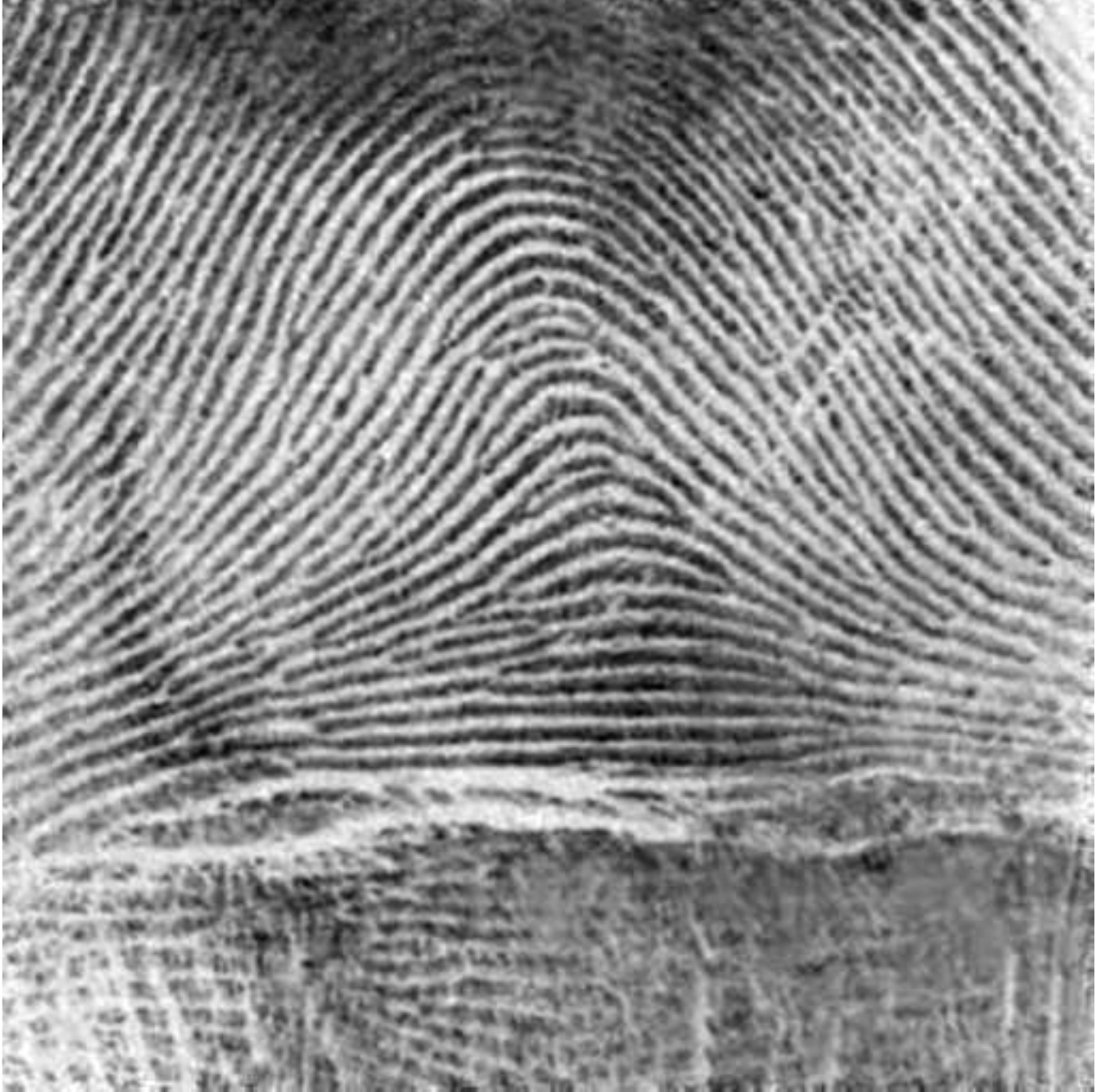}}
\subfigure[]{\includegraphics[width=.11\textwidth]{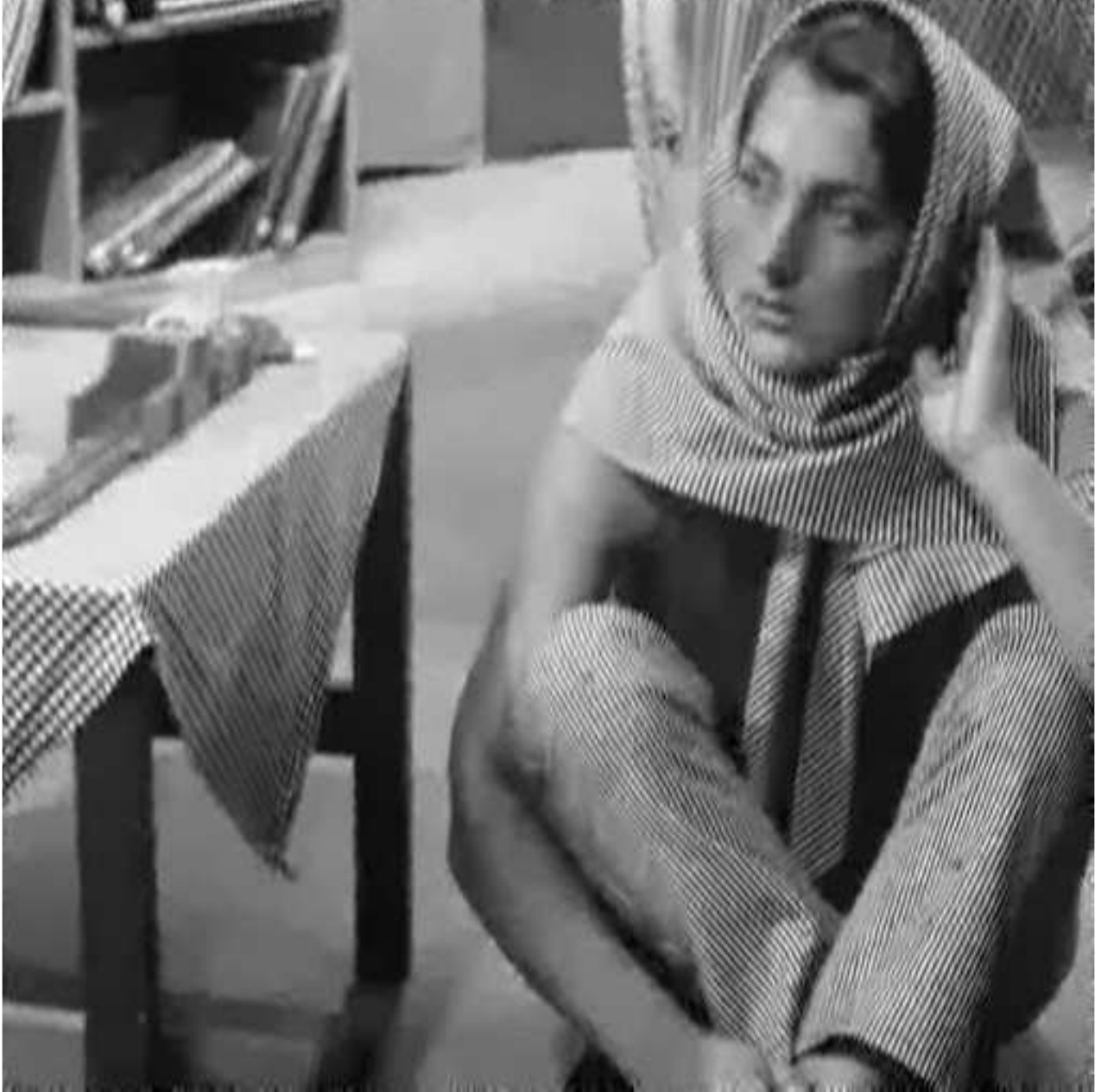}}
\subfigure[]{\includegraphics[width=.11\textwidth]{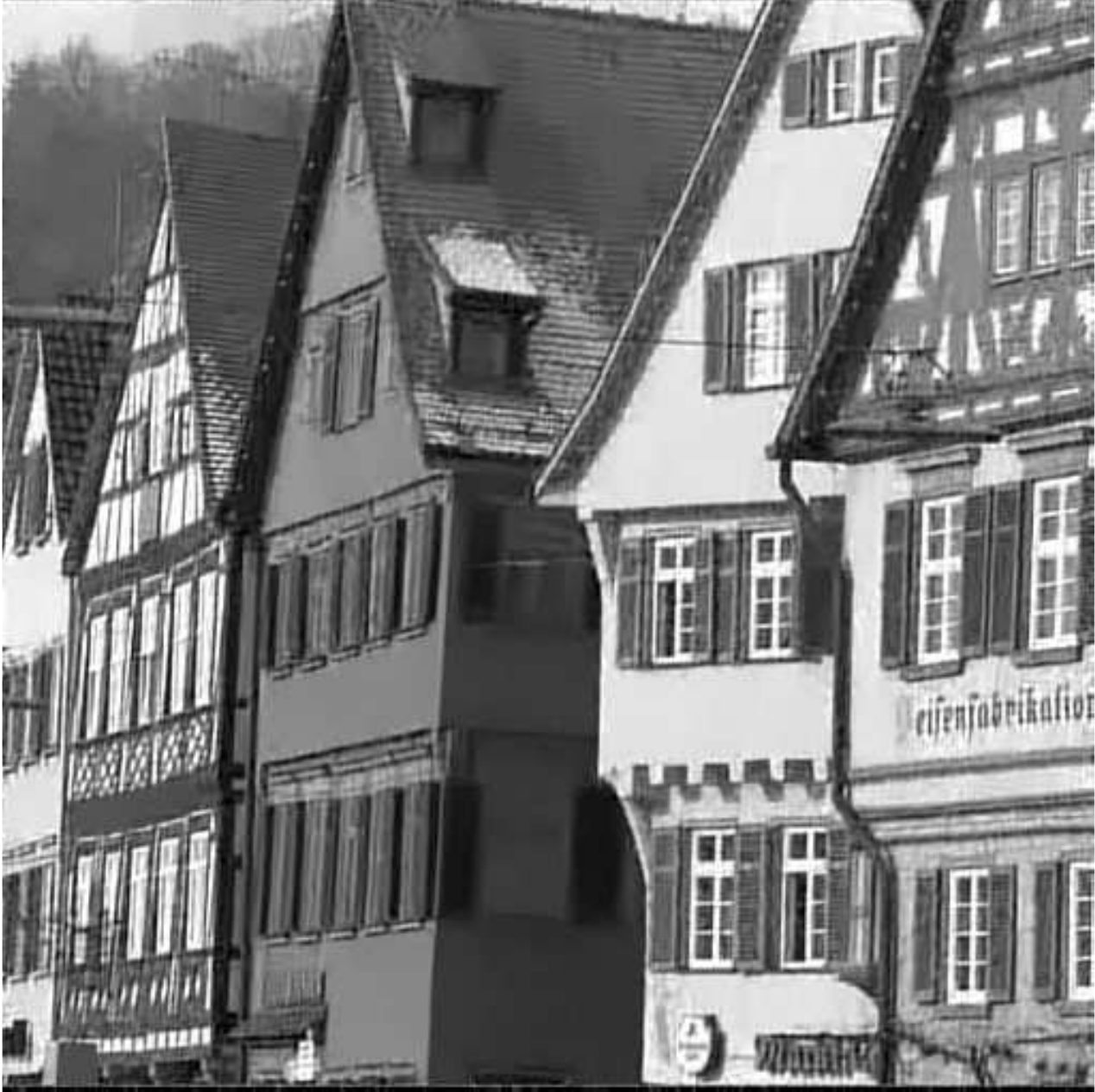}}\vskip -.05in
\end{center}
\vskip -.1in
\caption{CDP with $\delta=1.0\times 10^{-2}$. First row: PR; Second row: TVPR; Third row: ALGI; Fourth row: ALGII}
\label{cdp2}
\vskip -.1in
\end{figure}

\begin{figure}
\vskip -.25in
\begin{center}
\subfigure[]{\includegraphics[width=.11\textwidth]{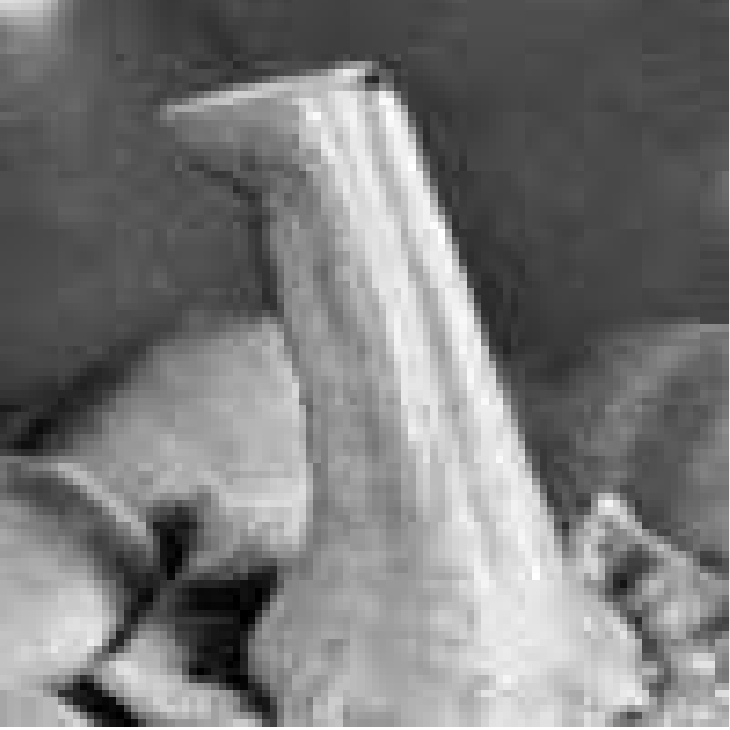}}
\subfigure[]{\includegraphics[width=.11\textwidth]{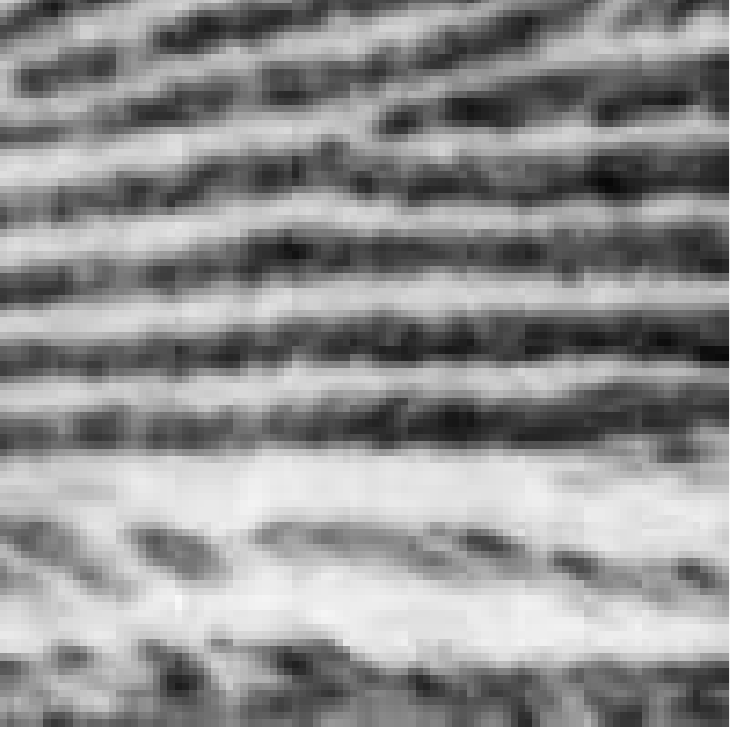}}
\subfigure[]{\includegraphics[width=.11\textwidth]{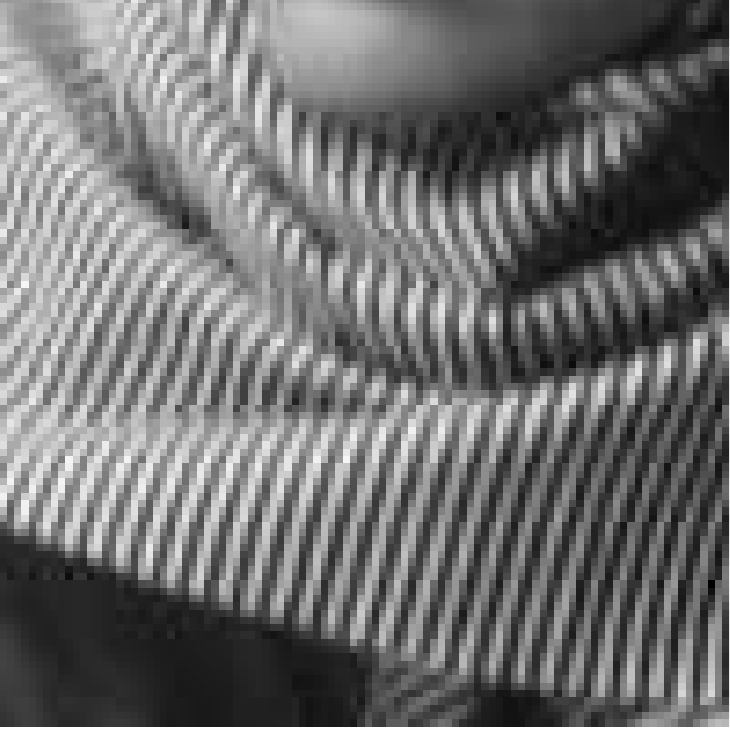}}
\subfigure[]{\includegraphics[width=.11\textwidth]{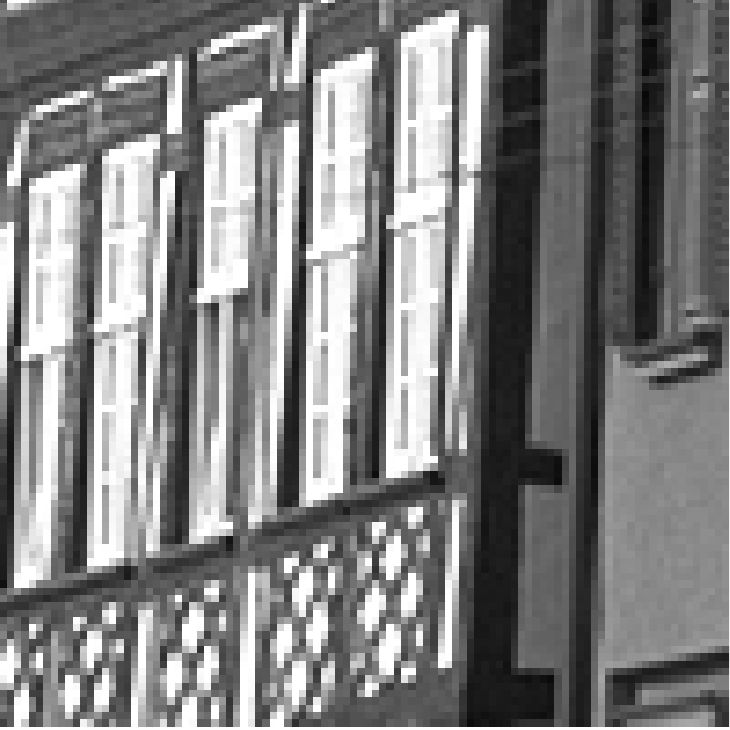}}\\\vskip -.05in
\subfigure[]{\includegraphics[width=.11\textwidth]{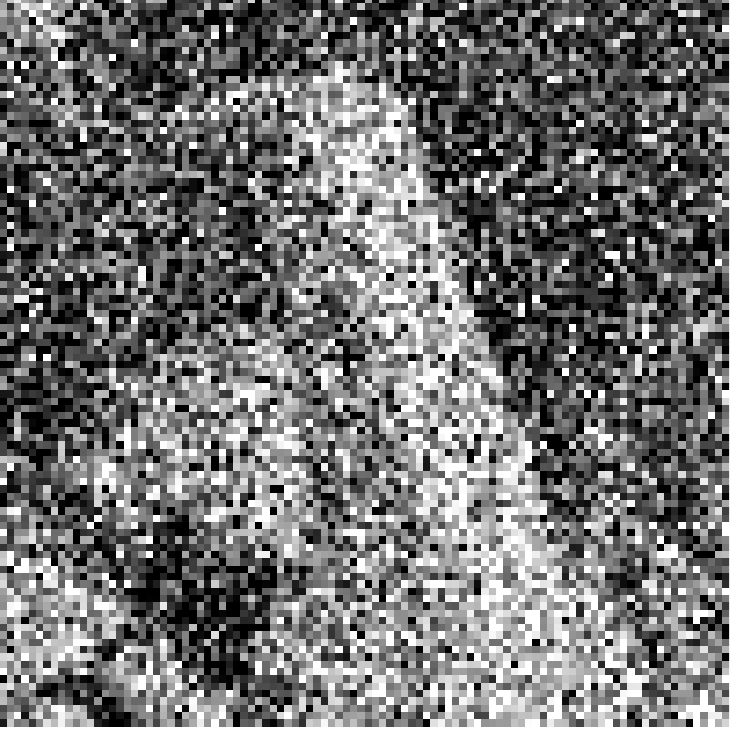}}
\subfigure[]{\includegraphics[width=.11\textwidth]{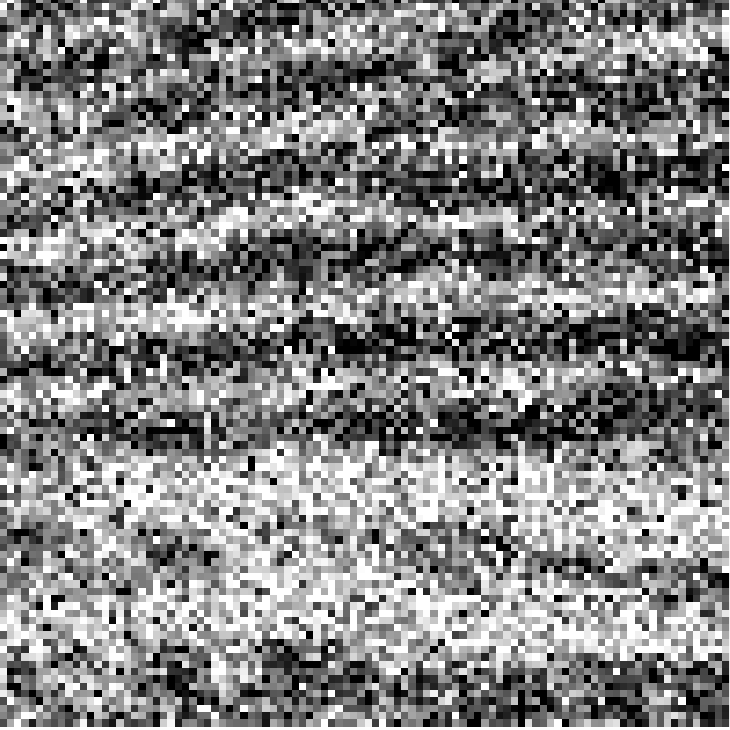}}
\subfigure[]{\includegraphics[width=.11\textwidth]{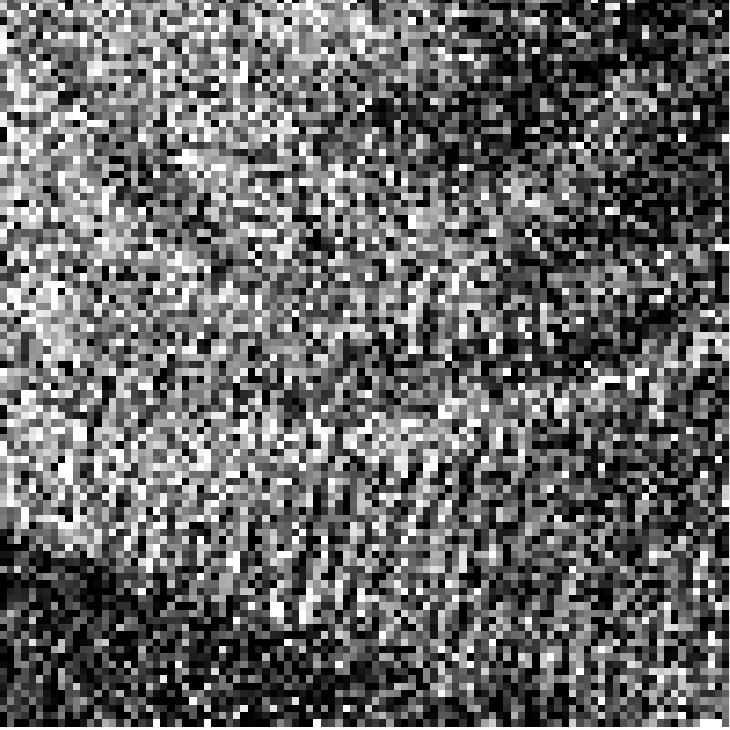}}
\subfigure[]{\includegraphics[width=.11\textwidth]{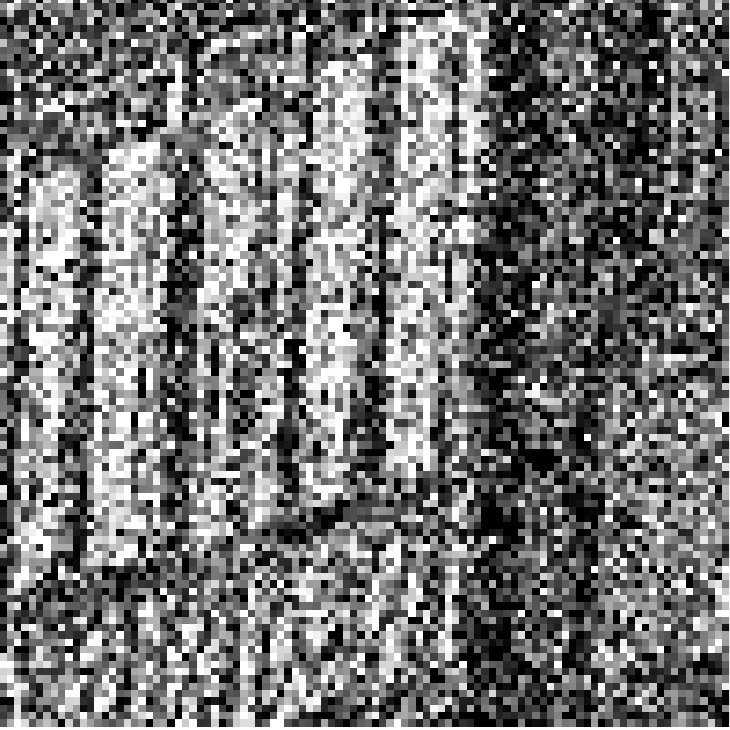}}\\\vskip -.05in
\subfigure[]{\includegraphics[width=.11\textwidth]{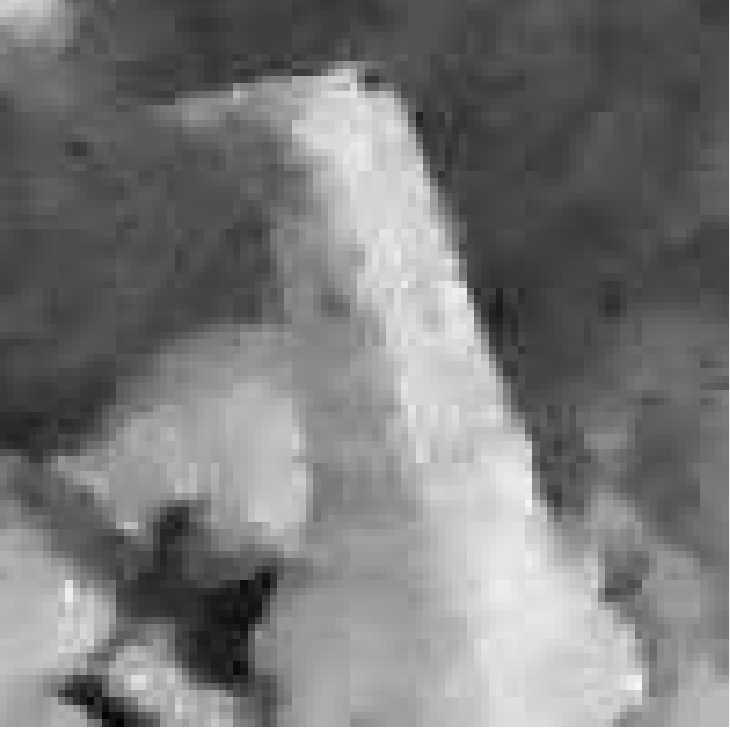}}
\subfigure[]{\includegraphics[width=.11\textwidth]{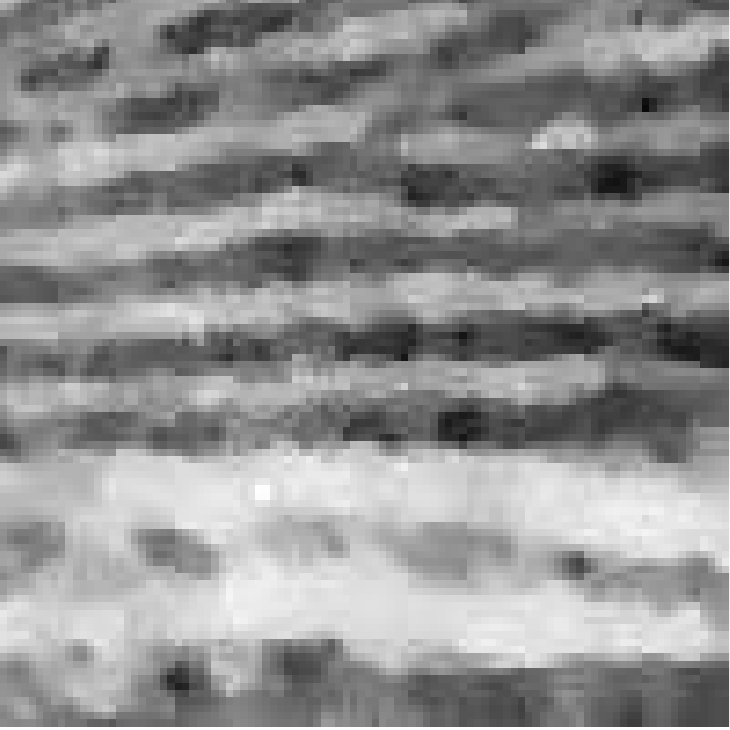}}
\subfigure[]{\includegraphics[width=.11\textwidth]{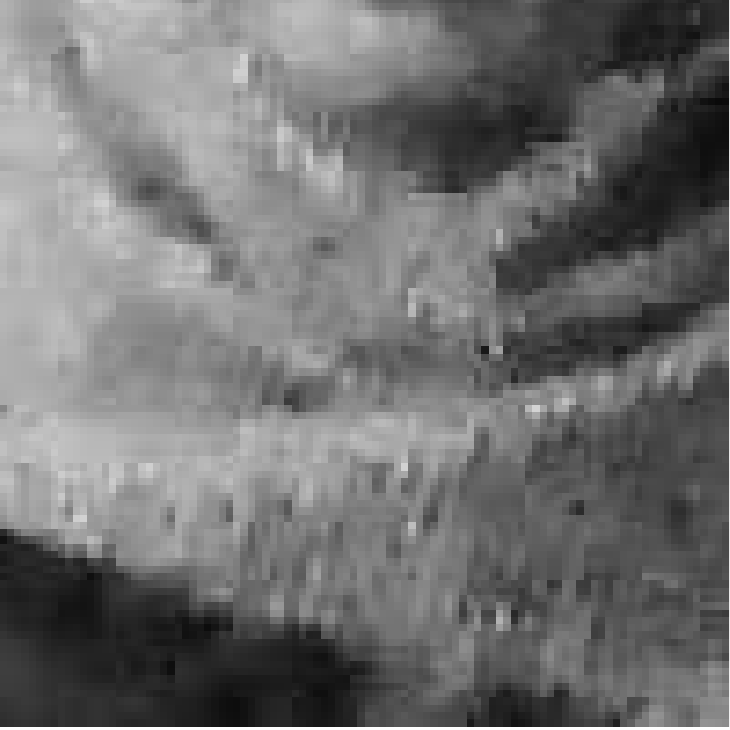}}
\subfigure[]{\includegraphics[width=.11\textwidth]{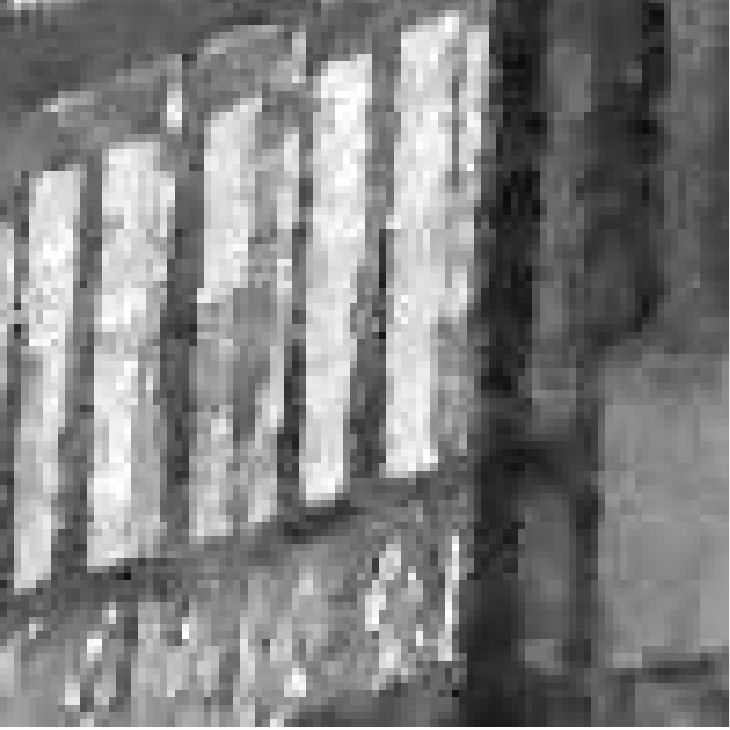}}\\\vskip -.05in
\subfigure[]{\includegraphics[width=.11\textwidth]{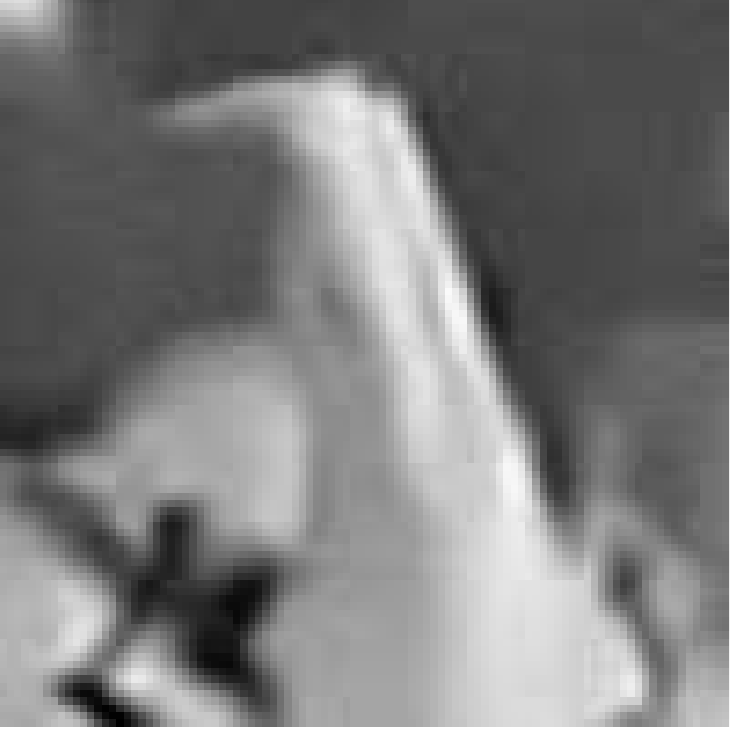}}
\subfigure[]{\includegraphics[width=.11\textwidth]{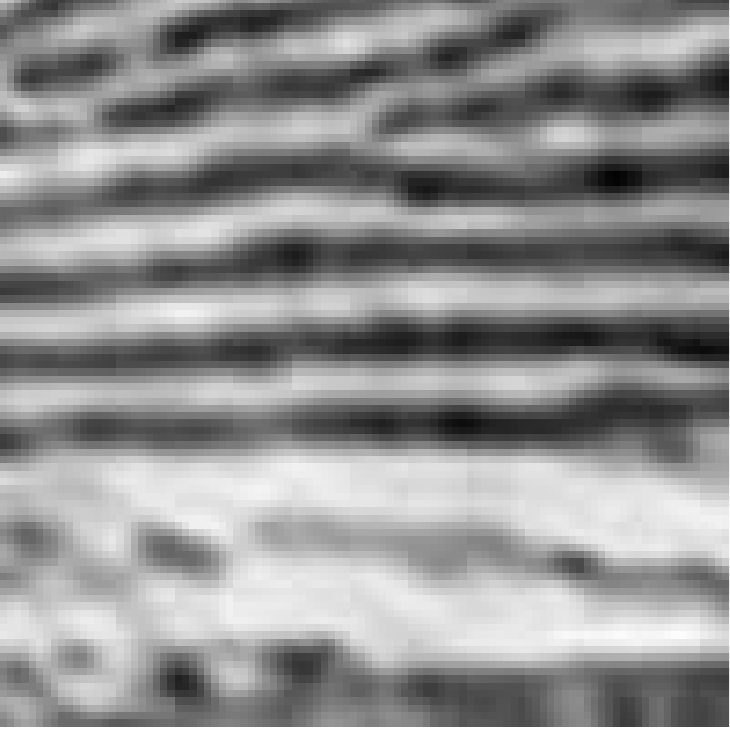}}
\subfigure[]{\includegraphics[width=.11\textwidth]{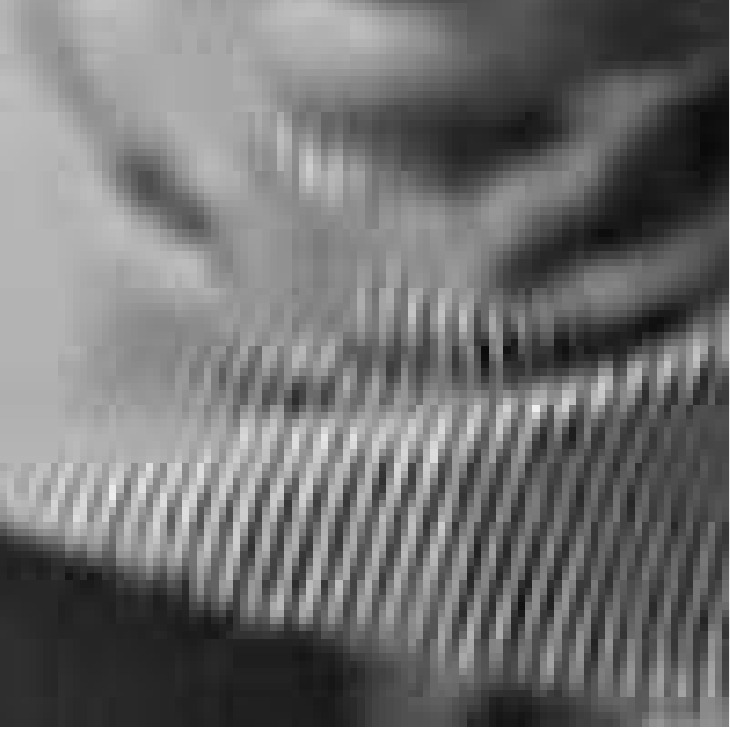}}
\subfigure[]{\includegraphics[width=.11\textwidth]{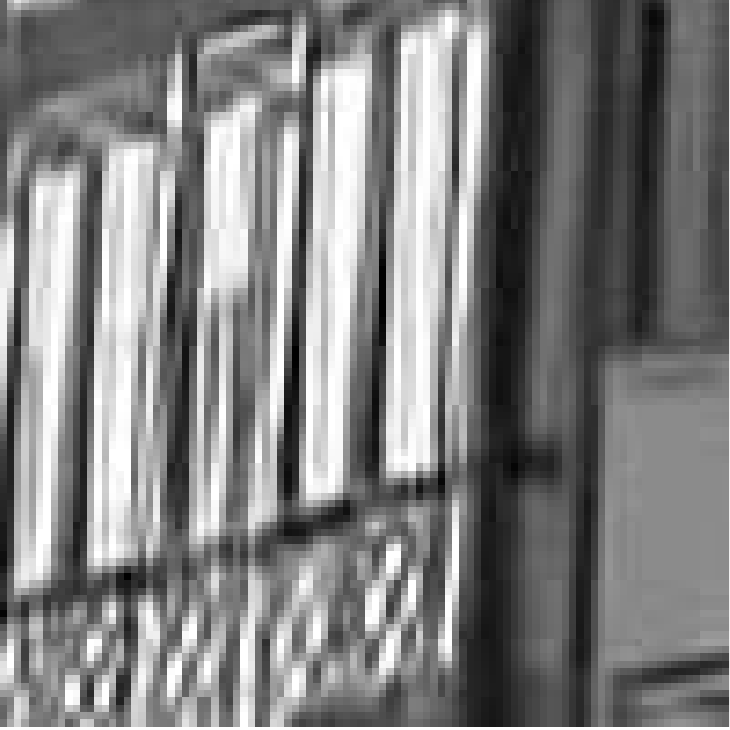}}\\\vskip -.05in
\subfigure[]{\includegraphics[width=.11\textwidth]{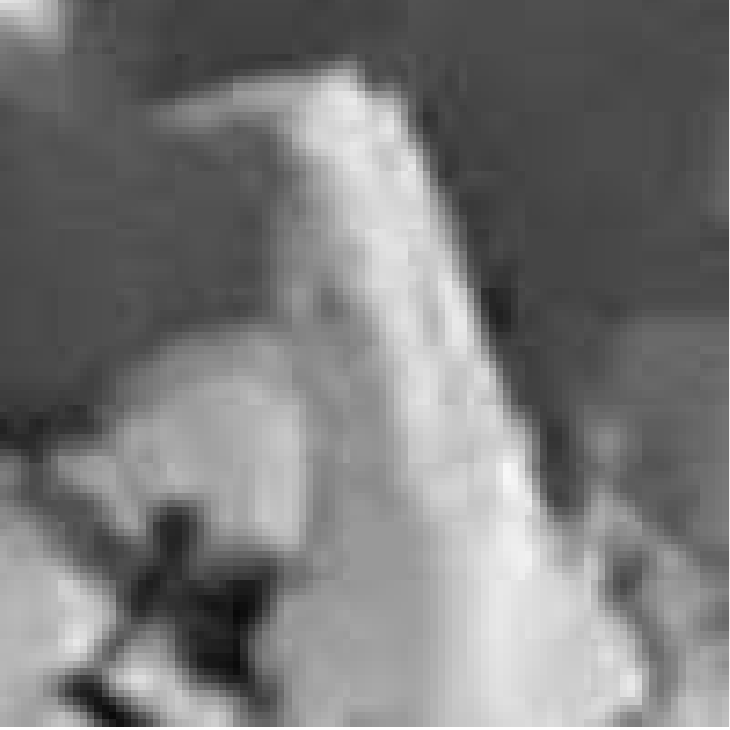}}
\subfigure[]{\includegraphics[width=.11\textwidth]{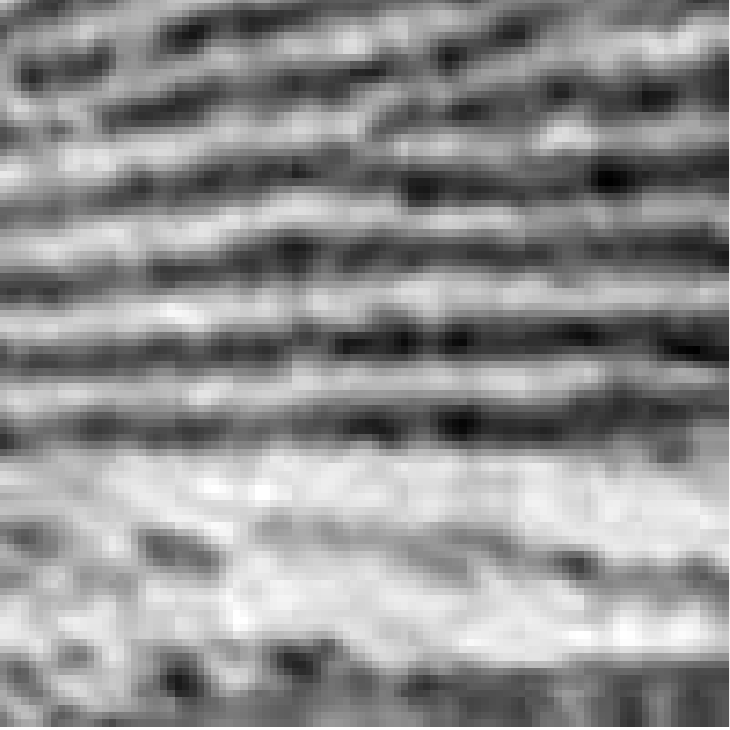}}
\subfigure[]{\includegraphics[width=.11\textwidth]{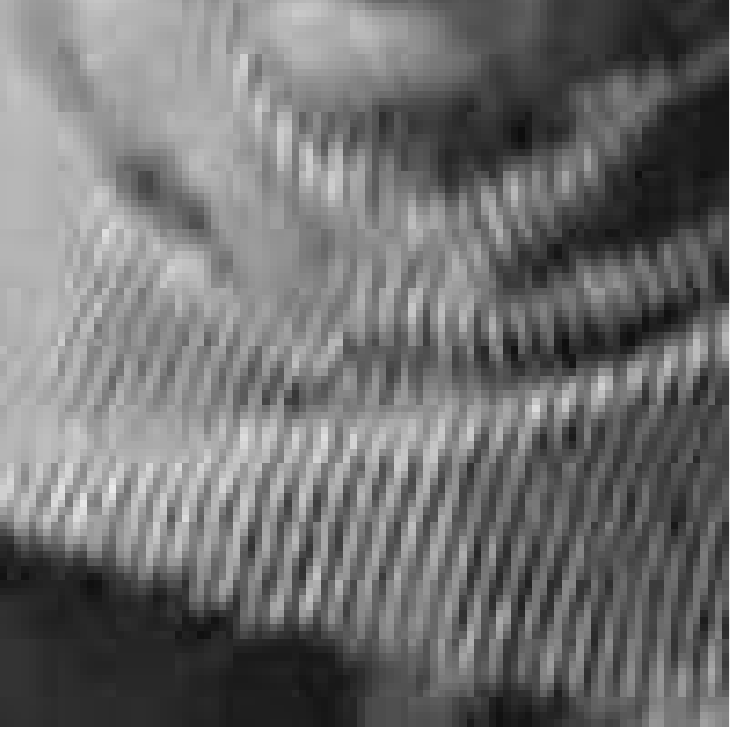}}
\subfigure[]{\includegraphics[width=.11\textwidth]{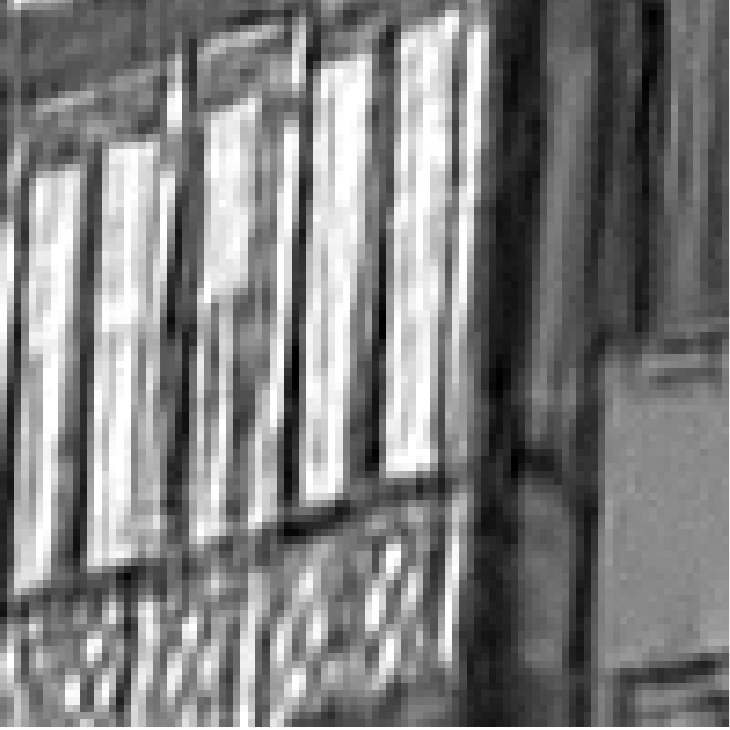}}\vskip -.05in
\end{center}
\vskip -.1in
\caption{Zoom-in images of Fig. \ref{cdp1} for CDP with $\delta=5.0\times 10^{-3}$. First row: Ground truth; Second row: PR; Third row: TVPR; Fourth row: ALGI; Fifth row: ALGII}
\label{cdp2-1}
\vskip -.2in
\end{figure}

\begin{figure}
\vskip -.4in
\begin{center}
\subfigure[]{\includegraphics[width=.11\textwidth]{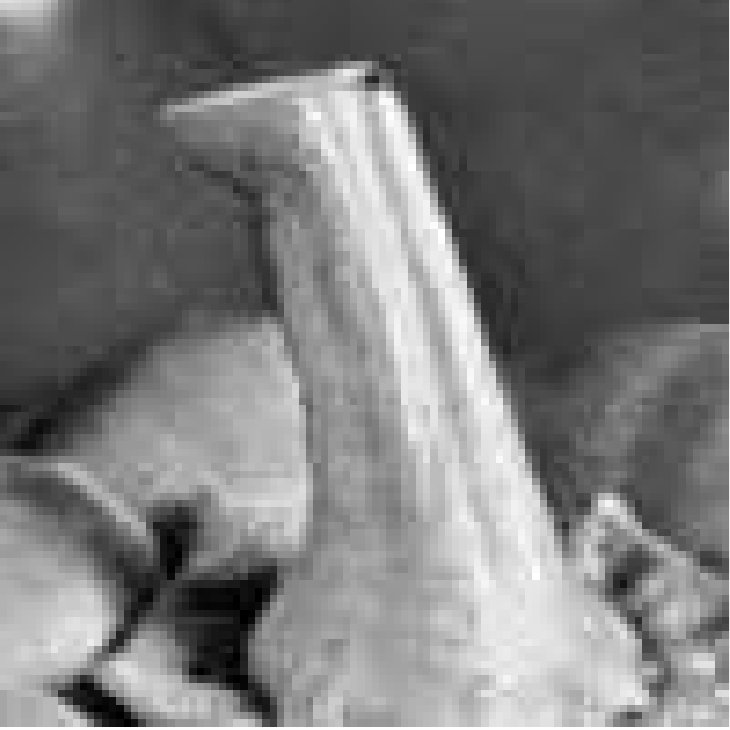}}
\subfigure[]{\includegraphics[width=.11\textwidth]{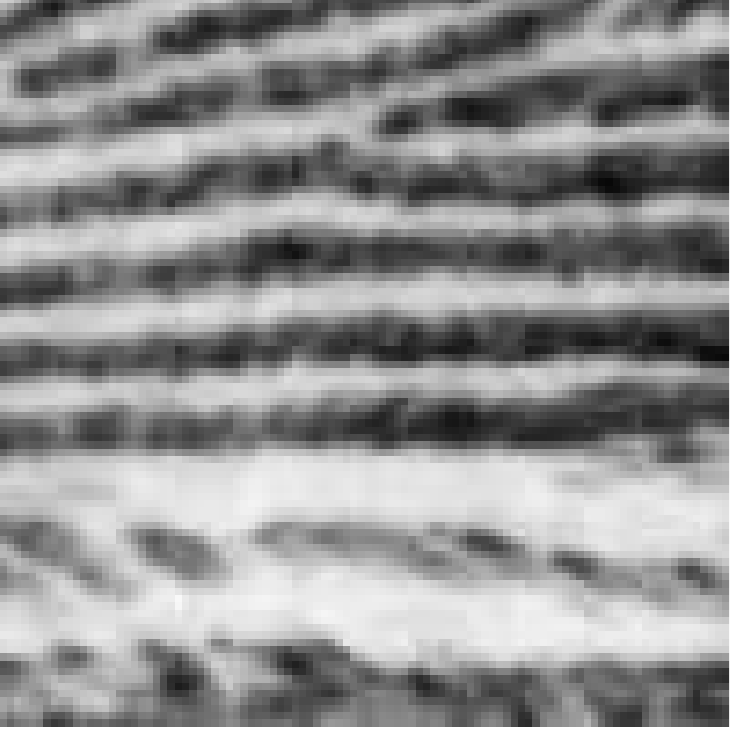}}
\subfigure[]{\includegraphics[width=.11\textwidth]{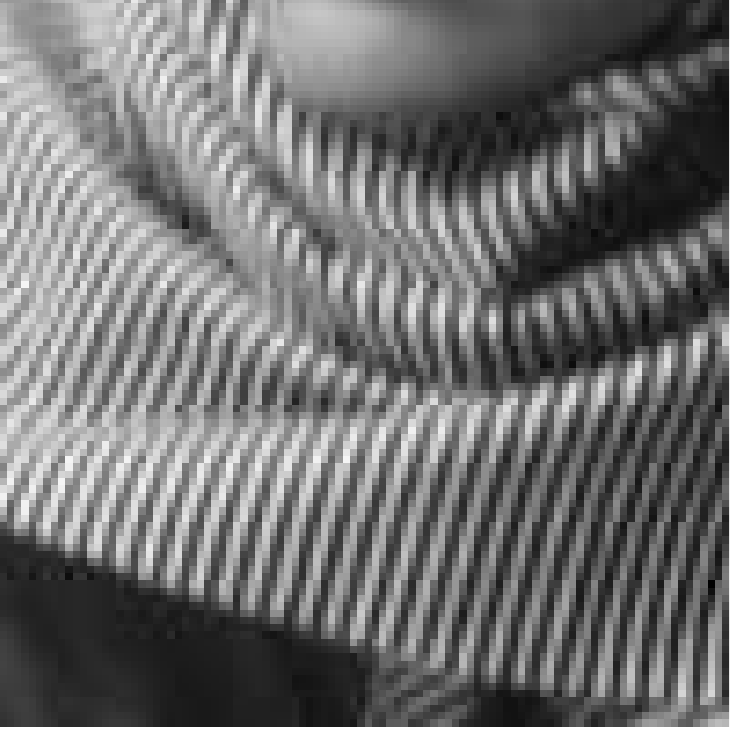}}
\subfigure[]{\includegraphics[width=.11\textwidth]{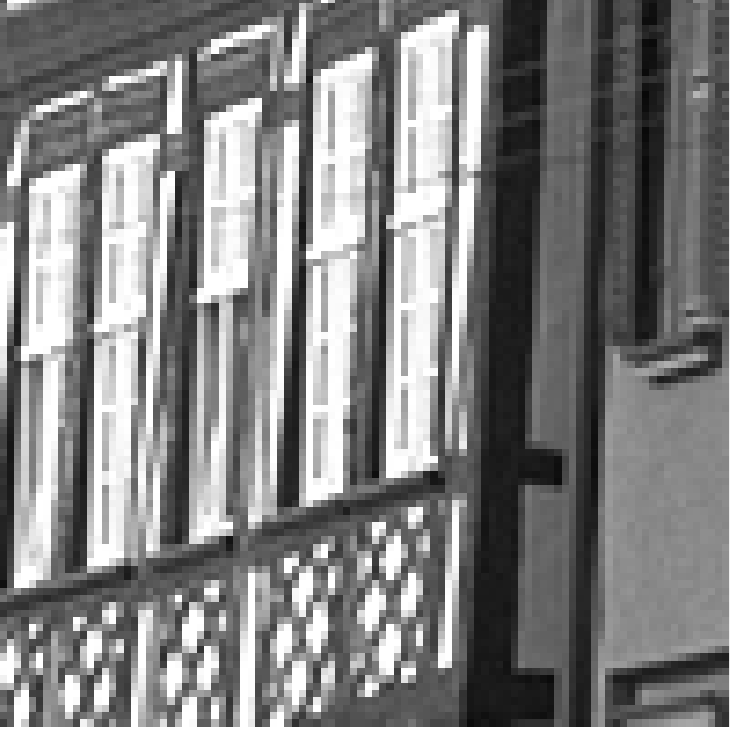}}\\\vskip -.05in
\subfigure[]{\includegraphics[width=.11\textwidth]{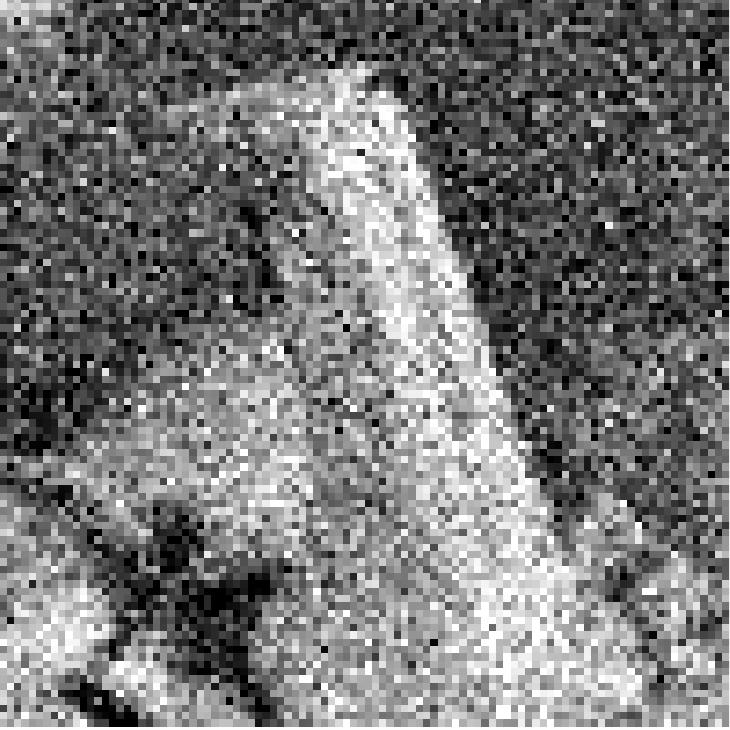}}
\subfigure[]{\includegraphics[width=.11\textwidth]{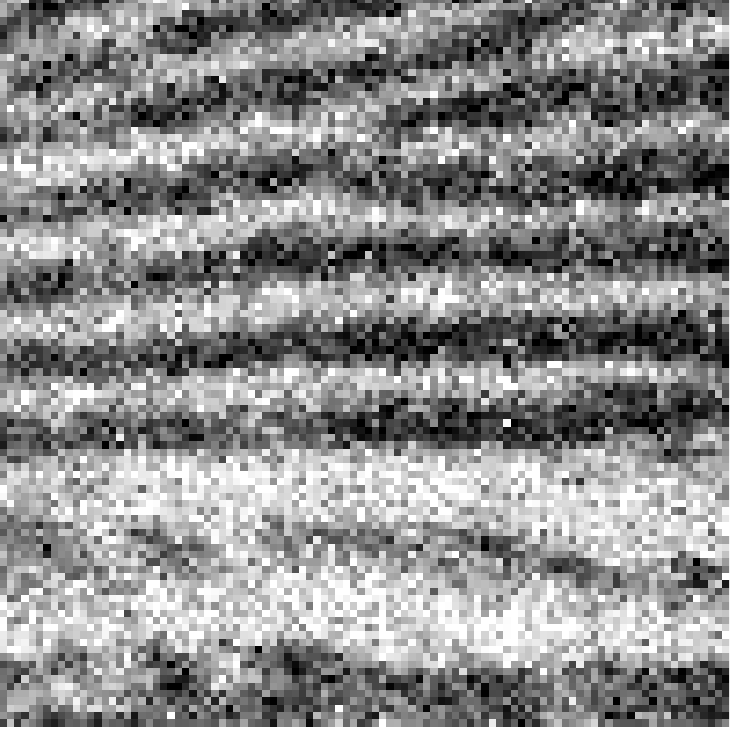}}
\subfigure[]{\includegraphics[width=.11\textwidth]{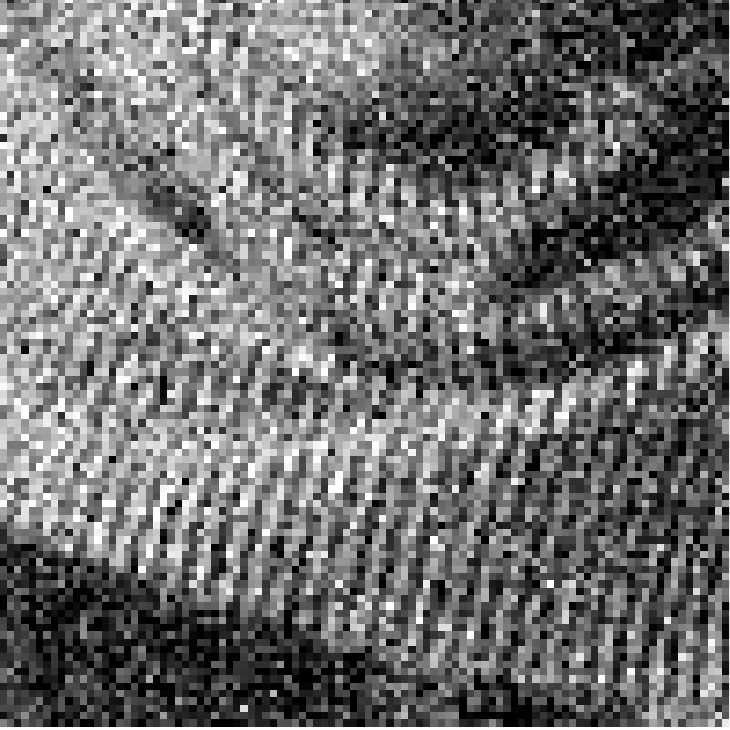}}
\subfigure[]{\includegraphics[width=.11\textwidth]{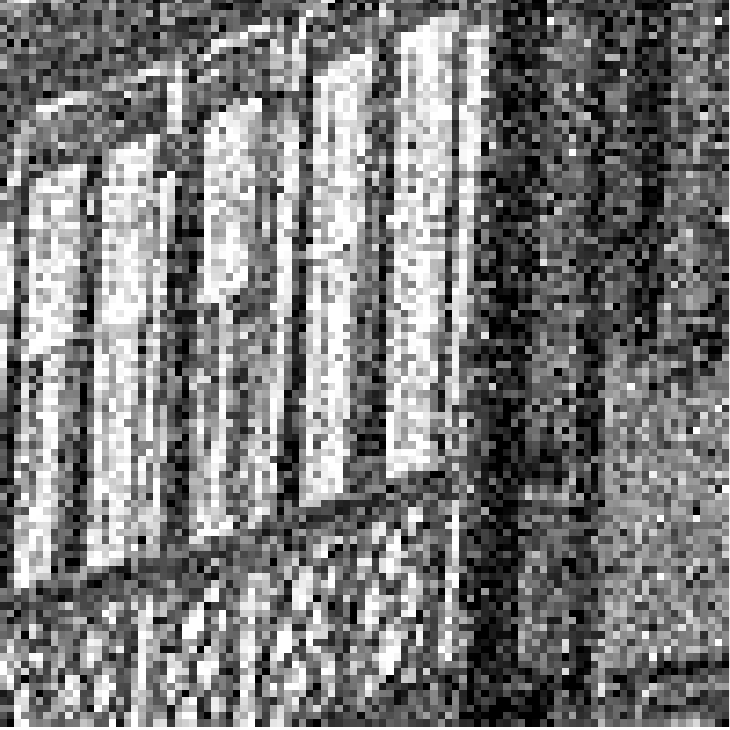}}\\ \vskip -.05in
\subfigure[]{\includegraphics[width=.11\textwidth]{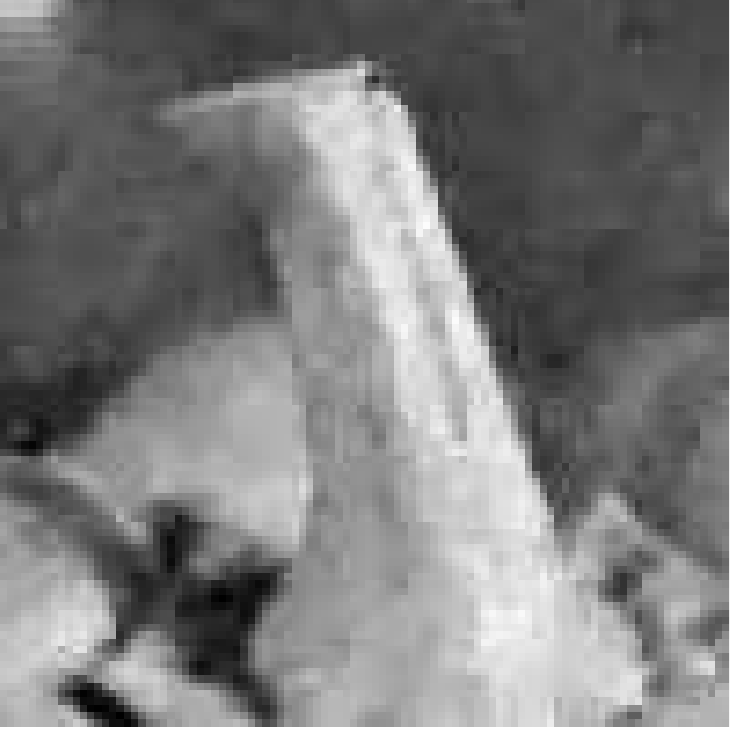}}
\subfigure[]{\includegraphics[width=.11\textwidth]{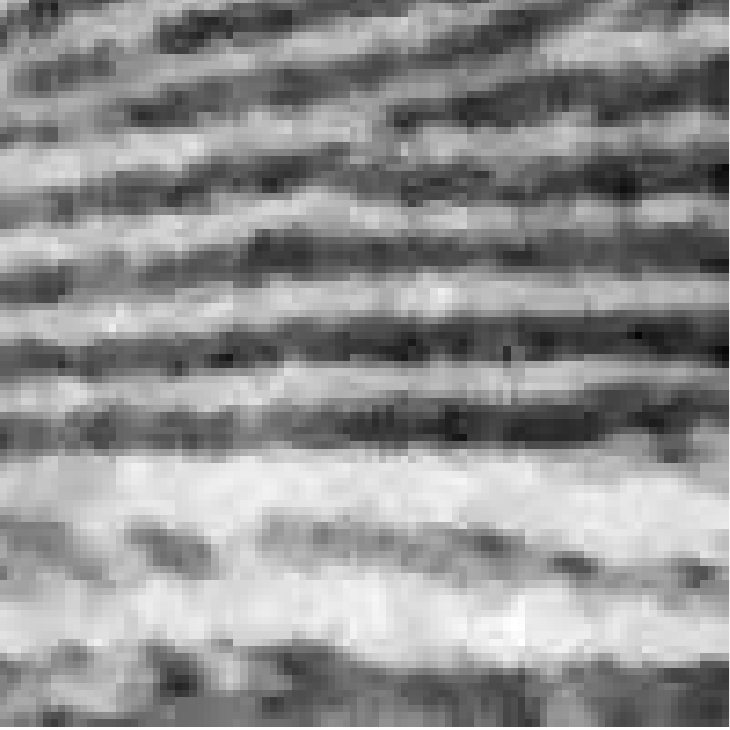}}
\subfigure[]{\includegraphics[width=.11\textwidth]{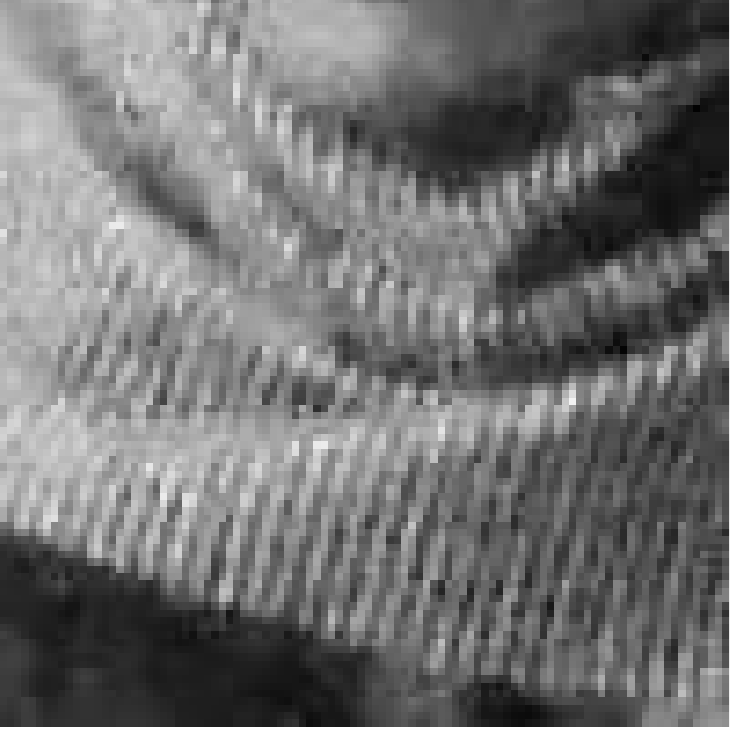}}
\subfigure[]{\includegraphics[width=.11\textwidth]{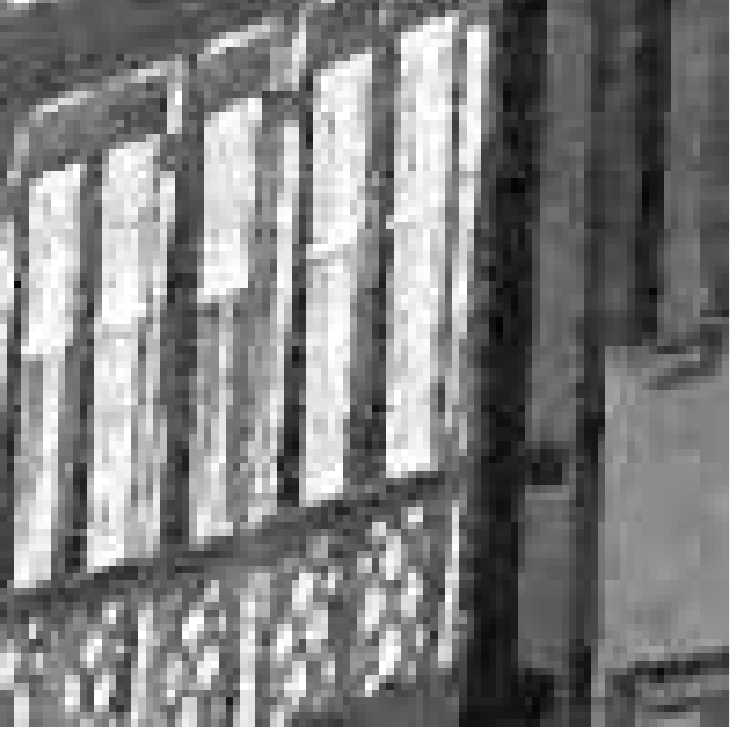}}\\\vskip -.05in
\subfigure[]{\includegraphics[width=.11\textwidth]{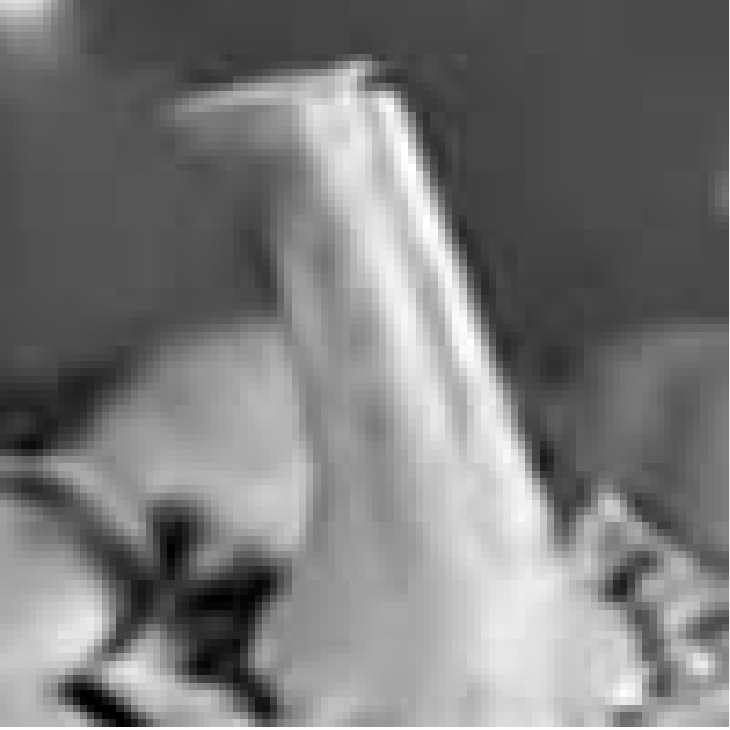}}
\subfigure[]{\includegraphics[width=.11\textwidth]{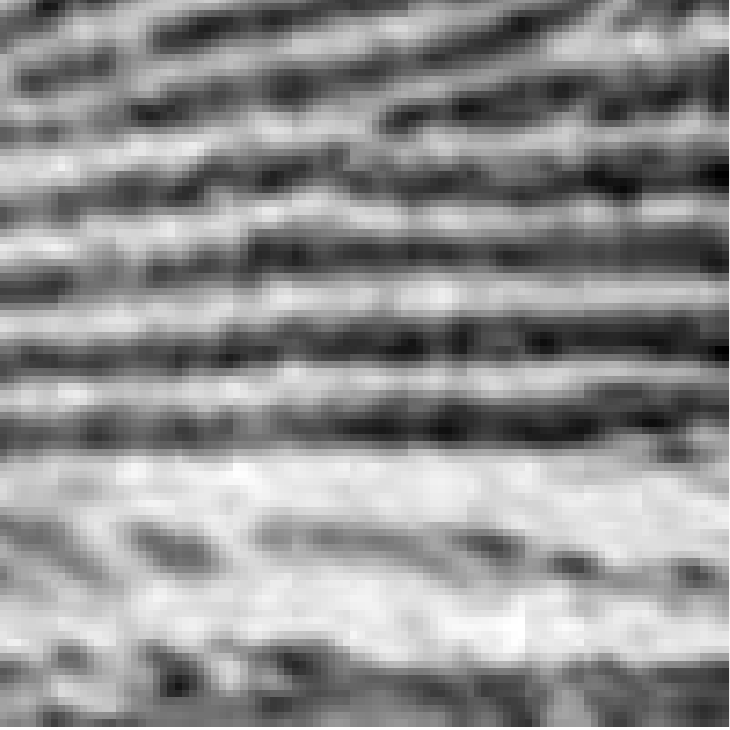}}
\subfigure[]{\includegraphics[width=.11\textwidth]{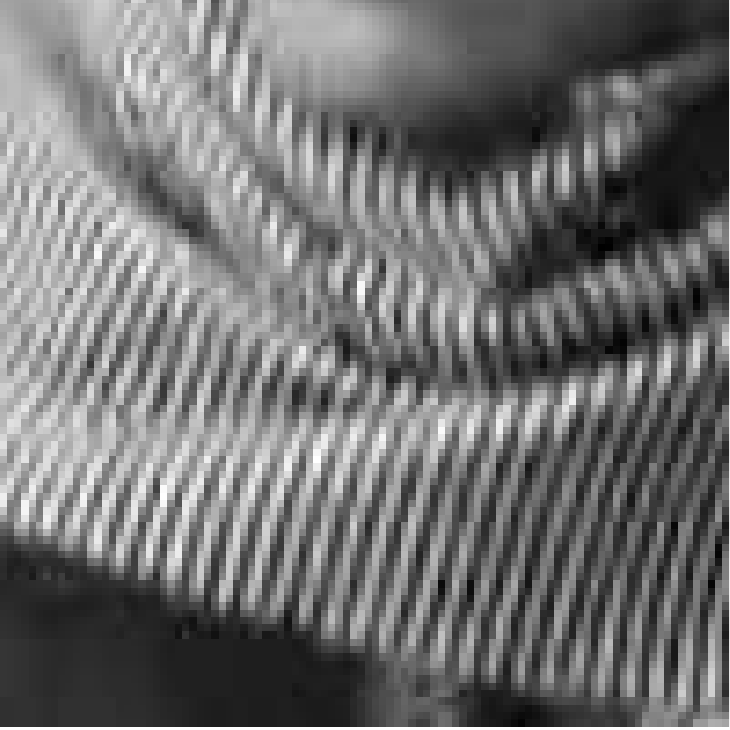}}
\subfigure[]{\includegraphics[width=.11\textwidth]{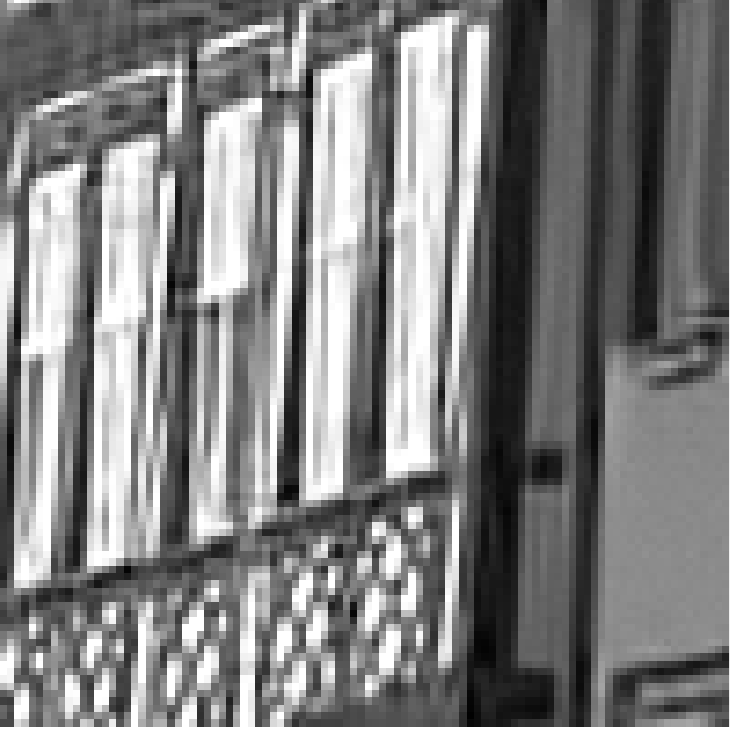}}\\\vskip -.05in
\subfigure[]{\includegraphics[width=.11\textwidth]{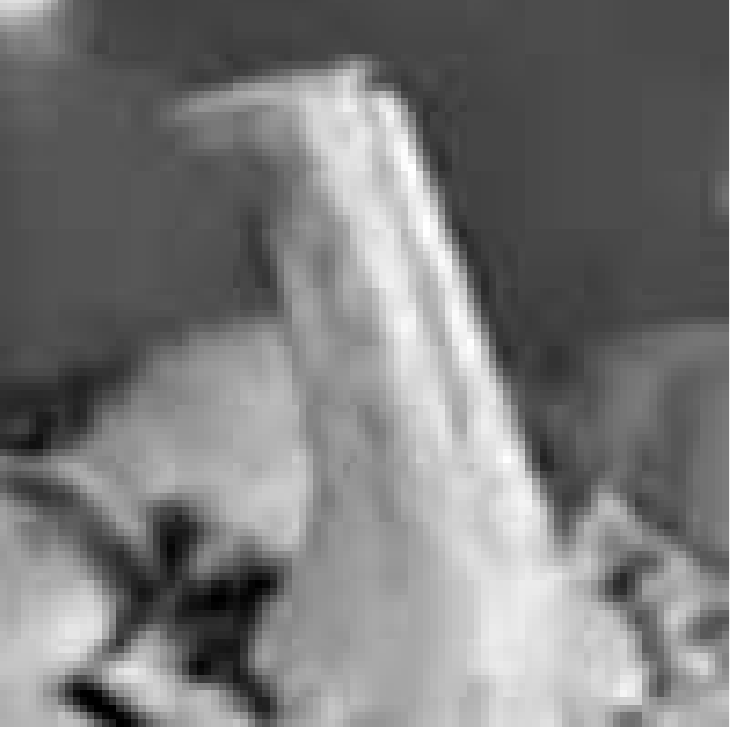}}
\subfigure[]{\includegraphics[width=.11\textwidth]{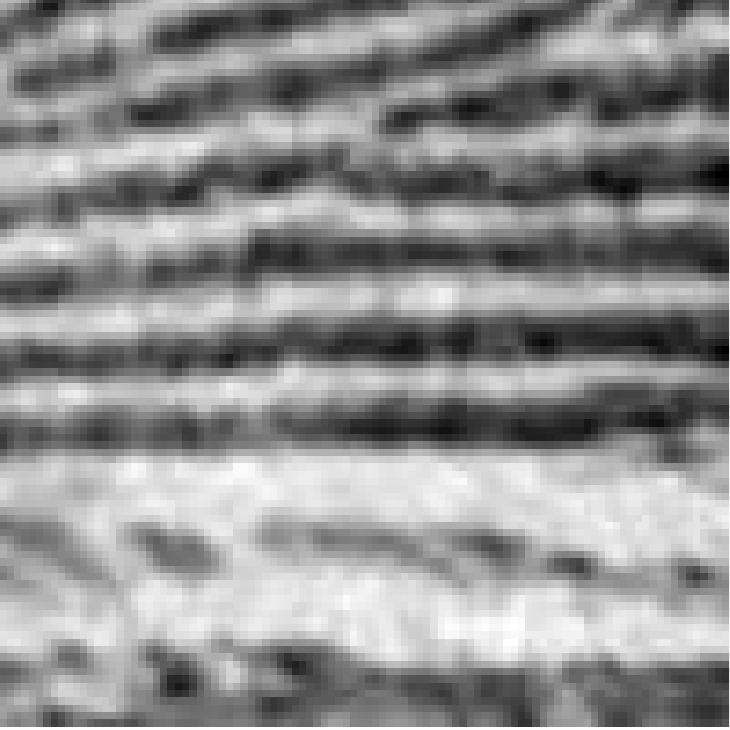}}
\subfigure[]{\includegraphics[width=.11\textwidth]{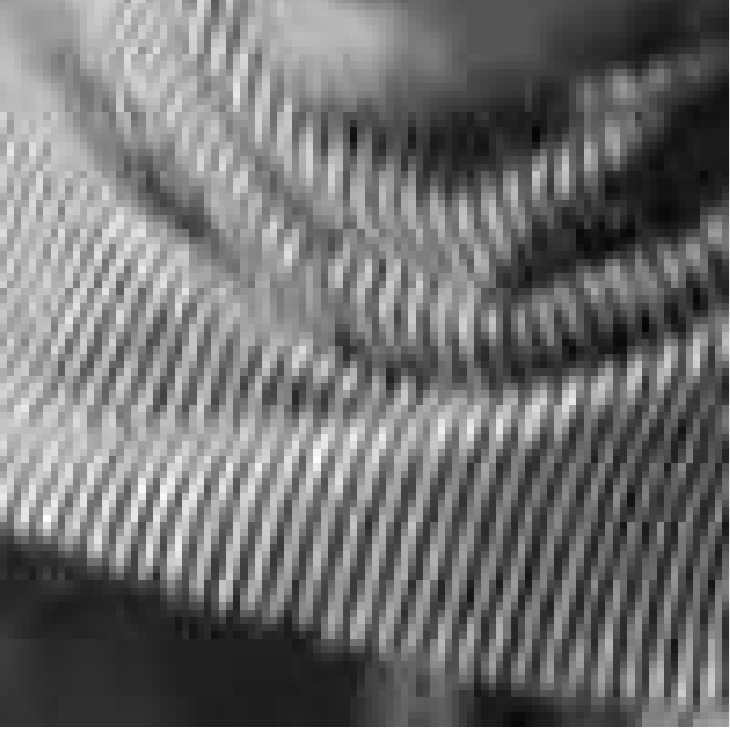}}
\subfigure[]{\includegraphics[width=.11\textwidth]{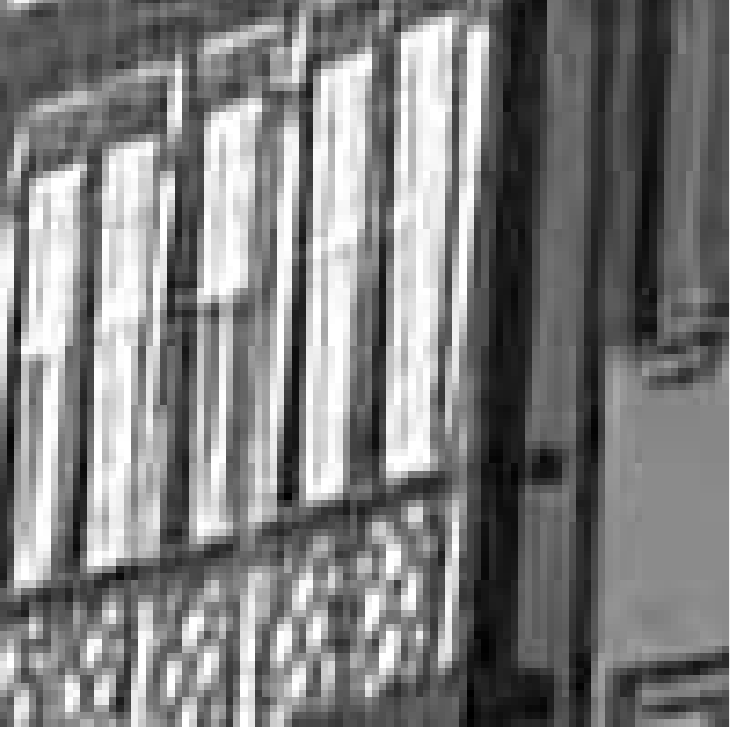}}\vskip -.05in
\end{center}
\vskip -.1in
\caption{Zoom-in images of Fig. \ref{cdp2} for CDP with $\delta=1.0\times 10^{-2}$. First row: Ground truth; Second row: PR; Third row: TVPR; Fourth row: ALGI; Fifth row: ALGII.}
\label{cdp2-2}
\vskip -.2in
\end{figure}

 Reconstructed results with different  noise levels are put in Fig. \ref{cdp1}-Fig. \ref{cdp2}, and  zoom-in parts of corresponding results are put in Fig. \ref{cdp2-1}-Fig. \ref{cdp2-2}. Readily one can see that ``PR'' generates very noisy results by observing the first rows of Fig. \ref{cdp1}-Fig. \ref{cdp2} and the second rows in the zoom-in results of Fig. \ref{cdp2-1}-Fig. \ref{cdp2-2}, where the edges, background and the repetitive structures are contaminated severely. By TV regularization method and our proposed dictionary learning methods, all the recovery images have sharper edges and cleaner background. However, by observing the zoom-in results in the third  rows of Fig. \ref{cdp2-1}-Fig. \ref{cdp2-2}, ``TVPR'' generates images with visible stair case artifacts, and some important texture information can not be kept, while ``ALGI'' and ``ALGII'' can both produce  high quality images and  one does not find any visible staircase artifacts in the zoom-in parts. By observing Fig. \ref{cdp2-1}(m), (q)  or Fig. \ref{cdp2-2}(m), (q), ``ALGI'' and ``ALGII'' generate  cleaner background than ``TVPR''. Moreover, by  Fig. \ref{cdp2-1}(n)-(p), (r)-(t)  or Fig. \ref{cdp2-2}(n)-(p), (r)-(t),  the textures are preserved pretty well compared with the results by ``TVPR''.  In order to qualify the improvement of our proposed methods, the corresponding SNRs are put in Table \ref{tab1-1} and Table \ref{tab1-2}, where the maximum in  the same row is marked in bold font. Inferred from these two tables, SNRs by ``TVPR'', ``ALGI'' and ``ALGII'' are about double of those by ``PR'', and SNRs by ``ALGI'' and ``ALGII'' increase averagely about 1.5dB and 2dB, respectively. The increase or improvement with the proposed methods compared with ``PR'' and ``TVPR'' is more obvious  on the image with more texture information \emph{e.g.} ``Fingerprint'' or  ``Barbara''  than  that  \emph{e.g.} ``Peppers'' with piecewise smooth features, which demonstrates that its advantage is to handle images with sophisticated texture features.

By observing the zoom-in parts, the results in the fourth rows of Fig. \ref{cdp2-1} and Fig. \ref{cdp2-2} seem more smoothing by ``ALGI'' than those in the fifth rows by ``ALGII''. On the other hand, the results by ``ALGII'' contain a bit more features than by those by ``ALGI''. In order to investigate the performance differences of ``ALGI'' and ''ALGII'', the sparsity of the coefficient matrix $\bm \alpha$ is put in Table \ref{tabSparse}, and one can readily see that the coefficient matrix is sparser by ``ALGI'' than by ``ALGII'', which possibly leads to those above behaviours. Anyway, both two algorithms produce comparable results with cleaner background and well-preserved textures. Hereafter, we only provide the results by ``ALGI'' since it has fewer parameters and seems more suitable for practical use.

\begin{table}
\begin{center}
\begin{spacing}{.9}
\begin{tabular}{|c||c|c|c|c|}
\hline
Images& PR &TVPR &ALGI&ALGII \\\hline\hline
Peppers&4.80&20.12&20.80&\bf 21.19  \\ \hline
Fingerprint&5.96&16.43&18.32&\bf 18.46\\ \hline
Barbara& 4.68&16.27&18.04&\bf 18.05\\ \hline
House&5.14&15.04&\bf 16.71&16.61\\ \hline\hline
Average& 5.15 &16.97&18.48&\bf 18.58\\
\hline
\end{tabular}
\end{spacing}
\end{center}
\vskip -.3in
\caption{SNRs for CDP with $\delta=5.0\times 10^{-3}$.}
\label{tab1-1}
\vskip -.2in
\end{table}

\begin{table}
\begin{center}
\begin{spacing}{.9}
\begin{tabular}{|c||c|c|c|c|}
\hline
Images& PR &TVPR &ALGI&ALGII \\\hline\hline
Peppers&9.65&23.03&24.23&\bf 24.27  \\ \hline
Fingerprint&11.31&20.14&21.96&\bf 22.06\\ \hline
Barbara& 9.61&18.67&\bf 21.99&21.44\\ \hline
House&10.11&18.48&\bf 20.11&19.60\\ \hline\hline
Average& 10.17 &20.08&\bf 22.07& 21.84\\
\hline
\end{tabular}
\end{spacing}
\end{center}
\vskip -.2in
\caption{SNRs for CDP with $\delta=1.0\times 10^{-2}$.}
\vskip -.15in
\label{tab1-2}

\end{table}
\begin{table}[h!]
\vskip 0in
\begin{center}
\begin{spacing}{.9}
\begin{tabular}{|c||c |c||c|c|}
\hline
\multirow{2}{*}{Images}&  \multicolumn{2}{c||}{$\delta=5.0\times 10^{-3}$}&\multicolumn{2}{c|}{$\delta=1.0\times 10^{-2}$}\\
\cline{2-5}
& ALGI &ALGII &ALGI&ALGII \\ \hline\hline
Peppers    &3.9\%&5.4\%  &6.9\%& 8.6\%  \\ \hline
Fingerprint&9.8\%&15.3\%  &20.5\%&24.3\%\\ \hline
Barbara    &4.9\%&7.7\%  &13.2\%&14.4\%\\ \hline
House       &8.6\%&14.7\% &20.8\%&22.3\%\\
\hline
\end{tabular}
\end{spacing}
\end{center}
\vskip -.2in
\caption{Sparsity $S(\bm\alpha)$  for CDP on real-valued images.}
\label{tabSparse}
\vskip -.3in
\end{table}

\subsubsection{Complex-valued Image}
\begin{figure}
\vskip -.in
\begin{center}
\subfigure[]{\includegraphics[width=.11\textwidth]{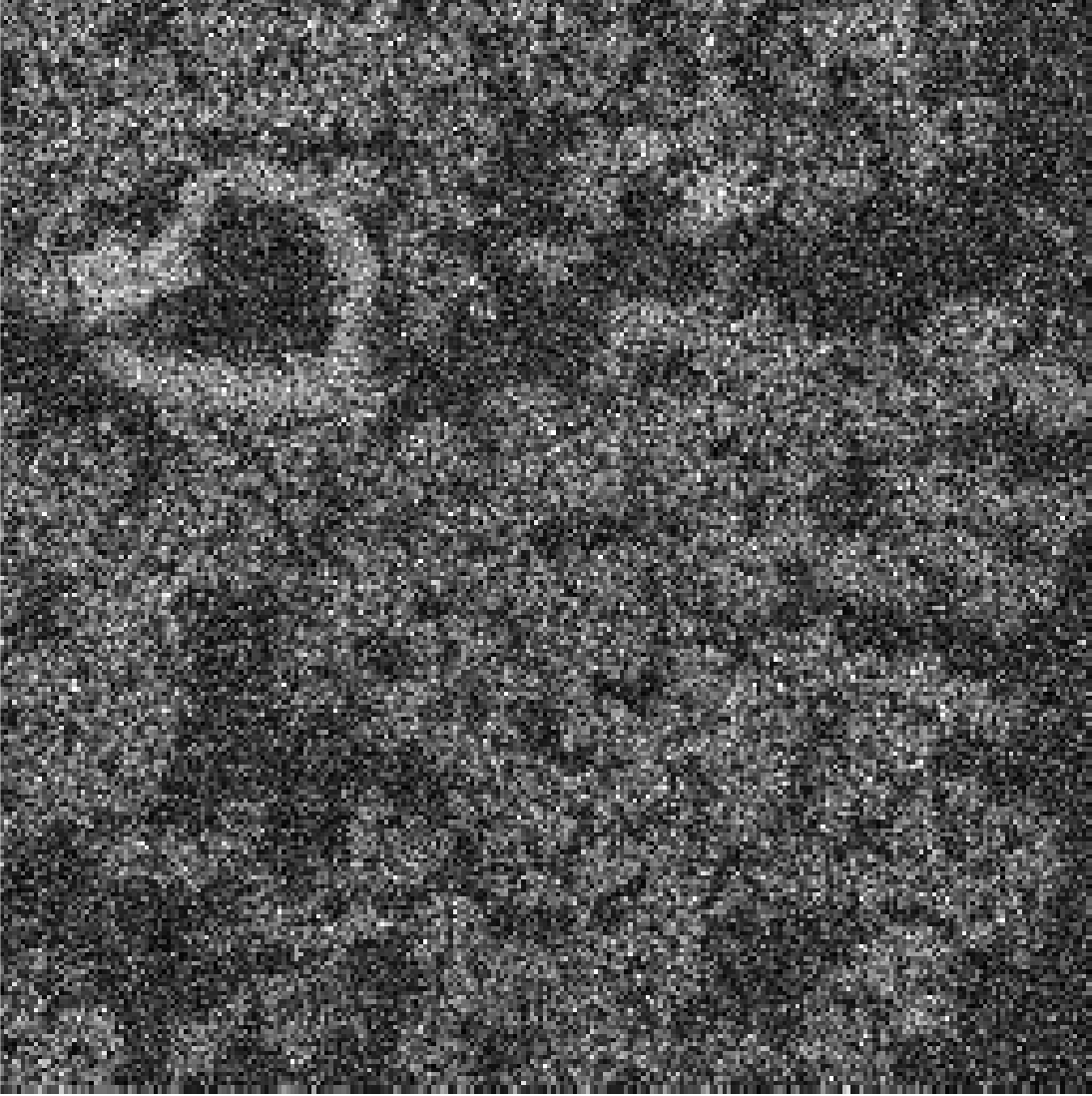}}
\subfigure[]{\includegraphics[width=.11\textwidth]{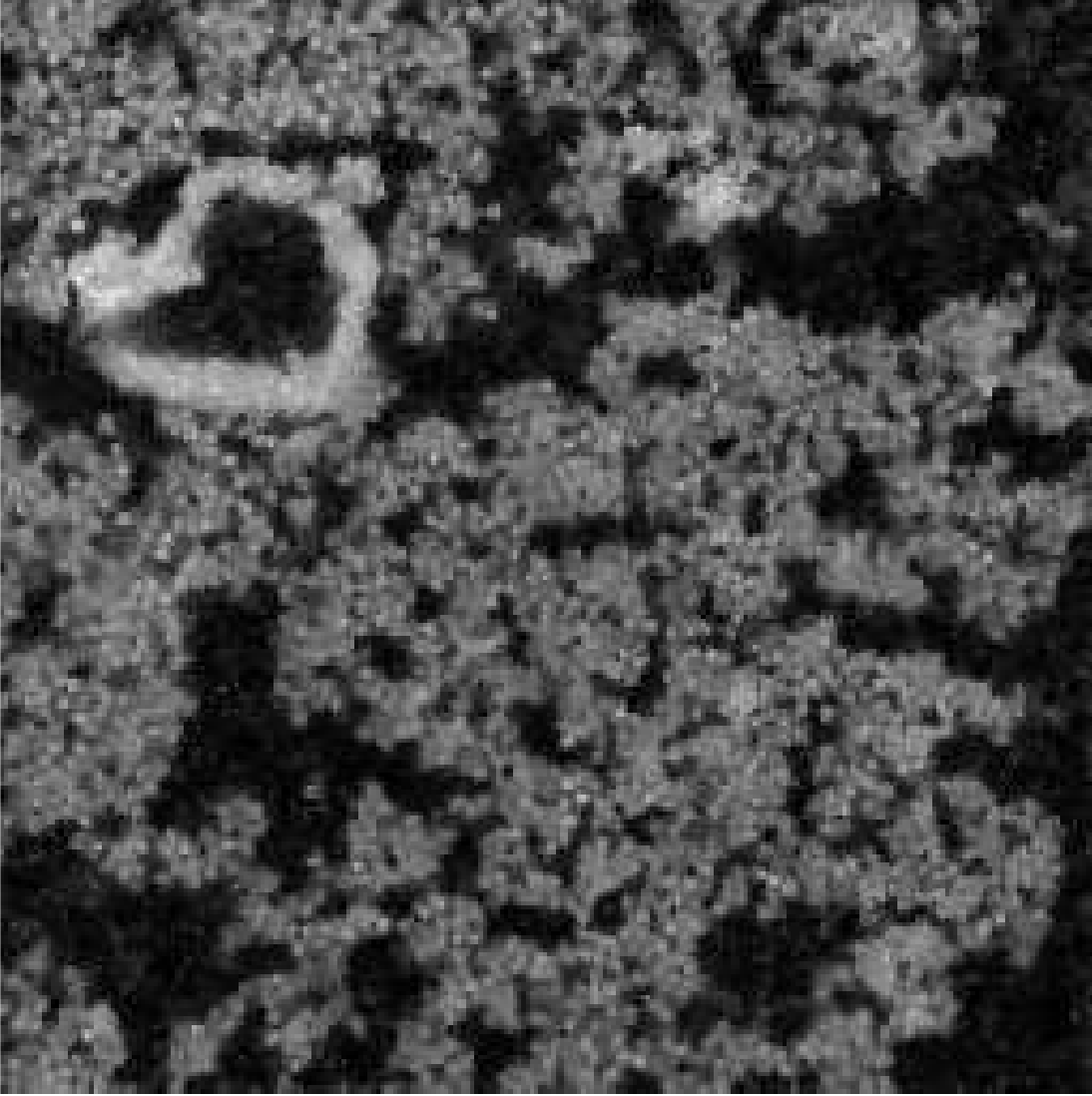}}
\subfigure[]{\includegraphics[width=.11\textwidth]{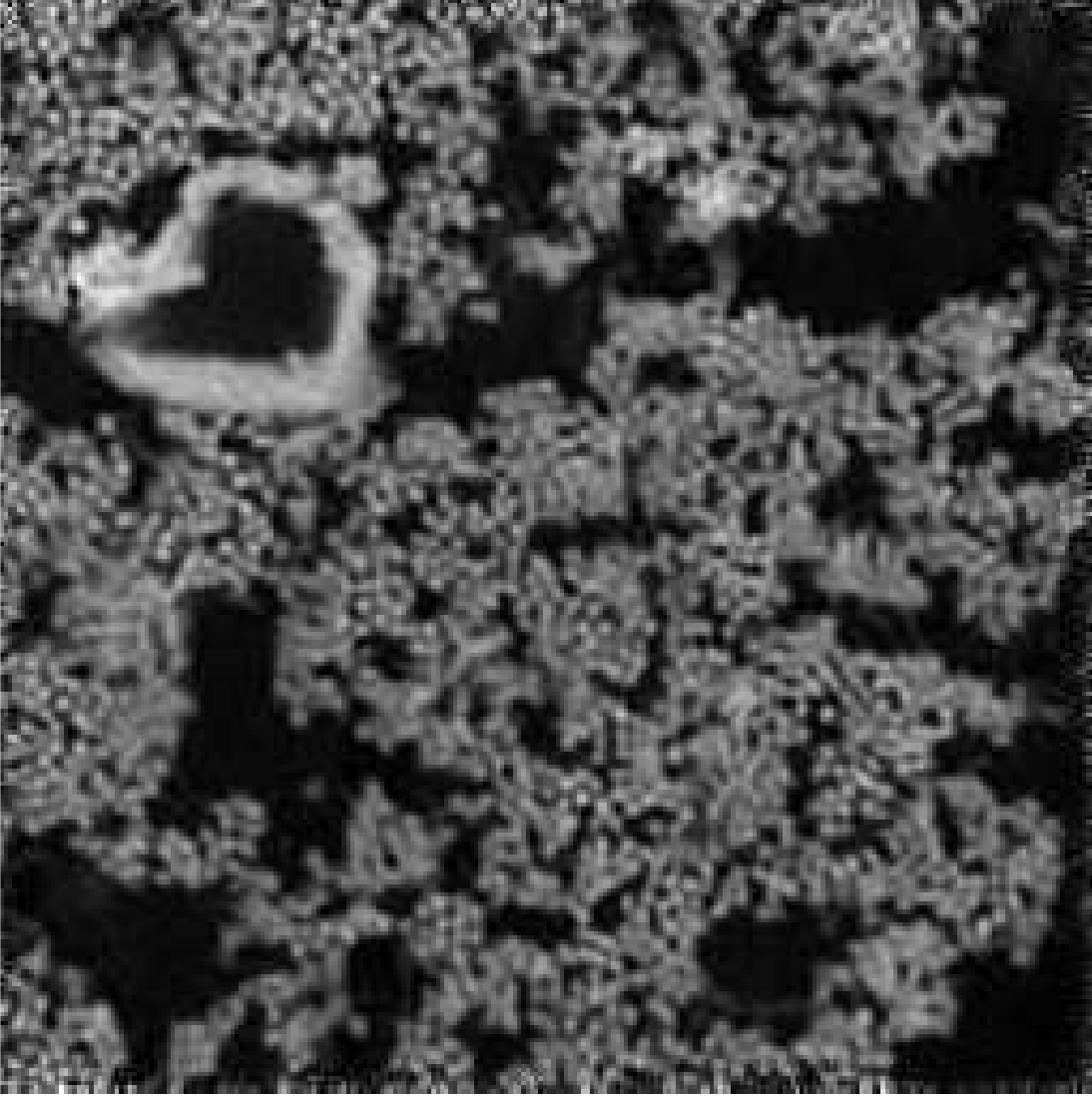}}
\subfigure[]{\includegraphics[width=.11\textwidth]{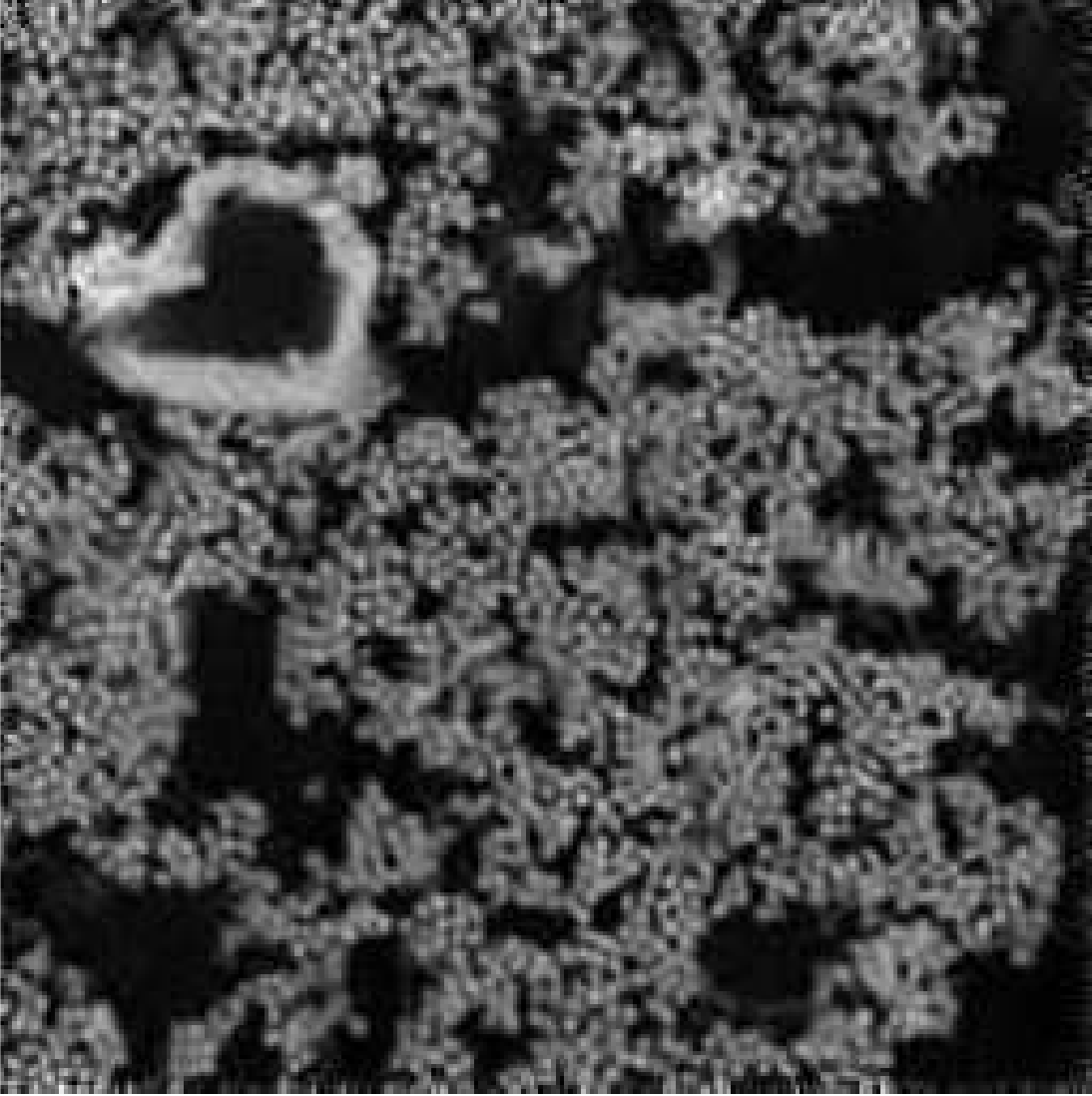}}\\\vskip -.05in
\subfigure[]{\includegraphics[width=.11\textwidth]{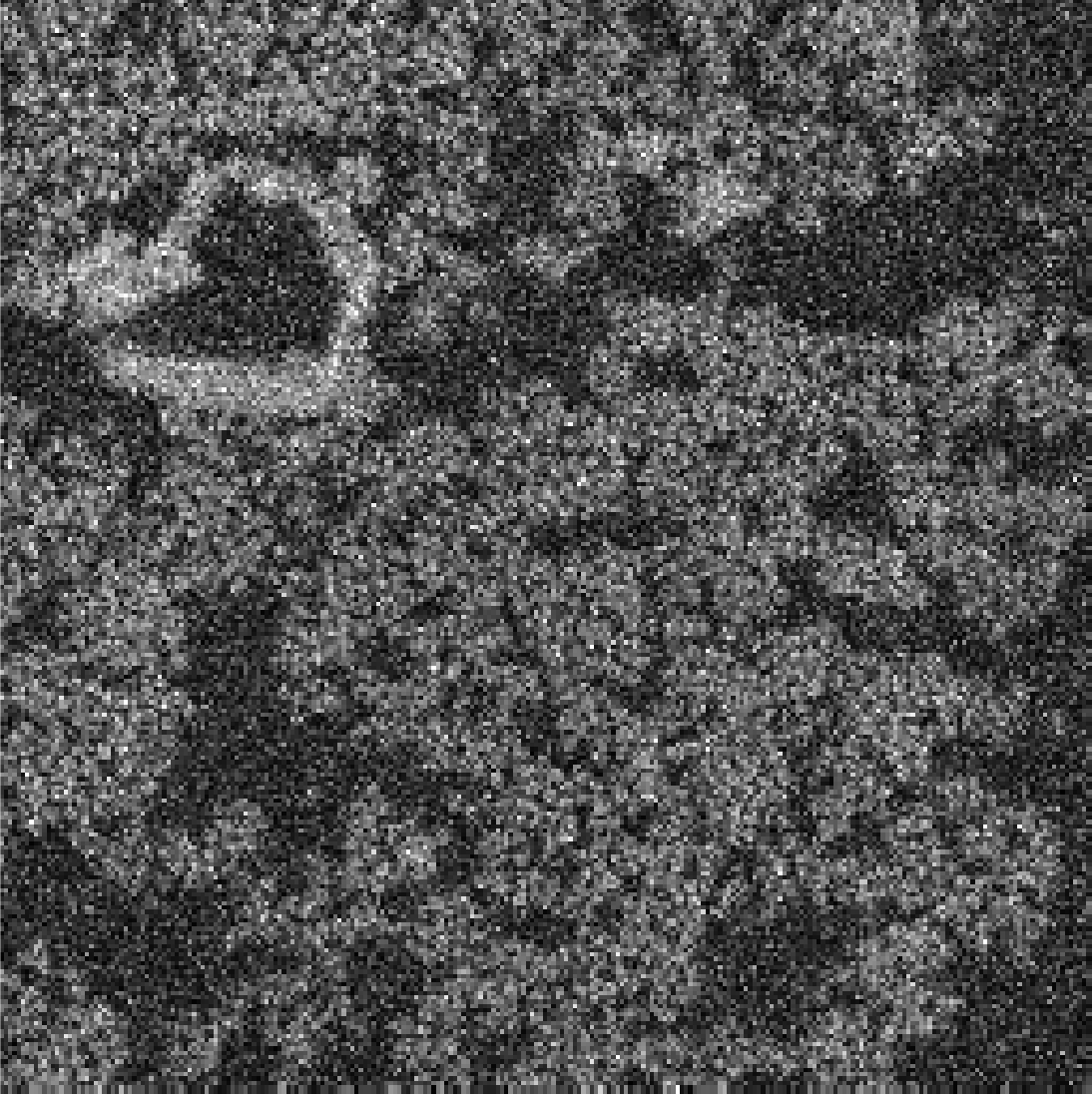}}
\subfigure[]{\includegraphics[width=.11\textwidth]{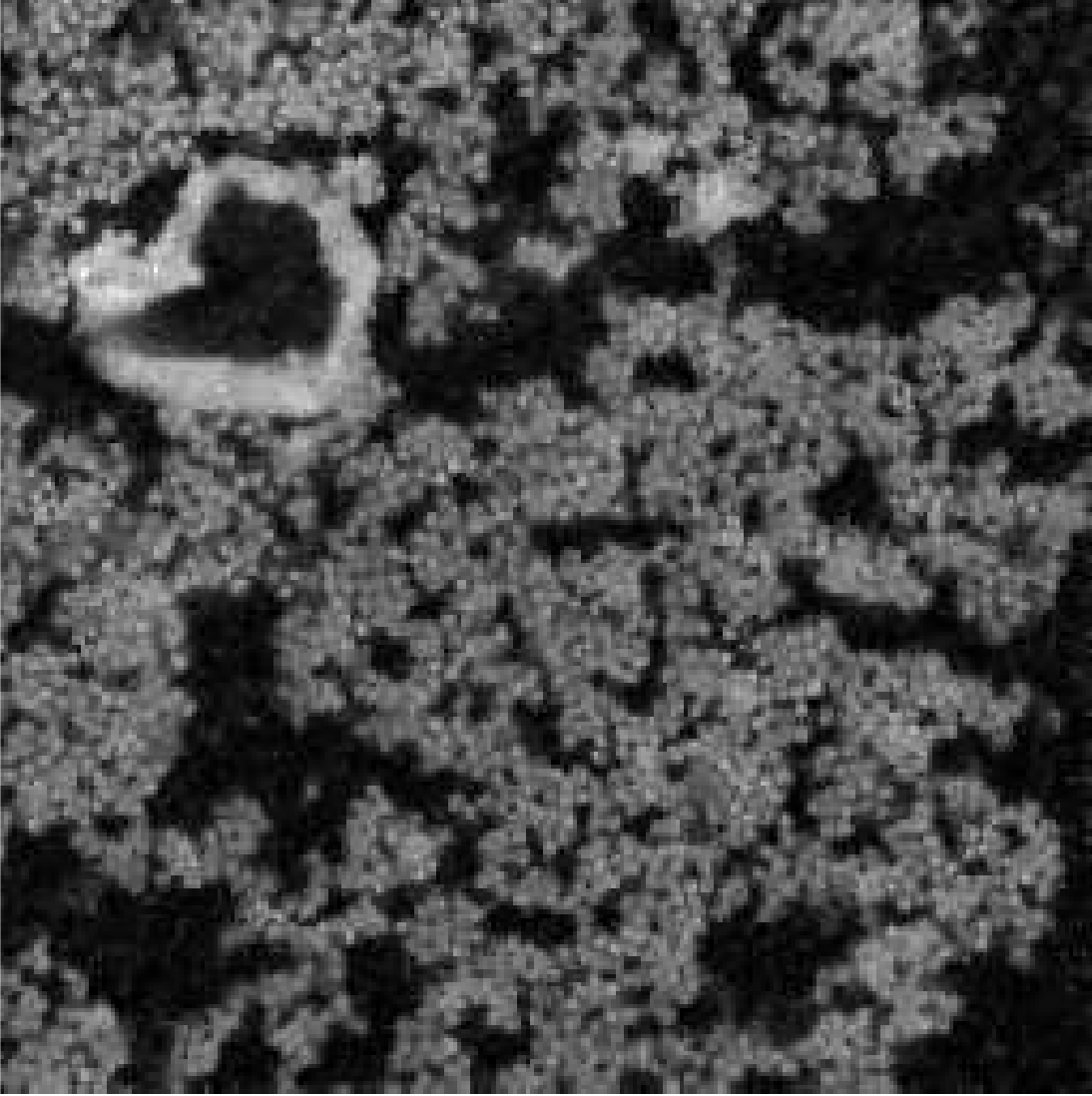}}
\subfigure[]{\includegraphics[width=.11\textwidth]{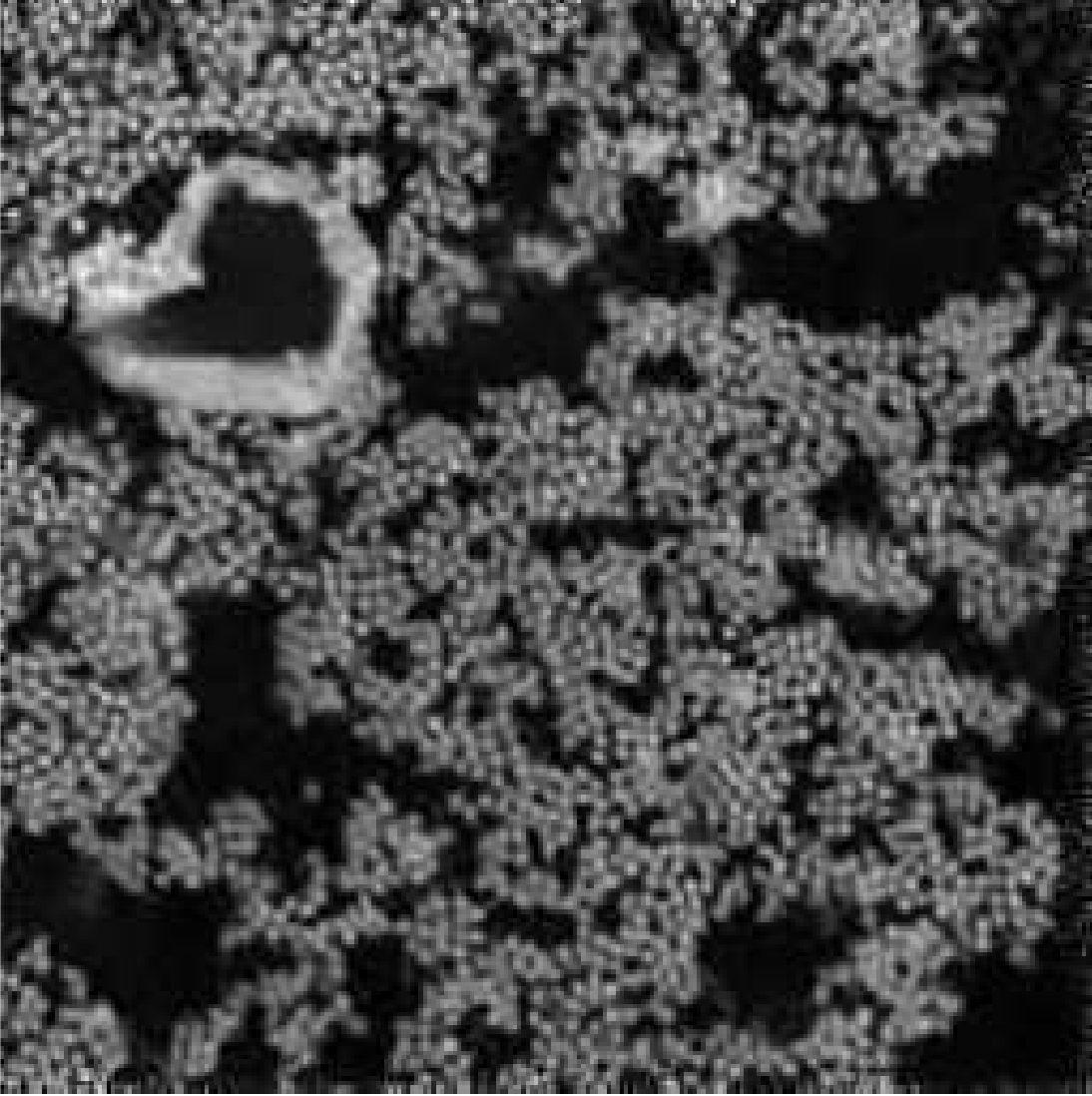}}
\subfigure[]{\includegraphics[width=.11\textwidth]{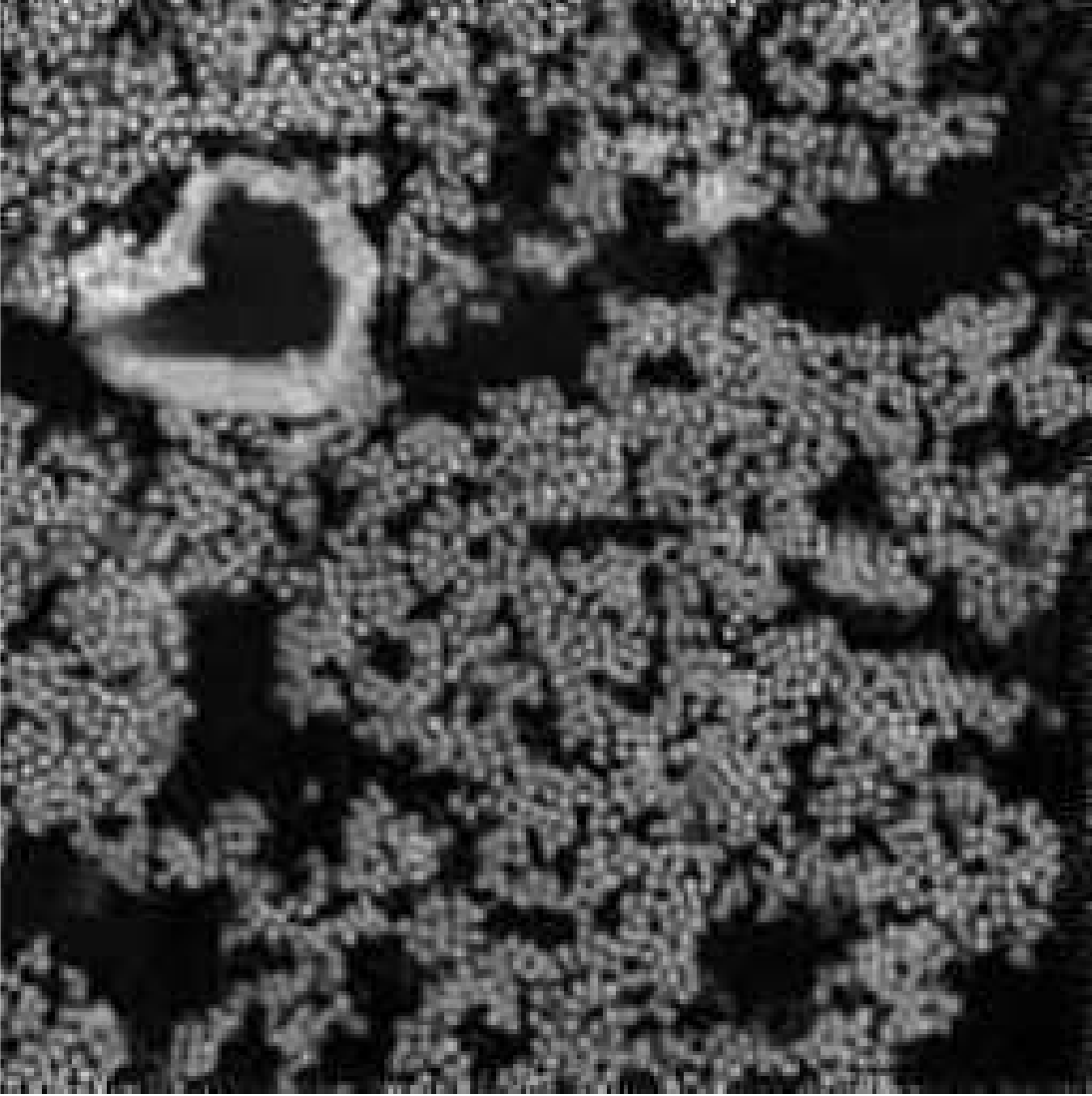}}\vskip -.05in
\end{center}
\vskip -.1in
\caption{CDP for complex-valued image  with $\delta\in\{0.08,0.1\}$. First row: $\delta=0.08$; Second row: $\delta=0.1$. From left to right: reconstructed images by PR in (a) and (e), TVPR in (b) and (f), ALGI in (c) and (g) and ALGI${}_{aniso}$ in (d) and (h).}
\label{cdp3}
\vskip -.2in
\end{figure}

More experiments  on the complex-valued image ``Goldballs'' in Fig. \ref{ground1}  for CDP  with peak level $\delta\in\{0.08, 0.1\}$ are performed, and please see the reconstructed results in Fig. \ref{cdp3}. Notice that for complex  valued images, we introduce two kinds of definitions for $L^0$ pseudo norm in \eqref{L0def}, and therefore we have two versions of PALM, \emph{i.e.} ``ALGI'' and ``ALGI${}_{aniso}$''. Set the parameters  $\eta=3.0\times 10^{-5}, \tau=1.0\times 10^{-3}$ for ``ALGI'' and $\eta=2.5\times 10^{-5}, \tau=8.0\times 10^{-4}$ for ``ALGI${}_{aniso}$''. Set iteration number in the outer loop as $T=50$.  By observing Fig. \ref{cdp3}, from noisy measurements one can only derive noisy results by ``PR''. By ``TVPR'', ``ALGI'' and ``ALGI${}_{aniso}$'', the recovery results are almost noise free and also have very clean background. ``TVPR'' can only recover the large scale structures at the left top corner, but can not keep the smaller repetitive structures. Obviously our proposed methods can produce better results, where both the large and small structures can be preserved well.  Related SNRs are put in Table \ref{tab3}, about 1dB, 2dB increases are gained by ``ALGI'', and ``ALGI${}_{aniso}$'' compared with ``TVPR''.

\begin{figure}
\vskip .15in
\begin{center}
\subfigure[]{\includegraphics[width=.07\textwidth]{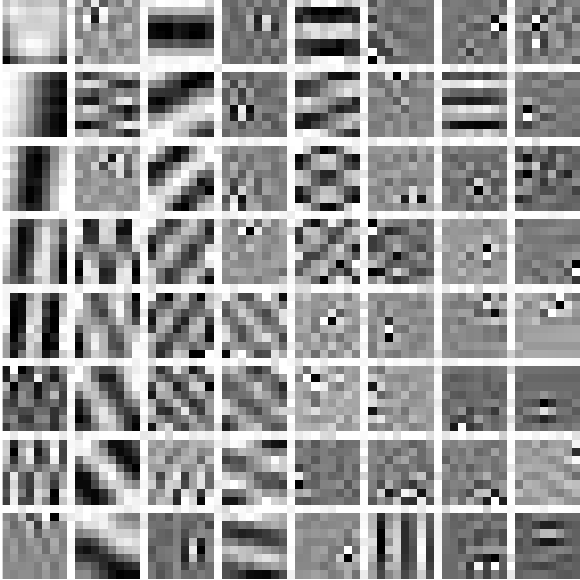}}~
\subfigure[]{\includegraphics[width=.07\textwidth]{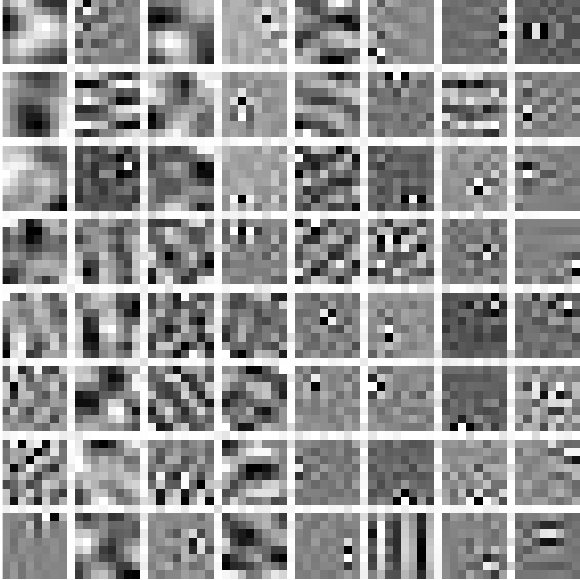}}~
\subfigure[]{\includegraphics[width=.07\textwidth]{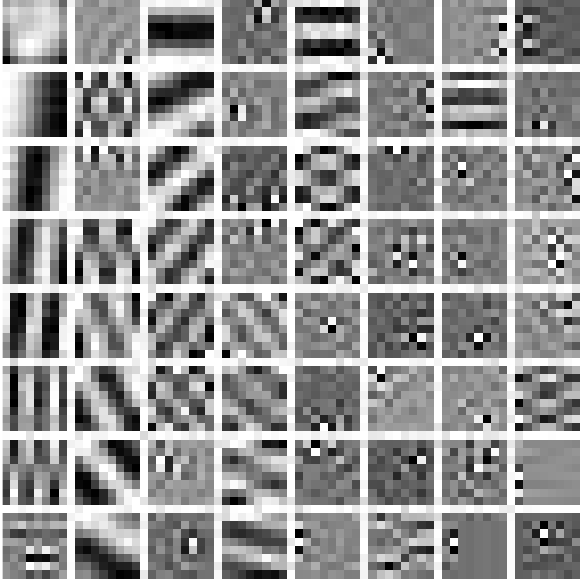}}~
\subfigure[]{\includegraphics[width=.07\textwidth]{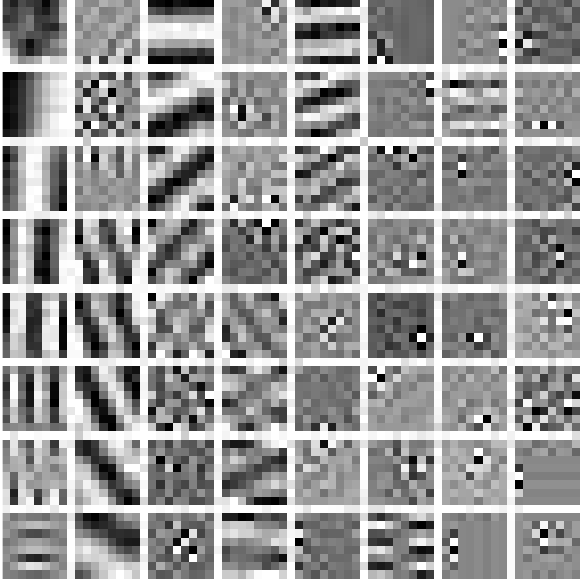}}\\\vskip -.05in
\subfigure[]{\includegraphics[width=.07\textwidth]{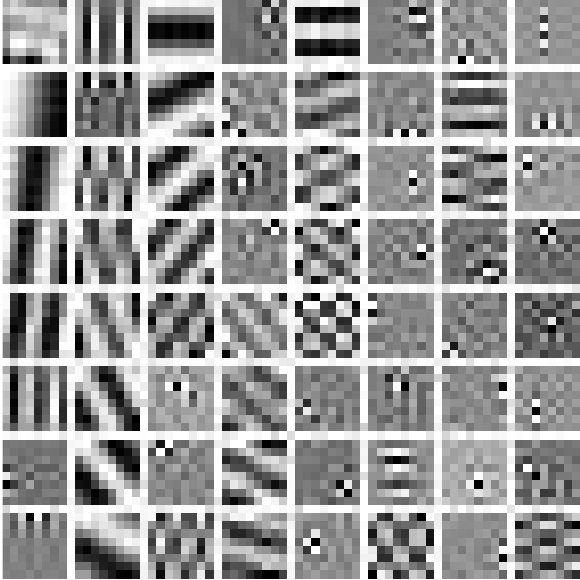}}~
\subfigure[]{\includegraphics[width=.07\textwidth]{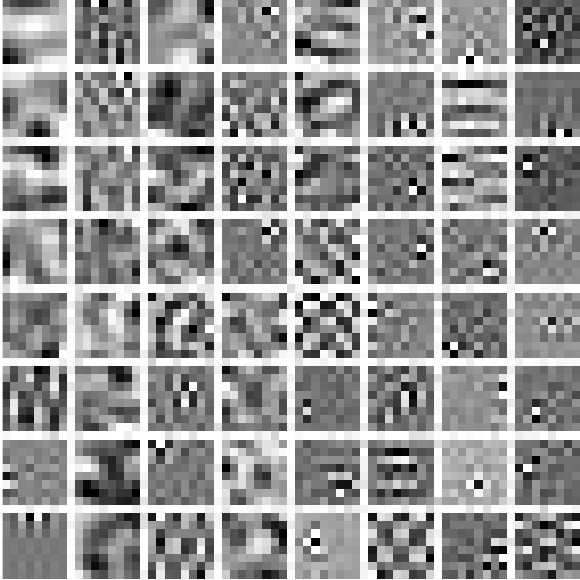}}~
\subfigure[]{\includegraphics[width=.07\textwidth]{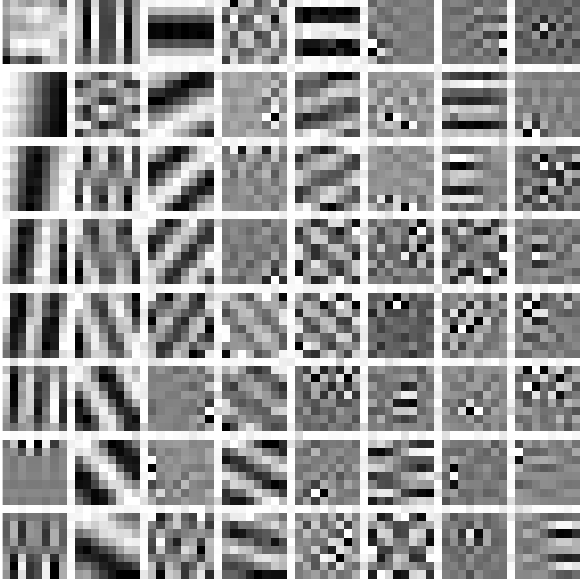}}~
\subfigure[]{\includegraphics[width=.07\textwidth]{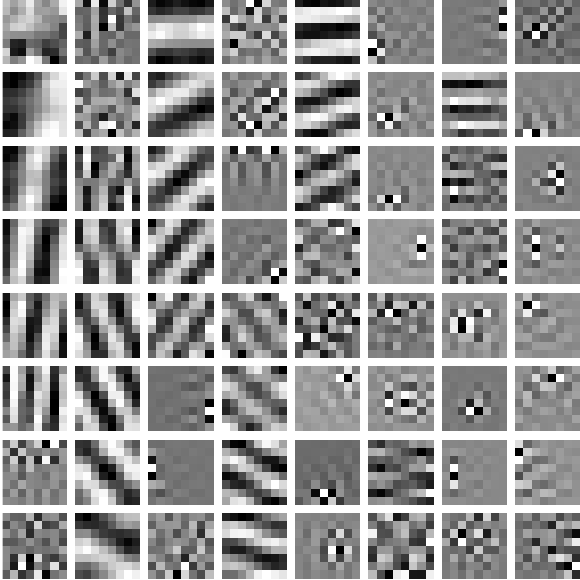}}\vskip -.05in
\end{center}
\vskip -.1in
\caption{Learned dictionaries for CDP corresponding to the results in Fig. \ref{cdp3}.  First row: $\delta=0.08$; Second row: $\delta=0.1$. From left to right: Real and complex parts of learned dictionaries by ALGI in the first two columns, and by ALGI${}_{aniso}$ in the third and fourth columns.}
\label{dic}
\vskip -.1in
\end{figure}

One also notices that ``ALG${}_{aniso}$'' produces the recovery results with higher accuracy than ``ALGI''.  In order to investigate the improvements by the anisotropic norm, we show the learned dictionaries in  Fig. \ref{dic}, and the sparsity of coefficient matrix $\bm \alpha$ in Table \ref{tabSparseC}.  Especially in Fig. \ref{dic}(d) and (h)  the imaginary parts of learned dictionary by ``ALGI${}_{aniso}$''  seem to have more features than those in Fig. \ref{dic}(b) and (f) by ``AlGI''. Meanwhile  the real parts of learned dictionaries by the anisotropic version algorithm are quite close to the imaginary parts, which is consistent with the similarity of structure of real and complex parts of ground truth in Fig. \ref{ground1}(f) and (g). Therefore, it produces better dictionaries by the anisotropic style in \eqref{L0def}. As a result, the sparsity is  strengthened, and one can observe that the sparsity levels $S(\Re(\bm\alpha))$ are much smaller by ``ALGI${}_{aniso}$'' than those by ``ALGI'' in Table \ref{tabSparseC}. It demonstrates that  anisotropic version algorithm  can help to  improve the image qualities  for  complex-valued images.


\begin{table}
\vskip .05in
\begin{center}
\begin{spacing}{}
\begin{tabular}{|c||c|c|c|c|}
\hline
$\delta$& PR &TVPR &ALGI&ALGI${}_{aniso}$ \\\hline\hline
0.08&2.91&10.66&11.52&\bf 12.72  \\ \hline
0.1&4.17&11.81&12.75&\bf 14.03\\
\hline
\end{tabular}
\end{spacing}
\end{center}
\vskip -.3in
\caption{SNRs of CDP for complex-valued image }
\label{tab3}
\vskip -.3in
\end{table}

\begin{table*}
\vskip -.1in
\begin{center}
\begin{spacing}{.95}
\begin{tabular}{|c||c|c|c||c|c|c|}
\hline
\multirow{2}{*}{Isotropic/Anisotropic $L^0$}&  \multicolumn{3}{c||}{$\delta=0.08$}&\multicolumn{3}{c|}{$\delta=0.1$}\\
\cline{2-7}
& $S(\bm\alpha)$ &$S(\Re(\bm\alpha))$&$S(\Im(\bm\alpha))$ & $S(\bm\alpha)$ &$S(\Re(\bm\alpha))$&$S(\Im(\bm\alpha))$ \\ \hline\hline
ALGI               &18.9\%&18.9\% &18.9\%&24.9\%&24.8\%&24.9\%\\ \hline
ALGI${}_{aniso}$   &20.5\%&20.5\% &0.25\%&25.9\%&25.9\%&0.16\%\\
\hline
\end{tabular}
\end{spacing}
\end{center}
\vskip -.3in
\caption{Sparsity  for CDP on complex-valued images with different $L^0$ pseudo norms.}
\label{tabSparseC}
\vskip -.1in
\end{table*}

%

\subsection{Ptychographic phase retrieval (PtychoPR) for complex-valued image}
\subsubsection{Fixed Sliding Distance}
\begin{figure}
\vskip -.25in
\begin{center}
\subfigure[]{\includegraphics[width=.11\textwidth]{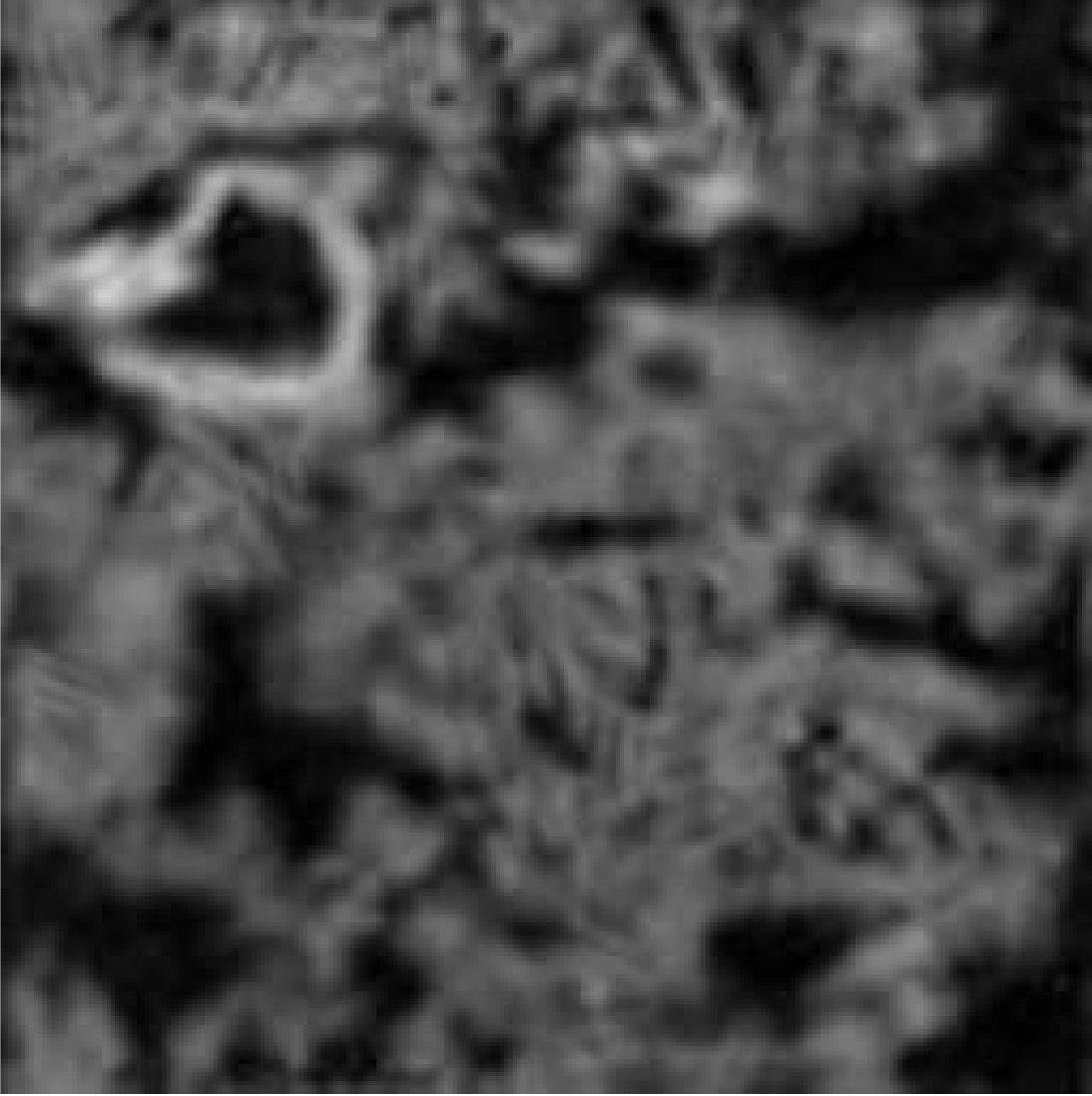}}
\subfigure[]{\includegraphics[width=.11\textwidth]{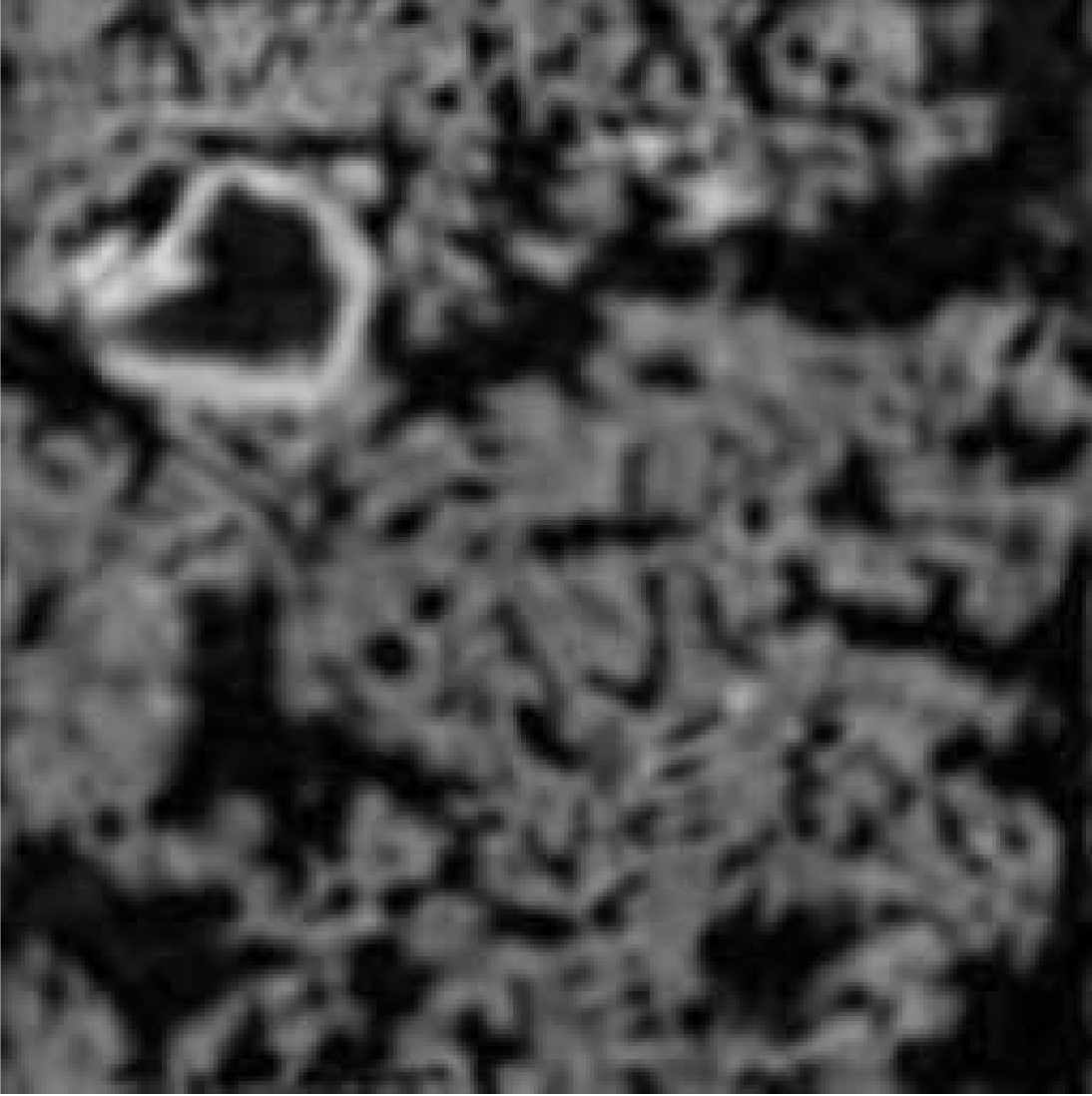}}
\subfigure[]{\includegraphics[width=.11\textwidth]{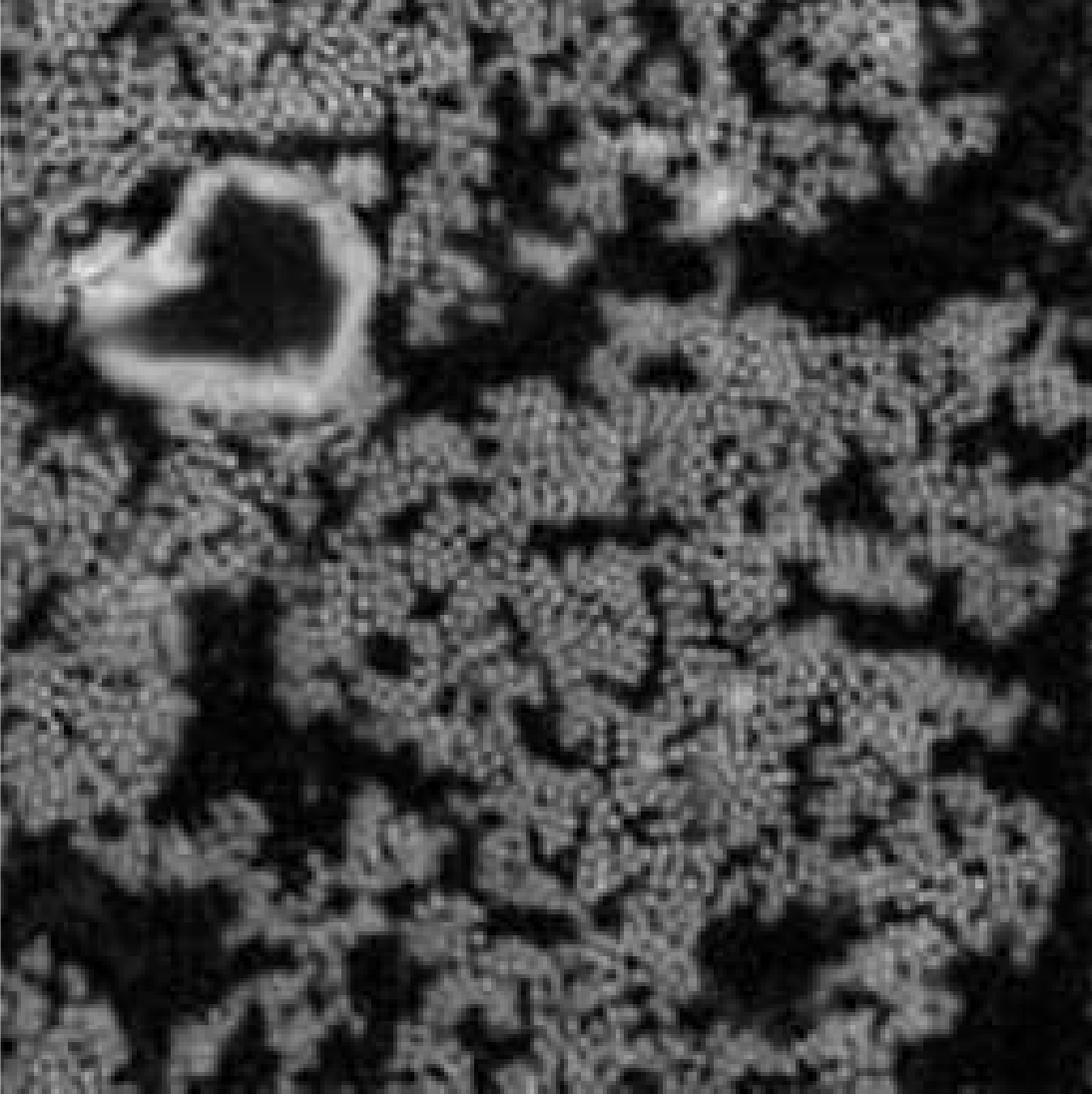}}
\subfigure[]{\includegraphics[width=.11\textwidth]{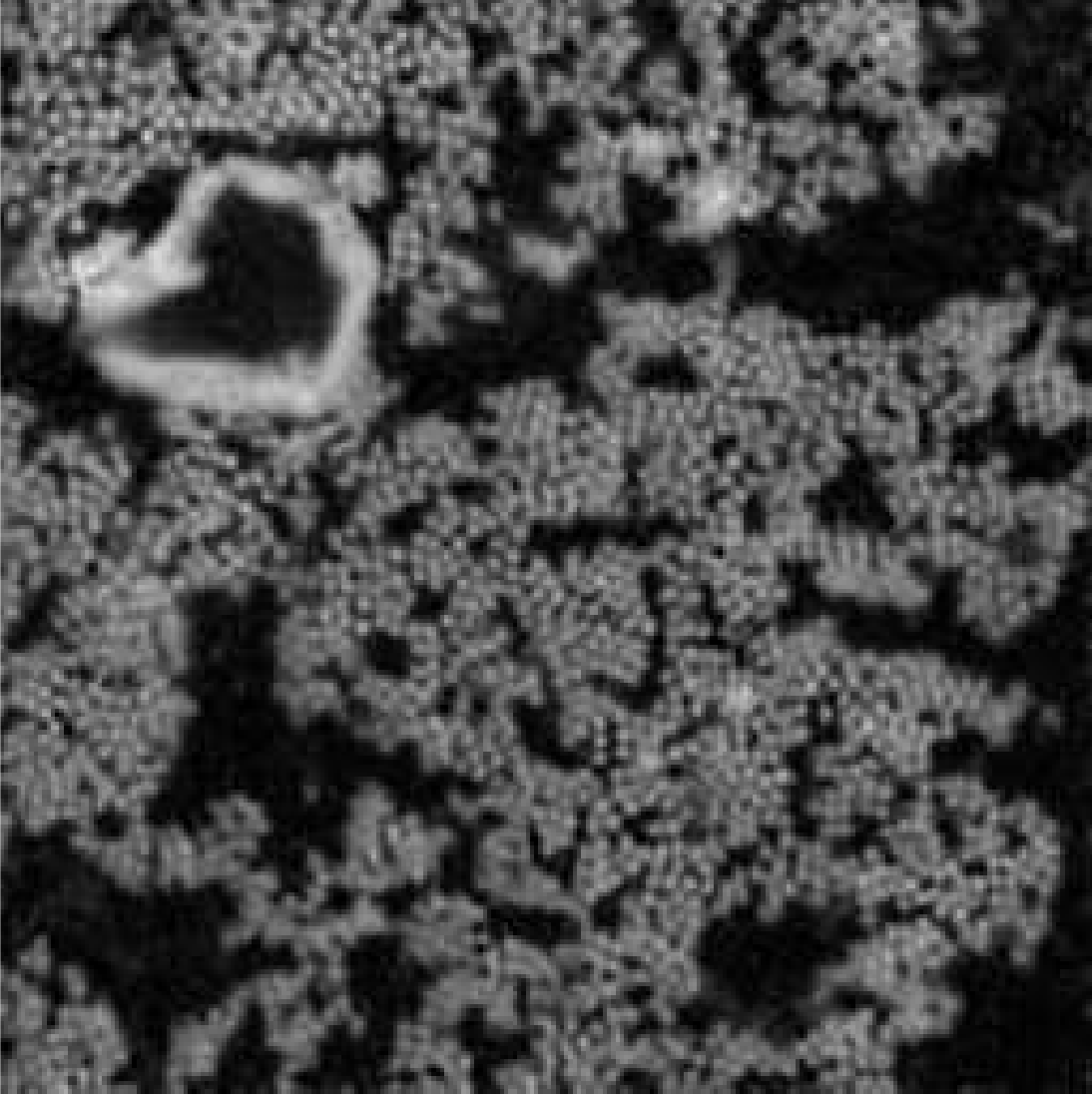}}\\\vskip -.05in
\subfigure[]{\includegraphics[width=.11\textwidth]{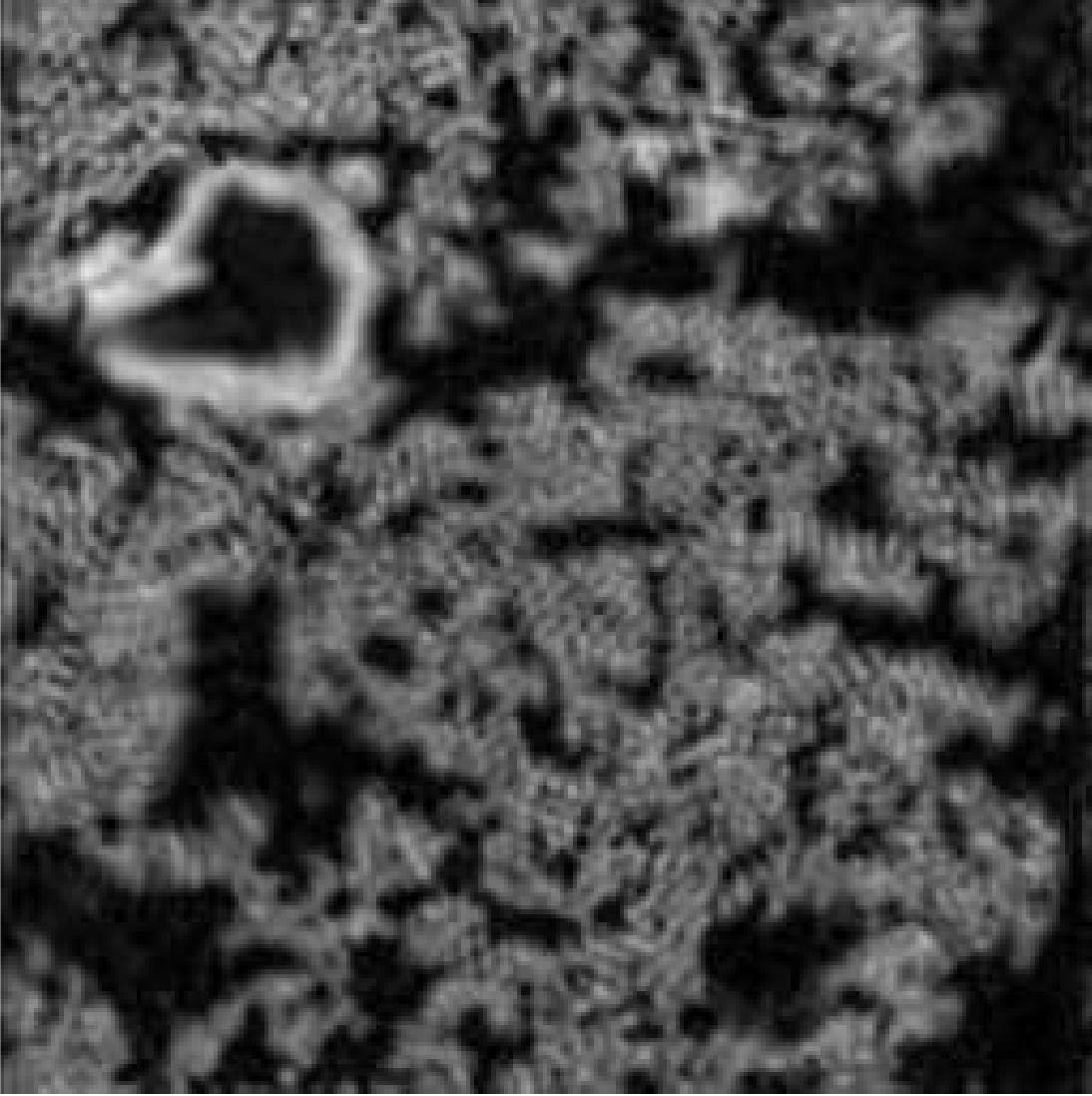}}
\subfigure[]{\includegraphics[width=.11\textwidth]{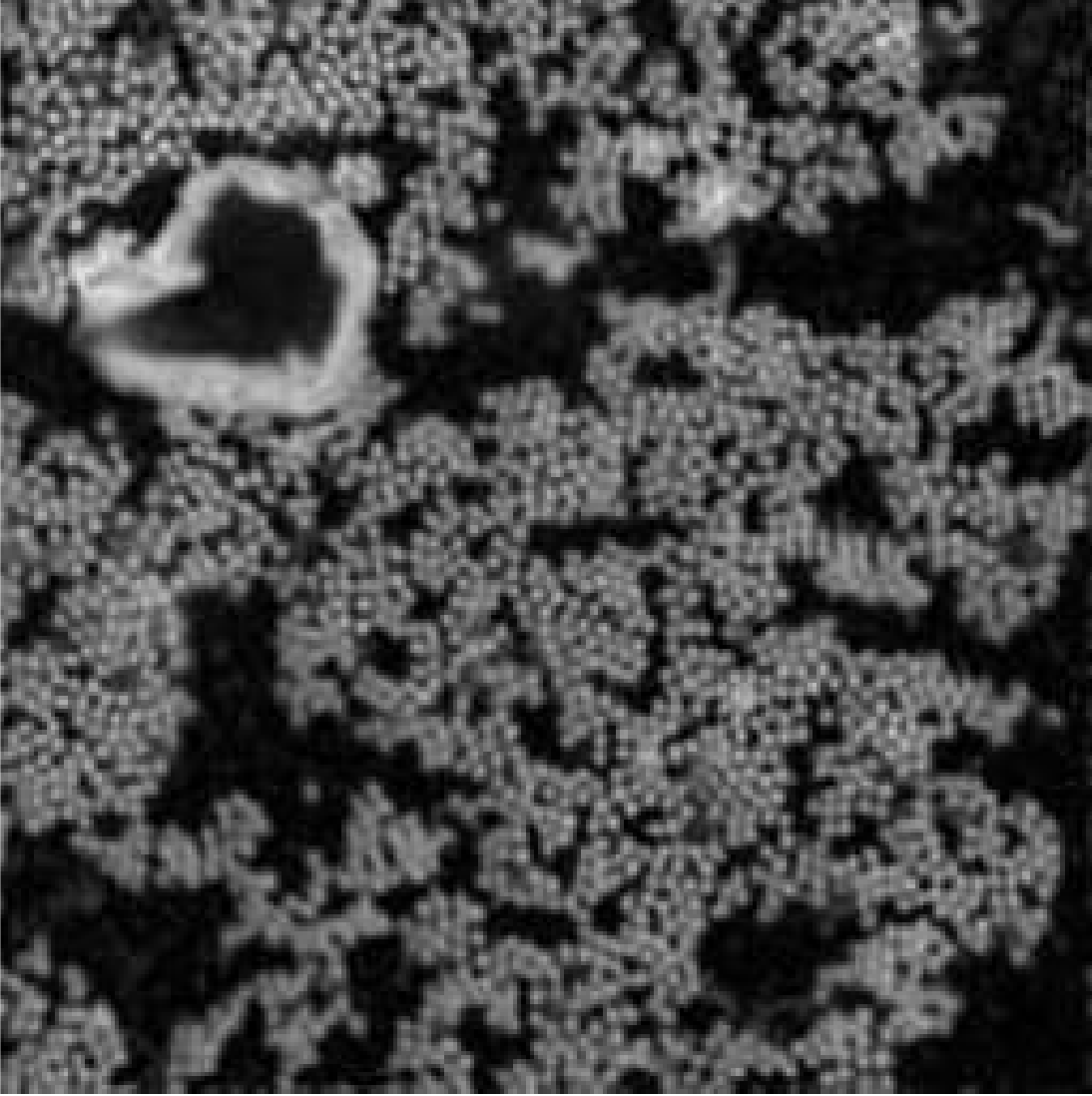}}
\subfigure[]{\includegraphics[width=.11\textwidth]{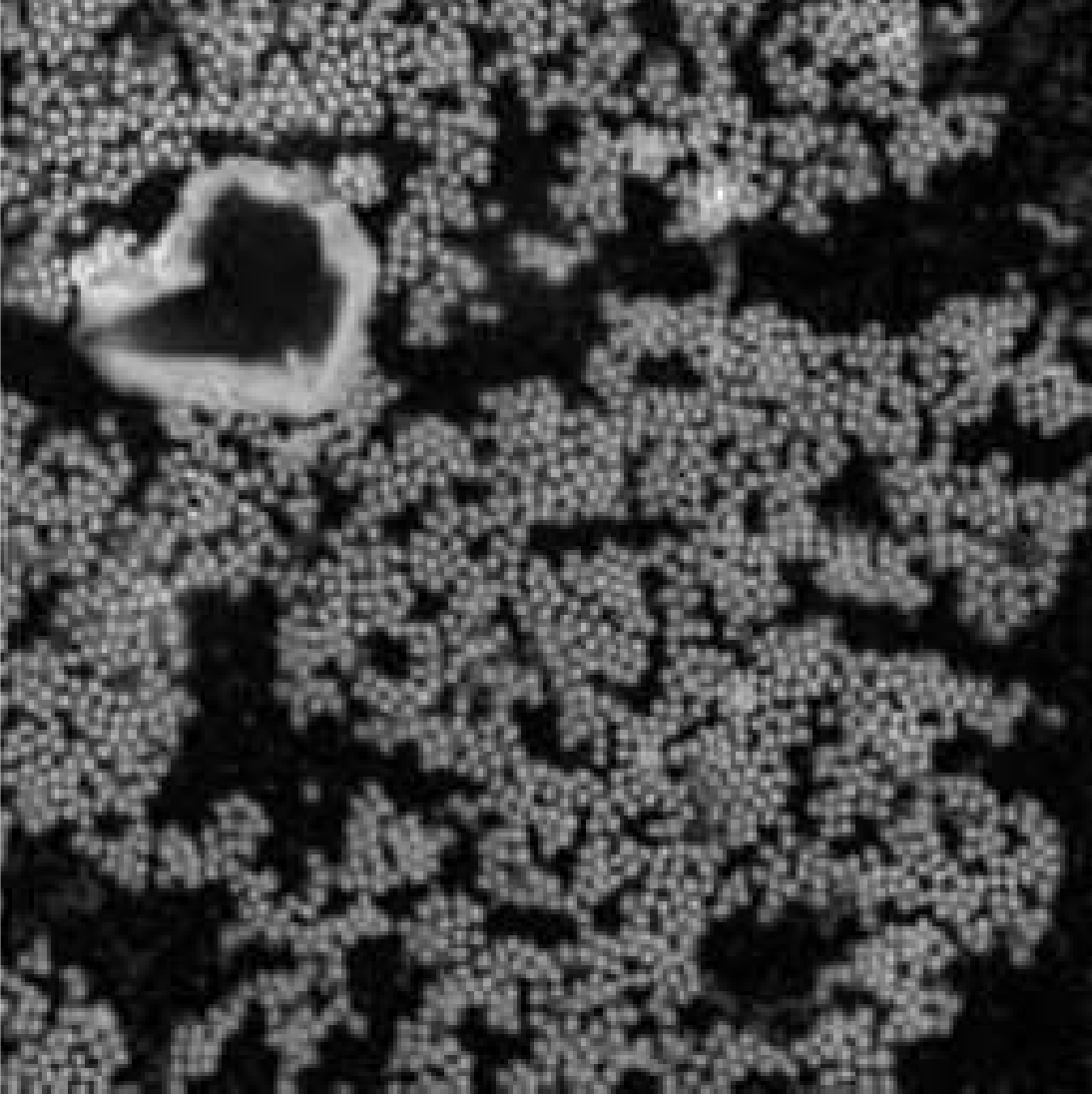}}
\subfigure[]{\includegraphics[width=.11\textwidth]{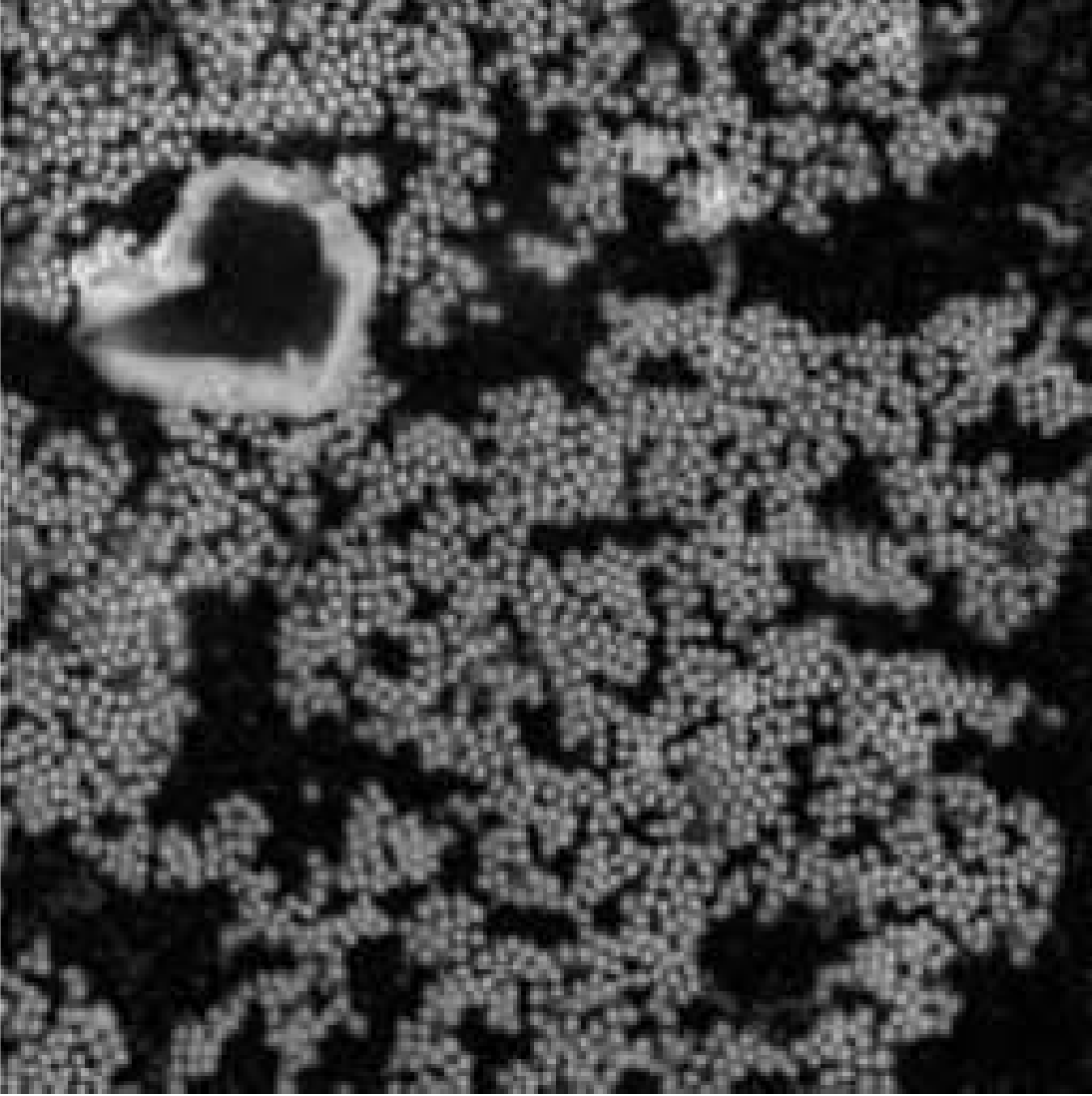}}\vskip -.05in
\end{center}
\vskip -.1in
\caption{PtychoPR for complex-valued image  with $\delta\in\{0.2,0.5\}$. First row: $\delta=0.2$; Second row: $\delta=0.5$. From left to right: reconstructed images by PR in (a) and (e), TVPR in (b) and (f), ALGI in (c) and (g) and ALGI${}_{aniso}$ in (d) and (h).}
\label{pty1}
\end{figure}

For PtychoPR, a zone plate lens and responding illumination mask are employed as \cite{marchesini2015alternating}.  Set sliding distance $\mathrm{SlidDist}=16$. The number of frames is $16\times16$  and the frame size is of $64\times 64$, and in such setting total number of measurements $m=16n$. Set peak level $\delta\in\{0.2, 0.5\}$. Here only two inner iterations for PALM is sufficient. Set outer iteration number $T=50$,  $\tau=5.0\times 10^{-3}, \eta=3.0\times 10^{-2}$ for $\delta=0.2$ and $\tau=8.0\times 10^{-3}, \eta=3.0\times 10^{-2}$ for $\delta=0.5$.   The reconstructed results are put in Fig. \ref{pty1}, and the corresponding SNRs are put in Table \ref{tab4}. Inferred from Fig. \ref{pty1}(a) and (e), it seems very blurry in the results of ``PR'' from noisy measurements, since the phaseless data are generated by structured deterministic illumination, and  corrupted  low frequency parts are worse than high frequency parts. By ``TVPR'', at higher noise level  $\delta=0.2$ in the first row of Fig. \ref{pty1}, it can produce results with sharp edges for large scale features and clean background, but can not preserve the smaller features at all. For the case with peak level $\delta=0.5$, ``TVPR'' can produce pretty good recovery results for both smaller and larger scale features. The proposed ``ALGI'' and ``ALGI${}_{aniso}$'' can recover the smaller features very well  especially at peak level $\delta=0.2$. Inferred from Table \ref{tab4}, SNRs are increased about 3.5dB, 4dB for ``ALGI'' and ``ALGI${}_{aniso}$'', respectively compared with ``TVPR'' at peak level $\delta=0.2$;  SNRs are increased about 1.5dB, 1.7dB for ``ALGI'' and ``ALGI${}_{aniso}$'' respectively compared with ``TVPR'' at peak level $\delta=0.5$. Such gains in SNRs implied our proposed algorithms can produce more accuracy results. Similarly to the previous subsection for complex-valued image, the anisotropic version algorithm  ``ALGI${}_{aniso}$'' has better performances than isotropic version algorithm.

%

\begin{table}
\begin{center}
\begin{spacing}{0.9}
\begin{tabular}{|c||c|c|c|c|}
\hline
$\delta$& PR &TVPR &ALGI&ALGI${}_{aniso}$ \\\hline\hline
0.2&6.85&8.91&12.35&\bf 12.84  \\ \hline
0.5&8.25&15.03&16.52&\bf 16.75\\
\hline
\end{tabular}
\end{spacing}
\end{center}
\vskip -.25in
\caption{SNRs of PtychoPR  for complex-valued image corresponding to Fig. \ref{pty1}  }
\label{tab4}
\vskip -.3in
\end{table}

\subsubsection{Variable Sliding Distance}
\begin{figure}
\vskip -.25in
\begin{center}
\subfigure[]{\includegraphics[width=.11\textwidth]{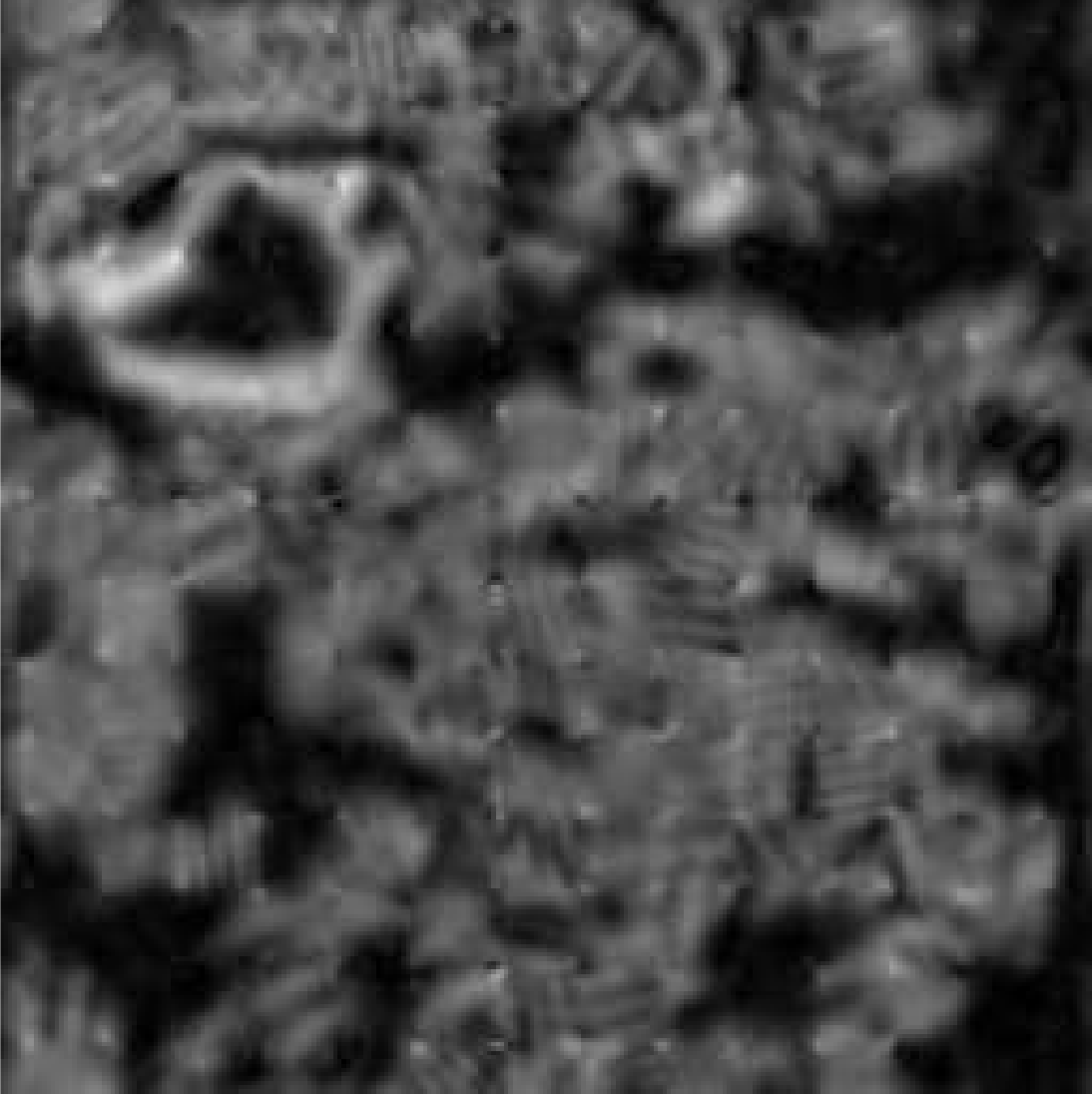}}
\subfigure[]{\includegraphics[width=.11\textwidth]{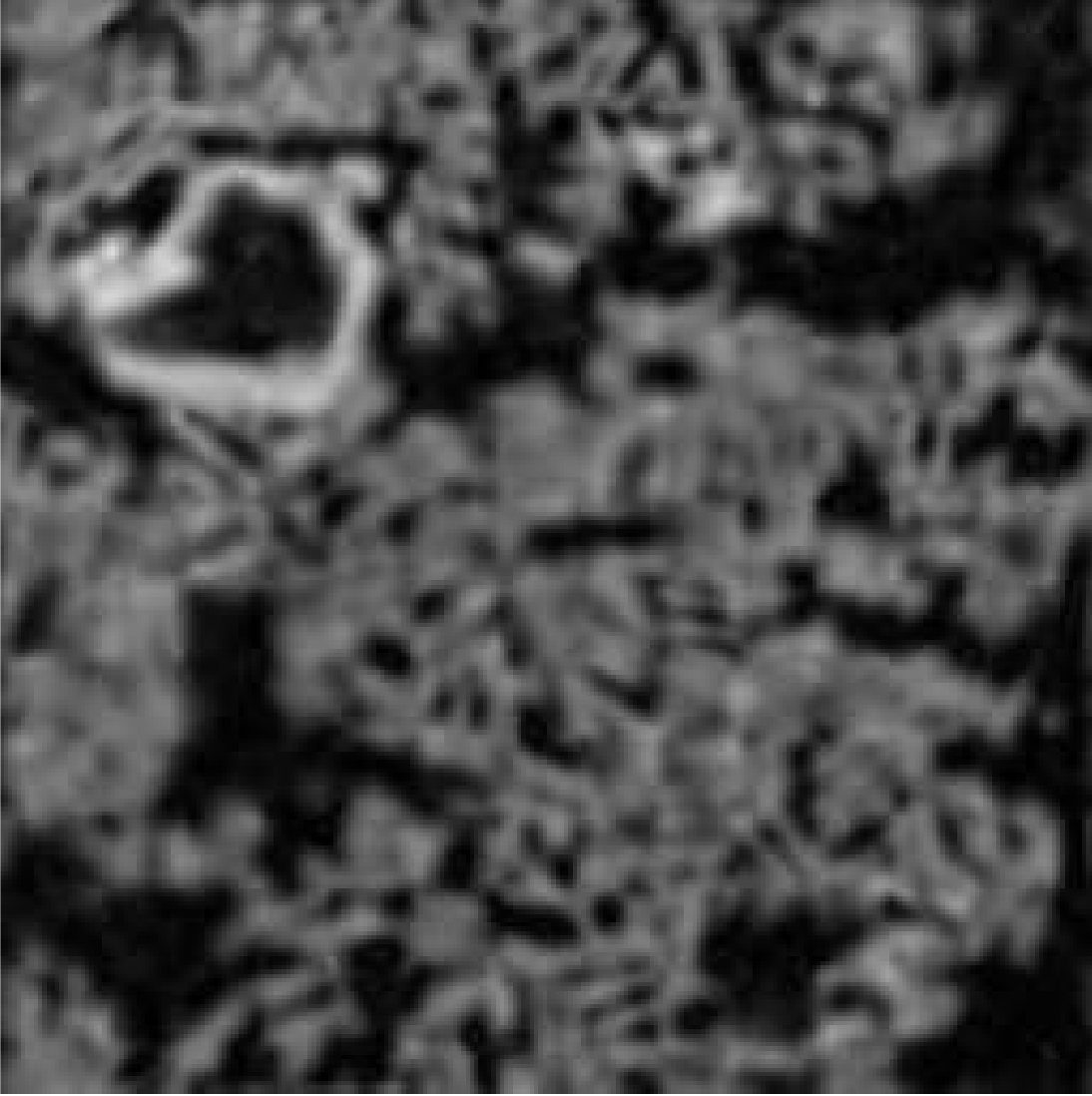}}
\subfigure[]{\includegraphics[width=.11\textwidth]{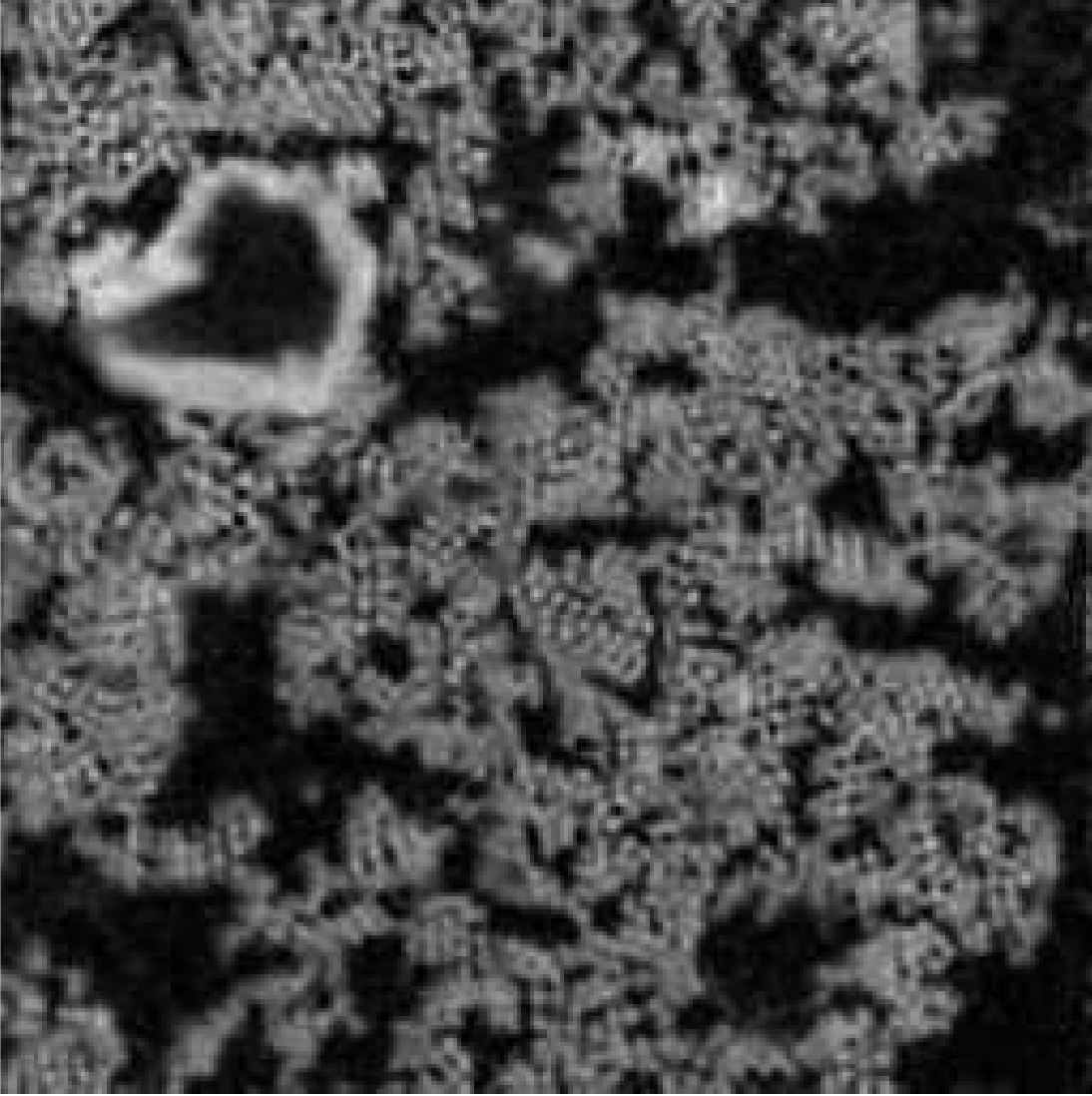}}
\subfigure[]{\includegraphics[width=.11\textwidth]{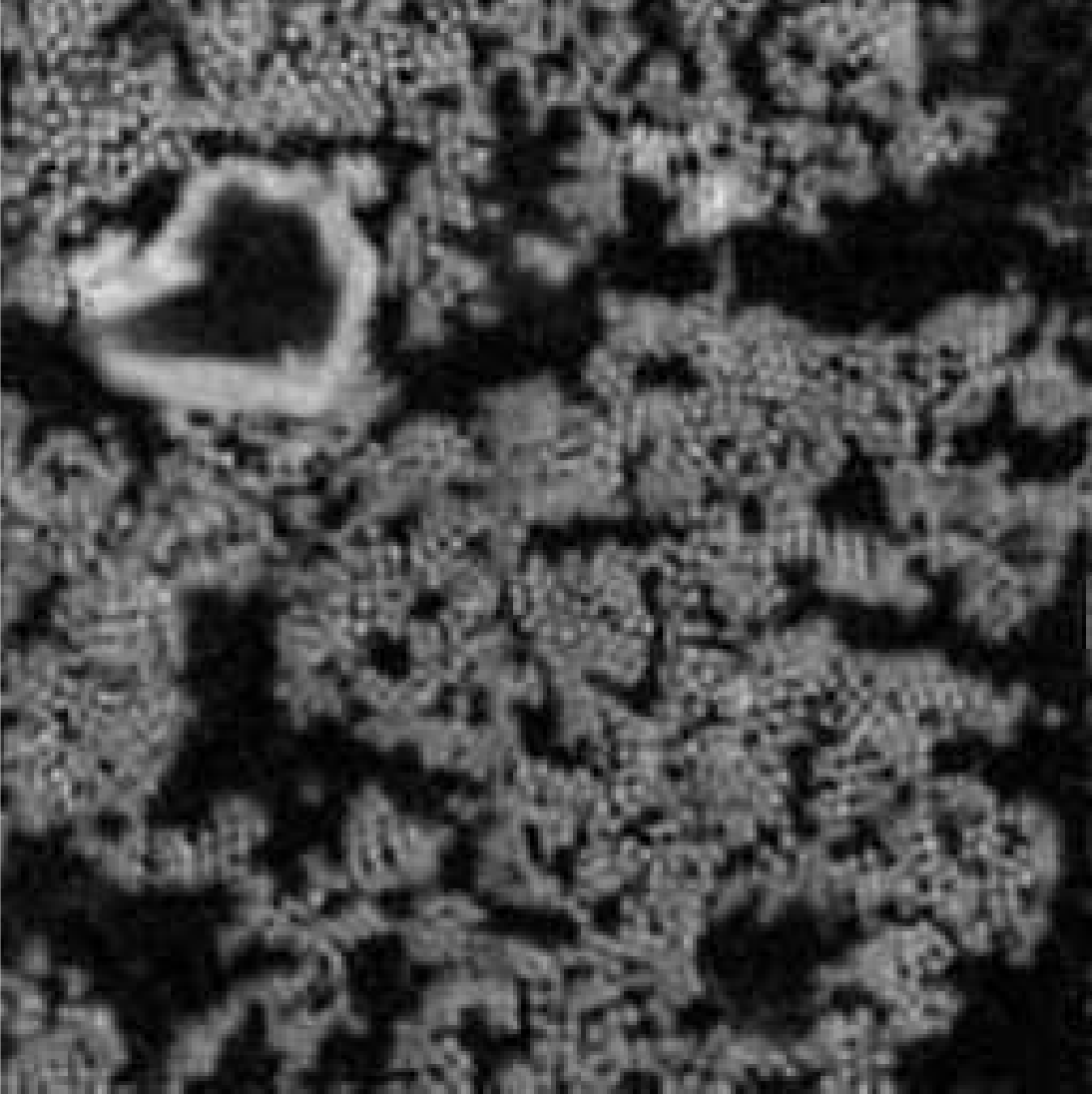}}\\\vskip -.05in
\subfigure[]{\includegraphics[width=.11\textwidth]{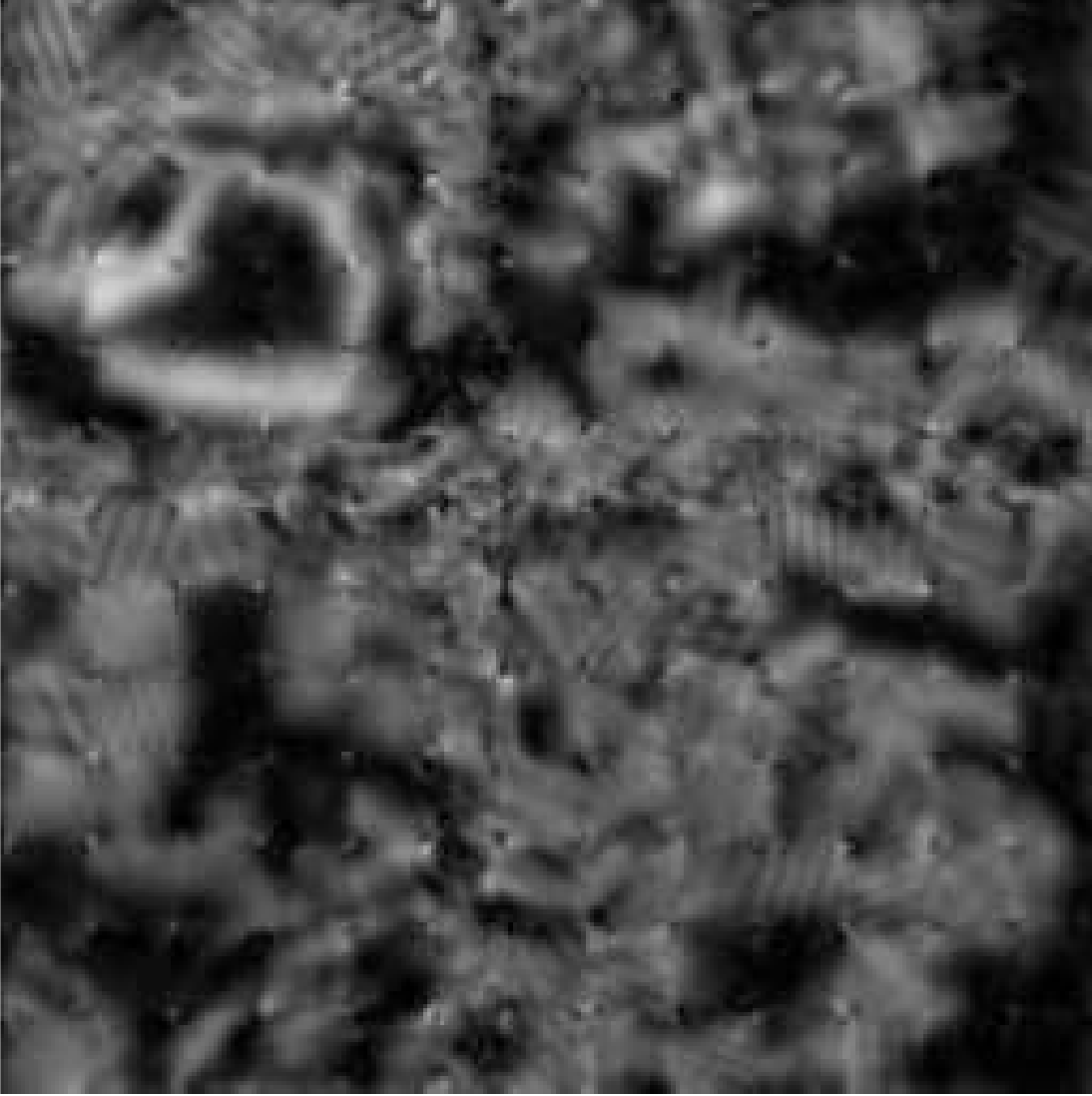}}
\subfigure[]{\includegraphics[width=.11\textwidth]{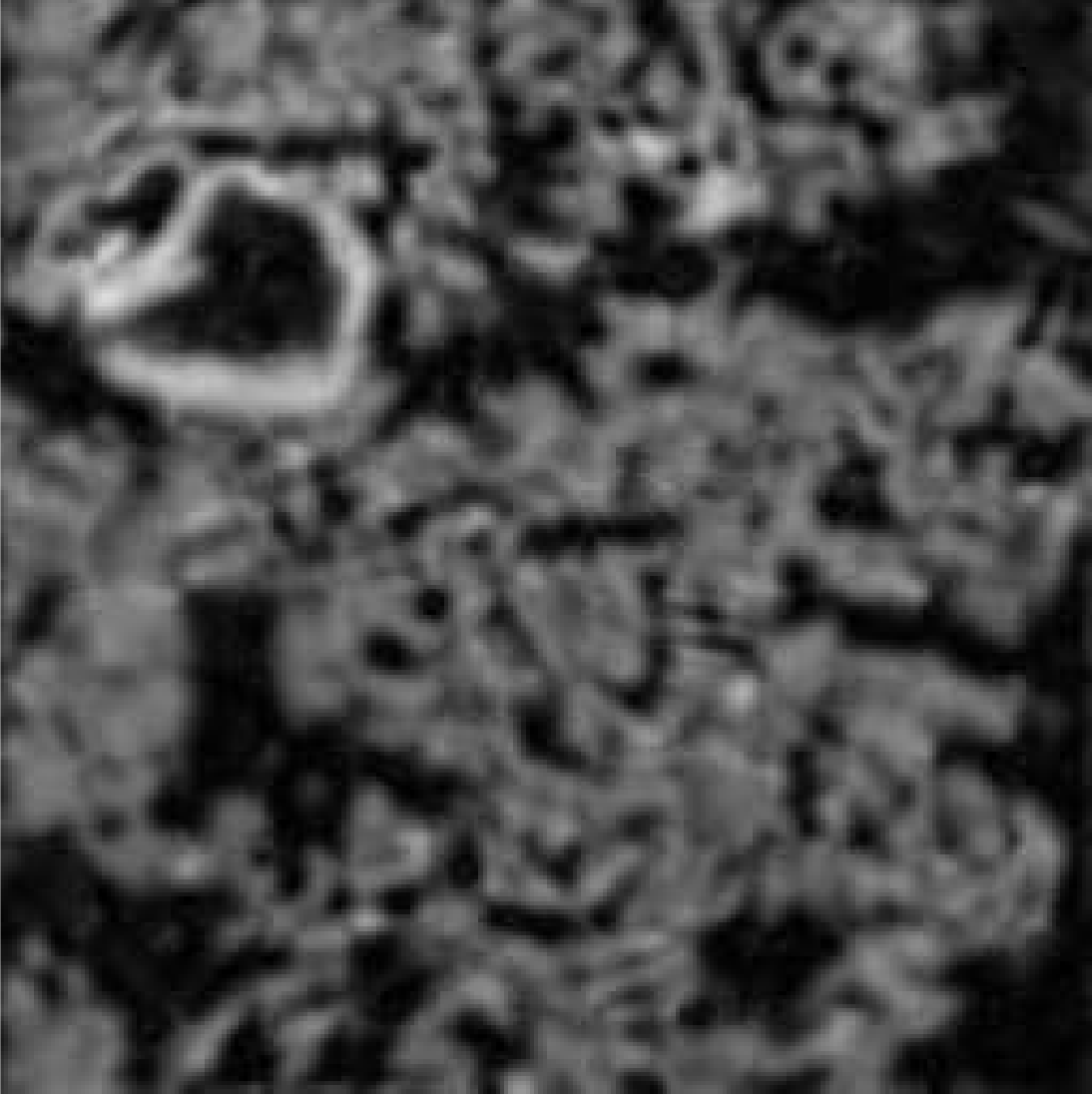}}
\subfigure[]{\includegraphics[width=.11\textwidth]{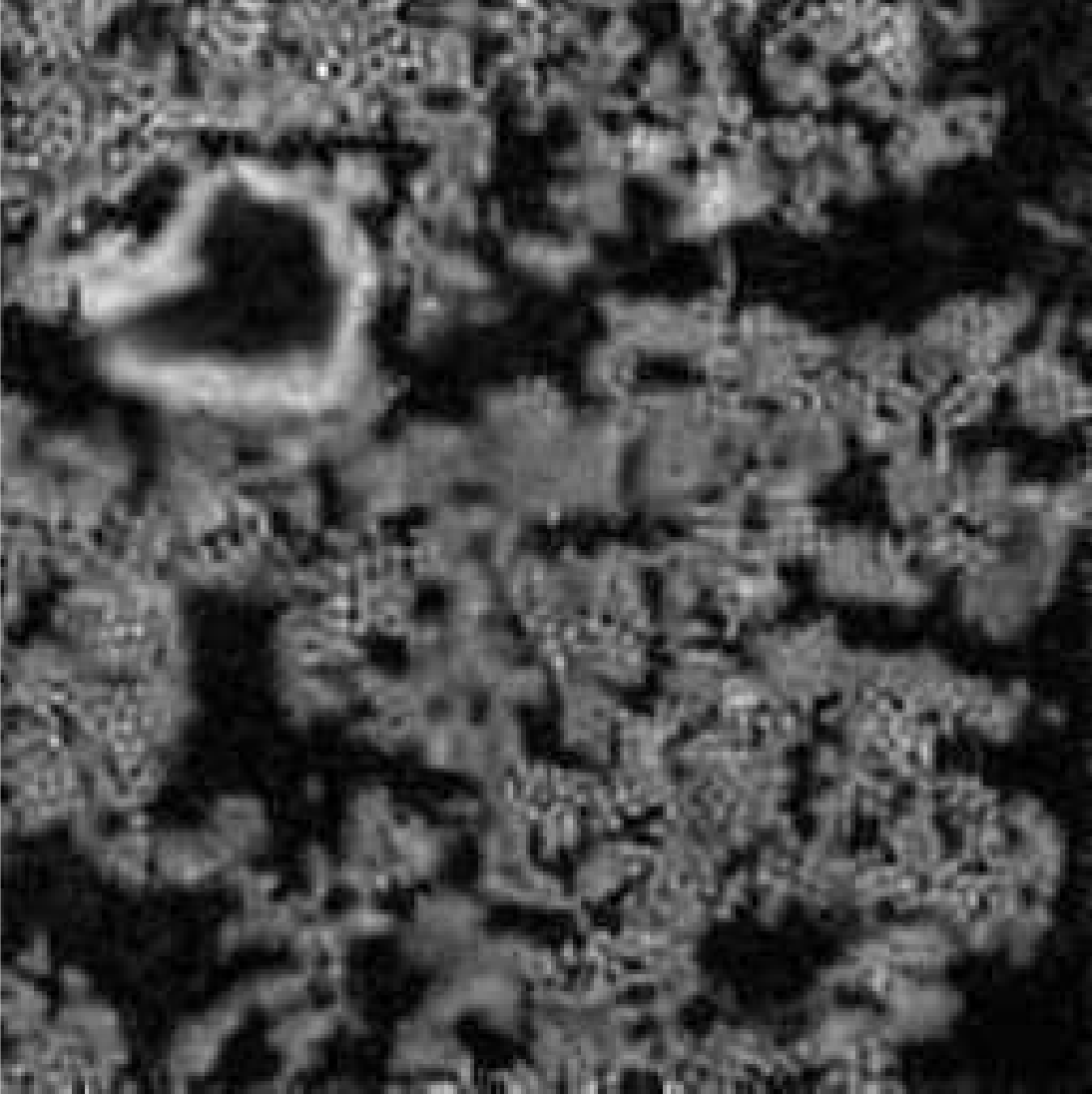}}
\subfigure[]{\includegraphics[width=.11\textwidth]{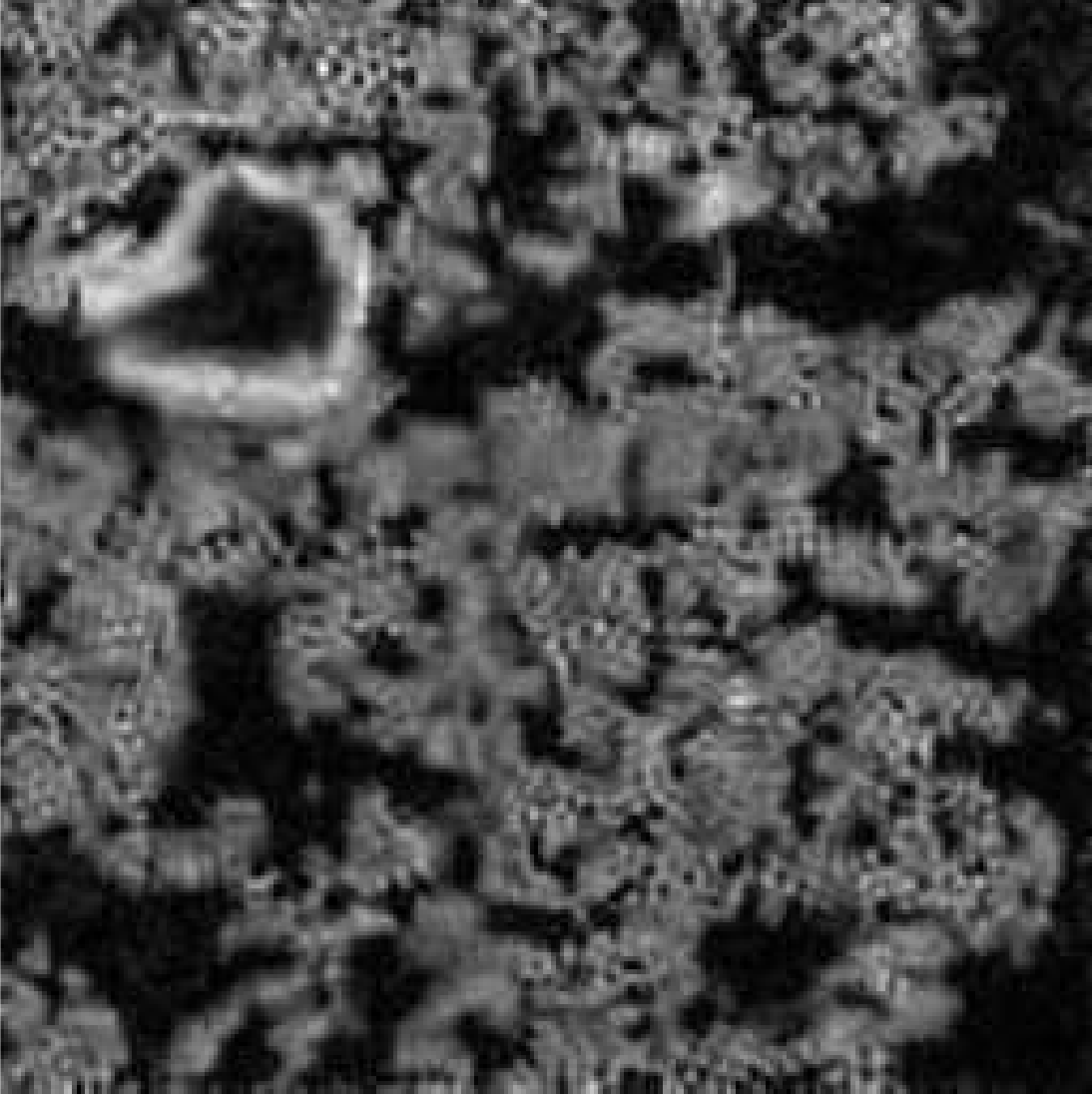}}\\\vskip -.05in
\subfigure[]{\includegraphics[width=.11\textwidth]{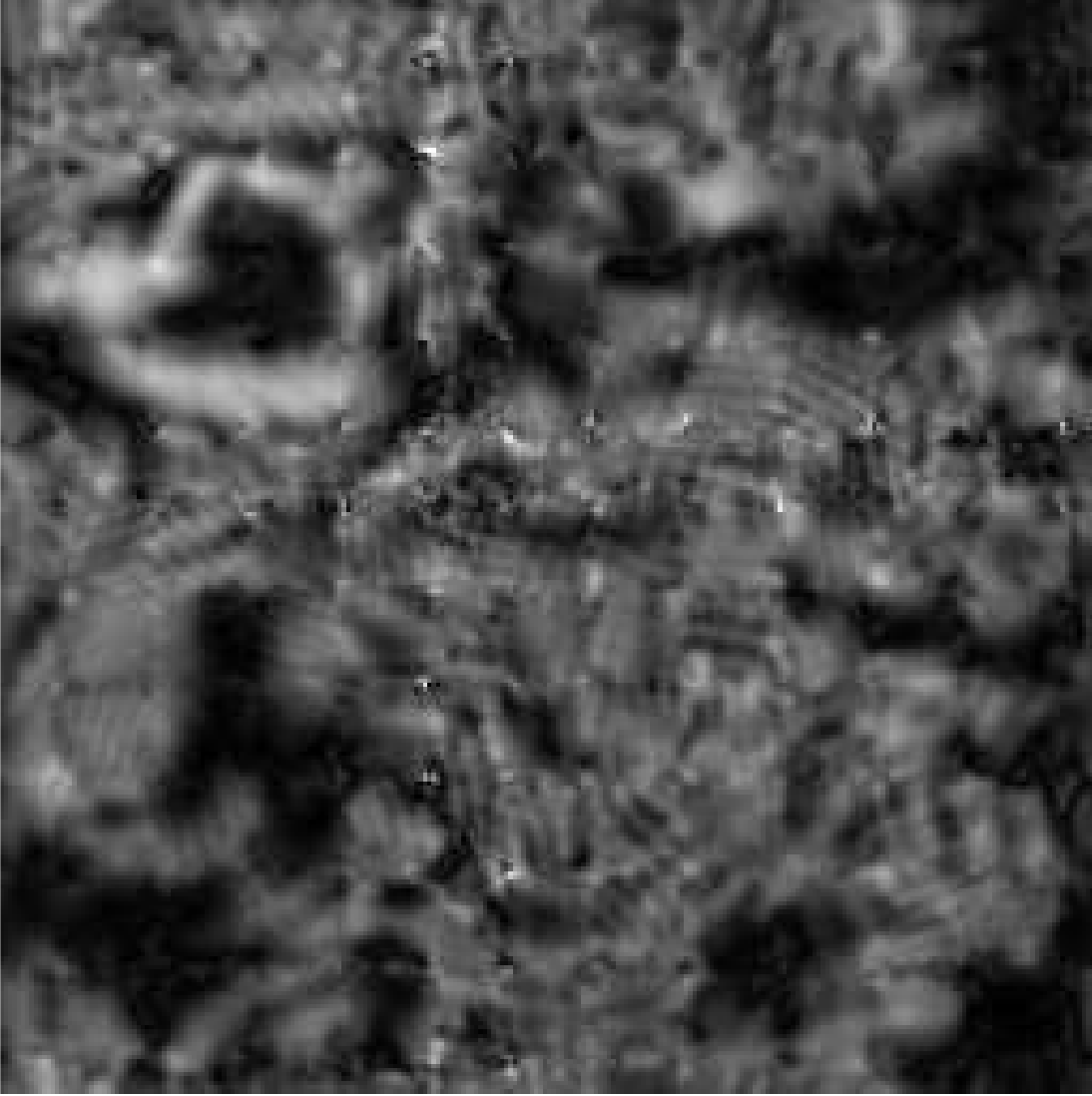}}
\subfigure[]{\includegraphics[width=.11\textwidth]{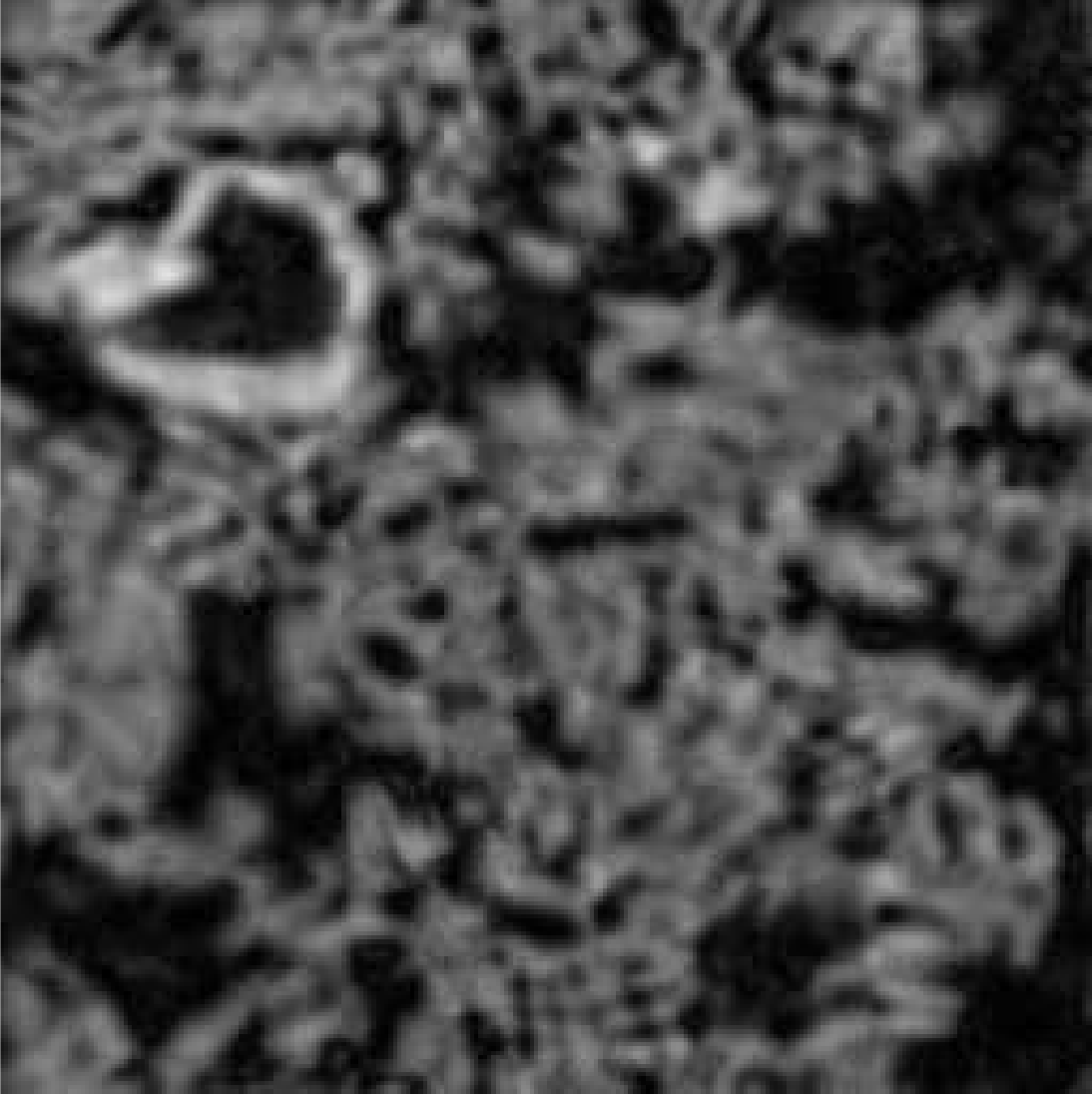}}
\subfigure[]{\includegraphics[width=.11\textwidth]{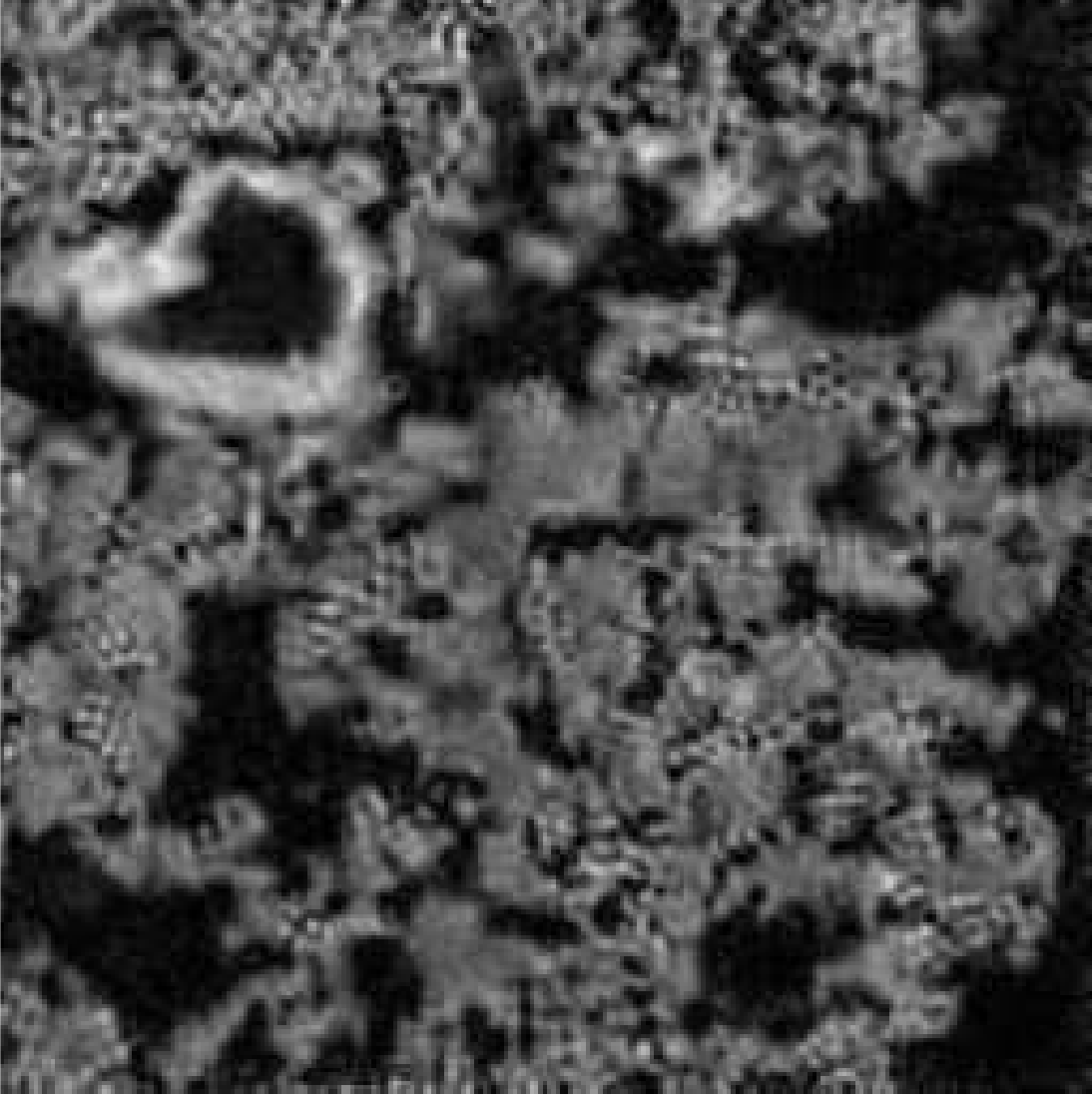}}
\subfigure[]{\includegraphics[width=.11\textwidth]{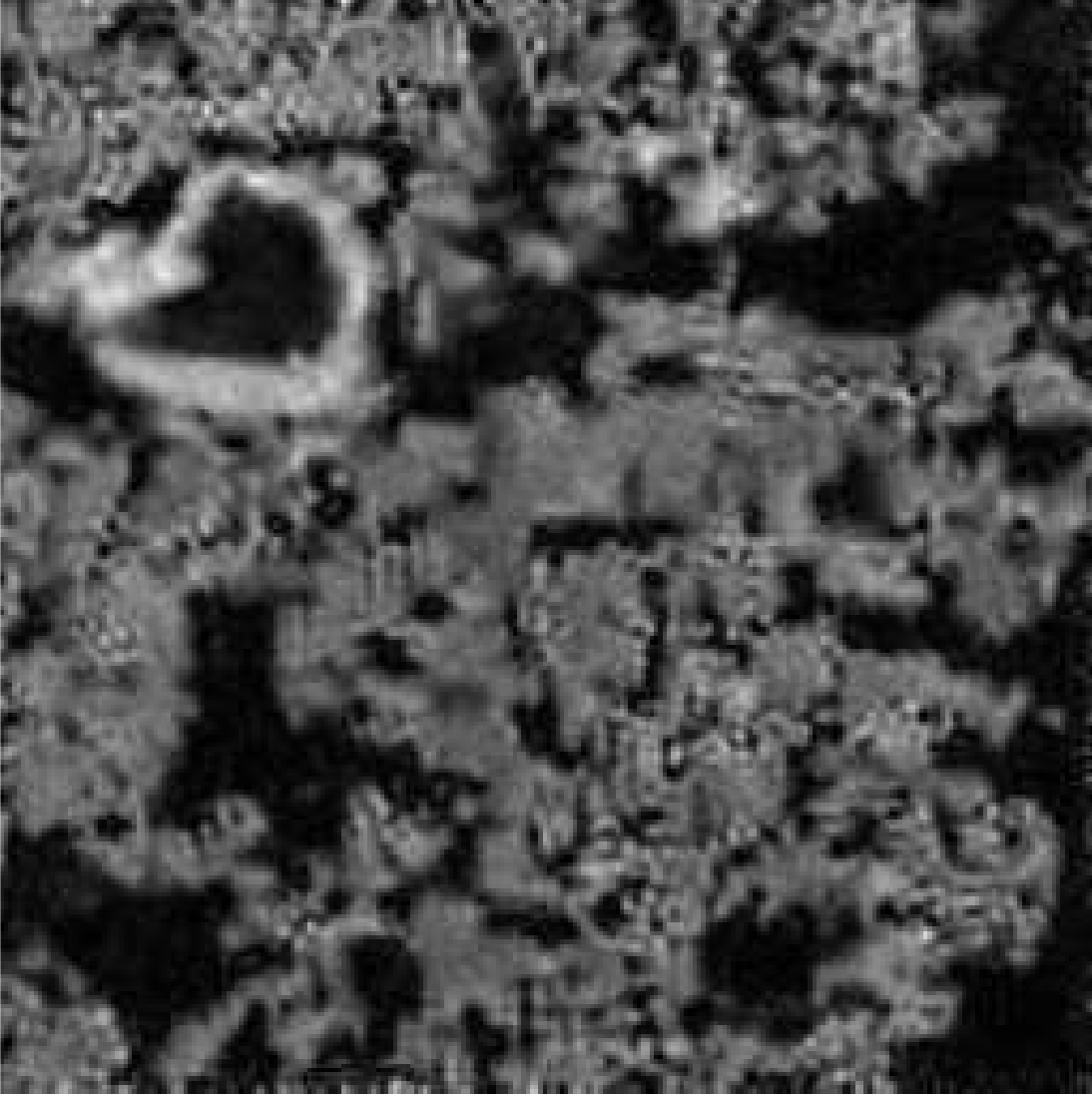}}\vskip -.05in
\end{center}
\vskip -.1in
\caption{PtychoPR for complex-valued image  with different sliding distances (peak level $\delta=0.2$). From first row to the last row: SlidDist$=18,20,22$ respectively. From left to right: reconstructed images by PR in (a), (e) and (i), TVPR in (b),(f) and (j), ALGI in (c), (g) and (k),  and ALGI${}_{aniso}$ in (d), (h) and (l).}
\label{pty2}
\vskip -.05in
\end{figure}

In order to further study the robustness of the proposed algorithms, more experiments with different sliding distances (collecting less data by increasing the sliding distances) are done. Set peak level $\delta=0.2$. Same parameters as the previous tests are used. Performances for PtychoPR with different number of measurements by increasing SlidDist are shown in Fig. \ref{pty2}, where $m/n=12.25, 9,$ and $7.56$ when the sliding distances $\mathrm{SlidDist}=18,20,22$ respectively. On can see that in the first column of Fig. \ref{pty2}, the results of ``PR'' are not only blurry, but also contain some visible structured artifacts.  ``TVPR'' can remove such artifacts, and recover some edges of large scale feature. Our proposed algorithms can further recover most of smaller features in Fig. \ref{pty2}(c) and (d). In an extreme case with $\mathrm{SlidDist}=22$ shown in the third row of Fig. \ref{pty2}, the proposed algorithms can also recover some parts of smaller features. We also put the corresponding SNRs in Table \ref{tab5}, where one can see the obvious gains at about 1.3dB, 1.6dB for ``ALGI'' and ``ALGI${}_{aniso}$'' respectively compared with ``TVPR''.

\begin{table}
\vskip -.15in
\begin{center}
\begin{spacing}{0.9}
\begin{tabular}{|c||c|c|c|c|}
\hline
SlidDist& PR &TVPR &ALGI&ALGI${}_{aniso}$ \\\hline\hline
18&6.50&8.67&10.83&\bf 11.23  \\ \hline
20&6.12&8.27&9.24&\bf 9.32\\ \hline
22&5.97&8.21&8.86&\bf 9.28\\ \hline\hline
Average&6.20&8.38&9.64& \bf9.94\\
\hline
\end{tabular}
\end{spacing}
\end{center}
\vskip -.2in
\caption{SNRs of PtychoPR with different sliding distances for complex-valued image corresponding to Fig. \ref{pty2}}
\label{tab5}
\vskip -.2in
\end{table}

\vskip -.2in
\subsection{Convergence Study}
\begin{figure}
\vskip -.1in
\begin{center}
\subfigure[]{\includegraphics[width=.12\textwidth]{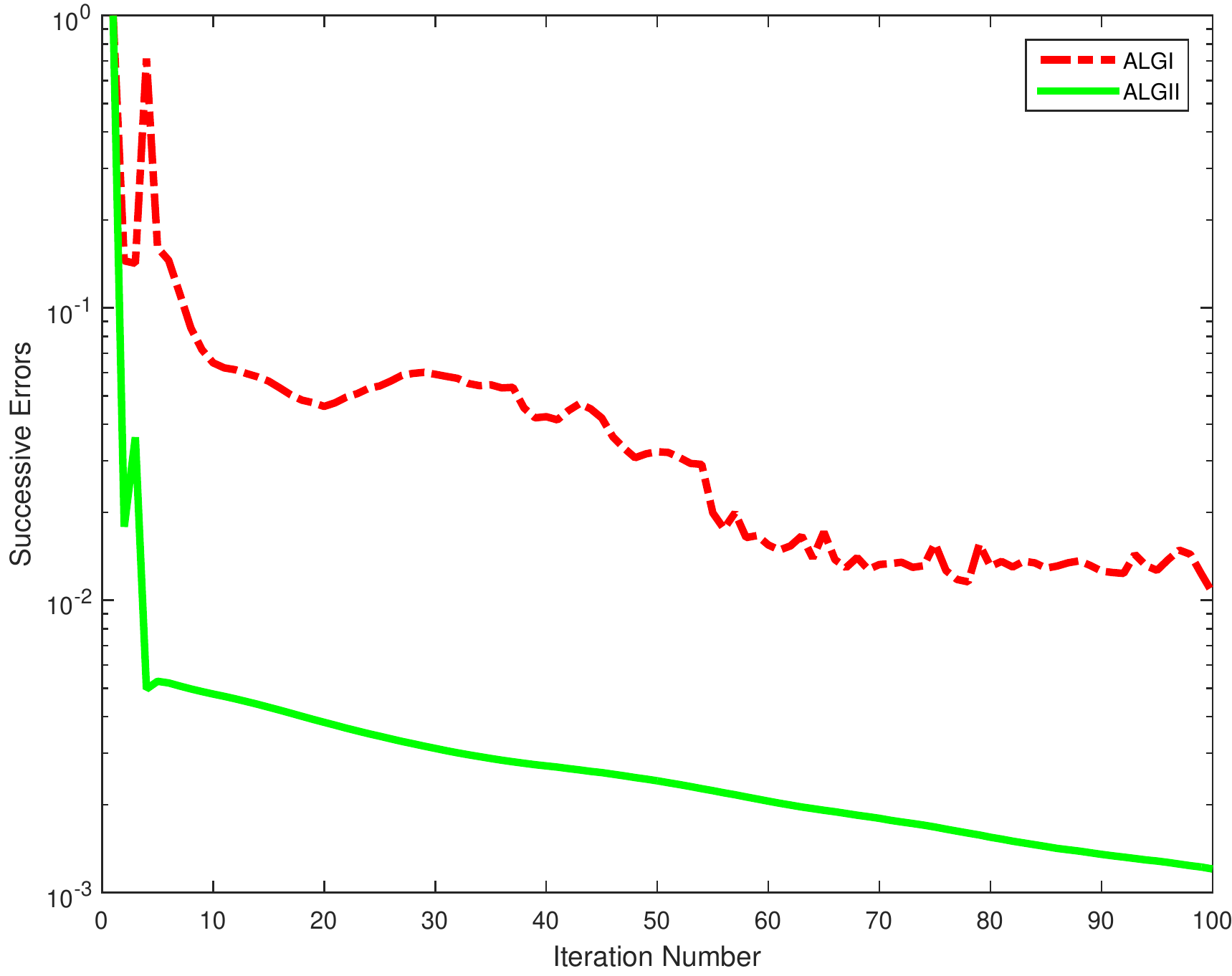}}\quad
\subfigure[]{\includegraphics[width=.12\textwidth]{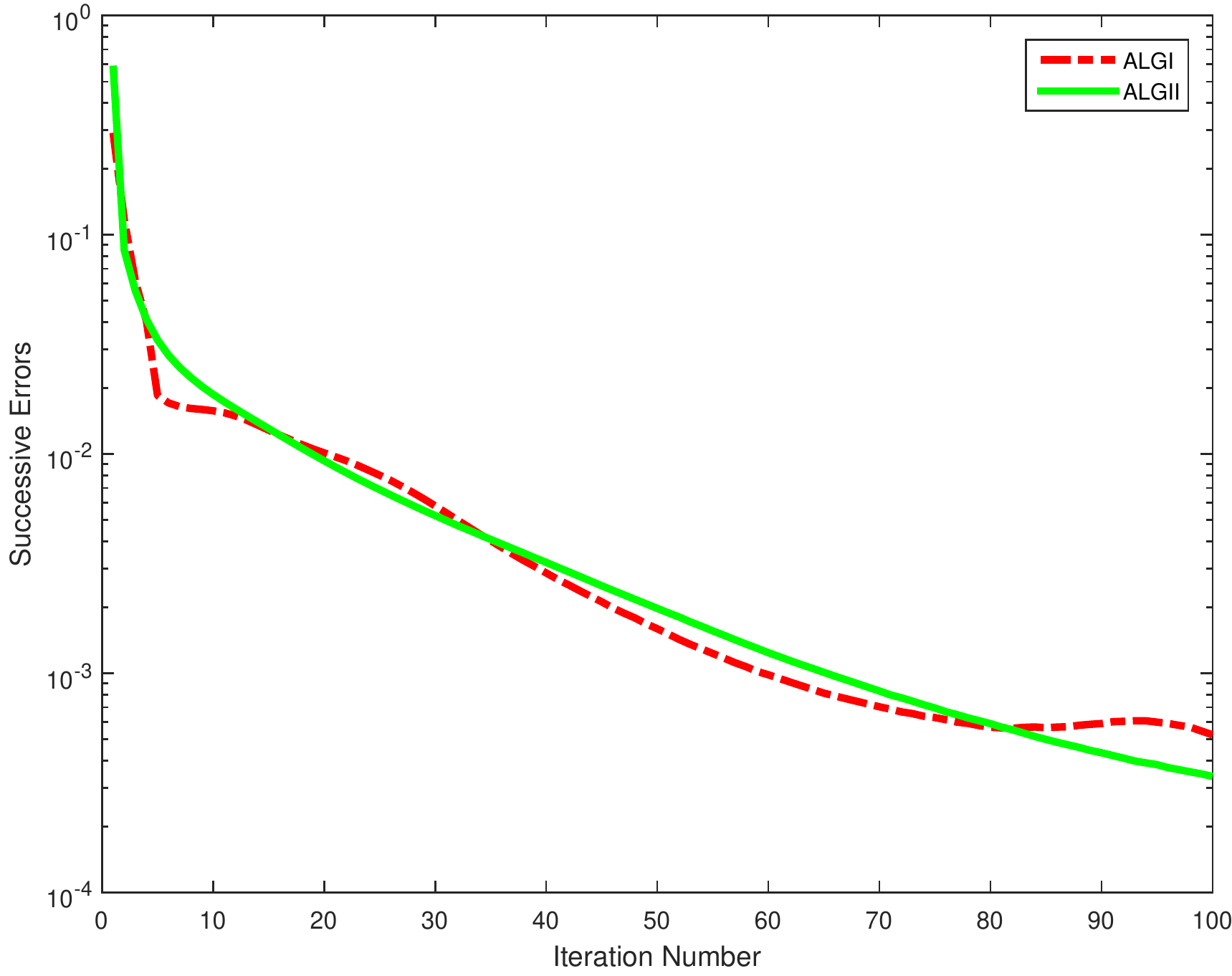}}\quad
\subfigure[]{\includegraphics[width=.12\textwidth]{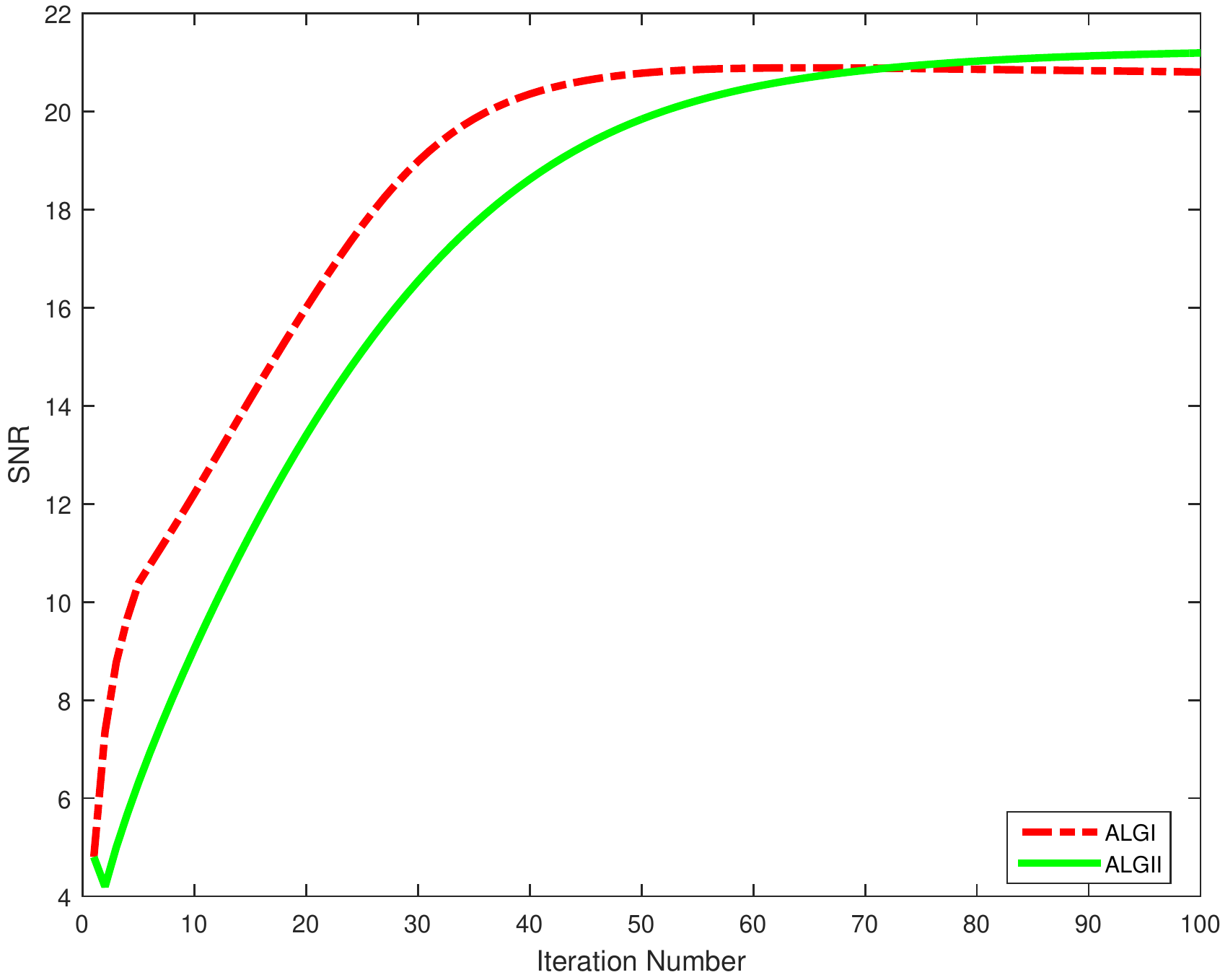}}\\
\subfigure[]{\includegraphics[width=.12\textwidth]{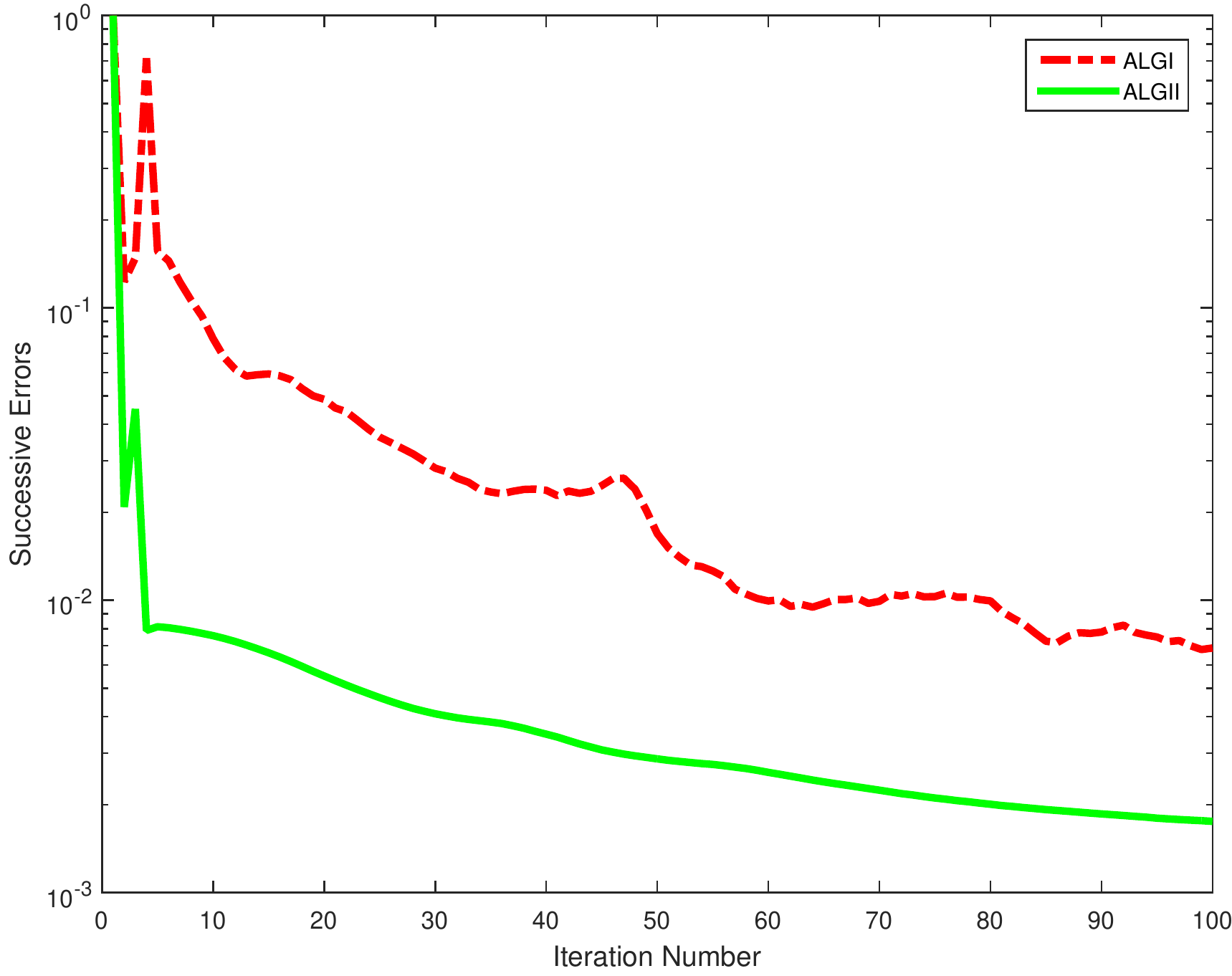}}\quad
\subfigure[]{\includegraphics[width=.12\textwidth]{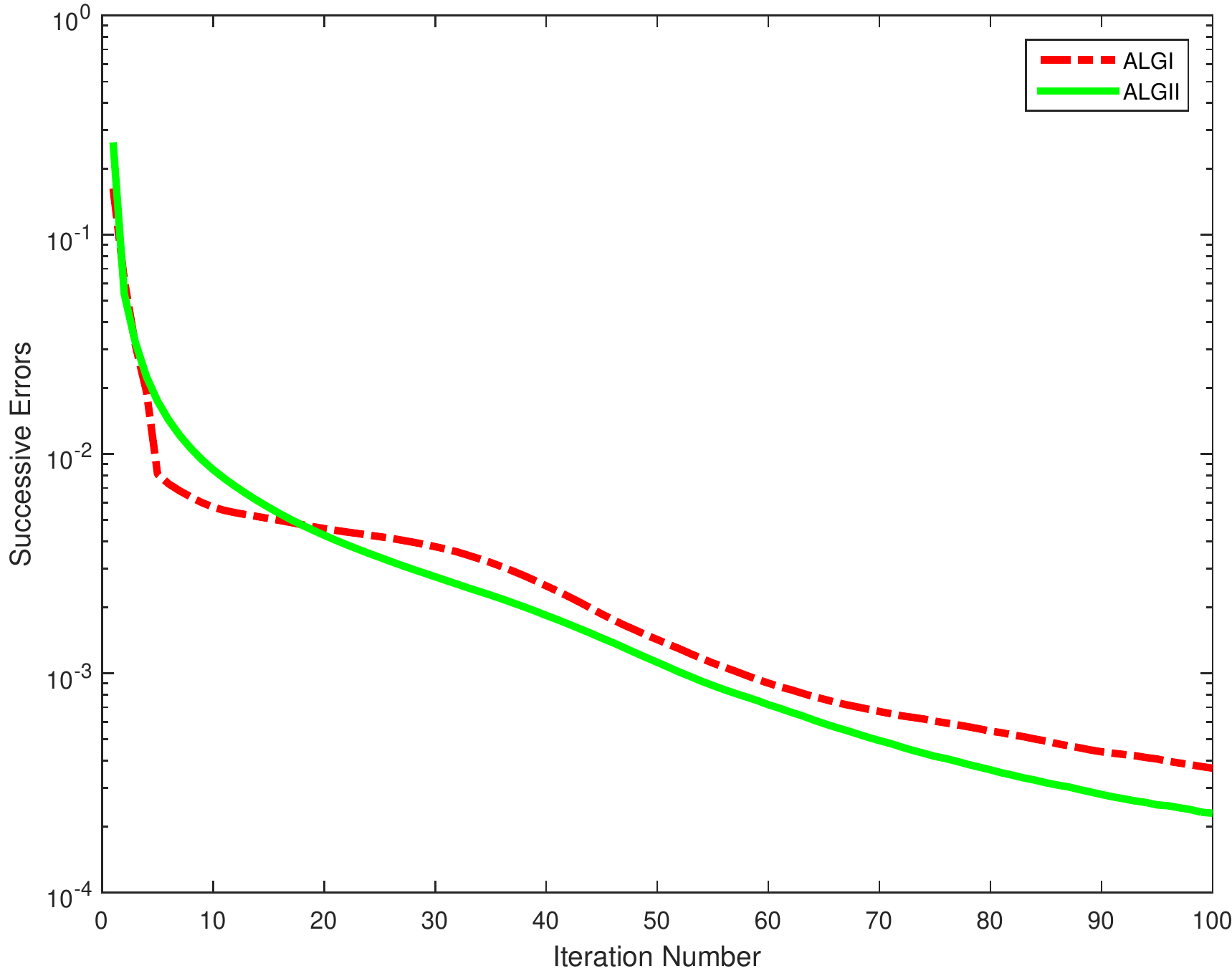}}\quad
\subfigure[]{\includegraphics[width=.12\textwidth]{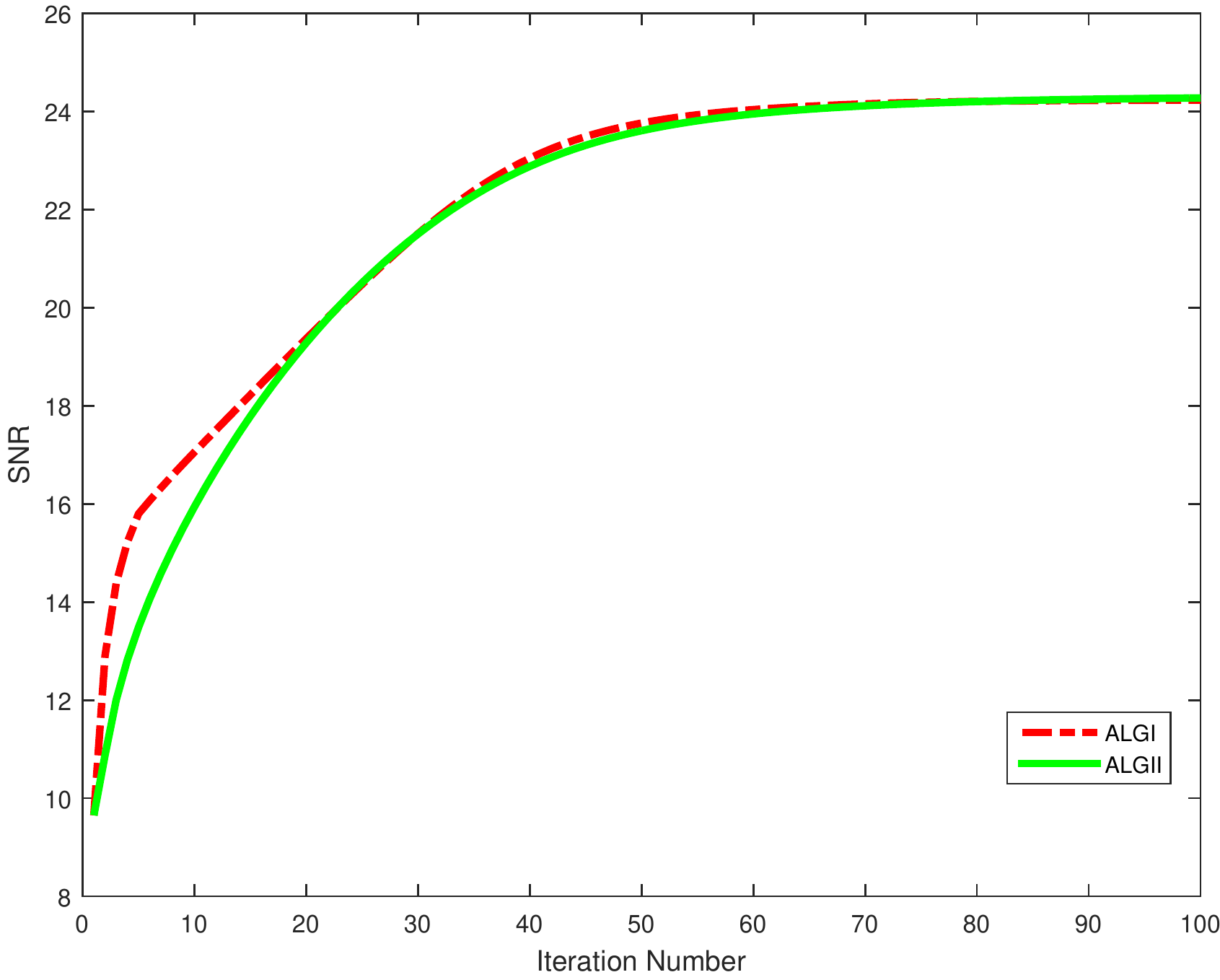}}
\end{center}
\vskip -.2in
\caption{Convergence histories for successive errors and SNRs  on ``Peppers''. First row: $\delta=5.0\times 10^{-3}$; Second row: $\delta=1.0\times 10^{-2}$. Histories of successive error of dictionaries, the iterative solution $u$ and SNRs from the left to right rows respectively.}
\label{fig_convergence}
\vskip -.2in
\end{figure}

To check the convergence of proposed algorithms, we monitor the histories of SNRs and the successive errors of iterative solution $u^k$  and $D^k$ \emph{w.r.t.} the iteration number $k$, which are defined as $\frac{\|u^k-u^{k-1}\|}{\|u^k\|},$ and $\frac{\|D^k-D^{k-1}\|}{\|D^k\|}.$
We show the histories of errors and SNRs in Fig. \ref{fig_convergence} for CDP on ``Peppers'', which implies the proposed algorithms are stable and convergent, which is consistent with the provided theories. It seems that the dictionary converges fast by ``ALGII'' than by  ``ALGI''. By inferred from Fig. \ref{fig_convergence} (c) and (f), SNR first increase, and then get stable, which demonstrates that the proposed algorithms are quite robust. Moreover, when the noise level increases, more iterations are needed. It is very interesting and important to give an optimal iteration condition, and we leave it as a future study.



\subsection{Performances \emph{w.r.t} Parameters}

\begin{figure}[h!]
\begin{center}
\subfigure{\includegraphics[width=.15\textwidth]{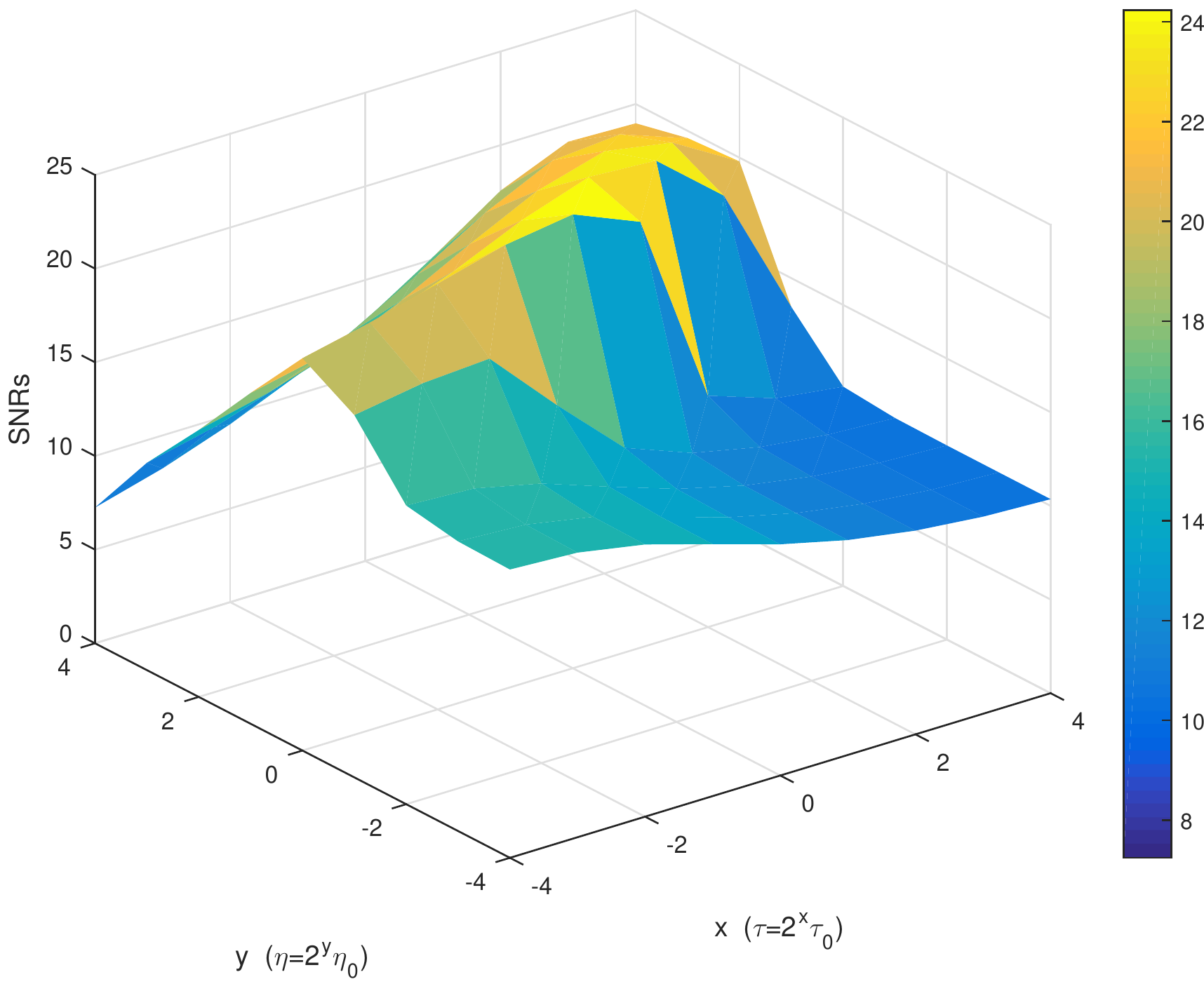}}
\end{center}
\vskip -.15in
\caption{SNRs for ALGI \emph{w.r.t} different parameters $\tau$ and $\eta$.}
\label{figPara}
\vskip -.15in
\end{figure}

Finally we provide some experiments to test the parameter impact, especially on the parameters $\tau$ and $\eta$ which balance  sparsity term and data fitting term.
We show the impact of ``ALGI'' on ``Peppers'' with $\delta=5.0\times 10^{-3}$. Particularly, we select
\[
\begin{split}
(\eta,\tau)\in &\{\eta^0\times2^{-x},\eta^0\times2^{-x+1},\cdots,\eta^0\times2^{x-1},\eta^0\times2^{x}\}\\
&\times \{\tau^0\times2^{-y},\tau^0\times2^{-y+1},\cdots,\tau^0\times2^{y-1},\tau^0\times2^{y}\},
\end{split}
\]
 with $x=y=4$, $\eta_0=8.0\times 10^{-6}, \tau_0=4.5\times 10^{-4}$,  and plot the corresponding SNRs  of the reconstructed images in Fig. \ref{figPara}.  One can see that too large or smaller parameters decrease the quality of reconstructed results. In the future, an automatic schemes for optimal parameter selection should be investigated for the best results.
\section{Conclusion}\label{secCon}
We propose a novel orthogonal dictionary  learning based model to denoise the phaseless measurements contaminated by Possion noise. Meanwhile,  two efficient algorithms are designed, \emph{i.e.} locally convergent AMM and  globally convergent PALM. The experiments demonstrate that both two algorithms are capable of  producing higher quality results, especially preserving the texture features  compared with the  TV  denoising method.  An incoherent dictionary based method and further speed up of the proposed algorithms should be investigated in the future.
\section*{ACKNOWLEDGMENT}
Dr. H. Chang was partially supported by China Scholarship Council (CSC), NSFC Nos. 11426165 and 11501413 and Innovation Project No. 52XC1605 of Tianjin Normal University. This work was also partially funded by the Center
for Applied Mathematics for Energy Research Applications, a joint ASCR-BES funded project within the Office of Science, US Department of Energy, under contract number DOE-DE-AC03-76SF00098.
\bibliographystyle{IEEEtran}

\appendix
One can obtain the following lemma to guarantee the boundedness of iterative solutions for Algorithm II immediately.
\begin{lemma}\label{lem2}
Assuming that there exist two positive constants $e^-,e^+$, such  that
\begin{equation}\label{ass1}
1/2<e^-\leq e^k\leq e^+\forall ~k,
\end{equation}
the sequence  $\{u^k,D^k,\bm\alpha^k\}$  is bounded.
\end{lemma}
\begin{IEEEproof}
Since the dictionary $D^k$ is orthogonal, one can readily get boundedness for it. Therefore we can set
$|(D^{k+1})^*R(u^{k+1})|\leq C_0,$ with a positive matrix $C_0$ independent with $k$, and the notation ``$\leq $'' denotes the pointwise relation for a matrix. By \eqref{solverIIalpha},
 we have
 \[
 \begin{split}
 |\bm\alpha^{k+1}|&\leq \big|1-\frac{1}{e^k}\big| |\bm\alpha^k|+\frac{1}{e^k}|(D^{k+1})^*R(u^{k+1})|\\
 &\leq \big|1-\frac{1}{e^+}\big| |\bm\alpha^k|+\frac{1}{e^-} C_0\\
 &\leq \big|1-\frac{1}{e^+}\big|^2 |\bm\alpha^{k-1}|+\big(|1-\frac{1}{e^+}|+1  \big) \frac{1}{e^-} C_0\\
 &\leq \big|1-\frac{1}{e^+}\big|^{k+1} |\bm\alpha^{0}|+\big( |1-\frac{1}{e^+}|^k+|1-\frac{1}{e^+}|^{k-1}\\
 &\qquad+\cdots+ |1-\frac{1}{e^+}|+1  \big) \frac{1}{e^-} C_0\\
 &=\big|1-\frac{1}{e^+}\big|^{k+1} |\bm\alpha^{0}|+\dfrac{1-|1-\frac{1}{e^+}|^{k+1}}{(1-|1-\frac{1}{e^+}|)|e^-}C_0.
 \end{split}
 \]
By \eqref{ass1}, $|1-\frac{1}{e^+}\big|<1,$ and {finally  one obtains the boundedness of  $\{\bm\alpha^k\}$.}
\end{IEEEproof}

Denote $Z=(u,D,\bm \alpha), Z^k=(u^k,D^k,\bm \alpha^k)$ with the norm $\|Z\|^2=\|u\|^2+\|D\|^2+\|\bm\alpha\|^2$.
The nonincrease of the objective functional can be derived as follows.
\begin{lemma}\label{energyDecay}
$\Upsilon(Z^{k+1})\leq \Upsilon(Z^k)-\frac{\lambda^+}{2}\|Z^{k+1}-Z^k\|^2.$
\end{lemma}
\begin{IEEEproof}
By  Lemma \ref{lem1} and the sufficient decrease property of proximal operator in Lemma 2 of \cite{bolte2014proximal}, one can readily prove it for the three iterative steps in Eqn. \eqref{eq1-1}-Eqn. \eqref{eq1-3}.
\end{IEEEproof}
One can also readily estimate the low bound of the subgradient  of $\mathcal H$  as follows.
\begin{lemma} \label{boundedGrad}
Denote three quantities for the subgradient as
\[
\begin{split}
&A_u^k=c^{k-1}(u^{k-1}-u^k)+\nabla_u\mathcal H(u^k,D^k,\bm\alpha^k)\\
&\qquad-\nabla_u\mathcal H(u^{k-1},D^{k-1},\bm\alpha^{k-1}),\\
&A_D^k=d^{k-1}(D^{k-1}-D^k)+\nabla_D\mathcal H(u^k,D^k,\bm\alpha^k)\\
&\qquad-\nabla_D\mathcal H(u^{k},D^{k-1},\bm\alpha^{k-1}),\\
&A_{\bm\alpha}^k=e^{k-1}({\bm\alpha}^{k-1}-{\bm\alpha}^k)+\nabla_{\bm\alpha}\mathcal H(u^k,D^k,\bm\alpha^k)\\
&\qquad-\nabla_{\bm\alpha}\mathcal H(u^{k},D^{k},\bm\alpha^{k-1}).
\end{split}
\]
For $A^k=(A_u^k,A_D^k,A_{\bm\alpha}^k)$, we have
\begin{equation}\label{eqGrad}
A^k\in \partial \Upsilon(u^k,D^k,\bm\alpha^k);
\end{equation}
and
\begin{equation}
\|A^k\|\leq (3\tilde\lambda+M)\|Z^k-Z^{k-1}\|.
\end{equation}
with $\tilde\lambda:=\max_k\max\{c^{k-1},d^{k-1},e^{k-1},\lambda_D^+,\lambda_D^-,\lambda_{\bm\alpha}^+\}$.
\end{lemma}
\begin{IEEEproof}
Eqn. \eqref{eqGrad} can be readily proved by  computing the first order optimal condition of \eqref{eq1-1} - \eqref{eq1-3}. We just estimate the bounded of $A_D^k$ as a example, and it is similar for the other two.
\[
\begin{split}
\|A_D^k\|&\leq d^{k-1}\|D^{k-1}-D^k\|+\|\nabla_D\mathcal H(u^k,D^k,\bm\alpha^k)\\
&\qquad-\nabla_D\mathcal H(u^{k},D^{k-1},\bm\alpha^{k-1})\|\\
&\leq d^{k-1}\|D^{k-1}-D^k\|+\|\nabla_D\mathcal H(u^k,D^k,\bm\alpha^k)\\
&\qquad-\nabla_D\mathcal H(u^{k},D^{k-1},\bm\alpha^{k})\|+\|\nabla_D\mathcal H(u^k,D^{k-1},\bm\alpha^k)\\
&\qquad-\nabla_D\mathcal H(u^{k},D^{k-1},\bm\alpha^{k-1})\|\\
&\leq d^{k-1}\|D^{k-1}-D^k\|+\lambda_D^+\|D^k-D^{k-1}\|+\lambda_D^-\|\bm\alpha^k-\bm\alpha^{k-1}\|\\
&\leq (d^{k-1}+\lambda_D^+)\|D^{k-1}-D^k\|+\lambda_D^-\|\bm\alpha^k-\bm\alpha^{k-1}\|.
\end{split}
\]
For $A_u^k, A_{\bm\alpha}^k$, we have
$$\|A_u^k\|\leq c^{k-1}\|u^{k-1}-u^k\|+M\|Z^k-Z^{k-1}\|,
$$
and
$$
\|A_{\bm\alpha}^k\|\leq (e^{k-1}+\lambda_{\bm\alpha}^+)\|{\bm\alpha}^{k-1}-{\bm\alpha}^k\|.
$$
By summing up the above three estimates, {pwe conclude to this lemma.}
\end{IEEEproof}

\end{document}